%
%


\magnification = 1200
\input amssym.def
\input amssym.tex
\input epsf

\def \qed {\hfill $\square$}
\def \R {\Bbb R}

\def \C {{\frak C}}

\def \ms {\medskip}
\def \sm {\setminus}

\def \dsp {\displaystyle}
\overfullrule=0pt
\def \diam{\mathop{\rm diam}\nolimits}
\def \dist{\mathop{\rm dist}\nolimits}

\def \Min{\mathop{\rm Min}\nolimits}

\def \length{\mathop{\rm length}\nolimits}
\centerline{\bf $C^{1+\alpha}$-REGULARITY}

\centerline{\bf FOR TWO-DIMENSIONAL ALMOST-MINIMAL SETS
IN $\R^n$}
\vskip 0.5cm
\centerline{ Guy David }
\vskip 1cm

\noindent
{\bf R\'{e}sum\'{e}.}
On donne une nouvelle d\'{e}monstration et une g\'{e}n\'{e}ralisation 
partielle du r\'{e}sultat de Jean Taylor [Ta] 
qui dit que les ensembles presque-minimaux 
de dimension 2 dans $\Bbb R^3$ (au sens d'Almgren) sont localement 
$C^{1+\alpha}$-equivalents \`{a}  des c\^ones minimaux. La 
d\'{e}monstration est plut\^ot \'{e}l\'{e}mentaire, mais utilise quand 
m\^eme un r\'{e}sultat local de s\'{e}paration de [D3] 
et une g\'{e}n\'{e}ralisation [DDT] 
du th\'{e}or\`{e}me du disque topologique de Reifenberg.
L'id\'{e}e directrice est encore que si $X$  est le c\^one sur un arc 
de petit graphe lipschitzien sur la sph\`{e}re, mais n'est pas 
contenu dans un disque, on peut utiliser un graphe de fonction harmonique 
pour d\'{e}former $X$ et diminuer substantiellement sa surface. Le r\'{e}sultat  
de s\'{e}paration local est utilis\'{e} pour se ramener \`{a} des 
unions de c\^ones sur des arcs de graphes Lipschitziens. 
Une bonne partie de la d\'{e}monstration se g\'{e}n\'{e}ralise 
\`{a}  des ensembles de dimension $2$ dans $\R^n$, mais dans ce cadre
notre r\'{e}sultat final de r\'{e}gularit\'{e} sur $E$ d\'{e}pend 
\'{e}ventuellement de la liste des c\^{o}nes minimaux obtenus par 
explosion \`{a} partir de $E$ en un point.

\bigskip \noindent
{\bf Abstract.}
We give a new proof and a partial generalization of Jean Taylor's result [Ta] 
that says that Almgren almost-minimal sets of dimension~$2$ in $\R^3$ 
are locally $C^{1+\alpha}$-equivalent to minimal cones. The proof is rather 
elementary, but uses a local separation result proved in [D3] 
and an extension of Reifenberg's parameterization theorem [DDT]. 
The key idea is still that if $X$ is the cone over an arc of small 
Lipschitz graph in the unit sphere, but $X$ is not contained in a disk, 
we can use the graph of a harmonic function to deform $X$ and diminish 
substantially its area. The local separation result is used to reduce to 
unions of cones over arcs of Lipschitz graphs. A good part of the proof 
extends to minimal sets of dimension $2$ in $\R^n$, but in this setting 
our final regularity result on $E$ may depend on the list of minimal cones 
obtained as blow-up limits of $E$ at a point.

 
\vskip 0.7cm

\noindent {\bf 1. Introduction}
\medskip

This paper can be seen as a continuation of [D3], 
and an attempt to give a slightly different proof of results of 
J. Taylor [Ta] 
and extend some of them to higher codimensions. In [D3] 
we described various sufficient conditions for an almost-minimal set 
$E$ of dimension~$2$ in $\R^n$ to be locally equivalent to a minimal 
cone through a biH\"older mapping.
Here we are interested in getting a more precise (typically, $C^{1,\alpha}$)
local equivalence with a minimal cone. A more precise argument seems 
to be needed to get these, because instead of proving that some 
quantities that measure the closeness to minimal cones stay small
(which may be obtained by compactness and limiting arguments), we shall
need to show that they decay at some definite speed.

Let us try to describe the scheme of the proof. As the reader
will see, the ideas are basically the same as in [Ta], 
even though we shall try to make the argument as modular as possible
and avoid currents. On the other hand, the proof will sometimes become 
more complicated because we want some of it to apply to two-dimensional 
sets in $\R^n$.

Let $E$ be a reduced almost minimal set (see Definition 1.10 below),
fix $x\in E$, and set
$$
\theta(r) = r^{-2} H^2(E\cap B(x,r))
\leqno (1.1)
$$
for $r$ small. Here and throughout the paper, $H^2$ denotes the 
two-dimensional Hausdorff measure.
It is important that if the gauge function
in the definition of almost-minimality  is small enough, then
$\theta(r)$ is almost nondecreasing, to the point that
$$
d(x) = \lim_{t\to 0^+} \theta(t)
\leqno (1.2)
$$
exists. Then we can also show that all the blow-up
limits of $E$ at $x$ are minimal cones with constant density $d(x)$.

For the local biH\"older equivalence of $E$ to a minimal cone, we do 
not need so much more than this, but for the $C^{1,\alpha}$ equivalence,
we shall also need to know that the distance in $B(x,r)$ from $E$ to a 
minimal cone goes to $0$ at a definite speed.
We shall introduce the density excess
$$
f(r) = \theta(r)-d(x) = \theta(r)-\lim_{t\to 0^+} \theta(t)
\leqno (1.3)
$$
and prove that under suitable conditions, it tends to $0$ 
like a power of $r$. See Theorem 4.5 for a differential inequality,
and Lemma 5.11 for the main example of ensuing decay.
Then we shall also show that $f(r)$ controls the Hausdorff distance 
from $E$ to minimal cones in small balls. [See Theorem 11.4 for the
main estimate, and Theorem 12.8 and Proposition~12.28 for consequences
when the gauge function is less than $Cr^b$.
It will then be possible to deduce the $C^{1,\alpha}$ equivalence 
from a variant [DDT] 
of Reifenberg's topological disk theorem [R1]. 
[See Theorem 1.15 and Corollary 12.25.]

The main ingredient in the proof of the fact that $\theta(r)$ is almost 
nondecreasing is a comparison of $E$ with the cone (centered at $x$)
over the set $E \cap \partial B(x,r)$, which gives an estimate
for $H^2(E\cap B(x,r))$ in terms of $H^1(E\cap \partial B(x,r))$.
Similarly, the key decay estimate for $f$ will rely on a comparison
with a more elaborate competitor with the same trace on $\partial B(x,r)$. 

More precisely, we shall first compare $E$ with the cone over 
$E \cap \partial B(x,r)$, but we shall also need to show 
that if $E \cap \partial B(x,r)$ is not close to $C \cap \partial B(x,r)$
for a minimal cone $C$, then we can also improve on the cone over 
$E \cap \partial B(x,r)$ and save a substantial amount of area.

In a way, the central point of the argument is the following 
observation. Let ${\cal C}_{T}$ denote a conical sector in the plane,
with aperture $T \in (0,\pi)$, and let $F : {\cal C}_{T} \to \R^{n-2}$ 
be a Lipschitz function with small norm, which is homogeneous 
of degree $1$. Thus the graph $\Gamma$ of $F$ is a cone in $\R^n$. 
Also suppose that $F(x) = 0$ on the two half lines that form the boundary 
of ${\cal C}_{T}$. The observation is that if $F \neq 0$, it is possible
to replace $\Gamma \cap B(0,1)$ with another surface $\Gamma'$ with
the same boundary, and whose area is substantially smaller.
If the Lipschitz constant is small enough, the main term in the
computation of the area of $\Gamma$ is the energy $\int |\nabla F|^2$,
and a simple Fourier series computation shows that when $F$ is 
homogeneous of degree $1$, it is so far from being harmonic in ${\cal C}_T$
that replacing it with the harmonic extension of its values on
$\partial({\cal C}_{T} \cap B(0,9/10))$ will produce a $\Gamma'$ with
substantially smaller area. Notice that by taking $T < \pi$,
we exclude the unpleasant case when $\Gamma$ is a plane but $F \neq 0$.
See Section 8 for the detailed construction.

A good part of our argument will consist in trying to reduce to this
simple situation. First we shall construct a net $g$ of simple curves in
$E\cap \partial B(x,r)$ (in Section 6). 
This will involve having some control on the topology of $E$, 
and this is the main place where we shall use some results from [D3]. 
Then we shall evaluate how much area we can win when we compare with
the cone over $E \cap \partial B(x,r)$, then with the cone over the net 
$g$, then with the cone over the simplified net $\rho$ obtained from
$g$ by replacing the simple curves with arcs of great circles. 
Finally, we shall also need to evaluate whether we can improve on
$\rho$ when its measure is different from $2rd(x)$.
At this point in the argument, and when $n > 3$, we apparently need to 
add an extra assumption (the full length property of Definition 4.10) 
on the minimal cones that we use to approximate $E$, or the blow-up 
limits of $E$ at $x$.

If we compare with J. Taylor's approach, maybe the main difference 
in the organization of the proof is that here we set aside (in [D3]) 
the topological information that we need to construct a competitor.
In ambient dimension $3$, this amounts to a separation property; in
higher dimensions, we shall use more of the biH\"older 
parameterization provided by [D3]. 
We use this information only once in the proof of the decay estimate, 
to simplify the construction of $\Gamma$ and the competitors. 
We also keep the final geometric part (the construction of a 
parameterization of $E$ once we know that it stays close to minimal cones) 
mostly separate from the main estimate,
and we import it from [DDT]. 

We shall try not to use too much geometric measure theory, but it is 
not clear that this attempt will be entirely successful. First of all, 
some of it is hidden in the rectifiability results for 
almost-minimal sets from [Al] and [DS] 
and the stability of these sets under limits [D1], 
that were used extensively in [D3]. 
We also use a generalization [DDT] 
of Reifenberg's parameterization theorem to situations where $E$ looks like 
cone of type $\Bbb Y$ or $\Bbb T$, both in [D3] 
to get a topological description of $E$ locally,
and more marginally here to say that estimates on density and Hausdorff 
distances to minimal cones can be translated into $C^{1,\alpha}$ estimates. 
Finally, we could not resist using the co-area theorem a few times.

One of the main motivations for this work was to understand 
the techniques of [Ta] 
and then apply them to a similar question for minimizers of the 
Mumford-Shah functional in $\R^3$; see [Le], 
where apparently more than the mere statement of the 
regularity result in [Ta] was needed. 

We also would like to know to which extent the $C^{1,\alpha}$ 
regularity result of [Ta] 
extends to two-dimensional sets in $\R^n$.
Here we shall only prove the $C^{1,\alpha}$ regularity of $E$
near $x$ under some assumption on the blow-up limits of 
$E$ at $x$, and many questions will remain unaddressed, in 
particular because we don't have a list of minimal cones. 

\medskip
Let us be more specific now and try to state our main result.
We shall work with almost-minimal (rather than minimal) sets in 
a domain. The extra generality is probably quite useful, 
and will not be very costly.
Our definition of almost-minimality will involve a gauge function, i.e.,
a nondecreasing function $h : (0,+\infty) \to [0,+\infty]$ such that 
$\lim_{\delta \to 0} h(\delta) = 0$. In  this definition, we are 
given an open set $U \i \R^n$ and a relatively closed subset $E$ of 
$U$, with locally finite Hausdorff measure, i.e., such that 
$$
H^2(E \cap B) < + \infty
\ \hbox{ for every compact ball } B \i U.
\leqno (1.4)
$$

We want to compare $E$ with competitors 
$F=\varphi_{1}(E)$, where $\{\varphi_{t}\}_{0 \leq t \leq 1}\,$, is a 
one-parameter family of continuous functions $\varphi_{t}: U \to U$,
with the following properties:
$$
\varphi_{0}(x) = x  \ \hbox{ for } x\in U,
\leqno (1.5)
$$
$$
\hbox{the function $(t,x) \to \varphi_{t}(x)$, from
$[0,1] \times U$ to $U$, is continuous,}
\leqno (1.6)
$$
$$
\varphi_{1} \hbox{ is Lipschitz}
\leqno (1.7)
$$
(but we never care about the Lipschitz constant in (1.7)) and, if we set 
$$
W_{t} = \{ x\in U \, ; \, \varphi_t(x) \neq x \}
\ \hbox{ and } 
\widehat W = \bigcup_{t\in [0,1]} \Big[ W_{t} \cup \varphi_t(W_{t}) 
\big],
\leqno (1.8)
$$
then
$$
\widehat W \hbox{ is relatively compact in $U \ $ and }
\diam(\widehat W) < \delta.
\leqno (1.9)
$$

\medskip
\proclaim Definition 1.10.
We say that the closed set $E$ in $U$ is an almost-minimal set 
(of dimension~$2$) in $U$, with gauge function $h$, 
if (1.4) holds and if 
$$
H^2(E \setminus F) \leq H^2(F \setminus E)
+ h(\delta) \, \delta^2
\leqno (1.11)
$$
for each $\delta > 0$ and each family $\{\varphi_t \}_{0 \leq t \leq 1}$ 
such that (1.5)-(1.9) hold, and were we set $F = \varphi_1(E)$.
We say that $E$ is reduced when
$$
H^2(E \cap B(x,r)) > 0 \ \hbox{ for $x\in E$ and $r>0$}.
\leqno (1.12)
$$

\medskip
This definition is of course inspired by Almgren's definition of 
``restricted sets" [Al], 
but we use a slightly different accounting.
We could use other definitions of almost-minimal sets; 
for instance, we could have replaced (1.11) with
$$
H^2(E \cap W_{1}) \leq H^2(\varphi_1(E \cap W_{1}))
+ h(\delta) \delta^2,
\leqno (1.13)
$$
and this would have yielded exactly the same class of 
almost-minimal sets with gauge $h$ (see Proposition 4.10 in [D3]). 
Also see Definition 4.1 in [D3] for a slightly smaller class of 
almost-minimal sets, with a definition closer to Almgren's.

Denote by $E^\ast$ the closed support of the restriction
of $H^2$ to $E$, i.e., set
$$
E^\ast = \big\{ x\in E \, ; \, H^2(E \cap B(x,r)) > 0 
\, \hbox{ for every } r>0  \big\}.
\leqno (1.14)
$$
Thus $E$ is reduced when $E = E^\ast$. It is fairly easy to see
that $E^\ast$ is closed in $U$, $H^2(E \setminus E^\ast)=0$, and $E^\ast$ is
almost-minimal with gauge $h$ when $E$ is almost-minimal 
with gauge $h$ (see Remark 2.14 in [D3]). 
Thus we may restrict our attention to reduced almost-minimal 
sets, as the other almost-minimal sets are obtained from the reduced
ones by adding sets of vanishing Hausdorff measure.

Even when we do not say this explicitly, the almost minimal sets
and minimal cones in this paper will be assumed to be reduced.

Even when $E$ is a reduced minimal set, it can have singularities, 
but the list of singularities is fairly small, especially when $n=3$.

First consider minimal cones, i.e., minimal sets that are also cones.
It was proved by E. Lamarle [La], 
A. Heppes [He], and J. Tayor [Ta], 
that there are exactly three types of (nonempty) 
reduced minimal cones in $\Bbb R^3$: the planes, the
sets of type $\Bbb Y$ obtained as unions of three half planes 
with a common boundary $L$ and that make $120^\circ$ angles along $L$, 
and sets of type $\Bbb T$ which are cones $T$ over the union of the 
edges of a regular tetrahedron, 
and centered at the center of the tetrahedron. 
See Figures 1.1 and 1.2.

\vskip 0.6cm  
\hskip 1.7cm  
\epsfxsize = 2.5cm 
\epsffile{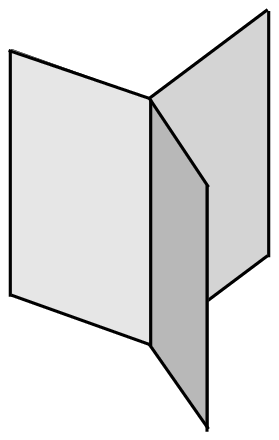}
\hskip 2.5cm 
\epsfxsize = 5.1cm \epsffile{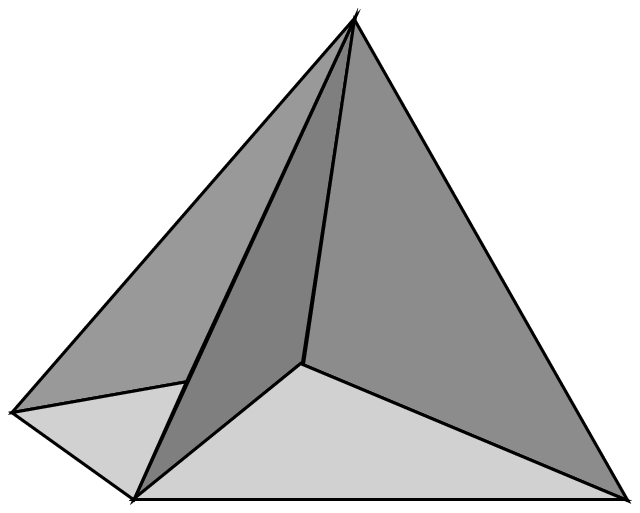}
\medskip
\noindent\hskip 1.2cm{\bf Figure 1.1.} A cone of type $\Bbb Y$ 
\hskip 1.9cm {\bf Figure 1.2.} A cone of type $\Bbb T$
\medskip 

When $n > 3$, the precise list of reduced minimal cones (centered 
at the origin) is not known, but they have the same structure as in $\R^3$. 
The set $K= E\cap \partial B(0,1)$ is composed of a finite collection
of great circles and arcs of great circles, which can only meet by sets
of three and with $120^\circ$ angles. See the beginning of Section 2
for a more precise description. 
We still have the $2$-planes (which correspond to a single circle) 
and the sets of type $\Bbb Y$ (where $K$ is the union of three half 
circles that meet at two antipodal points), and we call 
``sets of type $\Bbb T$" all the other nonempty reduced minimal cones. 
But for $n>3$ the list of sets of type $\Bbb T$ is larger, and not well 
investigated.

It turns out that if $E$ is a reduced almost-minimal set with a small 
enough gauge function $h$, then for $x\in E\cap U$ we have a reasonably good 
local description of $E$ near $x$. In particular, it is shown in [D3] 
(after [Ta] when $n=3$) 
that for $r$ small, $E$ coincides near $B(x,r)$ with the image of 
a reduced minimal cone by a biH\"older homeomorphism of space.

Our main statement here is roughly the same, except that we
require $h$ to be slightly smaller, add an assumption on the
list of blow-up limits of $E$ at $x$ when $n > 3$, but 
obtain a better regularity.

\medskip
\proclaim Theorem 1.15. 
Let $E$ be a reduced almost-minimal set in the open set $U \i \R^n$, 
with a gauge function $h$ such that $h(r) \leq C r^{\alpha}$ for $r$ 
small enough (and some choice of $C \geq 0$ and $\alpha >0$).
Let $x \in E$ be given. If $n > 3$, suppose in addition that 
some blow-up limit of $E$ at $x$ is a full length 
minimal cone  (See Definition 2.10). 
Then there is a unique blow-up limit $X$ of $E$ at $x$, and
$x+X$ is tangent to $E$ at $x$. In addition, there is a radius 
$r_0 > 0$ such that, for $0 < r < r_0$, there is a $C^{1,\beta}$ 
diffeomorphism $\Phi : B(0,2r) \to \Phi(B(0,2r))$, such that $\Phi(0)=x$ 
and $|\Phi(y) -x-y| \leq 10^{-2}r$ for $y\in B(0,2r)$,
and $E \cap B(x,r) = \Phi(X) \cap B(x,r)$.

\ms
Jean Taylor's result from [Ta] corresponds to $n=3$. 

See Definition 3.11 for the definition of blow-up limits, 
and Definition 2.10 for full length minimal cones. 
Let us just say here that every (reduced) minimal cone in $\R^3$
is full length. This is also the case of the explicitly 
known minimal cones in higher dimension. Thus we do not know 
whether every minimal cone is full length. See Section 14.

Roughly speaking, the full length condition
requires that if one deforms $X \cap \partial B(0,1)$ into a 
nearby net of arcs of great circles which is strictly longer, then the
cone over this deformation is far from being a minimal cone.

In the statement, $\beta > 0$ depends only on $n$, $\alpha$, and 
the full length constants for $X$ if $n>3$. The $C^{1,\beta}$ 
diffeomorphism conclusion means that $D\Phi$ and its inverse are 
H\"older-continuous with exponent $\beta$. 

Unfortunately, we have no lower bound on $r$, even when $n=3$
and $x$ is a point of type $\Bbb T$. [We could perhaps have hoped 
for lower bounds on $r$ that depend only on $h$ and the distances 
from $x$ to $\partial U$ and to the points of $E$ of ``higher types".]

\ms
Theorem 1.15 can be generalized in two ways. 
First, we can use gauge functions $h$ that do not decrease as fast 
as a power, but we shall always assume that
$$
h_1(r) = \int_0^r h(2t) {dt \over t} 
\ \hbox{ is finite for $r$ small enough},
\leqno (1.16)
$$
because we shall rely on a near monotonicity- 
and a parameterization result from [D3] where this is used. 

If we assume a suitably strong Dini condition (somewhat stronger than 
(1.16)), we can proceed as in Theorem 1.15, and get a decay
for the density excess $f(r)$ and the numbers $\beta_Z(x,r)$ that 
measure the good approximation of $E$ by minimal cones
which is sufficient to give the uniqueness of the tangent cone 
$X$ and a local $C^1$ equivalence to $X$. 

If we only assume (1.16), and this time assume that every
blow-up limit of $E$ at $x$ has the full length property, we still
get some decay of $f(r)$ and the $\beta_Z(x,r)$, but not enough 
to imply the $C^1$ equivalence of $E$ to a minimal cone or even the 
uniqueness of the tangent cone.
See Section 13 for a little more detail.

We can also replace the full length condition with a weaker
one, but this will lead to a slow decay of $f$ that
does not yield $C^1$ estimates. See Section 13 again. 

The plan for the rest of the paper is the following.
Section 2 is devoted to a rapid description of the minimal cones,
and the definition of the full length property. 

In Section 3 we define the density $\theta(r)$ and the density
excess $f(r)$, record simple consequences of
the almost monotonicity of density for almost minimal sets,
and say a few words about blow-up limits.

The main decay estimate for $f$ will come from a differential 
inequality that follows from a comparison inequality. 
This comparison inequality is stated in Section 4 (Theorem 4.5), 
where we also check that two technical assumptions are satisfied 
almost-everywhere. The proof takes up most of Sections 6-10.

In Section 5 we explain why Theorem 4.5 leads to a differential
inequality and a decay estimate for $f$.

The proof of Theorem 4.5 really starts in Section 6, 
and is performed on a single ball $B=B(x,r)$; the long-term goal
is to construct a nice competitor for $E$ in $B$, and thus
get an upper bound on $H^2(E\cap B)$.
In Section 6 we construct a net $g$ of simple curves $g_{j,k}$ 
in $\partial B$, that looks like the intersection of 
$\partial B$ with a minimal cone (see Lemma 6.11). 
The construction uses some separating properties of $E$.

In Section 7 we approximate the simple curves $g_{j,k}$
with Lipschitz curves $\Gamma_{j,k}$; this will be useful
because it is easier to find retractions on a neighborhood
of the cone over $\cup \Gamma_{j,k}$. Most of the section 
consists in checking that we diminish the length proportionally 
to the amount of $g_{j,k}$ that we replace with $\Gamma_{j,k}$.
The main ingredients are simple one-dimensional computations and 
the Hardy-Littlewood maximal function.

Section 8 explains how to replace the cone over a small Lipschitz 
curve like the $\Gamma_{j,k}$ with a Lipschitz surface with the
same trace on $\partial B$, but with smaller area.
The idea is to use the approximation of the minimal surface 
equation by the Laplace equation and observe that the cone
is the graph of a function which is homogeneous of degree $1$,
whence reasonably far from harmonic if it is not almost affine.
The estimate relies on simple Fourier series computations.

A first competitor is defined in Section 9. The main point is
to map most of $E \cap B$ to the finite union of Lipschitz
surfaces constructed in Section 8, without losing too much near
$\partial B$. Most of the area estimates for Theorem 4.5 are done 
at the end of the section.

Section 10 contains the final epiperimetry argument, which concerns 
the possible improvement of $E$ when it is equal, or very close to 
the cone over a union $\rho$ of geodesics (but yet is different from the 
approximating minimal cone $X$). This is the place where the geometry 
of $X$ plays a role, and where the full length property is used.

In Section 11, we switch to a slightly different subject, and show
that when $f(r)$ is very small, $E \cap B(0,C^{-1}r)$ is well approximated
by a minimal cone (but yet we lose a power in the estimate). See
Theorem 11.4 and Corollary 11.85 for the main statements. The proof uses 
the same competitor as in Section 9, and a slightly unpleasant construction
of transverse curves, which we use to control the way in which the
good approximating minimal cone that we get on almost every concentric 
sphere depends on the sphere.

Section 12 contains the consequence of the decay estimate 
(Theorem 4.5) and the approximation results (Theorem 11.4
and Corollary 11.85), when $h(r) \leq C r^b$, in terms of better 
approximation by cones and local $C^1$-equivalence to a minimal
cone. See Theorem~12.8, Proposition 12.28, and Corollary 12.25.
Theorem 1.15 is a special case of this.

Section 13 contains some technical improvements, concerning better 
estimates for Proposition 12.28 and Corollary 12.25, larger gauge functions 
(sufficiently large powers of ${\rm Log}(C/r)$), and weaker variants
of the full length property.

Finally the full length property is discussed again in Section 14,
where it is established for minimal cones in $\R^3$
and a few others.

\ms
The author wishes to thank T. De Pauw, J.-C. L\'{e}ger, P. Mattila,
and F. Morgan for discussions and help with J. Taylor's results,
and gratefully acknowledges partial support from the 
european network HARP. The pictures were done with the help of Inkscape.

\bigskip 
\noindent {\bf 2. Full length minimal cones.}
\medskip

We start with a rough description of the reduced 
minimal cones of dimension $2$ in $\R^n$, that we take 
from [D3]. We start with simple definitions. 

For us a cone centered at the origin will be a positive cone, 
i.e., a set $X$ such that $\lambda x \in X$ when $x\in X$ and
$\lambda \geq 0$. In general, a cone will just be a 
translation of a cone centered at the origin. 

A minimal set (in $\R^n$) is the same thing as an almost minimal 
set in $U = \R^n$ with the gauge function $h=0$. 
That is, $E$ is closed, $H^2(E\cap B) < + \infty$ for every 
ball $B$, and
$$
H^2(E \setminus f(E)) \leq H^2(f(E) \setminus E)
\leqno (2.1)
$$
for every Lipschitz function $f$ such that $f(x)=x$ out
of some compact set. The reader should maybe be warned that
the literature contains many other definitions of ``minimal sets"
or soap films.

By minimal cone, we simply mean a cone which is also
a minimal set in $\R^n$. We shall always assume that it is
reduced, even if we don't say it explicitly.

Observe that the $2$-planes, and the sets of type $\Bbb Y$
(three half $2$-planes with a common boundary $L$, and that make 
$120^\circ$ angles along $L$) are still minimal in any ambient 
dimension $n \geq 3$. All the other minimal cones will be called
cones of type $\Bbb T$.

The simplest example of a cone of type $\Bbb T$ is obtained by 
taking a regular tetrahedron in some affine $3$-space, and then 
letting $T$ be the cone over the union of the edges centered at 
the center of the tetrahedron. We get a minimal set, even if
$n > 3$ (because competitors could always be projected back 
to $\R^3$). 

When $n=3$, Lamarle [La], Heppes [He], and Taylor [Ta] 
showed that these are the only cones of type $\Bbb T$ that we can get.
When $n > 3$, the union of two orthogonal $2$-dimensional planes
is a minimal cone, and possibly this stays true if the two planes
are almost orthogonal. Also, the product of two one-dimensional sets
$Y$ lying in orthogonal two-planes is likely to be minimal.
And there may be lots of other, much wilder examples, but the truth is 
that the author does not know any. 

At least we know the following rough description of minimal
cones. Let $E$ be a (reduced) minimal cone of dimension $2$ in $\R^n$,
and suppose for convenience that $X$ is centered at the origin.
Set $K = E \cap \partial B(0,1)$. Then $K$ is a finite union of 
great circles and arcs of great circles $\C_j$, $j\in J$. 
The $\C_j$ can only meet when they are arcs of great circles, 
and only by sets of $3$ and at a common endpoint. Then they make angles of 
$120^\circ$ at that point. The arcs of circles have no free endpoints either 
(each endpoint is the end of exactly three arcs). 
In addition, there is a small constant $\eta_0 > 0$, that depends only
on $n$, such that
$$
\hbox{the length of each arc $\C_j$ is at least $10\eta_0$}
\leqno (2.2)
$$
and
$$\eqalign{
&\hskip 0.3cm \hbox{if $i, j \in J$, $i\neq j$, and $x\in \C_i$ is such that}
\dist(x,\C_j) \leq \eta_0, \hbox{ then $\C_i$ and $\C_j$}
\cr&\hbox{are both arcs of circles, and they have a common endpoint in }
\overline B(x,\dist(x,\C_j)).
}\leqno (2.3)
$$
Next set
$$
d(E) = H^2(E\cap B(0,1)) = {1 \over 2} \, H^1(K)
= {1 \over 2} \,\sum_{j\in J} \,\length(C_j),
\leqno (2.4)
$$
where the last equalities come from the fact that $E\cap B(0,1)$ is the
essentially disjoint union of the angular sectors over the arcs $\C_j$.
Recall from Lemma 14.12 in [D3] 
that there is a constant $d_T > {3 \pi \over 2}$ such that
$d(E) \geq d_T$ when $E$ is a minimal cone of type $\Bbb T$.

\ms

It will be easier to manipulate $K$ after we cut some of the
$\C_j$ into two or three shorter sub-arcs. For each $\C_j$
whose length is more than ${9 \pi \over 10}$, say, we
cut $\C_j$ into two or three essentially disjoint sub-arcs
$\C_{j,k}$ with roughly the same length. 
For the other arcs $\C_j$, those whose length is at most 
${9 \pi \over 10}$, we leave them as they are, i.e., decompose
them into a single arc $\C_{j,k}$. Altogether, we now have a decomposition
of $K$ into essentially disjoint arcs of circles $\C_{j,k}$,
$(j,k) \in \widetilde J$, such that
$$
10\eta_0 \leq \length(\C_{j,k}) \leq {9 \pi \over 10}
\ \hbox{ for } (j,k) \in \widetilde J.
\leqno (2.5)
$$
The $\C_{j,k}$ will be easier to use, because they
are geodesics of $\partial B$ that are determined by
their endpoints in a reasonably stable way.

Let us also distinguish two types of vertices. 
We denote by $V$ the collection of all the endpoints of the
various $\C_{j,k}$. Then let $V_0$ be the set of vertices
$x\in V$ that lie at the end of exactly three arcs $\C_{j,k}$.
Those are the vertices that were already present in
our initial decomposition of $K$ into $\C_j$. Finally set
$V_1 = V \setminus V_0$. Those are the endpoints that we
just added, and each one lies exactly in two arcs $\C_{j,k}$.
Notice that because of (2.3), if $(j,k)$ and $(j',k') \in \widetilde J$ 
are different pairs, then
$$
\dist(\C_{j,k},\C_{j',k'}) \geq \eta_0
\ \hbox{ or else $\ \C_{j,k}$ and $\C_{j',k'}$ have a common
endpoint $x\in V$.}
\leqno (2.6)
$$

We shall sometimes refer to this construction as the standard
decomposition of $K$. The precise choice of our decomposition
of the long $\C_j$ into two or three $\C_{j,k}$ does not
really matter, but we could make it unique by cutting $\C_j$
into as few equal pieces as possible.

\ms
Next we consider deformations of $K$, obtained as follows.
Let $\eta_1 < 10^{-1} \eta_0$ be given, and denote by
$\Phi(\eta_1)$ the set of functions $\varphi : V \to \partial B(0,1)$
such that
$$
|\varphi(x)-x| \leq \eta_1
\ \hbox{ for } x\in V.
\leqno (2.7)
$$

Let $\varphi \in \Phi(\eta_1)$ and $(j,k)\in \widetilde J$
be given. Denote by $a$ and $b$ the extremities of $\C_{j,k}$.
We denote by $\varphi_\ast(\C_{j,k})$ the geodesic of $\partial B(0,1)$
that goes from $\varphi(a)$ to $\varphi(b)$, where the uniqueness 
comes from (2.5) and (2.7). Finally set
$$
\varphi_\ast(K) = 
\bigcup_{(j,k) \in \widetilde J} \ \varphi_\ast(\C_{j,k})
\leqno (2.8)
$$
Thus we just deform $K$ as a net of geodesics, by moving
its vertices. Finally denote by $\varphi_\ast(X)$ the cone 
over $\varphi_\ast(K)$. 

Notice that by (2.5), (2.6), and (2.7), the 
geodesics $\varphi_\ast(\C_{j,k})$ all have lengths
larger than $9\eta_0$, and only meet at their 
vertices $\varphi(x), x\in V$, with angles 
that are still close to $2\pi/3$ when $x\in V_0$ and
to $\pi$ when $x\in V_1$. In particular,
$$
H^1(\varphi_\ast(K)) = 
\sum_{(j,k) \in \widetilde J}\length(\varphi_\ast(\C_{j,k})).
\leqno (2.9)
$$

\ms
\proclaim Definition 2.10.
Let $X$ be a (reduced) minimal cone centered at the origin.  
We say that $X$ is a full length minimal cone when there
is a standard decomposition of $K$ as above, an $\eta_1 < \eta_0/10$,
and a constant $C_1 \geq 1$, such that if $\varphi \in \Phi(\eta_1)$
is such that
$$
H^1(\varphi_\ast(K)) > H^1(K),
\leqno (2.11)
$$
then there is a deformation $\widetilde X$ of $\varphi_\ast(X)$ 
in $B(0,1)$ such that
$$
H^2(\widetilde X \cap B(0,1)) \leq H^2(\varphi_\ast(X)\cap B(0,1))
- C_1^{-1} [H^1(\varphi_\ast(K)) - H^1(K)].
\leqno (2.12)
$$

\ms
By deformation of $\varphi_\ast(X)$ in $B(0,1)$, we mean a 
set of the form $\widetilde X = f(\varphi_\ast(X))$, where
$f : \R^n \to \R^n$ is Lipschitz and such that 
$f(x)=x$ for $x\in \R^n\sm B(0,1)$ and $f(B(0,1)) \i B(0,1)$.

So we only require that when a small deformation $\varphi_\ast(K)$
is longer than $K$, it cannot be associated with a minimal cone,
and we require this with uniform ``elliptic" estimates.

Here we allowed ourselves to choose any standard decomposition of $K$,
but probably the required property does not depend on this choice.

Naturally, when $X$ is not centered at the origin, we say that it 
is a full length minimal cone when it is the translation of a 
full length minimal cone centered at the origin.


\ms
We shall see in Section 14 
that the minimal cones in $\R^3$ all have full length.
For instance, if $X$ is a plane, $K$ is a great circle,
and $H^1(\varphi_\ast(K)) > 2\pi$, it turns out that we can 
find a vertex $x \in V_1$ such that the angle of the two
geodesics $\varphi_\ast(\C_{j,k})$ that leave from $x$
is different from $\pi$ by at least 
$C^{-1} [H^1(\varphi_\ast(K)) - 2\pi]^{1/2}$, and then it 
is easy to find $\widetilde X$ as above.
The situation for cones of type $\Bbb Y$ or $\Bbb T$
is similar: the length excess of $\varphi_\ast(K)$
gives a lower bound on the wrongness of the angle at some vertex.

It is a priori possible that every minimal cone is full length,
but we really lack examples to discuss this seriously.
Notice however that if $X$ is a full length minimal cone, it cannot 
be embedded in a one-parameter family of minimal cones $X_t$ centered
at the origin such that
${\partial \over \partial t} \, H^2(X_t \cap B(0,1)) \neq 0$.
See Section 14 for a slightly longer discussion of these issues. 

\ms\noindent{\bf Remark 2.13.}
Some weaker version of Theorem 4.5 below will work with a
weaker (more degenerate) version of the full length property, where 
we replace $[H^1(\varphi_\ast(K)) - H^1(K)]$ in (2.12) with 
$[H^1(\varphi_\ast(K)) - H^1(K)]^N$. See 13.32.

\ms\noindent{\bf Remark 2.14.}
We shall not use the fact that our cones are minimal cones,
but only the description given at the beginning of this section
(that is, the description of $K$ as a union of great circles and
arcs of great circles, all the way up to (2.3). 
If $E$ is a cone that satisfies all these properties, 
we shall call it a \underbar{minimal-looking cone}.
This comes with an additional constant, the small constant 
$\eta_0$ in (2.2) and (2.3). Note that some minimal-looking 
cones are not minimal, for instance the cone over the union of
the edges of a cube in $\Bbb R^3$.

Essentially all the results proved in this paper have versions
for minimal-looking cones. See the end of Section 13 for a
slightly more detailed account.

\bigskip 
\noindent {\bf 3. Simple facts about density and limits.}
\medskip

Let $E$ be a reduced almost minimal set in $U \i \R^n$,
with gauge function $h$. We shall always assume that 
$$
h_1(r) = \int_0^{r} h(2t) {dt \over t} < +\infty
\ \hbox{ for $r$ small,}
\leqno (3.1)
$$
as in (1.16), in particular because (3.1) seems to be needed for 
the following result of almost-monotonicity for the density. 
There exist constants $\lambda \geq 1$ and $\varepsilon_0 > 0$, 
that depend only on the ambient dimension $n$, such that 
the following holds. Fix $x\in E$, and set 
$$
\theta(r) = r^{-2} H^2(E\cap B(x,r))
\leqno (3.2)
$$ 
as above. Let $R > 0$ be such that $B(x,R) \i U$, 
$h_1(R) < +\infty$, and $h(R) \leq \varepsilon_0$. Then
$$
\theta(r) \, e^{\lambda h_1(r)}
\hbox{ is a nondecreasing function of } r\in (0,R).
\leqno (3.3) 
$$
See for instance Proposition 5.24 in [D3]. 
Our main estimate will be an improvement of (3.3) when $h$ 
is small and $E$ is close to a full length minimal cone.
Before we come to that, let us mention some simple consequences 
of (3.3).

Notice that (3.1)-(3.3) imply the existence of the limit density
$$
d(x) = \lim_{r \to 0} \theta(r)
\leqno (3.4)
$$
for every $x\in E$. Then we can introduce the density excess
$$
f(r) = \theta(r) - d(x) = \theta(r)-\lim_{t\to 0^+} \theta(t)
\leqno (3.5)
$$
as in (1.3). Let us check that
$$
f(r) \leq f(s) + C h_1(s)
\ \hbox{ for } 0 < r < s < R
\leqno (3.6)
$$
when $B(x,2R) \i U$ and $h_1(2R) \leq \varepsilon_0/2$. 
Indeed, (3.1) yields $h(R) \leq \varepsilon_0$, so we can apply
(3.3) and get that
$$
\theta(r) \leq \theta(r) \, e^{\lambda h_1(r)}
\leq \theta(s) \, e^{\lambda h_1(s)}
\leq (1+Ch_1(s)) \, \theta(s)
\leq \theta(s) + C Ch_1(s)
\leqno (3.7)
$$
because $\theta(s) \leq C$ by the local Ahlfors-regularity
of $E$ (see Lemma 2.15 in [D3] in this context, 
but the property comes from [DS]); this is why  
we required that $B(x,2R) \i U$. Then (3.6) follows
by subtracting $d(x)$ from both sides.

When we let $r$ tend to $0$ in (3.7), we get that
$$
d(0) \leq \theta(s) + C h_1(s), 
\ \hbox{ or equivalently } \ f(s) \geq - C h_1(s),
\ \hbox{ for } \  0 < s < R.
\leqno (3.8)
$$

\ms
Let us also say a few words about blow-up limits.
For the definition of convergence, it is convenient to set
$$\eqalign{
d_{x,r}(E,F) &= r^{-1}\sup\big\{ \dist(y,E) \, ; \, 
y\in F \cap B(x,r) \big\}
\cr& \hskip 2cm
+ r^{-1}\sup\big\{ \dist(y,F) \, ; \, y\in E \cap B(x,r) \big\}
}\leqno (3.9) 
$$
when $E$ and $F$ are nonempty closed sets in a domain  $U$ that 
contains $B(x,r)$. By convention, 
$\sup\big\{ \dist(y,E) \, ; \, y\in F \cap B(x,r) \big\}=0$
when $F \cap B(x,r)$ is empty. 

\ms\proclaim Definition 3.10.
Let $\{ E_k \}$ be a sequence of subsets of $\R^n$,
and let $F$ be a closed set in the open set $U \i \R^n$.
We say that $\{ E_k \}$ converges to $E$ in $U$ when
$\lim_{k \to +\infty} d_{x,r}(E_k,F) = 0$ for every ball
$B(x,r)$ such that $\overline B(x,r) \i U$.

\ms
The definition would perhaps look more classical if we
replaced the balls $B(x,r)$ with compact subsets of $U$, but 
the result is the same because every compact subset of $U$
is contained in a finite union of balls $B(x,r)$ such that 
$\overline B(x,r) \i U$.

\ms\proclaim Definition 3.11.
A blow-up limit of $E$ at $x$ is a closed set $F$ that can be obtained
as limit in $\R^n$ of a sequence $\{ r_k^{-1}(E-x) \}$, where $\{ r_k \}$ 
is a sequence of radii such that $\lim_{k \to + \infty} r_k = 0$.

\ms
Since every ball in $\R^n$ is contained in a $B(0,\rho)$,
this means that $\lim_{k \to + \infty} d_{0,\rho}(E_k,F) = 0$
for every $\rho > 0$, where we set $E_k = r_k^{-1}(E-x)$, or equivalently 
that
$$
\lim_{k \to + \infty} d_{x, \rho r_k}(E,x+r_k F) = 0
\ \hbox{ for every $\rho > 0$.}
\leqno (3.12)
$$

The fact that the set $E_k$ naturally lives in $U_k = r_k^{-1}(U-x)$,
which may be smaller than $\R^n$, does not matter, because $U_k$
eventually contains every $\overline B(0,\rho)$.

Observe that $E$ always has at least one blow-up limit at $x$.
Indeed, by a standard compactness result on the Hausdorff distance
(see for instance Section 34 of [D2]), 
$$\eqalign{
&\hbox{every sequence $\{ r_k \}$ which tends to $0$ has a subsequence} 
\cr&\hskip 0.5cm
\hbox{for which the $r_k^{-1}(E-x)$ converge to some limit $X$.}
}\leqno (3.13)
$$
In addition, the blow-up limits are simpler, because of the following.

\ms\proclaim Proposition 3.14.
If $X$ is a blow-up limit of $E$ at $x$, then $X$
is a reduced minimal cone, with density
$$
H^2(X\cap B(0,1)) = \lim_{r \to 0} r^{-2} H^2(E \cap B(x,r)) = d(x).
\leqno (3.15)
$$

\ms
This is Proposition 7.31 in [D3], 
but let us rapidly say how it works.
The main point is a result of [D1], which says that 
the Hausdorff measure is lowersemicontinuous along sequences 
of quasiminimal sets. Then limits of reduced almost-minimal sets 
with a given gauge function $g$ are reduced almost-minimal sets 
with the same gauge function $g$.

Here we get the minimality of $X$ because $\{ r_k \}$ tends to $0$, 
so we can take $g(t) = h(rt)$ for arbitrarily small values of $r$.
Then 
$$
H^2(X\cap B(0,\rho)) = \lim_{r \to 0} r^{-2} H^2(E \cap B(x,r \rho)) 
= \rho^2 d(x)
\leqno (3.16)
$$
by (3.12) and because the Hausdorff measure goes to the limit
well; see Proposition 7.31 in [D3] 
again. So (3.15) holds, but also the density of $X$ is constant, 
and $X$ is a cone by Theorem~6.2 in [D3]. 
\qed

\ms
When $n=3$, there are only three types of (nonempty) reduced minimal 
cones: the planes, the sets of type $\Bbb Y$, and the cones over tetrahedron 
edges like $T$ in Figure 1.2. So Proposition 3.14 says that $d(x)$ can only take
the values $\pi$, $3\pi/2$, or $3 {\rm Argcos}(-1/3) \approx 1.82 \cdot \pi$,
depending on the common type of all the blow-up limits of $E$ at $x$.

When $n>3$, there are more possibilities. If $d(x)=\pi$, every blow-up 
limit of $E$ at $x$ is a plane and we say that $x$ is a $P$-point. 
If $d(x)=3\pi/2$, every blow-up limit is a set of type $\Bbb Y$, and
$x$ is called a $Y$-point. Otherwise, Lemma 14.12 in [D3] 
says that $d(x) \geq d_T$ for some $d_T > 3\pi/2$, and we say 
that $x$ is a $T$-point. We do not know much more in this case, 
there may even be many non mutually isometric blow-up limits of $E$ 
at $x$ (but they would all have the same density $d(x)$).
See [D3], Section 14 and the beginning of Section~16, 
for a little more detail.

\bigskip 
\noindent {\bf 4. The main comparison statement.}
\medskip

Here we want to state the main technical result of the
paper. Let $E$ be a reduced almost minimal set in $U \i \R^n$,
with a gauge function $h$ that satisfies (3.1), 
and let $x\in E$ and $r > 0$ be given. 

We shall make the following assumptions on $r$,
which will allow us to construct competitors and then
prove a differential inequality for $f$.
First we assume that for some small $\varepsilon > 0$,
$$
B(x,100r) \i U, 
\hskip0.2cm 
f(50r) = \theta(50r)-d(x) \leq \varepsilon,
\hskip0.2cm 
h_1(60r) = \int_0^{60r} h(2t) {dt \over t} \leq \varepsilon,
\leqno (4.1) 
$$
and there is minimal cone $X$ centered at the origin, such that
$$
d_{x,100r}(E,x+X) \leq \varepsilon. 
\leqno (4.2)
$$
We shall also use the standard decomposition of 
$K = X \cap \partial B(0,1)$ into arcs of great circles
$\C_{j,k}$, and the corresponding set of vertices $V$, as 
described in Section 2.

Our second assumption is a little more technical, but we shall 
see later that it is satisfied for almost every $r$.
We require that for every continuous nonnegative function $f$ on 
$\R^n$, 
$$\eqalign{
\lim_{\rho \to 0} \rho^{-1}
\int_{t\in (r-\rho,r)} \int_{E\cap \partial B(0,t)} f(z) \, dH^1(z) \, dt
= \int_{E\cap \partial B(0,r)} f(z) \, dH^1(z)
}\leqno (4.3)
$$
and that 
$$
\sup_{\rho > 0} \ \rho^{-1} \int_{t\in (r-\rho,r)}
H^1(E\cap \partial B(0,t)) \, dt < + \infty
\leqno (4.4)
$$
(the finiteness at $r$ of the corresponding maximal function).

\ms\proclaim Theorem 4.5. Let $\eta_1 < \eta_0/10$ be given.
If $\varepsilon > 0$ is small enough, depending 
on $\eta_1$, we have the following result. 
Let $E$ be a reduced almost minimal set in $U \i \R^n$,
with the gauge function $h$, and let $x\in E$ and $r > 0$ 
be such that (4.1)-(4.4) hold. Then there is a 
mapping $\varphi : V \to \partial B(0,1)$ such that (2.7) holds, 
$$
r H^1(\varphi_\ast(K)) 
\leq H^1(E\cap \partial B(x,r))
\leqno (4.6)
$$
and
$$\eqalign{
H^2(E \cap B(x,r))
\leq {r \over 2} &\, H^1(E\cap \partial B(x,r))
\cr&- 10^{-5} r \, [H^1(E\cap \partial B(x,r)) - r H^1(\varphi_\ast(K))]
+ 4 r^2 h(2r).
}\leqno (4.7)
$$
If in addition $X$ is a full length minimal cone with constants
$C_1$ and $\eta_1$, and such that $H^2(X\cap B(0,1)) \leq d(x)$, 
then 
$$\eqalign{
H^2(E \cap B(x,r))
\leq {r \over 2} &\, H^1(E\cap \partial B(x,r))
\cr&- \alpha r \, [H^1(E\cap \partial B(x,r)) - 2 d(x) r]
+ 4 r^2 h(2r),
}\leqno (4.8)
$$
where the small $\alpha > 0$ depends on $C_1$.

\ms
In (4.7) and (4.8), the term with square brackets corresponds
to an improvement over the estimate that leads to the monotonicity
property (3.3). More properties of $\varphi_\ast(K)$ and intermediate
objects will show up along the proof.

\ms
The proof of Theorem 4.5 will keep us busy for quite some time,
but let us already try to say how it will go. Of course
(4.7) and (4.8) will be obtained by comparing $E$ with a
competitor, what will be long is the construction of the 
competitor.

By homogeneity, we shall immediately reduce to the case when $x=0$ 
and $r=1$. We shall also assume that $H^1(E\cap \partial B(0,1))$
is not much larger than $2d(0)$, because otherwise (4.7) and (4.8)
hold more or less trivially.
First we shall find the points $\varphi(y)$, $y\in V$,
and at the same time a net of simple curves 
$g_{j,k} \i E \cap \partial B(0,1)$ that connect
them like the $\C_{j,k}$ do. That is, 
$g_{j,k}$ goes from $\varphi(y)$ to $\varphi(z)$,
where $y$ and $z$ denote the endpoints of $\C_{j,k}$.

For instance, if $X$ is a vertical cone of type $\Bbb Y$, 
the $g_{j,k}$ are six curves that connect two points near the
poles to three points near the equator.
The construction will be done in Section 6, 
and the proof will rely on separation properties that come
from the biH\"older description of $E$ that was obtained in [D3]. 

In Section 7 we shall replace the $g_{j,k}$ with 
Lipschitz graphs with small constants and the same endpoints,
the $\Gamma_{j,k}$. The construction is a simple maximal 
function argument, but we shall need to estimate 
$H^1(\Gamma_{j,k} \sm g_{j,k})$ in terms of how much length 
we win in the process. 

The reason why we prefer the $\Gamma_{j,k}$ is that it is easier
to modify the cone over a Lipschitz graph and diminish its area.
The idea is to approach the minimal surface equation with the
laplacian, and observe that a homogeneous function is rarely
harmonic. We shall improve the cone over each $\Gamma_{j,k}$
separately, and then glue the result to get a good competitor
for the cone over $\cup_{j,k} \Gamma_{j,k}$. This will save a 
substantial amount of area, unless the initial $\Gamma_{j,k}$
were already quite close to geodesics.

We do this modification in Section 8, 
and in Section 9 combine the three previous 
sections to get a competitor for $E$ and prove (4.7).

We also make sure that our competitor has a tip near
the origin that coincides with the cone over $\varphi_\ast(K)$,
so if $X$ is a full length minimal cone and 
$H^1(\varphi_\ast(K)) > 2 d(x)$, we can use Definition 2.10
to improve again our competitor near the origin, and prove (4.8).
So the proof of Theorem 4.5 will be completed at the end of
Section 10.

\bigskip
Let us now explain why some of our assumptions can be 
obtained easily. We start with (4.1) and (4.2).

\ms\proclaim Lemma 4.9.
Let $E$ be a reduced almost minimal set in $U \i \R^n$,
with a gauge function $h$ such that (3.1) holds, 
and let $x\in E$ be given.
For each $\varepsilon > 0$, we can find $r_0 > 0$ such that
for $0 < r < r_0$, (4.1) and (4.2) hold for some cone
$X$ which is a blow-up limit of $E$ at $x$.

\ms
Here $X$ is allowed to depend on $r$.
Notice that $H^2(X\cap B(0,1)) = d(x)$ (by (3.15)).
The full length property is not automatic, even though in
the context of Theorem 1.15, it will be enough to use this extra 
assumption once.

Let $x$ and $\varepsilon$ be given.
The fact that (4.1) holds for $r$ small follows directly
from (3.1) and (3.5). Suppose for a moment that
there are arbitrarily small radii $r$ such that (4.2) holds 
for no blow-up limit $X$. 
Let $\{ r_k \}$ be a sequence of such radii that tends to $0$.
By (3.13), we can replace $\{ r_k \}$ with a subsequence for which the 
$r_k^{-1}(E-x)$ converge to some limit $X$. Then $X$ is a blow-up 
limit of $E$ at $x$, and by (3.12) $d_{x,100r_k}(E,x+X)$ 
tend to $0$. This contradiction completes the proof of Lemma 4.9.
\qed

\ms
Let us also see why when $X$ is a plane or a cone with the 
density $d(x)$, we may even dispense with the density assumption 
in (4.1).

\ms\proclaim Lemma 4.10.
For each $\delta > 0$, we can find $\varepsilon >0$
such that if $E$ is a reduced almost minimal set in $U \i \R^n$,
$B(x,2r) \i U$, $h(3r) \leq \varepsilon$, and (4.2) holds for
some minimal cone $X$ centered at the origin and such that 
$H^2(X \cap B(0,1)) \leq d(x)$, then $f(r) \leq \delta$.

\ms
Indeed with these assumptions (verified with $\varepsilon/2$),
Lemma 16.43 in [D3], applied with $F = x+X$, implies that 
$$\eqalign{
H^2(E \cap B(x,r)) &\leq H^2(F \cap B(x,(1+\delta)r)) + \delta r^2
\cr& = H^2(X \cap B(0,(1+\delta)r)) + \delta r^2
\leq [(1+\delta)^2 d(x) + \delta] \, r^2,
}\leqno (4.11)
$$
hence $\theta(x,r) = r^{-2} H^2(E \cap B(x,r))
\leq d(x) + C \delta$ and $f(r) \leq C \delta$. 
The additional constant $C$ is unimportant, and Lemma~4.10 follows.
\qed

\ms
When (4.2) holds for some plane $X$, the condition that 
$H^2(X \cap B(0,1)) \leq d(x)$ is automatically satisfied,
because $\pi$ is the smallest possible density. 
This is not the case when $X$ is a cone of type $\Bbb Y$,
but at least, if $\varepsilon$ and $h_1(r)$ are small enough,
lemma 16.24 in [D3] gives a point of type $Y$ near $x$, 
to which we can apply Lemma 4.10. We don't know this when $X$
is of type $\Bbb T$.

We now turn to (4.3) and (4.4).

\ms\proclaim Lemma 4.12.
If $E$ is a reduced almost minimal set in $U \i \R^n$
and $B(x,r_0) \i U$, (4.3) and (4.4) hold for almost every 
$r \in (0,r_0)$.

\ms
We just need to know that $B(x,r_0) \i U$ to make sure
that $E$ is closed in $B(x,r_0)$ and $H^2(E \cap B(x,r_1)) < +\infty$
for $r_1 < r_0$. Now we just need to check that (4.3) and (4.4) hold 
for almost every $r \in (0,r_1)$, and take a countable union.

We want to apply the coarea formula, and for this we
need to know that $E$ is rectifiable. This is verified
in Section 2 of [D3] (see below (2.24) there). 
The reader should not worry that rectifiability 
is ``only" proved for generalized quasiminimal sets in [D3]; 
the point of introducing generalized quasiminimal sets was specifically 
to make sure that our almost-minimal sets would automatically be generalized 
quasiminimal sets. Compare Definitions~2.10 and 4.3 in [D3] 
(or read the few lines just above Lemma 4.7 there), and also see 
Proposition~4.10 there if you want to use Definition 4.8 in [D3]. 

So we are allowed to apply the coarea formula 
(Theorem~3.2.22 in [Fe]) 
on the rectifiable set $E\cap B(0,r_1)$. We apply it
to the $C^1$ mapping $g: z \to |z|$, and integrate against 
a bounded nonnegative Borel-measurable function $f$. 
The level sets are the $g^{-1}(t) = E\cap \partial B(0,t)$,
$0 < t < r_1$, and we get that
$$
\int_{E \cap B(0,r_1)} f(z) \, J_g(z) \, dH^2(z)
= \int_{t\in (0,r_1)} \bigg\{ \int_{E\cap \partial B(0,t)}
f(z) \, dH^1(z) \bigg\}\, dt,
\leqno (4.13)
$$
where $J_g$ is the appropriate $1$-jacobian (or gradient size) of $g$ on 
$E \cap B(0,r_1)$, and the fact that the inside integral 
$\int_{E\cap \partial B(0,t)} f(z) \, dH^1(z)$ is a measurable 
function of $t$ is part of the theorem.

For the sake of Lemma 4.12, we just need to know that
$J_g(z) \leq 1$, which is essentially obvious because $g$
is $1$-Lipschitz, but later on we shall need more information
that we give now.

Since $E$ is rectifiable, it has an approximate tangent plane 
$P(x)$ at $x$ for $H^2$-almost every $x\in E$. Incidentally,
we even know here that $P(x)$ is a true tangent plane
(but this won't really matter), because $E$ is locally Ahlfors-regular 
of dimension 2. See for instance Lemma 2.15 in [D3] for the 
local Ahlfors regularity and Exercise 41.21 in [D2] 
for the true tangent planes. We only need to compute $J_g(x)$ at points $x$ 
where $P(x)$ exists. Notice also that, essentially by definition, $J_g(x)$ is 
almost-everywhere the same as if $E$ were equal to $P(x)$. Thus 
$$
J_g(x)=\cos\alpha(x)
\ \hbox{ almost-everywhere on $E$},
\leqno (4.14)   
$$
where $\alpha(x) \in [0,\pi/2]$ denotes the angle of the radius 
$[0,x]$ with the tangent plane to $E$ at $x$ (that is, the smallest 
angle with a vector of the plane). When $n=3$, we could also define
$\alpha(x)$ by saying that $\pi/2 - \alpha(x)$ is the non oriented 
angle of $[0,x]$ with a unit normal to $P(x)$.

Return to the proof of Lemma 4.12. 
Set $\widetilde f(t) = \int_{E\cap \partial B(0,t)} f(z) \, dH^1(z)$
for $t\in (0,r_1)$. This is the inside integral in the 
right-hand side of (4.13), we already know that it is measurable,
and we also get that $\widetilde f$ is integrable, because (4.13) 
says that
$$
\int_{t\in (0,r)} \widetilde f(t) \, dt
= \int_{E \cap B(0,r_1)} f(z) \, J_h(z) \, dH^2(z)
\leq ||f||_{\infty} \, H^2(E \cap B(0,r_1)) < +\infty.
\leqno (4.15)
$$
Almost every $r \in (0,r_1)$ is a Lebesgue density point 
for $\widetilde f$, which means that
$$
\lim_{\rho \to 0} \ {1 \over 2 \rho} \, 
\int_{(r-\rho,r+\rho)} 
|\widetilde f(t) - \widetilde f(r)| \, dt = 0.
\leqno (4.16)
$$
Observe that (4.3) (for this single $f$) follows from
(4.16) and the triangle inequality, but we shall still
need some uniformity to control all the continuous functions at the 
same time. Let us take care of (4.4) in the mean time.

Take $f=1$, and denote by $b$ the corresponding function
$\widetilde f$. Thus $b(t) = H^1(E\cap \partial B(0,t))$.
As before, $b$ is integrable on $(0,r_1)$. Denote by $b^\ast$ its
(non centered) Hardy-Littlewood maximal function; thus
$$
b^\ast(r) = \sup \big\{ {1 \over |I|} \, \int_I b(t) dt \, ; \,
I \i (0,r_1) \hbox{ is an interval  that contains } r \big\}.
\leqno (4.17)
$$
for $r \in (0,r_1)$.
The Hardy-Littlewood maximal theorem (on page 5 of [St]) 
says that $b^\ast(r) < +\infty$ for almost-every $r \in (0,r_0)$.
Obviously (4.4) holds when $b^\ast(r) < +\infty$, so we
can safely return to (4.3).

Now let $C(\overline B(0,r_1))$ denote the set of continuous 
functions on $\overline B(0,r_1)$, and select a countable family 
$\{ f_i \}$, $i\in \Bbb N$, in $C(\overline B(0,r_1))$, which is 
dense in the $sup$ norm. We know that for almost every $r\in (0,r_1)$,
(4.17) holds, and so does (4.16) for every $\widetilde f_i$.

Fix $r\in (0,r_1)$ like  this.
We know that (4.4) holds for  $r$, and we just need to check
that (4.3), or even (4.16), holds for any given $f$. Set
$a(\rho) = \int_{(r-\rho,r+\rho)} 
|\widetilde f(t) - \widetilde f(r)| \, dt$
for $\rho > 0$ small (so that $(r-\rho,r+\rho) \i (0,r_1)$).
We just need to show that $\lim_{\rho \to 0} \rho^{-1} a(\rho) = 0$.

Let $\varepsilon > 0$ be given. Let $i$ be such that 
$||f-f_i||_\infty \leq \varepsilon$. Then 
$$\leqalignno{
a(\rho) &\leq \int_{(r-\rho,r+\rho)}
\Big[ \, |\widetilde f(t) - \widetilde f_i(t)| + 
|\widetilde f_i(t) - \widetilde f_i(r)| + 
|\widetilde f_i(r) - \widetilde f(r)| \Big] \, dt
\cr&
\leq \int_{(r-\rho,r+\rho)}
\big[ |\widetilde f_i(t) - \widetilde f_i(r)| + 2 ||f-f_i||_\infty
\, H^1(E\cap \partial B(0,t)) \big] \, dt 
& (4.18)
\cr&
\leq 2\varepsilon \rho + 2 \varepsilon \int_{(r-\rho,r+\rho)} 
H^1(E\cap \partial B(0,t)) \, dt
}
$$
for $\rho$ small enough, by (4.16) for $\widetilde f_i$
and because $||f-f_i||_\infty \leq \varepsilon$.
Since
$$
\int_{(r-\rho,r+\rho)} 
H^1(E\cap \partial B(0,r)) \, dt
= \int_{(r-\rho,r+\rho)} b(t)
\leq 2 \rho b^\ast(r)
\leqno (4.19)
$$
by definition of $b$ and $b^\ast$, (4.18) shows that
$a(\rho) \leq [2+ 4 b^\ast(r)] \, \varepsilon \rho$ 
for $\rho$ small enough. Thus $a(\rho)/\rho$ tends to $0$
and (4.16) holds for $f$, as needed.

This completes our proof of Lemma 4.12.
\qed

\bigskip 
\noindent {\bf 5. Differential inequalities and decay for 
the density excess $f$}
\medskip

In this section we show how the estimate (4.8) in the conclusion
of Theorem 4.5 can be used to derive a differential inequality
on $f$, which itself leads to a decay estimate for $f$.
Once we have all this, we shall use $f$ to control the geometry
of $E$, but this will not happen before Sections 11 and 12.

Here and in the rest of the paper, $E$ is a reduced two-dimensional 
almost-minimal set in the open set $U \i \R^n$, with a gauge function 
$h$ such that (3.1) holds. 

Let us first discuss how to recover 
functions like $\theta(r)$ from their derivative almost-everywhere. 
From now on, we shall work with $x=0$ to simplify our notation.

\medskip \proclaim Lemma 5.1. Suppose $0\in U$, and set
$D = \dist(0,\R^n \setminus U)$ and $v(r) = H^2(E \cap B(0,r))$
for $0<r<D$. Let $\mu$ denote the image by the radial projection 
$z \to |z|$ of the restriction of $H^2$ to $E$; thus $\mu$
is a positive Radon measure on $[0,D)$, and
$v(r) = \mu\big([0,r)\big) =
\int_{[0,r)} d\mu$ for $0<r<D$. Let $b(r)$ be
a positive $C^1$ function on $(0,D)$. Then
$$
b(y)v(y) - b(x)v(x) = \int_{[x,y)} b(r) \, d\mu(r) 
+ \int_{[x,y)} b'(r) v(r) \, dr
\ \hbox{ for } 0 < x < y < D.
\leqno (5.2)
$$
In addition, $v$ is differentiable almost-everywhere and
$$
b(y)v(y) - b(x)v(x) \geq \int_x^y [b(r) v'(r) + b'(r) v(r)] \, dr
\ \hbox{ for } 0 < x < y < D.
\leqno (5.3)
$$

\medskip
In order to prove (5.2), compute
$I = \int_{x \leq s \leq t < y} b'(s) ds d\mu(t)$ in two different ways.
If we integrate in $s$ first, we get that
$I = \int_{[x,y)} [b(t)-b(x)] d\mu(t)$, while integrating in $t$ first
yields $I = \int_{[x,y)} b'(s) (v(y)-v(s)) ds$. Then
$$\eqalign{
\int_{[x,y)} b(r) d\mu(r) + &\int_{[x,y)} b'(r) v(r) dr
= I + b(x) \int_{[x,y)} d\mu(t) -I + v(y) \int_{[x,y)} b'(s) ds
\cr & 
= b(x) [v(y)-v(x)] + v(y)[b(y)-b(x)] = v(y)b(y)-v(x)b(x),
}\leqno (5.4)
$$
so (5.2) holds. The fact that $v'(r)$ exists for almost every $r<D$ is
a standard fact about nondecreasing functions, and so is the fact that 
the absolutely continuous part of $d\mu$ is $v'(r) dr$. In particular, 
$v'(r)dr \leq d\mu$ and (5.3) follows from (5.2) because $b$ is nonnegative.
Lemma 5.1 follows. \qed

\medskip
Let us immediately record what this means in terms of
$\theta(r) = r^{-2} H^2(E\cap B(0,r)) = r^{-2} v(r)$.

\medskip \proclaim Lemma 5.5. Still assume that $0\in U$ and
set $D = \dist(0,\R^3 \setminus U)$. The function $\theta$ is
differentiable almost-everywhere on $(0,D)$, with
$$
\theta'(r) = r^{-2} v'(r) - 2 r^{-3} v(r) = r^{-2} v'(r) - 2 r^{-1} 
\theta(r) \ \hbox{ almost-everywhere,}
\leqno (5.6)
$$
and 
$$
\theta(y) - \theta(x) \geq \int_x^y \theta'(r)\, dr
\ \hbox{ for } 0 < x < y < D.
\leqno (5.7)
$$

\medskip
Indeed it is easy to see that if $v$ is differentiable at $r$, then
$\theta$ is differentiable at $r$ too, with the formula in (5.6).
Then (5.7) is just the same as (5.3). \qed

\medskip
Let us also check that (with the notation of Lemma 5.1)
$$
v'(r) \geq  H^1(E \cap \partial B(0,r)) 
\ \hbox{ for almost every } r \in (0,D).
\leqno (5.8)
$$

Pick any $r_1 \in (0,D)$ and apply the co-area formula
(4.13). [We are in the same situation as in Lemma 4.12,
where the assumption was that $r_0 \leq D$, except for 
the fact that now $x=0$.] 
Take for $f$ the characteristic function
of $A = B(0,r)\setminus B(0,r-\rho)$, with
$0 < \rho \leq r < r_1$. We get that 
$$\eqalign{
\int_{E \cap A}\, J_g(z) \, dH^2(z)
&= \int_{t\in (0,r_1)} \bigg\{ \int_{E\cap \partial B(0,t)}
f(z) \, dH^1(z) \bigg\}\, dt
\cr&
= \int_{t\in [r-\rho,r)}  H^1(E\cap \partial B(0,t)) \, dt
}\leqno (5.9)
$$
and, since $J_g(z) \leq 1$ in (4.13), we get that
$$
v(r)-v(r-\rho) = H^2(E \cap A) 
\geq \int_{t\in [r-\rho,r)}  H^1(E\cap \partial B(0,t)) \, dt.
\leqno (5.10)
$$
Now fix $r$, divide (5.10) by $\rho$, and let $\rho$ tend
to $0$. The left-hand side tends to $v'(r)$ if $v'(r)$ 
exist. The right-hand side tends to $H^1(E \cap \partial B(0,t))$
when $r$ is a Lebesgue point for the function
$t \to H^1(E \cap \partial B(0,t))$. 
This function is locally integrable by (5.9) (recall that
its measurability is part of the co-area formula), 
so the Lebesgue density theorem says that almost every $r\in (0,r_1)$
is a Lebesgue point, and (5.10) yields (5.8).

\ms
Let us now see how to deduce a differential inequality and 
a decay estimate from the conclusion of Theorem 4.5.

\ms\proclaim Lemma 5.11.
Let $0 < x < y < D=\dist(0,\R^3 \setminus U)$ be given, and
suppose that for some $\alpha \in (0,1/3)$, (4.8) holds
for  almost every $r \in (x,y)$. Then $f(r) = \theta(r) - d(0)$ 
is differentiable almost-everywhere on $(x,y)$, with
$$
r f'(r) \geq {4\alpha \over 1-2\alpha} \, f(r) - 24 h(2r) 
\ \hbox{ for almost every } r\in (x,y),
\leqno (5.12)
$$
and then
$$
f(x) \leq (x/y)^a f(y) + C x^a \int_x^y r^{-a-1} h(2r) dr,
\leqno (5.13) 
$$
where we set $a={4\alpha \over 1-2\alpha} \, $.

\ms
Indeed, we already know from Lemma 5.5 that $f$ is differentiable 
almost-everywhere on $(0,D)$, and (5.6) says that
$$
f'(r)= \theta'(r) = r^{-2} v'(r) - 2 r^{-1} \theta(r)
\geq r^{-2} H^1(E \cap \partial B(0,r)) - 2 r^{-1} \theta(r)
\leqno (5.14)
$$
almost-everywhere, by (5.8). 

Set $\dsp X = {1 \over 2r} H^1(E\cap \partial B(0,r))$.
Thus $rf'(r) \geq 2X-2\theta(r)$.
But (4.8) (divided by $r^2$) says that 
$(1-2\alpha) X \geq \theta(r) - 2\alpha d(0) -4 h(2r)$.
Hence
$$\eqalign{
rf'(r) &\geq 2X-2\theta(r) 
\geq {2 \over 1-2\alpha}\, \big[\theta(r) - 2\alpha d(0) -4h(2r)\big] 
- 2\theta(r) 
\cr&
\geq {4 \alpha \over 1-2\alpha}  \, \big[\theta(r) - d(0) \big] 
- {8 h(2r) \over 1-2\alpha}
= {4 \alpha f(r) \over 1-2\alpha}  - {8 h(2r) \over 1-2\alpha} \, ,
}\leqno (5.15)
$$
which yields (5.12) because $\alpha <1/3$.

\ms
Now let us integrate (5.12) to prove (5.13).
Set $g(r)= r^{-a} f(r)$; then 
$$
g'(r) = r^{-a} f'(r) - a r^{-a-1} f(r) = r^{-a-1} [r f'(r) - a f(r)]
\geq - 24 r^{-a-1} h(2r)
\leqno (5.16)
$$
almost-everywhere on $(x,y)$, by (5.12). We claim that
$$
g(y) - g(x) \geq \int_x^y g'(r) dr 
\geq - 24 \int_x^y r^{-a-1} h(2r) dr.
\leqno (5.17)
$$
By (5.16), we just need to worry about the first inequality.
Notice that $g(r) = r^{-a} f(r) = r^{-a} \theta(r) - r^{-a} d(0)$.
The analogue for $-r^{-a} d(0)$ of the first inequality in (5.17)
is trivial, and for $r^{-a} \theta(r) = r^{-a-2}v(y)$, it is just 
equivalent to (5.3) with $b(r)=r^{-a-2}$. So the first part of (5.17)
follows from Lemma 5.1, and (5.17) holds.
Now
$$
f(x) = x^a g(x) \leq x^a g(y) + 24 x^a \int_x^y r^{-a-1} h(2r) dr
\leqno (5.18)
$$
by (5.17); (5.13) and Lemma 5.10 follow. \qed

\ms
The next two examples are here to explain why 
(5.13) is really a decay estimate. The computation for
the first one will be used in Section 13; Example 5.21 
is just here to give an idea of what happens with weaker
decay conditions. 

\ms
\noindent {\bf Example 5.19.} 
Let us check that 
$$
f(x) \leq x^a y^{-a}f(y) + C' x^b
\leqno (5.20)
$$
when $x$ and $y>x$ are as in Lemma 5.11 and
$h(r) \leq C r^b$ for some $b \in (0,a)$.
(Even if $h$ is much smaller, we shall not get better 
bounds than $f(x) \leq C x^a$ anyway.)
Indeed $\int_x^y r^{-a-1} h(2r) dr \leq C 2^b \int_x^y r^{b-a-1} dr
= C 2^b (a-b)^{-1} (x^{-a+b}-y^{-a+b}) \leq C 2^b (a-b)^{-1} x^{-a+b}$, 
hence (5.20) follows from (5.13).

\ms
\noindent {\bf Example 5.21.}
Suppose that $h(r) \leq C [\log(A/r)]^{-b}$ for some constants
$A, b > 0$. Then
$$
f(x) \leq (x/y)^{a} f(y) + C (x/A)^{a/2} 
+ C \big[\log\big({A \over 2x}\big)\big]^{-b}  
\leq C_{A,y} \, \big[\log\big({A \over 2x}\big)\big]^{-b} 
\leqno (5.22)
$$
when $x<y<A/3$ satisfy the assumptions of Lemma 5.11.
[However, we should remember that we have to assume (3.1)
to prove (5.13), so we cannot really take $h(r) = C [\log(A/r)]^{-b}$
with $b \leq 1$.]

Indeed, (5.13) yields $f(x) \leq (x/y)^{a} f(y) 
+ C x^a \int_x^y r^{-a-1} [\log({A\over 2r})]^{-b} dr$.
We cut the domain of integration $(x,y)$ into two region, $I_1$ where 
$\log({A\over 2r}) \geq {1 \over 2} \, \log({A\over 2x})$ and 
$I_2 = (x,y)\sm I_1$. Then
$$
\int_{I_1} r^{-a-1} \big[\log\big({A\over 2r}\big)\big]^{-b} dr
\leq 
\big[{1 \over 2} \, \log\big({A\over 2x}\big)\big]^{-b} \int_{I_1} r^{-a-1} dr
\leq C x^{-a} \big[\log\big({A\over 2x}\big)\big]^{-b}.
\leqno (5.23)
$$ 
On $I_2$ we have that
${A\over 2r} \leq \big({A\over 2x}\big)^{1/2}$, hence
$r \geq (Ax/2)^{1/2}$. Thus
$$\eqalign{
\int_{I_2} r^{-a-1} \big[\log\big({A\over 2r}\big)\big]^{-b} dr
&\leq \big[\log\big({A\over 2y}\big)\big]^{-b}
\int_{I_2} r^{-a-1} dr
\cr&\leq C \int_{I_2} r^{-a-1} dr
\leq C (Ax/2)^{-a/2},
}\leqno (5.24)
$$
We sum (5.23) and (5.24), multiply by $Cx^a$, add $(x/y)^{a} f(y)$
and get (5.22).

\ms
Even if $h$ goes very slowly to $0$, (5.13) still gives a definite 
decay, with the same sort of argument as above. For instance, 
we can pick any $z\in [x,y]$, and observe that
$$\eqalign{
\int_x^y r^{-a-1} h(2r) dr
&\leq \int_x^z r^{-a-1} h(2z) dr + \int_z^y r^{-a-1} h(2y) dr
\cr&\leq a^{-1} x^{-a} h(2z) + a^{-1} z^{-a} h(2y), 
}\leqno (5.25)
$$
which leads to
$$
f(x) \leq x^a y^{-a}f(y) + C h(2z) + C x^a z^{-a} h(2y).
\leqno (5.26)
$$
One can then optimize the choice of $z$ (depending on $x$), but anyway 
it is  easy to force the right-hand side of (5.26) to tend to $0$.  
But again, remember that we use (3.1) to prove (5.13).

It is amusing that (5.13) leads to some definite decay for $f$
in all cases, even when (3.1) fails. But this probably does
not mean that we could use a variant of the more delicate method 
here to get the results of [D3] under weaker assumptions on $h$.

\ms 
We end the section with a discussion of a weaker form of 
Lemma 5.11 when (4.8) is replaced with the weaker inequality (5.27).
This part will only needed in Section 13. 

In the good cases, which include all the cases when $n = 3$
and $y$ is small enough, we shall be able to prove (4.8).
In some other cases, for instance if the tangent cones to $E$ only
satisfy a weaker full length condition, 
we may be a little less lucky but still get that
$$
\theta(r) \leq {1 \over 2r} H^1(E\cap \partial B(0,r)) -
\alpha \Big[ {1 \over 2r} H^1(E\cap \partial B(0,r)) - d(0) \Big]_+^{N} 
+ C h(2r)
\leqno (5.27)
$$
for some $N > 1$. Here the $+$ means that we take the positive
part (so $A_+^N = 0$ for $A \leq 0$).
Let us record the estimates that this would yield; we shall
discuss the weaker full length conditions in Section 13.

\ms\proclaim Lemma 5.28.
Let $0 < x < y < D=\dist(0,\R^3 \setminus U)$ be given, and
suppose that (5.27) holds for some $\alpha \in (0,1/2)$, some
$N > 1$, and almost-every $r \in (x,y)$.
Then $f$ is differentiable almost-everywhere on $(x,y)$, and
$$
r f'(r) \geq 2\alpha \, f(r)_+^N - C h(2r) 
\ \hbox{ for almost every } r\in (x,y).
\leqno (5.29)
$$

\ms
As before, the differentiability of $f$ almost everywhere comes from
Lemma 5.5. Set $X = {1 \over 2r} H^1(E\cap \partial B(0,r))$ as above;
thus (5.14) still says that $rf'(r) \geq 2X-2\theta(r)$, and (5.27) 
says that
$$
\theta(r) - X \leq -\alpha \big[X-d(0)\big]_+^N + C h(2r)
\leqno (5.30) 
$$
We first prove (5.29) when $\theta(r) \leq X$. Then
$$\eqalign{
rf'(r) &\geq 2X-2\theta(r) \geq 2\alpha \big[ X - d(0) \big]_+^{N} - 2C h(2r)
\cr&\geq 2\alpha \big[ \theta(r) - d(0) \big]_+^{N} - 2C h(2r)
= 2\alpha f(r)_+^{N} - 2C h(2r)
}\leqno (5.31) 
$$
by definition of $f(r)$.  So (5.29) holds in this case. 

If $\theta(r) > X$, (5.30) says that $\theta(r) - X \leq Ch(2r)$.
Then $\big[X-d(0)\big]_+^N \geq \big[\theta(r)-d(0)\big]_+^N 
-Ch(2r)$ and the proof of (5.31) yields 
$rf'(r) \geq 2\alpha f(r)_+^{N} - C h(2r)$ again. This proves (5.29) 
and Lemma 5.28.
\qed

\ms\noindent{\bf Remark 5.32.}
As before, (5.29) and Lemma 5.4 yield a decay estimate,
but which is  not nearly as good as in the situation of Lemma 5.11.
For instance, we claim that if (5.27) holds for almost-every $r \in (0,y)$
for some choice of $y < D$, $\alpha \in (0,1/2)$, and 
$N>1$,  and if 
$$
h(r) \leq C [{\rm Log}({1 \over r})]^{-{N \over N-1}}
\ \hbox{ for $r$ small, }
\leqno (5.33) 
$$
then
$$
f(x) \leq C_1 \Big[{\rm Log}\Big({2y \over x}\Big)\Big]^{-{1 \over N-1}}
\ \hbox{ for }  0  < x < y,
\leqno (5.34) 
$$
for some  $C_1$ that also depends on $y$ and $f(y)$.
Set 
$$
\varphi(r) = C_1 \big[{\rm Log}\big({2y \over r}\big)
\big]^{-{1 \over N-1}} \ \hbox{ for } r < y. 
\leqno (5.35) 
$$
We choose $\varphi$ of this form because
$$
r\varphi'(r) = {C_1 \over N-1} \big[{\rm Log}\big({2y \over r}\big)
\big]^{-{N \over N-1}} = {C_1^{1-N} \over N-1}\, \varphi(x)^N,
\leqno (5.36) 
$$
which is similar to (5.29).

If $C_1$ is large enough, $\varphi(y) > f(y)$. Then suppose that
$f(r) \geq \varphi(r)$ for some $r<y$, and denote by $x$ the
supremum of such $r$. Notice that $f(x) \geq \varphi(x)$
because $\varphi$ is continuous, $f(r) = \theta(r)-d(0)$, and 
$r^2\theta(r)$ is nondecreasing. Hence $x<y$.  
For $\delta > 0$ small,
$$\eqalign{
\theta(x+\delta) - \theta(x) &\geq \int_x^{x+\delta} \theta'(r) \, dr
\geq \int_x^{x+\delta} \big[ 2\alpha \, f(r)_+^N - C h(2r) 
\big]\, {dr \over r}
\cr& \geq \int_x^{x+\delta} \Big\{ 2\alpha \, f(r)_+^N 
- C  [{\rm Log}({2y \over r})]^{-{N \over N-1}} \Big\}\, {dr \over r}
\cr&\geq 2\alpha\delta  x^{-1} f(x_+)_+^N 
- C \delta x^{-1} [{\rm Log}({2y \over x})]^{-{N \over N-1}}
+ o(\delta)
}\leqno (5.37) 
$$
by Lemma 5.5, because $\theta' =f'$, by (5.29), and by (5.33). 
Here again, $C$ may depend on $y$. At the same time,
$$
\theta(x+\delta) - \theta(x) = f(x+\delta)-f(x)
\leq \varphi(x+\delta)- \varphi(x)
= \delta \varphi'(x) + o(\delta)
\leqno (5.38) 
$$
because $f(x+\delta) < \varphi(x+\delta)$ and
$f(x) \geq \varphi(x)$, which yields
$$
2\alpha f(x_+)_+^N  - C [{\rm Log}({2y \over x})]^{-{N \over N-1}}
\leq x\varphi'(x) =  {C_1 \over N-1} \big[{\rm Log}\big({2y \over x}\big)
\big]^{-{N \over N-1}}
\leqno (5.39) 
$$
by (5.36).
If $C_1$ is larger than $2(N-1)C$, this yields 
$2\alpha f(x_+)_+^N \leq {C_1 \over 2(N-1)} 
\big[{\rm Log}\big({2y \over x}\big) \big]^{-{N \over N-1}}$.
In addition, $\varphi(x) \leq f(x) \leq f(x_+)$ by 
definition of $x$ and because $\theta(r)$ is nondecreasing
and the other functions are continuous. Thus
$2 \alpha \varphi(x)^N \leq {C_1 \over 2(N-1)} 
\big[{\rm Log}\big({2y \over x}\big) \big]^{-{N \over N-1}}$.
This contradicts (5.35) if $C_1$ is large enough. So there
is no $x$ as above, $f(r) < \varphi(r)$ for $r<y$, and (5.34) 
follows.

Even if we take $h$ much smaller than in (5.33), the differential inequality 
in (5.29) will not give a much better decay than the one in (5.34).
For instance, $\varphi$ in (5.35) solves (5.29) with $h=0$
when $C_1$ is small , and it does not decay much near $0$.

And unfortunately (5.34) does not give enough control to allow us
to prove that $E$ is locally $C^1$-equivalent to a minimal cone.
See Section 13 for a rapid discussion of the weaker full length 
conditions.

\bigskip 
\noindent {\bf 6. Separation properties and the construction of a net of curves}
\medskip

The next few sections will be devoted to the proof
of Theorem 4.5. So we fix a reduced almost minimal set $E$
in the open set $U \i \R^n$, with gauge function $h$, 
and we let $x\in E$ and $r>0$ be such that (4.1) and (4.2)
hold. [We shall not use (4.3) and (4.4) before Section~9.] 
 
Since our statement is invariant under translations and dilations,
we can assume that $x=0$ and $r=1$, which will simplify the notation.

Our ultimate goal is to construct a nice competitor for $E$
in $B(0,1)$, which in particular will be better than the cone
over $E \cap \partial B(0,1)$ unless $E \cap \partial B(0,1)$
is already very nice.  In this section we concentrate on finding
a net of curves $g_{j,k} \i E \cap \partial B(0,1)$; see
Lemma 6.11. 

We already have a constant $\varepsilon > 0$ in play
(in (4.2)), and our construction will use another small parameter 
$\tau>0$, also to  be chosen later. Although we shall not always repeat 
this, all the estimates below hold only when $\tau$ is small enough, and then
$\varepsilon$ is small enough, depending on $\tau$. 

We set $B = B(0,1)$ and $\partial B = \partial B(0,1)$ to save
notation. Observe before we start that
$$
d(0)-C\varepsilon \leq \theta(r) \leq d(0)+C\varepsilon 
\ \hbox{ for } 0  < r \leq 50,
\leqno (6.1)
$$
by (4.1), (3.6), and (3.8).
The next lemma will allow us to restrict our attention
to the case when
$$
H^1(E \cap \partial B) \leq 2(1+\tau)\, d(0).
\leqno (6.2)
$$

\ms
\proclaim Lemma 6.3. If (6.2) fails, the conclusions of
Theorem 4.5 hold, and in particular
(4.8) holds for every $\alpha \in (0,1/3)$.

\ms
We start with (4.8), which is the most important.
If (6.2) fails,
$$\eqalign{
H^2(E\cap B(x,r)) - {r \over 2}  H^1&(E\cap \partial B(x,r)) +
\alpha r \Big[ H^1(E\cap \partial B(x,r)) - 2d(x)r \Big] 
\cr& 
= \theta(1) - {1-2\alpha \over 2} \, H^1(E\cap \partial B) 
- 2\alpha d(0)
\cr&
\leq \big[ d(0) + C\varepsilon \big] 
- (1-2\alpha) (1+\tau)\, d(0) - 2\alpha d(0)
\cr& = C\varepsilon - (1-2\alpha) \, \tau \, d(0) < 0
}\leqno (6.4)
$$
because $x=0$ and $r=1$, by (6.1), and if $\varepsilon$ is 
small enough; (4.8) follows.

\ms
Let us also check (4.6) and (4.7).
We take $\varphi(x) = x$ for $x\in V$. Recall that
for Theorem 4.5, $K = X \cap \partial B$, where $X$ still
denotes the minimal cone in (4.2), so we get that
$\varphi_\ast(K) = K = X \cap \partial B$. 
We need to  control $D = H^2(X\cap B(0,1))$, 
so we apply Lemma~16.43 in [D3], 
with $E = X$, $F=E$, $x=0$, $r=1$, and $\delta = \tau/3$.
The assumptions are satisfied by (4.1) and (4.2),
and we get that
$$\eqalign{
D &= H^2(X\cap B(0,1)) \leq H^2(E\cap B(0,1+\delta))+ \delta
\cr&
= (1+\delta)^2 \theta(0,1+\delta) + \delta
\leq (1+ 3 \delta) d(0)  = (1+ \tau) \, d(0) 
}\leqno (6.5)
$$
by (6.1), because $d(0) \geq \pi$, and if $\varepsilon$ is small enough.
The same  lemma, with $E$ and $X$ exchanged, also yields
$$
\theta(1) = H^2(E\cap B(0,1))
\leq H^2(X\cap B(0,1+\delta))+ \delta
= (1+\delta)^2 D + \delta
\leqno (6.6)
$$
hence 
$$
D \geq (1+\delta)^{-2} \theta(1) - \delta 
\geq (1+\delta)^{-2} [d(0)-C\varepsilon] - \delta
\geq (1- 3 \delta) \, d(0) = (1- \tau) \, d(0),
\leqno (6.7)
$$
again by (6.1) and because $d(0) \geq \pi$. Observe that (6.5) and 
(6.7) do not rely on the failure of (6.2).
Then (4.6) holds because 
$$
H^1(\varphi_\ast(K)) = H^1(K) = 2D  
\leq  2(1+\tau)\, d(0) \leq H^1(E \cap \partial B)
\leqno (6.8)
$$
because $X$ is the cone over $K$,
by (6.5) and because (6.2) fails. Next (4.7) holds because
$$\leqalignno{
H^2(E\cap B(x,r)) - {r\over 2}  &H^1(E\cap \partial B(x,r)) +
10^{-5}  r \Big[ H^1(E\cap \partial B(x,r)) 
- r H^1(\varphi_\ast(K)) \Big] 
\cr&
= \theta(1) - { 1 - 2 \cdot 10^{-5} \over 2} \, H^1(E\cap \partial B)
 - 2 \cdot  10^{-5} D  
\cr&
\leq \big[ d(0) + C\varepsilon \big] 
- (1 - 2 \cdot 10^{-5})(1+\tau) \, d(0) - 2 \cdot  10^{-5} D 
& (6.9)
\cr&
= C\varepsilon - \tau d(0) + 2 \cdot 10^{-5}(1+\tau) \, d(0)
- 2 \cdot  10^{-5} D
\cr&
\leq C\varepsilon - \tau d(0) + 2 \cdot 10^{-5} d(0)
[(1+\tau) - (1- \tau)] < 0
}
$$
(as for (6.4)) because $x=0$ and $r=1$, by (6.8) and (6.1), 
because (6.2) fails, and by (6.7). This  completes the proof of Lemma 6.3.
\qed

\ms
So we shall now assume that (6.2) holds.
Recall from (4.2) that we have a minimal cone $X$
centered at the origin, such that
$$
d_{0,100}(E,X) \leq \varepsilon.
\leqno (6.10)
$$
We shall again use a standard decomposition of 
$K = X \cap \partial B$ into geodesic arcs $\C_{j,k}$, 
$(j,k) \in \widetilde J$, as described in Section 2. 
We still denote by $V = V_0 \cup V_1$ the collection
of endpoints of the $\C_{j,k}$, with $V_0$ corresponding to
the original vertices (those who belong to three $\C_{j,k}$)
and $V_1$ to the added vertices (where only two $\C_{j,k}$
end, and make $180^\circ$ angles). We want to draw curves 
$g_{j,k}$ in $E \cap \partial B$ with the same structure, as in 
the following lemma.

\medskip
\proclaim Lemma 6.11. 
We can find points $\varphi(y)$, $y\in V$, and simple arcs
$g_{j,k}$ in $E \cap \partial B$, so that 
$$
\hbox{if $y$ and $z$ denote the endpoints of $\C_{j,k}$,
the endpoints of $g_{j,k}$ are $\varphi(y)$ and $\varphi(z)$,}
\leqno (6.12)
$$
$$
\hbox{the arcs $g_{j,k}$ are disjoint, except for their endpoints,}
\leqno (6.13)
$$
$$
|\varphi(y)-y| \leq C \tau 
\ \hbox{ for } y \in V,
\leqno (6.14)
$$
and
$$
\dist(z,\C_{j,k}) \leq C \tau 
\hbox{ for } z\in g_{j,k} \, .
\leqno (6.15)
$$

\ms
The proof will rely on a local description of $E$ as a biH\"older 
image of a minimal cone, which we take from [D3], 
plus a little bit of topology. 

Let us apply Lemma 16.19 in [D3] if $0$ is a point of type $P$, 
Lemma 16.25 in [D3] 
if $0$ is a $Y$-point, and Lemma 16.56 in [D3] 
if it is a $T$-point. Each time we take for $\tau$ the same
constant as in (6.2). [There is nothing subtle here, we are just saving 
some notation, and anyway we shall take $\tau$ small.] 
By (4.1), the lemma applies if $\varepsilon$ is small enough.
We get that $B$ is a biH\"older ball for $E$,
of the same type as $0$, and with the constant $\tau$.
By Definition 15.10 in [D3], 
this means that there is a reduced minimal cone $X$ centered 
at the origin, and a biH\"older mapping $f : B(0,2) \to \R^n$,
with the following properties:
$$
|f(y)-y| \leq \tau \ \hbox{ for } y\in B(0,2),
\leqno (6.16) 
$$
$$
(1-\tau) \, |y-z|^{1+\tau} \leq |f(y)-f(z)| 
\leq (1+\tau) \, |y-z|^{1-\tau}
\ \hbox{ for } y,z \in B(0,2),
\leqno (6.17) 
$$
$$
B(0,2-\tau) \i f(B(0,2)),
\leqno (6.18) 
$$
and
$$
E \cap B(0,2-\tau) \i f(X\cap B(0,2)) \i E.
\leqno (6.19) 
$$
In addition, the proof in [D3] says that we can take 
the same cone $X$ as in (6.10); the point is we can use
$Z(0,3) = X$ when we apply Corollary 15.11 in [D3], 
and the proof of [DTT] allows us to take $X=Z(0,3)$ above.  

In fact we shall only be interested in the restriction of
$f$ to $X\cap B(0,2)$, and in particular we do not need to know
(6.18) (or (6.17) when $y$ or $z$ lies out of $X$). The precise
biH\"older estimate will not be used either, but we need to know
that $f$ is a homeomorphism.
We shall use $f$ to find a set $G_1 \i E \cap \partial B$
with a good separation property; in the special case of
dimension $n=3$, and in particular if we used the ``Mumford--Shah" 
definition of almost minimal sets, we could obtain the separation 
property more easily. See Remark 6.46. 

\ms
First fix an arc $\C_{j,k}$, denote by $x$ and $x'$ the two endpoints
of $\C_{j,k}$, and set 
$$
\C'=\C'_{j,k} = \big\{ z \in \C_{j,k} \, ; \, 
\dist(z,\{ x, x' \}) \geq 5 \tau  \big\}
\leqno (6.20) 
$$
(a slightly smaller arc). If $\tau$ is smaller than
$\eta_0/10$, (2.2) says that $\C'$ is not empty.
Denote by $H$ the intersection of $\overline B(0,3/2)$ 
with the cone over $\C'$. Thus $H$ is a sector in a $2$-disk.
Then set 
$$
F = H \cap f^{-1}(E\cap \partial B) = \big\{ x\in H \, ; \, 
\vert f(x) \vert = 1 \big\},
\leqno (6.21) 
$$
because $f(H) \i E$ since $H \i X$ and by (6.19). Notice that
$$
\dist(z,\C') = \vert\vert z \vert - 1 \vert 
\leq \tau
\hbox{ for } z\in F,
\leqno (6.22)
$$
just because $\vert f(z) \vert = 1$ and by (6.16).
Also, $F$ is a level set of the restriction of $\vert f \vert$ to $H$,
so it separates $0$ from ${3 \over 2}\C'$ in $H$, just because (6.16)
says that $\vert f(0) \vert < 1$ and that $\vert f(y) \vert > 1$ for 
$y\in {3 \over 2}\C'$.

Now we use some topology. The set $H$ is compact, connected, 
locally connected, and simply connected. The compact set $F \i H$
separates the two connected pieces $\{ 0 \}$ and ${3 \over 2} \C'$. 
By 52.III.1 on page 335 of [Ku], 
there is a compact connected set $F_1 \i F$ such that
$$
\hbox{$F_1$ separates $0$ from ${3 \over 2} \C'$ in $H$.}
\leqno (6.23) 
$$

The author wishes to thank A. Ancona for telling him 
about this separation theorem, and is happy to share the 
statement and reference. We shall use it this way for the moment,
but in fact, since $H$ is a planar domain (a topological disk) 
we could manage with a less subtle theorem, 
for instance Theorem~14.3 on p.~123 in [Ne].  
We shall say a few words about this later, in Remark 6.45. 

For each $y \in \C'$, set 
$I_y = \big\{ ty \, ; \, 1/2 \leq t \leq 3/2 \big\}$. 
If $I_y$ did not meet $F_1$, we could use it to connect 
$0$ (or equivalently ${1 \over 2} \C'$) to 
${3 \over 2} \C'$ in $H\setminus F_1$.
So we can find $y^\sharp \in F_1 \cap I_y$.
Observe that $\vert y^\sharp-y \vert \leq
\vert\vert y^\sharp \vert - 1 \vert \leq \tau$
by (6.22), hence 
$\vert f(y^\sharp)- y \vert \leq \vert f(y^\sharp)- y^\sharp \vert
+ \vert y^\sharp - y \vert \leq 2\tau $ by (6.16).

Set $G_{j,k} = f(F_1)$. Thus $G_{j,k}$ is a connected subset of
$E \cap \partial B$ (by (6.19) and (6.21)), 
$$
\dist(y,G_{j,k}) \leq 2\tau \hbox{ for } y \in \C'
\leqno (6.24) 
$$
(because $f(y^\sharp) \in G_{j,k}$), and 
$$
\dist(z,\C') \leq 2\tau \hbox{ for } z \in G_{j,k} \, ,
\leqno (6.25) 
$$
because $z=f(y)$ for some $y \in F_1 \i F$, then 
$\vert y'-y \vert \leq \tau$ for some $y' \in \C'$, by (6.22),
and then $\dist(z,\C') \leq \vert f(y)-y' \vert 
\leq \vert f(y)-y \vert + \vert y-y' \vert  \leq 2\tau$.
 
When we apply (6.24) to the endpoints of $\C'$
and recall that $\C'=\C'_{j,k}$ is a geodesic which is barely
shorter than $\C_{j,k}$, we get that
$$
H^{1}(G_{j,k}) \geq \length(\C_{j,k}) - 15 \tau
\leqno (6.26) 
$$
because $G_{j,k}$ is connected.
By the definition (6.20) and (2.3) in particular (and if $\tau$ is small 
compared to $\eta_0$), the various $\C'_{j,k}$ lie at distances greater than 
$10\tau$ from each other, and then (6.25) says that the $G_{j,k}$ are disjoint.

Recall that $X \cap \partial B(0,1)$ is the disjoint union
of the $\C_{j,k}$, so
$$
\sum_{j,k} \length(\C_{j,k}) 
= H^{1}(X \cap \partial B(0,1))
= 2 H^2(X \cap B(0,1))
\leqno (6.27) 
$$
because $X$ is a cone. Also recall from (6.5) and (6.7) that 
$$
|H^2(X \cap B(0,1))-d(0)| \leq \tau d(0)
\leqno (6.28) 
$$
if $\varepsilon$ is small enough, depending on $\tau$.
For each fixed pair $(j,k) \in \widetilde J$,
$$\eqalign{
H^{1}(G_{j,k}) 
&= \sum_{(j',k') \in \widetilde J} H^{1}(G_{j',k'}) -
\sum_{(j',k') \neq (j,k)} H^{1}(G_{j',k'})
\cr&
= H^{1}\big(\bigcup_{(j',k') \in \widetilde J} G_{j',k'}\big) - 
\sum_{(j',k') \neq (j,k)} H^{1}(G_{j',k'})
\cr&
\leq H^{1}(E \cap \partial B) -
\sum_{(j',k') \neq (j,k)} \big[\length(\C_{j',k'})-15 \tau]
\cr&
\leq 2(1+\tau)\, d(0) - \big[ 2 H^2(X \cap B(0,1))
- \length(\C_{j,k}) \big] + C \tau 
\cr&
\leq \length(\C_{j,k}) + C_1 \tau
}\leqno (6.29) 
$$ 
because the $G_{j',k'}$ disjoint and contained
in $E \cap \partial B$, then by (6.26), (6.2), (6.27),
and (6.28).
Here $C$ and $C_1$ depend on $d(0)$ and the number of arcs,
but this does not matter, and in addition these constants would be easy 
to estimate in terms of $n$ alone. Similarly,
$$\eqalign{
H^{1}\big([E \cap \partial B]\setminus\bigcup_{j,k}G_{j,k}\big)
&=H^{1}(E \cap \partial B) - \sum_{j,k} H^{1}(G_{j,k})
\cr& \leq 2(1+\tau)\, d(0) 
- \sum_{j,k} \big[\length(\C_{j,k})-15 \tau]
\cr&
\leq 2(1+\tau)\, d(0)  - 2 H^2(X \cap B(0,1)) + C \tau
\leq C_2 \tau 
}\leqno (6.30) 
$$
again by (6.26), (6.2), (6.27), and (6.28).

\ms
Let us now construct analogues of the $G_{j,k}$ near
the vertices $x \in V$. After this, we shall connect the two types
of sets and later simplify the net that we get. So let $x\in V$
be given. Set 
$$
Y(x) = X \cap \partial B(0,1) \cap \overline B(x,C_3\tau),
\leqno (6.31) 
$$
where $C_3  = 100 + C_1 + C_2$. If $\tau$ is 
small enough, $Y(x)$ is either composed of three small arcs of
circle that make $120^\circ$ angles at $x$ (if $x\in V_0$), or
of two small arcs that leave from $x$ in opposite directions (if
$x\in V_1$). Then let $H$ denote the intersection
of $\overline B(0,3/2)$ with the cone over $Y(x)$, and set
$F = H \cap f^{-1}(E \cap \partial B)$ (as in (6.21)).

By the same separation result as above, there is a compact 
connected set $F_1 \i F$ that separates $0$ from 
${3 \over 2}Y(x)$ in $H$ (see Remark 6.45 below if you prefer 
to use the separation result in [Ne]). 
Then set $G(x) = f(F_1)$. As before, $G(x)$ is a 
connected subset of $E \cap \partial B$,
$$
\dist(y,G(x)) \leq 2\tau  \hbox{ for } y \in Y(x)
\leqno (6.32) 
$$
as in (6.24), and 
$$
\dist(z,Y(x)) \leq 2\tau \hbox{ for } z \in G(x) \, ,
\leqno (6.33) 
$$
as in (6.25). In particular, (6.32) says that $G(x)$ has branches 
that reach out reasonably far from $x$. 

\ms
Now we want connect the $G_{j,k}$ to the $G(x)$; we shall do this
in little tubes near the extremities of the $\C_{j,k}$. Let
$\C_{j,k}$ be given, and let $x, x' \in V$ denote its extremities.
Set 
$$
\C'' = \C''_{j,k,x} = \big\{ y\in \C_{j,k} \, ; \, 
10 \tau \leq |x-y| \leq (C_3-10) \tau \big\}. 
\leqno (6.34) 
$$
For $y\in \C_{j,k}$, denote by $D(y)$ the hyperdisk centered at $y$, 
of radius $3\tau$, and that lies in the hyperplane perpendicular to 
$\C_{j,k}$ at $y$. Let us check that
$$
\hbox{we can find $y\in \C''$ such that $E \cap \partial B \cap D(y)$
has exactly one point.}
\leqno (6.35) 
$$

Denote by $T=T_{j,k}$ the union of the $D(y)$, $y\in \C_{j,k}$.
By (6.25) and (6.20), $G_{j,k} \i T$. By (6.24), $G_{j,k}$ has points
within $2 \tau $ from the endpoints of $\C'$, so if we set
$$
\C''' =  \big\{ z \in \C_{j,k} \, ; \, 
\dist(z,\{ x, x' \}) \geq 10 \tau \big\} \i \C',
\leqno (6.36) 
$$
then $G_{j,k}$ crosses the union of the $D(y)$, $y\in \C'''$ (by 
connectedness). Hence, if we define the natural projection 
$\pi : T \to \C_{j,k}$ by $\pi(z)=y$ when $z\in D(y)$, we get that
$$
\pi(G_{j,k}) \hbox{ contains $\C'''$.}
\leqno (6.37) 
$$
Observe that $\pi$ is $(1+10\tau)$-Lipschitz on $T$
(by simple geometry), so (6.37) alone implies that
$H^{1}(G_{j,k}) \geq (1+10\tau)^{-1} H^{1}(\C''')
\geq (1+10\tau)^{-1} [\length(\C_{j,k}) - 21\tau]$.
We want to say that if (6.35) fails, we get a somewhat larger lower 
bound on $H^{1}(G_{j,k})$, that contradicts (6.29).

Put the lexicographic order on $\R^n$, and let $A$ denote the set of 
points $z\in G_{j,k}$ such that $\pi(z) \neq \pi(z')$ for $z'\in G_{j,k}$ 
strictly smaller than $z$ (thus, $z$ is the first point in $G_{j,k}$
with the projection $\pi(z)$).
First, $A$ is a Borel set. Indeed, $A$ is the 
intersection of the sets $A_m$, $1 \leq m \leq n$, of points $z\in G_{j,k}$ 
such that $\pi(z) \neq \pi(z')$ whenever $z' \in G_{j,k}$ is such that
$z'_i=z_i$ for $i<m$ and $z'_m < z_m$. And each $A_m$ is itself
a countable intersection of a countable union of closed sets $A_{m,p,q}$
where $\vert \pi(z) - \pi(z') \vert \geq 2^{-p}$ for
$z'\in G_{j,k}$ such that $z'_i=z_i$ for $i<m$ and $z'_m \leq z_m - 
2^{-q}$.

Next, $\pi(A)$ contains $\C'''$, by (6.37) and because $G_{j,k}$
is closed. That is, for each $\xi \in \C'''$, we can find points
$z\in G_{j,k}$ such that $\pi(z)=\xi$, and there is a first one
(get the smallest coordinate $z_j$ one $j$ at a time). So
$$\eqalign{
H^{1}(A) &\geq (1+10\tau)^{-1}  H^{1}(\pi(A))
\geq (1+10\tau)^{-1} H^{1}(\C''')
\cr& \geq (1+10\tau)^{-1} [\length(\C_{j,k}) - 21\tau]
}\leqno (6.38) 
$$
(by (6.36)) as before. 

Set $A' = G_{j,k} \setminus A$ and $A'' = E \cap \partial B
\setminus \bigcup_{(j',k')\in \widetilde J} \  G_{j',k'}$. 
If (6.35) fails, then for every $y\in \C''$, we can find at least two points 
$z \in E \cap \partial B \cap D(y)$. At most one of them lies in $A$,
because $\pi$ is injective on $A$; hence the other one lies in $A' \cup A''$, 
or in some other $G_{j',k'}$, $(j',k') \neq (j,k)$.
This last option is impossible, because if $z\in G_{j',k'}$, then
$\dist(z,\C'_{j',k'}) \leq 2\tau$ by (6.25), and $z$ lies out of $T$ 
(again if $\tau$ is small enough compared to $\eta_0$, and
by (6.20) and (2.3)). So we found $z\in A' \cup A''$ such that $\pi(z)=y$.
In other words, $\pi(A' \cup A'') \supset \C''$, hence
$$\eqalign{
H^{1}(A' \cup A'') &\geq (1+10\tau)^{-1} H^{1}(\pi(A' \cup A''))
\geq (1+10\tau)^{-1} H^1(\C'')
\cr&
\geq (1+10\tau)^{-1} [C_3 - 21] \,\tau
\geq (78+ C_1 + C_2)\, \tau
}\leqno (6.39) 
$$
(because $C_3 = 100+C_1+C_2$). Then
$$\eqalign{
H^1(G_{j,k}) &\geq H^1(A) + H^1(A')
\geq H^1(A) + H^1(A'\cup A'') - H^1(A'')
\cr&
\geq (1+10\tau)^{-1}  [\length(\C_{j,k}) - 21\tau]
+ (78+ C_1 + C_2)\, \tau - H^1(A'')
\cr&
\geq (1+10\tau)^{-1} [\length(\C_{j,k}) - 21\tau]
+ (78+ C_1 + C_2)\, \tau - C_2 \tau
\cr&
> \length(\C_{j,k}) + C_1 \tau
}\leqno (6.40) 
$$
by (6.38), (6.39), and (6.30), because $\length(\C_{j,k}) < \pi$,
and if $\tau$ is small enough.
This contradiction with (6.29) proves (6.35).

\ms
For each choice of $j$, $k$, and an endpoint $x$ of $\C_{j,k}$,
choose $y \in \C''$ as in (6.35), and denote by
$z_{j,k,x}$ the only point of $E \cap \partial B \cap D(y)$.
By the proof of (6.35), we already know that $z_{j,k,x} \in G_{j,k}$.
It also lies in $G(x)$, for the same sort of reasons: $G(x)$ is 
connected and lies close to $Y(x)$ (by (6.33)); 
by (6.32) it contains points of $T$ on both sides
of $D(y)$, so it meets $D(y)$, and then (6.35) says
that this happens at $z_{j,k,x}$.

\ms
The union of all the $G_{j,k}$ and the $G(x)$ is a subset
of $E \cap \partial B$ which is already rather nice, but we
shall need to simplify it to get the curves promised in Lemma~6.11.
We first fix $x\in V$ and simplify $G(x)$. 

Let us assume that $x\in V_0$; the other case will be simpler.
Denote by $z_1$, $z_2$, and $z_3$ the three $z_{j,k,x}$ 
that correspond to the $\C_{j,k}$ that touch $x$. 
Since $G(x)$ is connected, with $H^1(G(x)) \leq
H^1(E \cap \partial B) < +\infty$ (by (6.2)), there is a simple arc
$g_{1}$ in $G(x)$ that goes from $z_1$ to $z_2$. See for instance
Proposition 30.14 on page 188 of [D2]. 

For $1 \leq i \leq 3$, $z_i$ lies in the disk $D(y_i)$
for some $y_i$ in the corresponding $\C_{j,k}$, as in (6.35).
Let $T(x)$ denote the connected component of $x$ in 
$\big\{ z\in \R^n \, ; \, \dist(z,X\cap \partial B) \leq 3\tau
\big\} \setminus [D(y_1) \cup D(y_2) \cup D(y_3)]$. 
Thus $T(x)$ is a little tube centered at $x$ with a $Y$-shape.
Since $G(x) \i E$, (6.10) says that $G(x)$
can only meet $\partial T(x)$ on the $D(y_i)$, hence precisely
at the points $z_i$ (by (6.35)). Now $g_{1}$ is simple and
goes from $z_1$ to $z_2$, so its interior does not meet 
$z_1$ and $z_2$, and it can only meet $z_3$ at a single point.
Recall also that $g_1 \i G(x)$, which is contained in a slightly
longer tube (as in (6.33)), so $g_1$ cannot lie outside
of $T(x)$ and connect $z_1$ to $z_2$. Thus $g_1$ lies in the
interior or $T(x)$, except for its endpoints $z_1$ and $z_2$,
and perhaps $z_3$ that could be touched once.

If $x$ lies in $V_1$ instead, we construct $g_1$ as before,
except that $T(x)$ looks like a tube and there is no $z_3$, 
and we stop here.

Return to the case when $x\in V_0$. As before, there is a simple
arc $g'_2$ in $T(x)$ that goes from $z_3$ to $z_1$. This arc
meets $g_1$ (maybe immediately, if $z_3 \in g_1$, and maybe also
at $z_1$ or $z_2$; this does not matter). Let $\varphi(x)$ denote
the first point of $g_1$ that we meet when we run along $g'_2$, and 
denote by $g_2$ the arc of $g'_2$ between $z_3$ and $\varphi(x)$.
Our substitute for $G(x)$ will be $g(x)= g_1 \cup g_2$, which we see
as the union of three disjoint simple arcs in $T(x)$
that connect $\varphi(x)$ to the three $z_j$. [We don't care if
some of these arcs are reduced to one point.]

Now we define the $g_{j,k}$. Denote by $x_1$ and $x_2$
the two extremities of $\C_{j,k}$, and set
$z_i = z_{j,k,x_i}$ for $i=1,2$. Also denote by $y_i$
the point of $\C''$ such that $z_i$ is the only point of
$E \cap \partial B \cap D(y_i)$.

Since $G_{j,k}$ connects $z_1$ to $z_2$ and $H^1(G_{j,k})<+\infty$, 
we can find a simple arc $g$ in $G_{j,k}$, that goes from $z_1$
to $z_1$. Notice that $g \i G_{j,k} \i T$, where $T$ is as above
(below (6.35)). Let $T'$ denote the component of 
$T \setminus [D(y_1) \cup D(y_2)]$
that lies between these two disks; (6.10) says that $g$ can only meet 
$\partial T'$ on the disk $D(y_i)$, hence at the $z_i$. 
As before, $g$ is simple, so it only meets $\partial T'$ twice, 
at $z_1$ and $z_1$, and it is otherwise contained in the interior
of $T'$ (because we know from (6.25) that it is contained in $T$, and
otherwise it could not connect $z_1$ to $z_2$). 

Our path $g_{j,k}$ is obtained by following the arc of 
$g(x_1)$ between $\varphi(x_1)$ and $z_1$, then $g$, and then
the arc of $g(x_2)$ from $z_2$ to $\varphi(x_2)$. 

Observe that the $T(x)$, $x\in V$, are disjoint, that
$T(x)$ only meets the tube $T'$ associated to $\C_{j,k}$ 
when $x$ is one of the extremities of $\C_{j,k}$ (and then
the intersection is contained in the disk $D(y_{j,k,x})$ that 
contains $z_{j,k,x}$), and that the various
$T'$ are also disjoint. Then the arcs $g_{j,k}$ are simple,
and $g_{j,k}$ only meets $g_{j',k'}$ when $\C_{j,k}$ and
$\C_{j',k'}$ have a common endpoint $x$, in which case
$g_{j,k}\cap g_{j',k'}=\{ \varphi(x) \}$.

\ms
We are now ready to prove Lemma 6.11.
We already defined the $\varphi(x)$, $x\in V$, and the curves 
$g_{j,k}$. We have (6.12) by construction, we just checked (6.13),
and (6.14) holds because $\varphi(x) \in g(x) \i T(x)$, which itself 
is contained in $B(x,(C_3+3) \tau)$ (because the $y_i$ lie in $\C''$,
see (6.34), (6.35), and the definition of $T(x)$ somewhat below
(6.40)). 

So we just need to check (6.15). 
Recall that each $g_{j,k}$ is composed of an arc $g \i G_{j,k}$,
which stays within $2 \tau$ from $\C_{j,k}$ by (6.25), 
an arc of $g(x_1)$ which is contained in $T(x_1) \i B(x_1,(C_3+3)\tau)$,
and an arc of $g(x_2) \i T(x_2) \i B(x_2,(C_3+3) \tau)$.
Thus $g_{j,k}$ stays within $(C_3+3) \tau)$ from $\C_{j,k}$,
as needed. This completes our proof of Lemma 6.11.
\qed

\ms
Notice that
$$
H^1(g_{j,k}) \geq \length(\C_{j,k}) - C \tau 
\leqno (6.41) 
$$
because $\C_{j,k}$ is a geodesic with length at most
${9 \pi \over 10}$ (by (2.5)), and because by (6.14)
its endpoints are within $C \tau$ from the endpoints of 
$g_{j,k}$. On the other hand,
$$
\sum_{j,k} H^1(g_{j,k}) = H^{1}\big(\bigcup_{j,k} g_{j,k}\big) 
\leq H^1(E \cap \partial B) \leq 2(1+\tau)\, d(0)
\leqno (6.42)
$$
because the $g_{j,k}$ are disjoint and by (6.2), so
$$
H^1(g_{j,k}) \leq \length(\C_{j,k}) + C_4 \tau
\leqno (6.43)
$$
by the proof of (6.29), and then 
$$
H^{1}\big([E \cap \partial B]\setminus\bigcup_{j,k}g_{j,k}\big)
\leq C_5 \tau
\leqno (6.44) 
$$
by the proof of (6.30).

\smallskip
We end this section with two remarks on the construction of the
connected sets $G_{j,k}$ and $G(x)$.

\ms
\noindent{\bf Remark 6.45.}
We could also use Theorem 14.3 in [Ne], 
instead of the stronger separation result from [Ku], 
to construct the $G_{j,k}$ and analogues of the $G(x)$. 

We start with $G_{j,k}$. We still work in the intersection $H$
of $\overline B(0,3/2)$ with the cone over $\C' = \C'_{j,k}$ 
(see near (6.21)), and need to know that if the compact set 
$F = H \cap f^{-1}(E \cap \partial B)$ separates $0$ from 
${3 \over 2}\C'$ in $H$, then there is 
a compact connected set $F_1 \i F$ that still separates 
$0$ from ${3 \over 2}\C'$ in $H$ (as in (6.23)). 

Theorem 14.3 in [Ne] gives this, 
but it is stated when $H$ is the
two-dimensional sphere. We just need a small argument to reduce to 
that case. Let $\widetilde H$ be a second copy of $H$, with
the copy $\widetilde F$ of $F$ in it. Glue $\widetilde H$ to 
$H$, by identifying every point of $\partial H$
(two arc segments from the origin to the extremities of 
${3 \over 2}\C'$, plus ${3 \over 2}\C'$ itself) to its
copy in $\widetilde H$. This gives a sphere $S = H \cup \widetilde H$
and a compact set $\widehat F = F \cup \widetilde F \i S$.

Let us check that $\widehat F$ separates $0$ from ${3 \over 2}\C'$
in $S$. If not, there is a curve $\zeta \i S \setminus \widehat F$ that goes 
from $0$ to ${3 \over 2}\C'$. If $\zeta \i H$, we get a contradiction immediately 
with the fact that $F$ separates $0$ from ${3 \over 2}\C'$ in $H$.
Otherwise, let $\widehat \zeta$ denote the union of $\zeta$ and its
symmetric image $\widetilde \zeta$ (obtained by exchanging the copies
$H$ and $\widetilde H$); then $\widehat \zeta$ is connected because
$\zeta$ meets $\widetilde \zeta$ at the origin, and $\widehat \zeta$ 
does not meet $\widehat F$ because $\widehat F$ is symmetric. 

Set $\zeta' = \widehat \zeta \cap H$. It is connected too,
because if $O_1$ and $O_2$ are disjoint open sets in $H$ that cover
$\zeta'$, then the symmetric extensions $\widehat O_1$ and 
$\widehat O_2$ are disjoint, open in $S$, and cover 
$\widehat \zeta$, a contradiction. Finally, 
$\zeta'$ does not meet $F$ because $\widehat\zeta$ 
does not meet $\widehat F$. Since $\zeta'$ contains 
$0$ and meets ${3 \over 2}\C'$, we get a contradiction. So $\widehat F$ 
separates $0$ from ${3 \over 2}\C'$ in $S$.

By Theorem 14.3 in [Ne], 
we get a connected piece $F' \i \widehat F$ that separates 
$0$ from ${3 \over 2}\C'$ in $S$. Its closure is still connected and
contained in $\widehat F$ (because $\widehat F$ is compact),
so we can assume that $F'$ is compact.

Denote by $\widetilde F'$
the symmetric copy of $F'$, and set $\widehat F' = F' \cup \widetilde F'$.
Then $\widehat F'$ is also contained in $\widehat F$ (because
$\widehat F$ is symmetric), and it also separates 
$0$ from ${3 \over 2}\C'$ in $S$. In addition $F'$ meets $\partial H$, 
because otherwise it would not separate $0$ from ${3 \over 2}\C'$ in 
$S$, so $F'$ meets $\widetilde F'$, and $\widehat F'$ is connected. 

Set $F_1 = \widehat F' \cap H$. This set is contained in $F$
(because $\widehat F' \i \widehat F$, and it is connected too; 
the argument is the same as for $\zeta'$ above. It separates 
$0$ from ${3 \over 2}\C'$ in $H$ because $F'$ separates in $S$,
so it satisfies all the desired properties.

\ms
Now we construct substitutes for the $G(x)$, $x\in V$.
If $x\in V_1$, the set 
$Y(x) = X \cap \partial B(0,1) \cap \overline B(x,C_3 \tau)$ 
in (6.31) is composed of two short arcs of great circles
that leave from $x$ in opposite directions, the piece of cone $H$
is a again a simple topological disk, and we can use 
Theorem 14.3 in [Ne] as above.

When $x\in V_0$, $Y(x)$ has a fork, and we shall proceed differently.
Each branch of $Y(x)$ comes from a $\C_{j,k}$, and we use (6.35) to
find a point $y$ in the branch such the disk $D(y)$ only meets 
$E \cap \partial B$ once, at a point $z_{j,k,x}$ of $G_{j,k}$.
This  way, we get three points $y_{j,k,x}$, one in each branch
of $Y(x)$, which we decide to call $y_1$, $y_2$, and $y_3$. 
Also call $z_i$ the only point of $E \cap \partial B \cap D(y_i)$.
All we need to do is find a connected set $G(x) \i T(x)$ that connects 
the three $z_i$, because after this we notice that since $D(y_i)$ 
only meet $E \cap \partial B$ once, each $z_i$ lies in the corresponding 
$G_{j,k}$, and we can proceed as before.

First we want to connect $z_1$ to $z_2$. Let $Y_i$ denote
the arc of $Y(x)$ that goes through $y_i$, and let
$H$ be the intersection of $\overline B(0,3/2)$ with the cone 
over $Y_1 \cup Y_2$. Then set $F = H \cap f^{-1}(E \cap \partial B)$, 
and notice that $F$ separates $0$ from ${3 \over 2}Y$ in $H$, again 
because it is a level set. Since $H$ is a simple topological disk, 
we can apply Theorem 14.3 in [Ne] as above, 
and find a compact connected set $F_1 \i F$ that separates 
$0$ from ${3 \over 2}Y$ in $H$. Denote by $y'_1$ and $y'_2$
the endpoints of $Y_1$; notice that $F_1$ contains a point in 
$I_1 = [0,y'_1]$, because otherwise we could use $I_1$ to connect 
$0$ to ${3 \over 2}Y$. Similarly, $F_1$ meets $I_2 = [0,y'_2]$.

Now set $G = f(F_1)$; $G$ is compact, connected, and contained in 
$E \cap \partial B \cap T(x)$. It contains points in $f(I_1)$
and $f(I_2)$, so it crosses the two disks $D(y_i)$,
and hence contains $z_1$ and $z_2$. Similarly, there is a 
compact connected set $G' \i E \cap \partial B$ that contains
$z_1$ and $z_3$, and we can take $G(x) = G \cup G'$.
Thus we can use Theorem 14.3 in [Ne] 
instead of 52.III.1 in [Ku]. 

\ms
\noindent{\bf Remark 6.46.}
When $n = 3$, we can get the connected sets $G_{j,k}$ and $G(x)$
(and then proceed as above) with a little less information than 
the full parameterization $f$ of (6.16)-(6.19).

Let $X$ be as in (6.10), and denote by $W_i$, $1 \leq i \leq 2$, 
$3$, or $4$ (depending on the type of $X$) the connected components 
of $\overline B \setminus X_{100\varepsilon}$, where 
$X_{100\varepsilon} = \big\{ y\in \R^3 \, ; \, \dist(y,X) 
\leq 100\varepsilon \big\}$ denotes a closed $100\varepsilon$-neighborhood 
of $X$. We know from (6.10) that the $W_i$ don't meet $E$. 
Now assume that 
$$
\hbox{the $W_i \cap \partial B$ lie in different connected 
components of $E \cap \partial B$.}
\leqno (6.47)
$$
[This is obviously the case if we know that $E$ separates the 
$W_i$ in $\overline B$, and the existence of the parameterization 
$f$ in (6.16)-(6.19) gives this information.] 

Then we can use a slightly different separation argument 
(also based on Theorem 14.3 in [Ne]) 
to produce the connected sets $G_{j,k}$ and $G(x)$. This time we work
directly in slices of $\partial B$ and find connected pieces of 
$E \cap \partial B$ that separate, instead of using $f^{-1}$ and working 
on a level set in a piece of cone.

For instance, when we need to construct $G(x)$ for some
$x\in V_0$, we work in the topological disk 
$H = \partial B \cap B(x,10^{-1})$, notice that $E\cap H$ separates 
the three main regions of $H$ from each other, and use the same argument 
as in Remark 6.45 to find connected sets $G$ and $G'$ in $E \cap H$ 
that separate them too.

It is a little sad that (still when $n=3$) the author was not able to
find a simple proof of the needed separating property that would not 
use the parameterisation from [D3], because it is the only place 
in the argument where we use the geometric part of [D3] and [DDT] 
heavily. On the other hand, if we want to get all the way to the 
regularity property stated in Theorem 1.15, this does not make  a 
huge difference, because arguments of the same type 
(i.e., an extension of the Reifenberg parameterization theorem) 
are needed there anyway. 

Let us also mention that in the special case of ``Mumford-Shah"
minimal and almost minimal sets (see Section 18 of [D3]), 
(6.47) is a rather easy consequence of the definitions 
(see the proof of Theorem 1.9 in [D3], and in particular 
the proof of (18.31), near the end of Section 18), 
So at least we get a slightly simpler proof in this case.

\bigskip
\noindent {\bf 7. Construction of a Lipschitz graph and length estimates} 
\ms

The construction of this section will be applied, after a rotation, 
to each of the curves $g_{j,k}$, $(j,k) \in \widetilde J$, 
that were produced in Lemma 4.11, but the section is fairly independent, 
and we use slightly different notation.

Here again $B$ is the unit ball in $\R^n$, and we give ourselves 
a simple rectifiable curve $\gamma$ in $\partial B$. We assume that
$$
9\eta_0 \leq \length(\gamma) \leq {10\pi \over 11},
\leqno (7.1)
$$
where $\eta_0 > 0$ is as in (2.2).

The first inequality is mostly here for convenience or normalization.
The second one will be needed in some estimates; the point is to make 
sure that we are far enough from the situation of non uniqueness for 
geodesics, and incidentally $10 \pi/11$ could be be replaced with 
any constant strictly smaller than $\pi$.

Denote by $a$ and $b$ the extremities of $\gamma$. 
We assume that
$$
\length(\gamma) \leq \dist_{\partial B}(a,b) + \tau_1,
\leqno (7.2)
$$
where we denote by $\dist_{\partial B}$ the geodesic distance on 
$\partial B$, and $\tau_1$ is a positive constant that we can take 
as small as we want.

We also assume, for simplicity, that $a$ and $b$ lie on the horizontal 
$2$-plane $P$ through the origin, and that 
$$
\dist(z,P) \leq \tau_1 \ \hbox{ for }z \in \gamma.
\leqno (7.3)
$$
This does not cost us too much, because (7.1) and (7.2) imply
the analogue of (7.3), with $\tau_1$ replaced with a constant 
$\tau'_1$ which tends to $0$ with $\tau_1$. 
We won't need to bother, though, because (7.3)
will come for free when $\gamma$ is one of the $g_{j,k}$.

\smallskip
We want to construct a Lipschitz curve $\Gamma$ on $\partial B$,
with the same extremities $a$ and $b$, whose Lipschitz constant 
is at most $\eta$, and which has a big intersection with $\gamma$. 
Here $\eta > 0$ is given in advance, small enough, and we are allowed 
to choose $\tau_1$ very small, depending on $\eta$. Later on, we shall
connect this section with the previous one, and we shall be allowed 
to take $\tau$ small, depending on $\eta$ and $\tau_1$.
It will be important to estimate $H^1(\gamma \setminus \Gamma)$, 
for instance, in terms of 
$$
\Delta L = \length(\gamma)-\dist_{\partial B}(a,b),
\leqno (7.4)
$$
and not just $\tau_1$ or $\tau$. See (7.30)-(7.32) 
below for the main properties of $\Gamma$.
Let us first prove some estimates, and then construct $\Gamma$. 

Let $z : I \to \partial B$ denote a parameterization of $\gamma$
by arc-length. Here $I$ is a compact interval, and 
$|I| = \length(\gamma)$. We write
$z(t) = (z_1(t),z_2(t),v(t))$, with $v(t) \in \R^{n-2}$.
Notice that $|v(t)| \leq \tau_1$, by (7.3). Then set
$w(t) = (1-|v(t)|^2)^{1/2} = (z_1^2(t)+z_2^2(t))^{1/2} \geq 1-\tau_1$,
and write
$$
z(t) = \big(\cos\theta(t)w(t), \sin\theta(t)w(t),v(t)\big)
\leqno (7.5)
$$ 
for $t\in I$. We take a continuous determination of
$\theta$, which is easy because $(z_1(t),z_2(t))$ does not
vanish. Notice that $w$, and then $\theta$, are Lipschitz. Also set 
$$
l = \length(\gamma) = |I| \ \hbox{ and } \ d=\dist_{\partial B}(a,b).
\leqno (7.6)
$$ 
Thus $\Delta L = l-d$. 
We shall assume that $\gamma$ roughly runs counterclockwise on average, 
so that
$$
\int_{I} \theta'(t) dt = \dist_{\partial B}(a,b) = d \, ;
\leqno (7.7)
$$
otherwise, we could always parameterize $\gamma$ backwards and get
the same results.

\medskip
\proclaim Lemma 7.8. If $\tau_1$ is small enough,
$$
\int_I |v'(t)|^2 dt \leq 14 \, \Delta L.
\leqno (7.9)
$$

\medskip
First differentiate (7.5) to get that
$$\eqalign{
z'(t) = \theta'(t)& w(t) \big(-\sin\theta(t),\cos\theta(t),0 \big)
\cr&+ w'(t) \big(\cos\theta(t),\sin\theta(t),0\big)
+ \big(0,0,v'(t)\big).
}\leqno (7.10)
$$ 
The three pieces are orthogonal, so
$$
|z'(t)|^2 = \theta'(t)^2 w(t)^2 + |w'(t)|^2 + |v'(t)|^2. 
\leqno (7.11)
$$
Recall that $w(t)^2 + |v(t)|^2 = 1$, so
$w(t) w'(t) = - \langle v'(t),v(t) \rangle$. Also,
$z$ is a parameterization by arc-length, so $|z'(t)| = 1$ 
almost-everywhere, i.e., 
$$
\theta'(t)^2 w^2(t) = 1 - |w'(t)|^2 - |v'(t)|^2
= 1 - {\langle v'(t),v(t) \rangle^2 \over w(t)^2} - |v'(t)|^2.
\leqno (7.12)
$$
Set $\lambda = 101/100$ and recall that $|v(t)| \leq \tau_1$
by (7.3), so $w^{-2}(t) = [1-|v(t)|^2]^{-1} \leq 1+ \lambda |v(t)|^2$ 
if $\tau_1$ is small enough. Hence, ignoring the non-positive 
middle term in (7.12),
$$
\theta'(t)^2 
\leq \big[1+\lambda |v(t)|^2 \big] (1 - |v'(t)|^2) 
\leq 1 + \lambda |v(t)|^2 - |v'(t)|^2
\ \hbox{ almost-everywhere on } I.
\leqno (7.13)
$$

We may assume (by translation) that $I$ starts at the origin. Thus
$I=[0,l]$, with $l = \length(\gamma)$ as above. Observe that 
$v(0)= v(l)=0$ because $a$ and $b$ lie on $P$. 
We can extend $v$ to $[-l,l]$ so that it is odd, and then
write $v$ as a sum of sines. That is, write 
$$
v(t) = \sum_{k \geq 1} c_k\sin(\pi kt/l)
\leqno (7.14)
$$
(with vector-valued coefficients $c_k$ when $n > 3$). Then
$$
||v||_{L^2(I)}^2 = {l \over 2} \sum_{k\geq 1} \, |c_k|^2
\leqno (7.15)
$$
and
$$
||v'||_{L^2(I)}^2 = {l \over 2}  \sum_{k\geq 1} (\pi k/l)^2 |c_k|^2
\geq \Big({\pi \over l}\Big)^2 ||v||_{L^2(I)}^2 
\geq {121 \over 100} \, ||v||_{L^2(I)}^2 \, ,
\leqno (7.16)
$$
by (7.1) and (7.6). We now integrate (7.13) and get that
$$
\int_I \theta'(t)^2 dt
\leq |I| + \lambda ||v||_{L^2(I)}^2 - ||v'||_{L^2(I)}^2
\leq l - \Big(1 - {100\lambda \over 121}\Big) ||v'||_{L^2(I)}^2
\leq l - {1\over 7} \,||v'||_{L^2(I)}^2 \, .
\leqno (7.17)
$$
Notice that 
$$
d^2 = \Big\{ \int_I \theta'(t) dt \Big\}^2
\leq l \int_I \theta'(t)^2 dt
\leq l^2 - {l\over 7} \int_I |v'(t)|^2 dt
\leqno (7.18)
$$
by (7.7), Cauchy-Schwarz, (7.6), and (7.17). Hence
$$
\int_I |v'(t)|^2 dt \leq {7 \over l} \, (l^2-d^2)
= 7 \, {l+d \over l} \, \Delta L \leq 14 \, \Delta L
\leqno (7.19)
$$
(see below (7.6)). Lemma 7.8 follows. \qed  

\medskip
We also need some control on $\theta'$. Set 
$$
f(t)= 1 + 2|v(t)|^2 - \theta'(t). 
\leqno (7.20)
$$
Observe that $\theta'(t)^2 \leq 1 + 2 |v(t)|^2$ almost-everywhere, by 
(7.13), so $f(t) \geq 0$. On the other hand,
$$\eqalign{
\int_I f(t) dt &= l + 2 \int_I |v(t)|^2 dt - \int_I \theta'(t) dt
= l - d + 2 \int_I |v(t)|^2 dt 
\cr&
= \Delta L + 2 \int_I |v(t)|^2 dt
\leq \Delta L + 2 \int_I |v'(t)|^2 dt 
\leq 30 \, \Delta L
}\leqno (7.21)
$$
by (7.7), the line below (7.6), (7.16), and (7.19).
So $f(t)$ is often small, and hence $\theta'(t)$ is rarely 
smaller than $1/2$, say.

\medskip
We are now ready to construct our $\eta$-Lipschitz curve $\Gamma$.
It will be obtained by modifying $z$ on open 
intervals where some maximal functions are large. First denote by $v^\ast$ 
the (non centered) Hardy-Littlewood maximal function of the derivative 
of $v$; 
thus
$$
v^\ast(t) = \sup \Big\{ {1 \over |J|} \int_{J} |v'(s)| ds 
\, ; \, J \hbox{ is an interval contained in $I$ such that } t\in J \Big\}.
\leqno (7.22)
$$
The Hardy-Littlewood maximal theorem (see for instance [St], p.5) 
says that $||v^\ast||_{2}^2 \leq C ||v'||_{2}^2 \leq C \Delta L$,
by (7.19). Set $Z_1 =  \{ t\in I \, ; \, v^\ast(t)> \eta/4 \}$.
Thus 
$$
|Z_1| \leq 16 \eta^{-2} ||v^\ast||_{2}^{2} 
\leq C \eta^{-2}  \Delta L.
\leqno (7.23)
$$
Similarly set 
$$
f^\ast(t) = 
\sup \bigg\{ {1 \over |J|} \int_{J} f(s) \, ds 
\, ; \, J \hbox{ is an interval contained in $I$ such that } t\in J \bigg\}
\leqno (7.24)
$$
and $Z_2 =  \{ t\in I \, ; \, f^\ast(t) > 1/2 \}$. Then
$$
|Z_2| \leq 2 ||f^\ast||_{{\rm Weak} L^1} \leq C ||f||_1 \leq C \Delta L,
\leqno (7.25)
$$
by the Hardy-Littlewood maximal theorem again and (7.21). Finally set
$Z = Z_1 \cup Z_2$; thus 
$$
|Z| \leq C \eta^{-2}  \Delta L.
\leqno (7.26)
$$

It is easy to see (just from the definition of $Z_1$ and $Z_2$ with
maximal functions) that $Z$ is open in $I$. We want to keep $z$ as it 
is on $I \setminus Z$, and on the open intervals that compose $Z$,
replace $\gamma$ with arcs of geodesics.

First observe that if $t\in I\setminus Z$, $t' \in I$, and $J$ denotes the 
interval with endpoints $t$ and $t'$, then
$$
|v(t') - v(t)| 
\leq \int_J |v'(u)|\, du
\leq |t'-t| \, v^\ast(t) \leq \eta |t'-t|/4,
\leqno (7.27)
$$
by (7.22) and because $t\notin Z_1$. Similarly, if in addition $t' \geq t$, 
$$\eqalign{
\theta(t') - \theta(t) &=
\int_J [f(t) + \theta'(t)] \, dt - \int_J f(t)dt
= \int_J [ 1 + 2v(t)^2] dt - \int_J f(t)dt
\cr&
\geq (t'-t) - (t'-t)f^\ast(t) \geq (t'-t)/2
}\leqno (7.28)
$$
by (7.20), (7.24), and because $t\notin Z_2$. If instead $t' < t$, the 
same argument yields $\theta(t) - \theta(t') \geq (t-t')/2$. Altogether,
$$
\theta(t') - \theta(t) \geq (t'-t)/2 \ \ \hbox{ when 
$0  \leq t \leq t' \leq l$, and $t$ and $t'$ do not both lie in $Z$.}
\leqno (7.29)
$$

So we want to keep the part of $\gamma$ that corresponds to
$I \setminus Z$, and replace $z(Z)$ with arcs of geodesics.
Recall that $Z$ is open in $I$, so it is an at most countable union 
of disjoint intervals $I_j$ that are open in $I$. 
Call $a_j$ and $b_j$ the endpoints of $I_j$, and denote by
$\rho_{j}$ the arc of geodesic on $\partial B$ that goes from
$z(a_j)$ to $z(b_j)$. For each $j$, we replace the arc of $\gamma$ 
between $z(a_j)$ and $z(b_j)$ with $\rho_{j}$.
This gives a new curve $\Gamma$. First observe that
$$
\Gamma \hbox{ has the same endpoints as } \gamma,
\leqno (7.30)
$$
because even if the endpoints of $I$ lie in $Z$, we keep
$z(0)$ and $z(l)$ as endpoints of the corresponding geodesics. Also,
$$
H^1(\Gamma \setminus \gamma) \leq H^1(\gamma \setminus \Gamma) 
\leq C \eta^{-2} \Delta L;
\leqno (7.31)
$$
the first inequality holds because the geodesic $\rho_{j}$ 
is never longer than the  arc of $\gamma$ that it replaces, so that
$H^1(\Gamma) \leq \length(\Gamma) \leq \length(\gamma) = H^1(\gamma)$ 
(because $\gamma$ is simple), while the second inequality comes from 
(7.26) and the fact that $z$ is $1$-Lipschitz.
To end this section, we need to check that 
$$
\Gamma \hbox{ is a Lipschitz graph with constant }\leq \eta,
\leqno (7.32)
$$
and at the same time we shall explain what we mean by this
(see (7.42) or (7.44)). Let us first check that
$$
|v(b_j) - v(a_j)| \leq \eta (b_j-a_j)/4
\ \hbox{ and } \ 
\theta(b_j) - \theta(a_j) \geq (b_j-a_j)/2.
\leqno (7.33)
$$
Recall that $I_j$ is a connected component of $Z$, which is open in $I$.
So $a_j \in I \setminus Z$, except perhaps if $a_j = 0$. Similarly, 
$b_j\in I\setminus Z$, unless it is the final endpoint of $I$.
Both things cannot happen at the same time, because
$$
|b_j-a_j| \leq |Z| \leq C \eta^{-2} \Delta L 
\leq C \eta^{-2} \tau_1 < |I| 
\leqno (7.34)
$$ 
by (7.26), (7.2), (7.4), (7.1), and if $\tau_1$ is small enough. 
So $a_j$ or $b_j$ lies in $I\setminus Z$, we can apply (7.27) and
(7.29) to them, and we get (7.33).

We want to estimate the average slope of $\gamma$ between $z(a_j)$
and $z(b_j)$. Denote by $\pi$ the projection on the horizontal plane.
For the next estimate, it is convenient to use complex notations; then
$\pi(z(t)) = w(t) e^{i\theta(t)}$ by (7.5), and
$$\eqalign{
|\pi(z(b_j))-\pi(z(a_j))| 
&= \big|w(b_j)e^{i\theta(b_j)} - w(a_j)e^{i\theta(a_j)} \big|
\cr&= |e^{i\theta(a_j)}| \, 
\big|w(b_j)e^{i[\theta(b_j)-\theta(a_j)]} - w(a_j)\big|
\cr&= \big|w(b_j)[e^{i[\theta(b_j)-\theta(a_j)]} -1]
+ [w(b_j) - w(a_j)] \big|
\cr&\geq |w(b_j)| \big| e^{i(\theta(b_j)-\theta(a_j))} - 1 \big|
-|w(b_j)- w(a_j)|.
}\leqno (7.35)
$$ 
Recall that $|v(b_j)| \leq \tau_1$ by (7.3), so
$w(b_j) \geq 1-\tau_1$. Also, $|b_j-a_j|$ is as  small as we want, 
by (7.34), and $|\theta(b_j)-\theta(a_j)|$ is at most twice larger, 
for instance because (7.13) says that $|\theta'| \leq 2$. 
Then $|e^{i(\theta(b_j)-\theta(a_j))} - 1| \geq {9 \over 10} 
\, |\theta(b_j)-\theta(a_j)| \geq {9 \over 20} \, (b_j-a_j)$,
by (7.33). Finally, $|w(b_j)-w(a_j)| \leq |v(b_j)-v(a_j)| 
\leq \eta (b_j-a_j)/4$, because $x \to (1-x^2)^{1/2}$ is $1$-Lipschitz
near the origin, and by (7.33) again. Altogether, 
$$
|\pi(z(b_j))-\pi(z(a_j))| 
\geq {9 \over 20} \, (1-\tau_1) (b_j-a_j) - \eta (b_j-a_j)/4
\geq {8 \over 20} \, (b_j-a_j)
\leqno (7.36)
$$
if $\tau_1$ and $\eta$ are small enough. Then (7.33) yields
$$
|v(b_j) - v(a_j)| \leq \eta (b_j-a_j)/4 \leq {5 \eta \over 8} \,
|\pi(z(b_j))-\pi(z(a_j))|.
\leqno (7.37)
$$

Recall that $|z(b_j))-z(a_j)| \leq |b_j-a_j|$ is as  small as we want, 
so the variation of the unit tangent vector to the geodesic $\rho_j$ 
between $z(a_j))$ and $z(b_j))$ is also as small as we want. Then
the slope of that tangent is less than $3\eta/4$.

The reader is probably already convinced that (7.32) will easily 
follow from this, but let us complete the argument brutally. Let us
even parameterize $\Gamma$ with a function 
$\widetilde z : I \to \partial B$, as we did for $\gamma$ in (7.5).
We need to define functions $\widetilde\theta$, $\widetilde v$, 
$\widetilde w$, and $\widetilde z$ on $I$, so that
$$
\widetilde z(t) = \big(\cos\widetilde\theta(t) \, \widetilde w(t), 
\sin\widetilde\theta(t) \, \widetilde w(t),\widetilde v(t)\big)
\leqno (7.38)
$$ 
and, since $\widetilde z(t) \in \partial B$,
$\widetilde w(t)^2 + |\widetilde v(t)|^2 = 1$.

We keep $\widetilde\theta(t)= \theta(t)$, $\widetilde v(t)= v(t)$, 
$\widetilde w(t)=w(t)$, and $\widetilde z(t)=z(t)$ for 
$t\in I \setminus Z$.
Let $I_j$ be a component of $Z$, and again denote its extremities 
by $a_j$ and $b_j$. We keep $\widetilde\theta(a_j)= \theta(a_j)$
and $\widetilde\theta(b_j)= \theta(b_j)$, even if $a_j$ or $b_j$
lies in $Z$, and we define $\widetilde\theta$ so that it is affine
on $I_j$. Recall that the geodesic $\rho_j$ from $z(a_j)$ to $z(b_j)$
is short, with a small slope; then, given 
$\theta \in [\theta(a_j),\theta(b_j)]$, there is a unique point
$\xi = \xi(\theta) \in \rho_j$ such that $\pi(\xi) = e^{i\theta}|\pi(\xi)|$.
In other words, the fact that $\widetilde z(t) \in \rho_j$ and 
that we want (7.38) with a given $\widetilde\theta(t) \in 
[\theta(a_j),\theta(b_j)]$ determine $\widetilde z(t)$,
$\widetilde v(t)$, and $\widetilde w(t)$ uniquely, and
these functions are smooth on $I_j$.

Also recall that because of (7.37), the slope of the tangent vector 
to $\rho_j$ stays smaller than $3\eta/4$. Then
$$
|\widetilde v'(t)| \leq {3\eta \over 4} \, 
|(\pi \circ \widetilde z)'(t)| \ \hbox{ on $I_j$.}
\leqno (7.39)
$$
But (7.38) yields
$$
|(\pi \circ \widetilde z)'(t)|^2 = 
\widetilde \theta'(t)^2 \widetilde w(t)^2 + \widetilde w'(t)^2
\leq \widetilde\theta'(t)^2 \widetilde w(t)^2 
+|\widetilde v'(t)|^2 |\widetilde v(t)|^2 \widetilde w(t)^{-2}
\leqno (7.40)
$$ 
by the proof of (7.11), and because $\widetilde w(t) \widetilde w'(t) 
= - \langle \widetilde v'(t),\widetilde v(t) \rangle$
(since $|\widetilde v|^2+|\widetilde w|^2 =1$).
By (7.3), the two extremities of $\rho_j$ lie within $\tau_1$ of $P$,
so $\rho_j$ itself lies within $\sqrt 2 \tau_1$ of $P$,
$|\widetilde v(t)|^2 \leq 2 \tau_1^2$, $|\widetilde w(t)|^2 \geq 
1-2\tau_1^2$, and hence
$$\eqalign{
|\widetilde v'(t)|^2 
&\leq \Big({3\eta \over 4}\Big)^2 \, |(\pi \circ \widetilde z)'(t)|^2
\leq \Big({3\eta \over 4}\Big)^2 \, 
\big[\widetilde\theta'(t)^2 \widetilde w(t)^2 +
|\widetilde v'(t)|^2 |\widetilde v(t)|^2 \widetilde w(t)^{-2}\big]
\cr&
\leq \Big({3\eta \over 4}\Big)^2 \, \widetilde\theta'(t)^2 
+ \Big({3\eta \over 4}\Big)^2 \, 
{2 \tau_1^2 \over 1-2\tau_1^2} \, |\widetilde v'(t)|^2
\leq \Big({3\eta \over 4}\Big)^2 \, \widetilde\theta'(t)^2 
+ 10^{-2} |\widetilde v'(t)|^2
}\leqno (7.41)
$$ 
by (7.39) and (7.40), and if $\eta$ and $\tau_1$ are small enough. 
This yields
$$
|\widetilde v'(t)| \leq {4\eta\over 5} \; \widetilde\theta'(t)
\ \hbox{ for } t \in I_j
\leqno (7.42)
$$
(recall that $\widetilde\theta'(t) > 0$ by (7.33)
and because $\widetilde\theta$ is affine on $I_j$).

\smallskip
We are ready to check that
$$
|\widetilde v(t')-\widetilde v(t)|
\leq {4\eta \over 5} \, [\widetilde\theta(t') -\widetilde\theta(t)]
\ \hbox{ for } t,t' \in I \hbox{ such that } t \leq t'.
\leqno (7.43)
$$
When $t$, $t'$ lie in a same $I_j$, this follows from (7.42).
Next suppose that $t \in I_j$ and $t'\in I_k$, with $k \neq j$. Then
$b_j$ and $a_k$ both lie in $I\setminus Z$ (because $Z$ is open), and
$$\eqalign{
|\widetilde v(t')-\widetilde v(t)|
&\leq |\widetilde v(t')-\widetilde v(a_k)|
+ |\widetilde v(a_k)-\widetilde v(b_j)|
+ |\widetilde v(b_j)-\widetilde v(t)|
\cr&
\leq {4\eta \over 5} (\widetilde\theta(t')-\widetilde\theta(a_k)) 
+ {\eta \over 2} \, (\widetilde\theta(a_k)-\widetilde\theta(b_j)) 
+ {4\eta \over 5} (\widetilde\theta(b_j)-\widetilde\theta(t))
\cr&\leq {4\eta \over 5} (\widetilde\theta(t')-\widetilde\theta(t))
}\leqno (7.44)
$$
by (7.42), (7.27), and (7.29), as needed. 
The case when $t$, or $t'$, or both, lie in $I\setminus Z$ is similar. 

So (7.43) holds. We can take this, or rather the fact that 
$\widetilde v(t)$ is an $\eta$-Lipschitz function of $\widetilde\theta(t)$, 
as a strong definition of (7.32). [That is, the altitude of a point of $\Gamma$
is an $\eta$-Lipschitz function of its projection on $\partial B \cap P$.] 

Notice that (7.43) implies the other reasonable definition 
of (7.32), i.e., that 
$$
\widetilde v(t) \hbox{ is an $\eta$-Lipschitz function of }
\pi(\widetilde z(t)).
\leqno (7.45)
$$ 
Indeed $|(\pi \circ \widetilde z)'(t)| \geq \widetilde\theta'(t)
\,\widetilde w(t) \geq {5 \over 6} \, \widetilde\theta'(t)$
by the first half of (7.40) and because $|\widetilde w(t)|^2 \geq 
1 - 2\tau_1^2$. Then $\widetilde\theta(t)$ is a ${6 \over 5}$-Lipschitz 
function of $\pi(\widetilde z(t))$, and (7.45) follows from (7.43).

\bigskip
\noindent {\bf 8. Improvement of the cone over a small Lipschitz graph} 
\medskip

In this section we consider a Lipschitz graph $\Gamma$ with small 
constant contained in the unit sphere, and we try to find a surface 
in the unit ball, with the same boundary as the cone over $\Gamma$, 
but with a smaller area.

Later on, this construction will be applied to the Lipschitz curves
$\Gamma_{j,k}$ that were constructed in the previous section, which 
themselves will come from the curves $g_{j,k}$ that we found in 
Section 6. The general idea is to find successive improvements 
of our almost-minimal set $E$, but the way we reduce to Lipschitz 
graphs will be explained later.

So we are given a curve $\Gamma$ in the unit sphere
$\partial B$. We assume for simplicity that its extremities $a$ and
$b$ lie in the horizontal $2$-plane. Let us even assume that
$$
a=(1,0,0) \ \hbox{  and } \ b=(\cos T,\sin T,0) \hbox{  for some }
T \in [8\eta_0, 10\pi/11],
\leqno (8.1)
$$
where again $\eta_0 > 0$ is as in (2.2).
We avoid the case when $a$ and $b$ are nearly 
antipodal for the same sort of reasons as in Section 7.
We also assume that $\Gamma$ is a Lipschitz graph with
constant at most $\eta$, by which we mean that we can find an
$\eta$-Lipschitz function $v : [0,T] \to \R^{n-2}$, with
$v(0)= v(T) = 0$, such that if we set
$$
\hbox{$w(t) = (1-|v(t)|^2)^{1/2}$ and }
z(t) = \big(\cos t \, w(t), \sin t \, w(t),v(t)\big)
\ \hbox{ for $t\in [0,T]$,}
\leqno (8.2)
$$ 
$z$ is a parameterization of $\Gamma$ by $[0,T]$.

\medskip\noindent{\bf Remark 8.3.}
The curve $\Gamma$ of Section 7 satisfies these requirements. Indeed,
$\length(\Gamma) \leq \length(\gamma) \leq 10\pi/11$
by (7.31) and (7.1); $\length(\Gamma) \geq \dist_{\partial B}(a,b)
\geq \length(\gamma) - \tau_1 \geq 9\eta_0-\tau_1 \geq 8\eta_0$ 
by (7.2), (7.1), and if $\tau$ is small enough; finally
the representation (8.2) follows from (7.43) and (7.38): 
set $\widetilde\theta(t) = u$ in (7.38) (this change of variable is all right 
because $\widetilde\theta$ is strictly increasing), and observe that 
$v(u) = \widetilde v(t)$ is an $\eta$-Lipschitz function of
$u$, by (7.43).

\medskip
As in Section 7, we can take $\eta$ as small as we want, and we shall use 
this smallness to approximate minimal surfaces with graphs of harmonic functions. 
We shall try to estimate various gains in terms of $\length(\Gamma)-T$ 
and the equivalent quantity $||v'||_2^2$.

Define a homogeneous function $F$ on a sector with aperture $T$ by
$$
F(r \cos t, r \sin t) = {r \, v(t) \over w(t)}
\ \hbox{ for $r \geq 0$ and $t \in [0,T]$.}
\leqno (8.4)
$$
Also set 
$$
D_T = \{ (r\cos t, r\sin t) \, ; \, r \in (0,1) \hbox{ and }
t \in (0,T) \}, 
\leqno (8.5)
$$
and then denote by $\Sigma_F$ the graph of $F$ over $\overline D_T$.
Notice that $\Gamma \i \Sigma_F$, by (8.2) and (8.4).

We want to find a new function $G$, also defined on $\overline D_T$,
such that
$$
G(r \cos t, r \sin t) = F(r \cos t, r \sin t)
\ \hbox{ for $t \in \{ 0,T \}$, and for $9/10 \leq r \leq 1$,}
\leqno (8.6)
$$
and whose graph $\Sigma_G$ has a smaller area. 
[We required that $G(r \cos t, r \sin t) = F(r \cos t, r \sin t)$
for $r \geq 9/10$ to make sure that we do not change anything
near $\Gamma$.] The natural thing to do is to take the harmonic 
function with the same boundary values as $F$, but we shall modify 
this a little near $\partial B$ (to keep it Lipschitz) and
near the center, to allow further modifications.
In fact, we shall require that
$$
G(r \cos t, r \sin t) = 0
\ \hbox{ for $0 \leq r \leq 2\kappa$ and $0 \leq t \leq T$,}
\leqno (8.7)
$$
where $\kappa$ is a small positive constant that does not even depend 
on $n$.

\ms\proclaim Lemma 8.8.
We can find a Lipschitz function $G : \overline D_T \to \R^{n-2}$
such that (8.6) and (8.7) hold,
$$
|\nabla G(z)| \leq C ||v'||_{\infty} \leq C \eta
\ \hbox{ for almost-every } z\in D_T,
\leqno (8.9)
$$
and 
$$
H^{2}(\Sigma_F) - H^{2}(\Sigma_G)
\geq 10^{-4}  \int_0^T |v'(t)|^2 dt
\geq 10^{-4} [\length(\Gamma) - T].
\leqno (8.10)
$$

\ms
A good part of the proof is common with the paper [D3], 
which we shall quote for some of the computations. Set 
$$
f(t) = F(\cos t,\sin t)= {v(t) \over w(t)}
= v(t) (1-|v(t)|^2)^{-1/2}
\ \hbox{ for $0 \leq t  \leq T$.}
\leqno (8.11)
$$
Note that $||v||_\infty \leq \eta T$ because
$v$ is $C\eta$-Lipschitz and vanishes at $0$, hence
$w(t) \geq 1 - \eta T$ for $0 \leq t  \leq T$, and then
$$
\hbox{$f$ is $2 ||v'||_\infty$-Lipschitz, and
$f(0)=f(T) = 0$.}
\leqno (8.12)
$$

A first attempt for $G$ is to take the harmonic extension
of $f$ to $D_T$ (with Dirichlet conditions on the boundary),
which we define as follows. First we use the fact that 
$f(0) = f(T) = 0$ to write $f$ as a sum of sines, i.e.,
$$
f(t) = \sum_{k \geq 1} \beta_k \sin(\pi k t / T)
\leqno (8.13)
$$
for $0 \leq t  \leq T$. (See (13.4) in [D3].) 
Notice that
$$
{\pi^2 \over 2T} \, \sum_{k \geq 1} k^2 |\beta_k|^2 
= \int_{0}^T  |f'(t)|^2 dt < 4 \eta^2 T
\leqno (8.14)
$$
(see (13.5) in [D3] and use (8.12)), 
so $\sum_{k \geq 1} |\beta_k| < +\infty$ by Cauchy-Schwarz, and the series 
in (8.13) even converges pointwise. We set
$$
G_1(\rho \cos t, \rho \sin t) =
\sum_{k \geq 1} \beta_k  \rho^{\pi k / T} \sin(\pi k t / T)
\leqno (8.15)
$$
for $t \in \{ 0, T \}$ and $0 \leq \rho \leq 1$, as in (13.7) in [D3], 
and where we also get the normal convergence of the series from (8.14).
Then $G_1(\cos t, \sin t) = f(t)$, so
$$
G_1(z) = F(z) \hbox{ on $\overline D_T \cap \partial B(0,1)$,}
\leqno (8.16)
$$
and 
$$
G_1(\rho \cos t, \rho \sin t) = 0
\ \hbox{ for $t \in \{ 0, T \}$ and $0 \leq \rho \leq 1$}
\leqno (8.17)
$$
trivially. Direct computations using the expansions in (8.13) 
and (8.15) and the orthogonality of sines and cosines show that
$$
{1 \over 2} \int_0^T |f'(t)|^2 dt
\leq \int_{D_T} |\nabla F|^2
\leq \int_0^T |f'(t)|^2 dt,
\leqno (8.18)
$$
as in [D3], Lemma 13.9, and  
$$
\int_{D_T} |\nabla G_1|^2 
\leq {2T/\pi \over 1+(T/\pi)^2} \int_{D_T} |\nabla F|^2
\leq {220 \over 221} \int_{D_T} |\nabla F|^2
\leqno (8.19)
$$
as in (13.19) in [D3], 
and where the last inequality holds because $T \leq 10\pi/11$
and ${\lambda \over (1+\lambda^2)}$ is an increasing function
of $\lambda \in [0,1]$. 

Incidentally, it is fortunate for us that even when $n > 3$, 
what we need later is really 
$|\nabla G_1|^2 =\big| {\partial G_1 \over \partial x} \big|^2
+\big| {\partial G_1 \over \partial y} \big|^2$, which is easier to 
compute, rather than for instance the square of the operator norm of 
the differential $DG_1$ (acting from $\R^2$ to $\R^{n-2}$), 
which may be different.
The reader may check that in polar coordinates, one may also compute
$|\nabla G_1|^2$ as $\big| {\partial G_1 \over \partial r} \big|^2
+\big| {1 \over r} {\partial G_1 \over \partial \theta} \big|^2$.

\smallskip
We unfortunately need some additional control on $\nabla G_1$. 
The radial and tangential derivatives of $G_1$ at 
$z = (\rho \cos t, \rho \sin t)$ are
$$
{\partial G_1 \over \partial \rho}
(\rho \cos t, \rho \sin t) 
= \sum_{k \geq 1} \beta_k \, {\pi k \over T \rho} 
\, \rho^{k\pi /T} \sin(\pi k t/T)
\leqno (8.20)
$$
and
$$
{1 \over \rho} {\partial G_1 \over \partial t} 
(\rho \cos t, \rho \sin t) = 
\sum_{k \geq 1} \beta_k \, {\pi k \over T \rho} 
\, \rho^{k\pi /T} \, \cos(\pi k t/T), 
\leqno (8.21) 
$$
which yields
$$\eqalign{
\int_0^T |\nabla G_1(\rho \cos t, \rho \sin t)|^2 dt
&= \int_0^T \Big|{\partial G_1 \over \partial \rho} \Big|^2 dt
+
\int_0^T \Big| {1 \over \rho} {\partial G_1 \over \partial t} \Big|^2 dt
\cr& \hskip -2cm
= T \sum_{k \geq 1} \Big| \beta_k \, {\pi k \over T \rho} \Big|^2
\, \rho^{2k\pi /T}
= {\pi^2 \over T} \sum_{k \geq 1} k^2 |\beta_k|^2 \, 
\rho^{-2} \, \rho^{2k\pi /T}
}\leqno (8.22) 
$$
as in [D3], (13.13) and its analogue for the tangential derivative. 
Since $T < \pi$, $\rho^{-2} \, \rho^{2k\pi /T} \leq 1$ for $\rho \leq 1$
and we immediately get that
$$
\int_0^T |\nabla G_1(\rho \cos t, \rho \sin t)|^2 dt
\leq {\pi^2 \over T} \sum_{k \geq 1} k^2 |\beta_k|^2
= 2 \int_0^T |f'(t)|^2 dt
\leqno (8.23)
$$
by (8.14). We can also get a pointwise estimate for $\rho < 1$
by adding (8.20) to (8.21) and applying Cauchy-Schwarz. We get that
$$\leqalignno{
|\nabla G_1(\rho \cos t, \rho \sin t)|  
&\leq {2 \pi \over T \rho} \sum_{k \geq 1} k \, |\beta_k| \, \rho^{k\pi /T}
\leq C \rho^{-1} \Big\{ \sum_{k \geq 1} k^2 |\beta_k|^2 \Big\}^{1/2}
\Big\{ \sum_{k \geq 1} \rho^{2k\pi /T} \Big\}^{1/2}
\cr&
\leq C \rho^{-1} ||f'||_2 \,\Big\{ \sum_{k \geq 1} \rho^{2k\pi /T} \Big\}^{1/2}
= C \rho^{-1} ||f'||_2 \, 
\Big\{ { \rho^{2\pi/T} \over 1 - \rho^{2\pi/T}} \Big\}^{1/2}
& (8.24)
\cr&
\leq C \rho^{{\pi \over T} -1} (1-\rho)^{-1/2} ||f'||_2 
\leq C (1-\rho)^{-1/2} ||f'||_2 
}
$$
by (8.14) and because $T \leq \pi$. 
The estimate diverges slightly near $\rho = 1$,
and unfortunately this is not surprising, because the harmonic 
extension of a Lipschitz function on the unit disk is sometimes
a little less than Lipschitz.

To ameliorate this, we pick a radius $r\in (0,1)$ and 
use a different function $G_2$ in the annular region 
$A = \overline D_T \cap \overline B(0,1) \sm B(0,r)$. 
That is, we keep
$$
G_2(z) = G_1(z) \ \hbox{ for } z\in \overline D_T \cap \overline B(0,r)
\leqno (8.25)
$$ 
but take
$$
G_2(\rho \cos t, \rho \sin t) = 
{ 1-\rho \over 1-r} \, G_1(r \cos t, r \sin t)
+ { \rho-r \over 1-r} \, G_1(\cos t, \sin t)
\leqno (8.26)
$$
when $(\rho \cos t, \rho \sin t) \in A$.
[We interpolate linearly the values of $G_1$ on the
circles $\partial B(0,r)$ and $\partial B(0,1)$.] Note that
the two definitions coincide on $\partial B(0,r)$, and that
$$
G_2(z) = G_1(z) = F(z) 
\ \hbox{ for } z\in \overline D_T \cap \partial B(0,1)
\leqno (8.27)
$$
by (8.16). The radial derivative on $A$ is such that
$$\eqalign{
\Big|{\partial G_2 \over \partial \rho}(\rho \cos t, \rho \sin t)\Big|
& = { 1 \over 1-r} \, 
\big| G_1(\cos t, \sin t)-G_1(r \cos t, r \sin t) \big|
\cr&
\leq  { 1 \over 1-r} \, \int_r^1 | \nabla G_1(x\cos t, x\sin t)| dx
\cr&
\leq C (1-r)^{-1} ||f'||_2 \int_r^1 (1-x)^{-1/2}dx
\leq C (1-r)^{-1/2} ||f'||_2
}\leqno (8.28)
$$
by (8.24). For the derivative in the tangential direction,
$$\eqalign{
\Big| {1 \over \rho} \,
{\partial G_2 \over \partial t}(\rho \cos t, \rho \sin t)\Big|
&\leq 
{r \over \rho} \, 
\Big| {\partial G_1 \over \partial t}(r \cos t, r \sin t) \Big|
+ {1\over \rho} \, 
\Big| {\partial G_1 \over \partial t}( \cos t, \sin t) \Big|
\cr& 
\leq  | \nabla G_1(r \cos t, r \sin t)| +  2|f'(t)|
\cr&
\leq C (1-r)^{-1/2} ||f'||_2 +  2|f'(t)|
}\leqno (8.29)
$$
because $G_1(\cos t, \sin t) = f(t)$ (by (8.11) and (8.27))
and by (8.24). Altogether, 
$$
|\nabla G_2
(\rho \cos t, \rho \sin t)|
\leq C (1-r)^{-1/2} ||f'||_2 + 2 |f'(t)| 
\leq C (1-r)^{-1/2} ||v'||_\infty
\leqno (8.30)
$$
on $A$, because $f$ is $2||v'||_\infty$-Lipschitz (by (8.12)). 
Recall from (8.25) that $G_2=G_1$ on $\overline D_T \cap \overline B(0,r)$;
there $|\nabla G_2(\rho \cos t, \rho \sin t)| \leq C (1-r)^{-1/2} ||f'||_2$
directly by (8.24), and so
$$
G_2 \hbox{ is $C (1-r)^{-1/2} ||v'||_\infty$-Lipschitz on 
$\overline D_T$,}
\leqno (8.31)
$$
because $G_2$ is continuous across $\partial B(0,r)$.
We do not mind the factor $(1-r)^{-1/2}$ much here,
because we shall soon take $r = 1- 10^{-6}$.

We also need to estimate $\displaystyle \int_A |\nabla G_2(z)|^2 dz$.
First,
$$\eqalign{
\int_A \Big|{\partial G_2 \over \partial \rho}(z)\Big|^2 dz
&\leq (1-r)^{-2} \int_{\rho = r}^1 \int_{t=0}^T 
\Big\{ \int_{x=r}^1 | \nabla G_1(x\cos t, x\sin t)| dx \Big\}^2 
\rho d\rho dt
\cr&
\leq (1-r)^{-1} \int_{t=0}^T  \Big\{
\int_{x=r}^1 | \nabla G_1(x\cos t, x\sin t)| dx \Big\}^2 dt
\cr&
\leq \int_{t=0}^T  \int_{x=r}^1
| \nabla G_1(x\cos t, x\sin t)|^2 dx dt
\leq 2 (1-r) \int_0^T |f'(t)|^2 dt
}\leqno (8.32)
$$
by the first part of (8.28), because the inside integral does not
depend on $\rho$, by Cauchy-Schwarz, Fubini, and (8.23)
(with $\rho =  x$ and integrated on $(r,1)$).
Similarly,
$$\eqalign{
\int_A \Big|{1 \over \rho} \, {\partial G_2 \over \partial t}(z)
\Big|^2 dz
&\leq  \int_{\rho = r}^1 \int_{t=0}^T 
\Big\{ | \nabla G_1(r \cos t, r \sin t)| + 2|f'(t)| \Big\}^2
\rho d\rho dt
\cr&
\leq  (1-r) \int_{t=0}^T 
\Big\{ | \nabla G_1(r \cos t, r \sin t)| + 2|f'(t)| \Big\}^2 dt
\cr&
\leq  2(1-r) 
\Big\{ \int_{t=0}^T |\nabla G_1(r \cos t, r \sin t)|^2 dt 
+ 4 \int_{t=0}^T |f'(t)|^2  dt \Big\}
\cr&
\leq 12 (1-r) \int_0^T |f'(t)|^2 dt
}\leqno (8.33)
$$
by (8.29), because the integrand does not depend on $\rho$,
and by (8.23) (with $\rho=r$). Altogether,
$$
\int_A |\nabla G_2(z)|^2 dz 
\leq 14 (1-r) \int_0^T |f'(t)|^2 dt
\leq 10^{-4} \int_0^T |f'(t)|^2 dt
\leqno (8.34)
$$
by (8.32) and (8.33), and because we take $r = 1 - 10^{-6}$.

\smallskip
We are now ready to define $G$. We keep 
$$
G(z) = F(z) \ \hbox{ for $z \in \overline D_T \setminus B(0,9/10)$,}
\leqno (8.35)
$$
as suggested by (8.6), and set
$$
G(z) = {9 \over 10} \, G_2 \big({10z\over 9}\big) 
\ \hbox{ for } z \in \overline D_T \cap B(0,9/10)\setminus 
B(0,3\kappa),
\leqno (8.36)
$$
where the small constant $\kappa$ will be chosen soon.
Notice that when $z \in \partial B(0,9/10)$,
$F(z) = {9 \over 10} \, F \big({10z\over 9}\big)
= {9 \over 10} \, G_1 \big({10z\over 9}\big)
={9 \over 10} \, G_2 \big({10z\over 9}\big)$
because $F$ is homogeneous of degree $1$ and by (8.27).
So our two definitions match on $\partial B(0,9/10)$. 
Recall from (8.7) that we need to take
$$
G(z) = 0 \ \hbox{ on } \overline D_T \cap B(0,2\kappa).
\leqno (8.37)
$$
The simplest way to make a continuous transition is to 
interpolate linearly (as before) and take
$$\eqalign{
G(\rho \cos t, \rho \sin t) &= 
{\rho - 2\kappa \over \kappa} \, G(3\kappa \cos t, 3\kappa \sin t) 
\cr& 
= {9 \over 10} \, {\rho - 2\kappa \over \kappa} \, 
G_1 \Big({30\kappa\over 9} \cos t, {30\kappa\over 9} \sin t \Big) 
}\leqno (8.38)
$$
for $0 \leq t \leq T$ and $2\kappa \leq \rho \leq 3\kappa$.
[We used (8.36) and (8.25) here.]

\ms
This completes our definition of $G$; let us now check that
it satisfies the requirements of Lemma 8.8. 

We need to know that $G(\rho \cos t, \rho \sin t) 
= F(\rho \cos t, \rho \sin t) = 0$ when $t \in \{ 0, T \}$,
to complete our proof of (8.6). But for such $t$,
$F(\rho \cos t, \rho \sin t) = 0$ directly by (8.4) and because
$v(0)=v(T)=0$ (see above (8.2)). Next, 
$G_1(\rho \cos t, \rho \sin t)  = 0$ by (8.17), 
and then $G_2(\rho \cos t, \rho \sin t)  = 0$ by the definition 
(8.25) or (8.26).
Then $G(\rho \cos t, \rho \sin t)  = 0$, because it is obtained
from $F$ or $G_2$ by (8.35), (8.36), (8.37) or (8.38).

Since (8.7) holds by (8.37), we turn to the Lipschitz condition (8.9). 
Since $G$ is continuous across $\partial B(0,9/10)$, 
$\partial B(0,3\kappa)$, and $\partial B(0,2\kappa)$, 
we just need to investigate the remaining
annular domains separately. Out of $B(0,9/10)$, 
$G = F$ by (8.35), so it is enough to check that
$$
F \hbox{ is $5||v'||_\infty$-Lipschitz on } \overline D_T.
\leqno (8.39)
$$
Recall that $v$ is $\eta$-Lipschitz and vanishes at the origin, so
$|v(t)| \leq ||v'||_\infty T \leq \eta T$, $w(t) \geq 1-\eta T$ by (8.2),
and $|w'(t)| \leq 2 ||v'||_\infty$, again by (8.2). Then (8.4) yields
$$\eqalign{
|\nabla F(r \cos t, r \sin t)|
&\leq \Big|{\partial F \over \partial r}(r \cos t, r \sin t)\Big|
+ {1 \over r} \Big|{\partial F \over \partial t}(r \cos t, r \sin t)\Big| 
\cr& 
\leq{|v(t)| \over w(t)}+{|v'(t)| \over w(t)}+{|v(t)||w'(t)|\over w(t)^2}
\leq 5 ||v'||_\infty
}\leqno (8.40)
$$
if $\eta$ is small enough, and as needed.

Return to $G$. In $B(0,9/10)\sm  B(0,3\kappa)$, we use (8.36) and (8.31). 
Since $G = 0$ in $B(0,2\kappa)$, we are just left with 
$B(0,3\kappa) \sm  B(0,2\kappa)$, where (8.38) yields
$$
\left| {1 \over \rho} \,{\partial G \over \partial t}
(\rho \cos t, \rho \sin t)\right|
\leq {3 \kappa \over \rho} \, 
\Big| \nabla G_1 \Big({30\kappa\over 9} \cos t, 
{30\kappa\over 9} \sin t \Big) \Big| 
\leq C ||f'||_2 \leq C ||v'||_\infty
\leqno (8.41)
$$
by (8.24) and (8.12), and 
$$
\left| {\partial G \over \partial \rho}(\rho \cos t, \rho \sin t)\right|
\leq \kappa^{-1} 
\Big| G_1 \Big({30\kappa\over 9} \cos t, {30\kappa\over 9} \sin t \Big) \Big| 
\leq C ||f'||_2 \leq C ||v'||_\infty
\leqno (8.42)
$$
this time by (8.17) and (8.24). So (8.9) holds.

\smallskip
We are left with the verification of (8.10). 
Since $F=G$ out of $B(0,9/10)$ by (8.35), 
$$
H^{2}(\Sigma_F) - H^{2}(\Sigma_G)
= H^{2}(\Sigma'_F) - H^{2}(\Sigma'_G),
\leqno (8.43)
$$
where $\Sigma'_F$ is the graph of the restriction of $F$
to $D' = D_T \cap B(0,9/10)$, and similarly for $\Sigma'_G$.
By the area formula, 
$$
H^{2}(\Sigma'_G) = \int_{D'} J(z) dz,
\leqno (8.44)
$$
where $J$ is the jacobian of the parameterization
$H: z \to (z,G(z))$, which we compute now. Fix $z\in D'$ where
$G$ is differentiable, denote by $e_1$ and $e_2$ the first  two vectors 
of the canonical basis of $\R^n$ (so $e_1$ and $e_2$ form a basis of 
the horizontal plane $P$), and set 
$v_1 = {\partial G \over \partial x}(z)$
and $v_2 = {\partial G \over \partial y}(z)$. Thus, if we set
$w_1 = DH(z)(e_1)$ and $w_2 = DH(z)(e_2)$, we get that
$w_i = e_i+v_i$ for $i=1,2$.
Now $J(z) = |w_1||w_2| \sin\alpha$, where $\alpha \in [0,\pi/2]$ is the
angle of $w_1$ with $w_2$, so 
$$\eqalign{
J(z)^2 &= |w_1|^2 |w_2|^2 \sin^2\alpha
= |w_1|^2 |w_2|^2 [1-\cos^2\alpha]
=  |w_1|^2 |w_2|^2 - \langle w_1,w_2 \rangle^2
\cr& \hskip-0.1cm
= (1+|v_1|^2)(1+|v_2|^2) - \langle v_1,v_2 \rangle^2
= 1 + |v_1|^2 + |v_2|^2 + |v_1|^2|v_2|^2 - \langle v_1,v_2 \rangle^2
\cr& 
\leq 1 + |v_1|^2 + |v_2|^2 + |v_1|^2|v_2|^2
\cr&
\leq 1 + |\nabla G(z)|^2 + |\nabla G(z)|^4 
\leq 1 + (1+C\eta) |\nabla G(z)|^2
}\leqno (8.45)
$$
where the last inequality  comes from (8.9).
Recall that $(1+u)^{1/2} \leq 1 + {u \over 2}$ for $u \geq 0$. 
Applying this with $u = J(z)^2 - 1$ yields
$$\eqalign{
H^{2}(\Sigma'_G) 
&= \int_{D'} J(z) \, dz = \int_{D'} (1+u)^{1/2} dz
\leq \int_{D'} (1+ {u \over 2}) \, dz
\cr&
= H^2(D') + {1 \over 2} \int_{D'} u \, dz
\leq H^2(D') + {(1+C\eta) \over 2} \int_{D'} |\nabla G(z)|^2 dz
}\leqno (8.46)
$$
by (8.44) and (8.45). Next
$$\eqalign{
\int_{D'} |\nabla G|^2 
&\leq \int_{D' \cap B(0,3\kappa)} |\nabla G|^2
+ \int_{D' \setminus B(0,3\kappa)} |\nabla G|^2
\cr&
\leq C \int_{D' \cap B(0,3\kappa)} ||f'||_2^2
+ \int_{D' \setminus B(0,3\kappa)} 
\big|\nabla\big[{9 \over 10} \, G_2(10z/9)\big]\big|^2
\cr&
\leq C \kappa^2 ||f'||_2^2 + {81 \over 100} \int_{D_T} |\nabla G_2|^2
}\leqno (8.47)
$$
by (8.41), (8.42), (8.36), and a linear change of variable.
We choose $\kappa$ so small that $C \kappa^2 \leq 10^{-4}$
in (8.47), and since
$$
\int_{D_T} |\nabla G_2|^2 
\leq \int_{D_T} |\nabla G_1|^2 + \int_{A} |\nabla G_2|^2 
\leq {220 \over 221} \int_{D_T} |\nabla F|^2
+ 10^{-4} ||f'||_2^2 
\leqno (8.48)
$$
by (8.25), (8.19) and (8.34), we get that
$$
\int_{D'} |\nabla G|^2 
\leq {81 \over 100} \, {220 \over 221} \int_{D_T} |\nabla F|^2
+ 2 \cdot 10^{-4} \, ||f'||_2^2 \, .
\leqno (8.49)
$$

We also need a lower bound for $H^{2}(\Sigma'_F)$. 
This time the area formula yields
$$
H^{2}(\Sigma'_F) = \int_{D'} J_F(z) dz
\geq \int_{D'} \big\{ 1 + |\nabla F(z)|^2 \big\}^{1/2}
\leqno (8.50)
$$
where the lower bound comes from the first two lines of (8.45).
Recall from (8.39) that $F$ is $5\eta$-Lipschitz
on $\overline D_T$. Also,
$$
(1+u)^{1/2} \geq 1 + {u \over 2} - {u^2 \over 8}
\geq 1 + {u \over 2} - {25 \eta^2 u \over 8}
\geq 1 + {u \over 2} \, (1-7\eta^2)
\leqno (8.51)
$$
for $0 \leq u \leq 25 \eta^2$; we apply this with
$u = |\nabla F(z)|^2$, integrate, and get that
$$
H^{2}(\Sigma'_F) \geq \int_{D'} \big\{ 1 + |\nabla F|^2 \big\}^{1/2}
\geq H^2(D') + {1-7 \eta^2 \over 2} \int_{D'}  |\nabla F|^2.
\leqno (8.52)
$$
Notice that $\int_{D'}  |\nabla F|^2 = {81 \over 100} \int_{D_T}  |\nabla F|^2$
because $F$ is homogeneous of degree 1; then 
$$\eqalign{
H^{2}(\Sigma_F) - &H^{2}(\Sigma_G)
= H^{2}(\Sigma'_F) - H^{2}(\Sigma'_G)
\cr&
\geq {1-7 \eta^2 \over 2} \int_{D'} |\nabla F|^2
- {1+C\eta \over 2} \int_{D'} |\nabla G|^2 
\cr&
\geq {81 \over 200} \, \Big[1-7 \eta^2 - (1+C\eta )\,{220 \over 221}\Big] 
\int_{D_T} |\nabla F|^2 - (1+C\eta ) 10^{-4}  ||f'||_2^2
\cr&
\geq {81 \over 200} \, {1 \over 300} \int_{D_T} |\nabla F|^2
- (1+C\eta ) 10^{-4} ||f'||_2^2
\cr&
\geq  \Big[ {81 \over 12 \cdot 10^4}-(1+C\eta )10^{-4} \Big] ||f'||_2^2
\geq 5 \cdot 10^{-4} ||f'||_2^2
}\leqno (8.53)
$$
by (8.43), (8.52), (8.46), (8.49), if $\eta$ is small enough,
and by (8.18). Since $f = v (1-|v|^2)^{-1/2}$ by (8.11) and 
$|v| \leq \eta T$ everywhere, we easily get that
$||f'||_2^2 \geq {1 \over 2} ||v'||_2^2$, and the first part
of (8.10) follows from (8.53).
So our proof of Lemma 8.8 will be complete as soon as we check that
$\length(\Gamma) - T \leq \int_0^T |v'(t)|^2 dt$

We differentiate in (8.2) or lazily use (7.11) with $\theta'=1$,
and get that  $|z'(t)|^2 = w(t)^2 + |w'(t)|^2 + |v'(t)|^2$.
Recall that $ww' = - \langle v,v' \rangle$ because
$w^2  + |v|^2 =1$, so $|w'| \leq |v||v'|/w$ and
$|z'(t)|^2 \leq 1 + |v'(t)|^2 [1+|v(t)|^2/w(t)^2] \leq 1 + 2 |v'(t)|^2$
(recall that $|v(t)| \leq \eta T$ and hence $w(t) \geq 1-\eta T$).
Hence $|z'(t)| \leq 1 + |v'(t)|^2$ (recall that
$(1+u)^{1/2} \leq 1 + u/2$ for $u \geq 0$), and
$$
\length(\Gamma) = \int_0^T |z'(t)| dt
\leq T + \int_0^T |v'(t)|^2 dt,
\leqno (8.54)
$$
as needed. Lemma 8.8 follows.
\qed

\bigskip 
\noindent {\bf 9. Retractions near Lipschitz graphs and the 
construction of a first competitor}
\medskip

We shall now use the previous sections to construct a first 
competitor for $E$ in $B(0,1)$. We still assume that (6.2) holds
(otherwise Lemma 6.3 gives the desired estimates), and we 
found in Section 6 a collection of curves $g_{j,k}$, 
$(j,k) \in \widetilde J$, with the properties described in Lemma 6.11.
Recall in particular that the extremities of $g_{j,k}$ are
$\varphi(a)$ and $\varphi(b)$, where $a$ and $b\in V$ are the 
endpoints of the corresponding $\C_{j,k}$.

Denote by $D_{j,k}$ the cone over the geodesic that goes from $a$
to $b$. Thus $D_{j,k}$ is a plane sector bounded by the half lines
$[0,a)$ and $[0,b)$. Also denote by $P_{j,k}$ the plane that contains
$D_{j,k}$.

Next apply to each arc $g_{j,k}$ the construction of Section 7.
More precisely, we should first apply a rotation ${\cal R}$ such that
${\cal R}(P_{j,k})$ is the horizontal plane $P$ of Section 7, and apply
the construction of Section 7 to $\gamma = {\cal R}(g_{j,k})$, but, in
order to save notation, we shall often pretend that ${\cal R}$ is 
the identity. Let us check that $\gamma$ satisfies our 
assumptions (7.1)-(7.3). 

First, the geodesic distance in $\partial B$ from $\varphi(a)$ to 
$\varphi(b)$ is $\dist_{\partial B}(\varphi(a),\varphi(b))
\geq \dist_{\partial B}(a,b) - C \tau \geq 9 \eta_0$ by
(6.14) and (2.5), so $\length(\gamma) \geq 9 \eta_0$.
On the other hand,
$\length(\gamma) = \length(g_{j,k}) = H^1 (g_{j,k})
\leq \length(\C_{j,k}) + C_4 \tau < {10 \pi \over 11}$ 
because $g_{j,k}$ is simple, by (6.43), and by (2.5), so (7.1) 
holds. We also get (7.2), because 
$\length(\gamma) \leq \length(\C_{j,k}) + C_4 \tau 
= \dist_{\partial B}(a,b) + C_4 \tau 
\leq \dist_{\partial B}(\varphi(a),\varphi(b)) + C \tau$
by (6.14) and if $\tau$ is small enough compared to $\tau_1$.
[See the discussion above (7.4) concerning the order in which we 
choose the constants.]
Finally, (7.3) follows from (6.15).

So we can apply the construction of Section 7 to $g_{j,k}$,
and we get a curve $\Gamma_{j,k}$ in $\partial B(0,1)$, which
is an $\eta$-Lipschitz graph over $P_{j,k}$ and has a large 
intersection with $g_{j,k}$. Denote by $\Gamma_{j,k}^\ast$
the cone over $\Gamma_{j,k}$. Then $\Gamma_{j,k}^\ast$
is the graph of a homogeneous Lipschitz function defined on
$D_{j,k}$. As was observed in Remark 8.3, we can apply the 
construction of Section 8 to the curve 
$\Gamma ={\cal R}(\Gamma_{j,k})$, where ${\cal R}$ is the same 
rotation as above, which sends $P_{j,k}$ to $P$. We get a 
new surface $\Sigma_{j,k}$, with the properties described near
Lemma 8.8. In particular, $\Sigma_{j,k}$ is the  graph over
$P_{j,k}$ of a $(C \eta)$-Lipschitz function $G_{j,k}$ defined on 
$D_{j,k}$ (see (8.9)), and since $\Gamma_{j,k}^\ast$ is the
graph of the homogeneous function $F$ in in (8.4)
(see the remark below (8.5)), 
$$
\hbox{$\Sigma_{j,k}$ coincides with $\Gamma_{j,k}^\ast$
out of $B(0,19/20)$}
\leqno (9.1)
$$
by (8.6) and because $F$ and $G_{j,k}$ are both Lipschitz 
with small constants.
Also recall that $\Sigma_{j,k}$ and $\Gamma_{j,k}^\ast$ are both
bounded by the two half lines through $\varphi(a)$ and $\varphi(b)$
that bound $D_{j,k}$ (see (8.6), (8.4), and the line above (8.2)), 
and that 
$$
H^{2}(\Sigma_{j,k} \cap B) 
\leq H^{2}(\Gamma_{j,k}^\ast \cap B) - 10^{-4}
[\length(\Gamma_{j,k}) - \dist_{\partial B}(\varphi(a),\varphi(b))],
\leqno (9.2)
$$
by (8.10) and the definition (8.1). Set
$$
g = \bigcup_{(j,k) \in \widetilde J} \ g_{j,k} \,  , \ 
\Gamma = \bigcup_{(j,k) \in \widetilde J} \ \Gamma_{j,k}  \, , \ 
\Gamma^\ast = \bigcup_{(j,k) \in \widetilde J} \ \Gamma_{j,k}^\ast  
\,  , \
\hbox{ and } 
\Sigma = \bigcup_{(j,k) \in \widetilde J} \ \Sigma_{j,k}  \, .
\leqno (9.3)
$$
We shall need to define a projection of a neighborhood
of $\Sigma$ onto $\Sigma$; which will allow us later to deform 
$E \cap \overline B(0,r)$ onto a subset of $\Sigma$. 

\ms\proclaim Lemma 9.4.
There is a neighborhood $U$ of $E\cap \overline B$ in $\overline B$
and a $50$-Lipschitz mapping $p_2 : U \to \Sigma \cap \overline B$,
such that 
$$
p_2(z) = z \ \hbox{ for } z\in \Sigma \cap \overline B \sm B(0,1/3)
\leqno (9.5)
$$ 
and 
$$
p_2(z) \in \Sigma \cap B(0,1/2) \ \hbox{ for } z\in U \cap B(0,1/3).
\leqno (9.6)
$$

\ms
We shall start with a first projection on a simpler set $R$. Set
$$
R_{j,k} = \big\{ z\in \R^n \, ; \, 
\dist(z,D_{j,k}) \leq 10^{-2} \eta_0 \dist(z,\partial D_{j,k})\big\}
\ \hbox{ and } \ 
R = \bigcup_{(j,k) \in \widetilde J} \ R_{j,k}  \, ,
\leqno (9.7)
$$
where $\partial D_{j,k}$ is the union of the two half lines that
bound $D_{j,k}$ (the half lines through $\varphi(a)$ and $\varphi(b)$,
with the notation above), and $\eta_0$ is as in Section 2. 
Thus the $R_{j,k}$ are small conic sectors
around the $D_{j,k}$. Observe that
$$
\Sigma_{j,k} \i R_{j,k}
\leqno (9.8)
$$
rather trivially, because $\Sigma_{j,k}$ is a small Lipschitz
graph over $D_{j,k}$ bounded  by the two half lines of $\partial D_{j,k}$.

We shall first project a small neighborhood of $\Sigma \cap \partial B$
onto $R \cap \partial B$. Set 
$$
W = \big\{ z\in \partial B \, ; \, \dist(z,V) < C_6 \tau \big\},
\leqno (9.9)
$$
where $V$ is our set of vertices from Section 2, and $C_6$
will be chosen soon, somewhat larger than the constants in Lemma 6.11.

We shall only need to define the first projection $p_0$ on 
$[R \cap \partial B] \cup W$, and on $R \cap \partial B$ we just take
$p_0(z)=z$. So we just need to take care of the $B_x = \partial B 
\cap B(\varphi(x),C_6 \tau)$, $x\in V$.
Observe that the $B_x$ are far from each other, by (2.5) and (2.6), so 
we can treat them separately. 

Fix $x \in V$. There are two or three arcs $\C_{j,k}$ that end at $x$.
For all the other arcs $\C_{j',k'}$, (2.6) says that $\C_{j',k'}$
lies further than $\eta_0/2$ away from $B_x$ (if $\tau$ is small 
enough), and then (6.14) says that $D_{j',k'}$ still lies further than 
$\eta_0/3$ away from $B_x$ too. So $R \cap B_x$ is simply composed of 
the two or three $R_{j,k}$ for which $x$ is an endpoint of $\C_{j,k}$.

Fix such a pair $(j,k)$. Near $B_x$, $R_{j,k}$ is a thin conical 
neighborhood in $\partial B$ of the geodesic $\rho_{j,k}$ that connects 
$\varphi(x)$ to some other point of $\varphi(V)$, and $p_0$ will
be easier to define first on a larger conical neighborhood $A_{j,k}$ of
$\rho_{j,k}$. For $z\in B_x \sm \{ \varphi(x) \}$, 
denote by $\theta_{j,k}(z)\in [0,\pi]$ the angle of $z-\varphi(x)$ 
with the tangent at $\varphi(x)$ of $\rho_{j,k}$ (pointing away from 
$\varphi(x)$). Then set $A_{j,k} = \big\{ z\in B_x \sm \{ \varphi(x) \} 
\, ; \, \theta_{j,k}(z) \leq {99 \over 100} \, {\pi \over 3} \big\}$.

We define $p_0$  on $A_{j,k}$ so that $p_0(z) = z$ when $z\in R_{j,k}$,
$p_0(A_{j,k}) \i B_x \cap R_{j,k}$, and $p_0(z) = \varphi(x)$ when
$\theta_{j,k}(z) = {99 \over 100} \, {\pi \over 3}$. 
For instance, project $z\in A_{j,k}$ onto $R_{j,k}$ along a vector
field in $\partial B$ which is roughly rotationally invariant around 
the tangent line to $\rho_{j,k}$ at $\varphi(x)$ and tangent to 
$\partial A_{j,k}$ along its conical boundary, as suggested by Figure 9.1. 

We can choose $p_0$ so that it is ${19 \over 10}$-Lipschitz
on $A_{j,k}$, say. Now observe that the two or three $A_{j,k}$
are disjoint, because the $\rho_{j,k}$ make angles larger
than $2\pi/3 - C \tau$ at $\varphi(x)$. So we get a clear 
definition of $p_0$ on $\cup_{j,k} A_{j,k}$.
We take $p_0(z) = \varphi(x)$ on $B_x \sm \cup_{j,k} A_{j,k}$;
then $p_0$ is continuous across the conical boundaries of the $A_{j,k}$, 
and we get a ${19 \over 10}$-Lipschitz map $p_0$ defined on $B_x$.

\hskip 2.2cm  
\epsfxsize = 7.5cm \epsffile{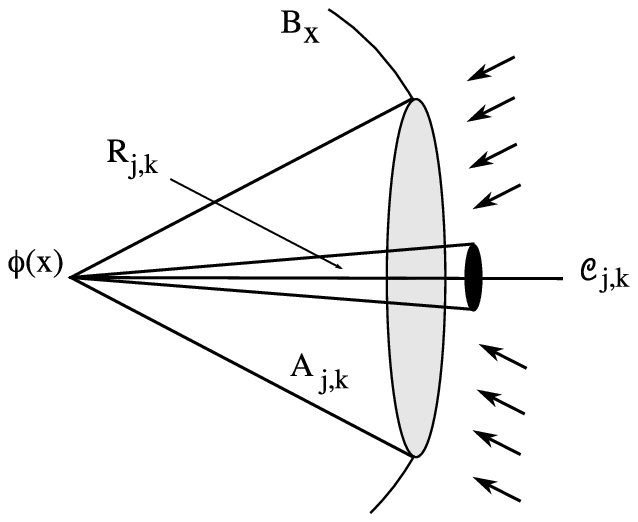}
\medskip
\noindent
{\bf Figure 9.1.} Project along the roughly invariant
vector field whose values on the vertical plane are suggested by the 
arrows. This is a picture in $\partial B$, which we assimilate to 
$\R^{n-1}$.
\medskip 

Recall that we decided to take $p_0(z)=z$ on
$R$; the two definitions fit on the $B_x$ 
(because $R \cap A_{j,k} = R_{j,k}$ and $R$ does not
meet $B_x \sm \cup_{j,k} A_{j,k}$), and then
$$
p_0 \ \hbox{ is $2$-Lipschitz on $[R \cap \partial B] \cup W$,}
\leqno (9.10)
$$
because the $B_x$ are so far from each other and 
$p_0(z) = z$ on $R \cap B_x$, which is not empty. Record that
by construction,
$$
p_0(z) = z \hbox{ for } z\in R \cap \partial B
\ \hbox{ and } \ p_0(z) \in R \hbox{ for } 
z\in [R \cap \partial B] \cup W.
\leqno (9.11)
$$

The second projection $p_1$ will be defined on 
$$
U = \big\{ z\in \overline B \, ; \, \dist(z,D) \leq C_7 \tau \big\},
\leqno (9.12)
$$
where $D$ is the union of the $D_{j,k}$, and $C_7$ will depend
on $\eta_0$, but not on $C_6$. First set 
$$\left\{\eqalign{
m(z) &= 0 \hskip2cm \hbox{ for } 0 \leq |z| \leq 10^{-1}\, ; \cr
m(z) &= |z| \hskip1.8cm \hbox{ for } |z| \geq 2\cdot 10^{-1} \, ; \cr
m(z) &= 2(|z|-10^{-1}) \ \, \hbox{ for } 
10^{-1} \leq |z| \leq 2\cdot 10^{-1}
}\right . \leqno (9.13)
$$
(so that $m$ is continuous and piecewise linear) and then
$$
p_1(z) = m(z) \, p_0(z/|z|)
\hbox{ for } z\in U.
\leqno (9.14) 
$$
Let us check that this is well defined.
First, we take $p_1(0)=0$, which is legitimate because 
$p_1(z) = 0$ for $z\in B(0,10^{-1}) \sm \{ 0 \}$ (since $m(z)=0$).
We also need to check that $p_0(z/|z|)$ is defined when
$|z| \geq 10^{-1}$. Set $w=z/|z|$; we want to check that
$w\in [R \cap \partial B] \cup W$. Observe that 
$\dist(w,D) \leq 10 C_7 \tau$ by (9.12) and because $D$ 
is a cone, and let $(j,k) \in \widetilde J$ be such that 
$\dist(w,D_{j,k}) \leq 10 C_7 \tau$. If 
$\dist(w,\partial D_{j,k}) \geq 1000 \eta_0^{-1} C_7 \tau$, 
then $w\in R_{j,k}$ by (9.7), and we are happy. 
Otherwise, recall that $\partial D_{j,k}$ is just the union of
two half lines through points $\varphi(x)$, with $x\in V$.
So we can find $x\in V$ such that 
$|w-\varphi(x)| \leq 2000 \eta_0^{-1} C_7 \tau$, and then
$x\in B_x \i W$, by (6.14) and if $C_6$ is large enough, depending
on $C_7$ and $\eta_0$.
Thus $w\in [R \cap \partial B] \cup W$ in all cases, and
$p_1$ is well defined.

Let us check that
$$
p_1 \hbox{ is $22$-Lipschitz.}
\leqno (9.15)
$$
For $z_1, z_2 \in U \sm B(0,10^{-1})$,
$$\leqalignno{
|p_1(z_1) - p_1(z_2 )| &\leq 
m(z_1) \, | p_0(z_1/|z_1|) - p_0(z_2/|z_2|)|
+ |m(z_1)-m(z_2)|\, |p_0(z_2/|z_2|)|
\cr&
\leq 2 m(z_1) \Big| {z_1 \over |z_1|}-{z_2 \over |z_2|} \Big|
+ |m(z_1)-m(z_2)| \leq 22 |z_1-z_2|
& (9.16)
}$$
by (9.10) and because the radial projection $z \to {z \over |z|}$ is
$10$-Lipschitz on $\overline B \sm B(0,10^{-1})$.
If instead $z_1 \in U \sm B(0,10^{-1})$ but $z_2 \in B(0,10^{-1})$,
we get that
$|p_1(z_1) - p_1(z_2 )| = |p_1(z_1)| \leq m(z_1) \leq 2 |z_1-z_2|$;
the case when $z_1, z_2 \in B(0,10^{-1})$ is trivial, and (9.15) 
follows.

Observe that
$$
p_1(z) \in R \cap \overline B
\hbox{ for } z\in U,
\leqno (9.17)
$$
just by (9.11) and because $R$ is a cone. We shall need to
know that the $R_{j,k}$ are essentially disjoint.
Let us even check that
$$
\dist(z, \bigcup_{(j',k') \neq (j,k)} R_{j',k'}) \geq 
\Min\Big\{ {\eta_0 |z|\over 10} \, ; \,
\dist(z,\partial D_{j,k}) \Big\}
\ \hbox{ for } z\in R_{j,k}.
\leqno (9.18)
$$
Since $R$, the $R_{j,k}$, and $\partial D_{j,k}$ are cones,
it is enough to check this when $z\in R_{j,k}\cap \partial B$. 
Recall from the definition (9.7) that  
$$
\dist(z,D_{j,k}) \leq 
10^{-2} \eta_0 \dist(z, \partial D_{j,k}) \leq 10^{-2} \eta_0.
\leqno (9.19)
$$
Denote by $a$ and $b$ the extremities of $\C_{j,k}$.
Then $\partial D_{j,k}$ is the union of the two half lines
through $\varphi(a)$ and $\varphi(b)$, and $D_{j,k}$ is the cone 
over the geodesic $\rho_{j,k}$ between $\varphi(a)$ and $\varphi(b)$.
Thus, if $\xi \in D_{j,k}$ minimizes the distance to $z$,
$$
\dist(z,\rho_{j,k})
\leq \dist(z,\xi/|\xi|)
\leq |z-\xi| + (1-|\xi|) \leq 2 |z-\xi|
= 2 \dist(z, D_{j,k})
\leqno (9.20)
$$
because $D_{j,k}$ is the cone over $\rho_{j,k}$ and $|z|=1$.
Now every point of $\rho_{j,k}$ is $C \tau$-close to $\C_{j,k}$,
because the two geodesics $\C_{j,k}$ and $\rho_{j,k}$ have
almost the same extremities, by (6.14), and because (2.5)
gives the stability of geodesics. So
$$\eqalign{
\dist(z,\C_{j,k}) &\leq \dist(z,\rho_{j,k}) + C \tau
\leq 2 \dist(z,D_{j,k}) + C \tau
\cr&
\leq 2 \cdot 10^{-2} \eta_0 + C \tau
\leq 3 \cdot 10^{-2} \eta_0
}\leqno (9.21)
$$
by (9.20) and (9.19). Let us prove (9.18) by contradiction,
and assume that we can find $(j',k') \neq (j,k)$ and 
$z' \in R_{j',k'}$ such that
$$
|z'-z| < \Min\Big\{ {\eta_0 \over 10} \, ; \,
\dist(z,\partial D_{j,k}) \Big\}.
\leqno (9.22)
$$
We want to get a contradiction. Set $w'=z'/|z'|$, and observe that 
$\dist(w',\C_{j',k'}) \leq 3 \cdot 10^{-2} \eta_0$
by the proof of (9.21). Then
$$\eqalign{
\dist(z,\C_{j',k'}) &\leq |z-z'| + |z'-w'| + \dist(w',\C_{j',k'})
\cr&\leq |z-z'| + (1-|z'|) + 3 \cdot 10^{-2} \eta_0
\cr&\leq 2 |z-z'| + 3 \cdot 10^{-2} \eta_0
\leq \eta_0/5
}\leqno (9.23)
$$
because $|z|=1$ and by (9.22). By (2.6), (9.21), and (9.23),
$\C_{j,k}$ and $\C_{j',k'}$ have a common 
endpoint $x$. Since $\C_{j,k}$ and $\C_{j',k'}$ are fairly 
short (by (2.5)) and leave from $x$ with large angles, 
we also have that $z$ and $z'$ lie in 
$B(x,\eta_0/2) \i B(\varphi(x),\eta_0)$. But in 
$B(\varphi(x),\eta_0)$ the situation is simple, 
$D_{j,k}$ and $D_{j',k'}$ are two half planes with a 
common boundary and that make an angle larger than $100^\circ$,
and $R_{j,k}$ and $R_{j',k'}$ are small conical neighborhoods
of $D_{j,k}$ and $D_{j',k'}$; then (9.22) fails, which proves (9.18).

\ms
Define $\pi_{j,k}$ on $\overline B \cap R_{j,k}$ by
the fact that $\pi_{j,k}(z) \in \Sigma_{j,k}$ and 
$\pi_{j,k}(z) - z \in P_{j,k}^\perp$. That is, if we 
denote by $\pi'_{j,k}$ the orthogonal projection on $P_{j,k} \, ,$
$\pi_{j,k}(z)$ is the point of the graph $\Sigma_{j,k}$ 
such that $\pi'_{j,k}(\pi_{j,k}(z)) = \pi'_{j,k}(z)$.
Trivially,
$$
\pi_{j,k}(z) = z \ \hbox{ for } z\in \overline B \cap \Sigma_{j,k}.
\leqno (9.24)
$$
Notice that
$$
|\pi'_{j,k}(z)-z| = \dist(z,D_{j,k}) 
\leq 10^{-2} \eta_0 \dist(z, \partial D_{j,k})
\leqno (9.25)
$$
for $z\in R_{j,k}$, by (9.7). Next 
$$
|\pi_{j,k}(z)-\pi'_{j,k}(z)|
\leq C \eta \dist(\pi'_{j,k}(z),\partial D_{j,k})
\leq 2C \eta \dist(z, \partial D_{j,k})
\leqno (9.26)
$$
because $\Sigma_{j,k}$ is a $C\eta$-Lipschitz graph
that contains $\partial D_{j,k}\cap \overline B$,
and by (9.25). Thus
$$
|\pi_{j,k}(z)-z| \leq 2 \cdot 10^{-2} 
\eta_0 \dist(z, \partial D_{j,k})
\ \hbox{ for } z\in R_{j,k} \, ,
\leqno (9.27)
$$ 
by (9.25) and (9.26). We want to define $\pi$ on 
$\overline B \cap R$ by setting
$$
\pi(z) = \pi_{j,k}(z) \ \hbox{ when } 
z\in \overline B \cap R_{j,k}
\leqno (9.28)
$$ 
so we need to verify that the different definitions 
match on the intersections. We first check that
when $(j',k') \neq (j,k)$,
$$
R_{j,k} \cap R_{j',k'} \cap \overline B
\i \partial D_{j,k} \cap \partial D_{j',k'} \cap \overline B
\i \Sigma \cap \overline B.
\leqno (9.29)
$$ 
Indeed, if $z \in R_{j,k} \cap R_{j',k'} \cap \overline B$,
(9.18) says that $z \in \partial D_{j,k}$, and by symmetry 
$z \in \partial D_{j',k'}$ too. The second inclusion 
follows from (8.6), (8.4), and the line above (8.2).

Now (9.29) says that $\pi_{j,k}(z) = \pi_{j',k'}(z) = z$
on $R_{j,k} \cap R_{j',k'} \cap \overline B$, and our
various definitions of $\pi$ match.

We'll need to know that
$$
\pi \hbox{ is $2$-Lipschitz on $\overline B \cap R$.}
\leqno (9.30)
$$
To prove (9.30), let $z, z' \in \overline B \cap R$ be given, 
and let $(j,k)$ and $(j',k')$ be such that $z \in R_{j,k}$ 
and $z' \in R_{j',k'}$. If $(j',k')=(j,k)$, we easily
get that $|\pi(z)-\pi(z')| = |\pi_{j,k}(z)-\pi_{j,k}(z')| 
\leq (1+C\eta) |z-z'|$, by definition of $\pi_{j,k}$ and because 
$\Sigma_{j,k}$ is a $C\eta$-Lipschitz graph over $D_{j,k}$.
So we may assume that $(j',k') \neq (j,k)$. Then
$$\eqalign{
|\pi(z)-\pi(z')| 
&\leq |z-z'| + |\pi_{j,k}(z)-z| + |\pi_{j',k'}(z')-z'| 
\cr& 
\leq |z-z'| + 2 \cdot 10^{-2} \eta_0 \,
\big[\dist(z, \partial D_{j,k}) + \dist(z', \partial D_{j',k'}) \big]
}\leqno (9.31)
$$
by (9.27).
Without loss of generality, we may assume that
$\dist(z, \partial D_{j,k}) \geq \dist(z', \partial D_{j',k'})$,
and then (9.31) yields
$$
|\pi(z)-\pi(z')| \leq |z-z'| + 4 \cdot 10^{-2} \eta_0 \,
\dist(z, \partial D_{j,k}).
\leqno (9.32)
$$
Recall from (9.18) that 
$$
|z-z'| \geq 
\Min\Big\{ {\eta_0 |z|\over 10} \, ; \,
\dist(z,\partial D_{j,k}) \Big\}.
\leqno (9.33)
$$
If $|z-z'| \geq \dist(z,\partial D_{j,k})$, (9.32) says
that $|\pi(z)-\pi(z')| \leq (1+4 \cdot 10^{-2} \eta_0)|z-z'|$,
and we are happy. Otherwise, $|z-z'| \geq {\eta_0 |z|\over 10}$
and (9.32) says that
$|\pi(z)-\pi(z')| \leq |z-z'| + 4 \cdot 10^{-2} \eta_0 \, |z|
\leq (1+4 \cdot 10^{-1}) \, |z-z'|$; (9.30) follows.

\ms
It could be that $\pi$ sends some points of
$\overline B \cap R$ slightly out of $\overline B$, so we
compose it with the standard retraction $r$ on the sphere,
defined by $r(z)=z$ when $z\in \overline B$ and $r(z)=z/|z|$
otherwise. Let us check that
$$
r \circ \pi(z) \in \Sigma \ \hbox{ for } z\in \overline B \cap R,
\leqno (9.34)
$$
where $\Sigma$ still denotes the union of the $\Sigma_{j,k}$.
Let $z\in \overline B \cap R$ be given, and let $(j,k)\in \widetilde J$
be such that $z\in R_{j,k}$. If $r \circ \pi(z) = \pi_{j,k}(x)$,
we are happy because $\pi_{j,k}(x) \in \Sigma_{j,k}$ by definition
of $\pi_{j,k}$. Otherwise, $r \circ\pi_{j,k}(z) \neq \pi_{j,k}(z)$
(by (9.28)), so $\pi_{j,k}(z)$ lies out of $\overline B$. 
Denote by $\xi$ its orthogonal projection on $P_{j,k}$,
and observe that $\xi \in \overline B$, because it is also the
the orthogonal projection of $z$ on $P_{j,k} \,$, and 
$z \in \overline B$.

Next, $|\pi_{j,k}(x)-\xi| \leq C \eta$ because 
$\Sigma_{j,k}$ is the graph over $D_{j,k}$ of a 
$C\eta$-Lipschitz function $G_{j,k}$ which is equal to $0$
on the two half lines that bound  $D_{j,k}$ (see (8.9) and
(8.6)). Hence $|\xi| \geq |\pi_{j,k}(x)| - |\pi_{j,k}(x)-\xi|
\geq 1 - C \eta$.
But on the region where $9/10 \leq |\xi| \leq 1$, 
(8.6) and (8.4) say that $G_{j,k}$ is a homogeneous function,
so $r \circ \pi_{j,k}(z) \in \Sigma_{j,k}$ because 
$\pi_{j,k}(z) \in \Sigma_{j,k}$. So (9.34) holds.

\ms
We are now ready to prove Lemma 9.4. We take 
$$
p_2 = r \circ \pi \circ p_1 
\leqno (9.35)
$$
with $p_1$ defined by (9.14). We know that $p_1$ is
defined on $U$ and take values in $\overline B \cap R$
(see (9.14) and (9.17)), so $p_2$ also is defined on $U$,
by (9.28), and takes value in $\Sigma \cap \overline B$, 
by (9.34). It is $44$-Lipschitz because $p_1$ is $22$-Lipschitz
by (9.15), $\pi$ is $2$-Lipschitz by (9.30), and $r$ is $1$-Lipschitz.
We shall check soon that
$$
U \hbox{ is a neighborhood of $E\cap \overline B$ in $\overline B$,}
\leqno (9.36)
$$
but let us take care of the other properties first.
We start with (9.5). Let 
$z\in \Sigma \cap \overline B \sm B(0,1/3)$
be given, and let $(j,k)$ be such that $z\in \Sigma_{j,k}$.
Then $z\in R_{j,k}$, by the definition (9.7) of $R_{j,k}$, because
$\Sigma_{j,k}$ is a $C\eta$-Lipschitz graph over $D_{j,k} \cap 
\overline B$ which contains $\partial D_{j,k} \cap \overline B$
(see (8.9) and (8.6)), and if $\eta$ is small enough.
So $p_0(z/|z|) = z/|z|$, by (9.11), and the definition
(9.14) says that $p_1(z)=z$, because $m(z) = |z|$ 
(by (9.13) and because $z\in \overline B \sm B(0,1/3)$).
By (9.28) and (9.24), $\pi(p_1(z)) = \pi(z)=z$, and then
$p_2(z) = r(z) = z$ by (9.35) and because $z\in \overline B$.
So (9.5) holds.

Now  we check (9.6). For $z\in U \cap B(0,1/3)$,
$p_1(z) \in R \cap B(0,1/3)$ because $m(z) \leq |z|$ 
and $p_0$ takes values in $R$ (by (9.13) and (9.11)).
Then $\pi(p_1(z)) \in \Sigma \cap B(0,1/2)$ by (9.28),
the fact that $\pi_{j,k}$ takes values in $\Sigma_{j,k}$,
and (9.27). Finally $p_2(z) = \pi(p_1(z)) \in \Sigma \cap B(0,1/2)$,
by (9.35), and  (9.6) follows.

So we are left with (9.36) to prove. Let $x\in E \cap \overline B$
be given. Recall from (4.2) that $\dist(x,X) \leq 100\varepsilon$,
where $X$ is also the minimal cone which was used in Section 6
to construct the $g_{j,k}$. Thus 
$\dist(x,|x|\C_{j,k}) \leq 100\varepsilon$ for some $(j,k)$.
Denote by $a$ and $b$ the endpoints of $\C_{j,k}$; since
$D$ contains the geodesic from $\varphi(a)$ to $\varphi(b)$,
(6.14) implies that $\dist(x,D) \leq 100\varepsilon+C\tau
< C_7 \tau$ (if $C_7$ is large enough), which implies
(9.36) (compare with (9.12)).
This completes our proof of Lemma 9.4.
\qed

\ms
Our next task is to use the mapping $p_2$ to construct
competitors for $E$ in $\overline B$.

First let us use the Whitney extension theorem to extend 
$p_2$ to $\R^n$ in a Lipschitz way. Since we can compose
$p_2$ with $r$ without changing its values on $U$
(recall that $p_2(U) \i \overline B$), we can do the extension
so that $p_2(z) \in \overline B$ for all $z$.

Then let $\xi > 0$ be very small, to be chosen later, 
define $\psi$ by 
$$\left\{\eqalign{
\psi(t) &= 1 \ \hbox{ for } 0 \leq t \leq 1-\xi \, , \cr
\psi(t) &= 0 \ \hbox{ for } t \geq 1 \, , \cr
\psi \ \  &\  \hbox{is affine on } [1-\xi,1],
}\right . \leqno (9.37)
$$
and set
$$
f_1(z) = \psi(|z|) \, p_2(z) + [1-\psi(|z|)] \, z
\ \hbox{ for } z\in \R^n.
\leqno (9.38)
$$
Observe that $f_1$ is Lipschitz, and that 
$$
f_1(z) = z \hbox{ for } z\in \R^n\setminus B
\ \hbox{ and } \ f_1(B) \i \overline B
\leqno (9.39)
$$
by (9.38), because $\psi(t) = 0$ for $t \geq 1$,
and because $p_2(\overline B) \i \overline B$.

Our first competitor for $E$ will be $F_1 = f_1(E)$.
Notice that $f_1$ satisfies the conditions (1.5)-(1.9) 
that define the admissible deformations, with 
$\widehat W \i \overline B$
and where we set $f_t(z) = t f_1(z)+(1-t)z$
for $t \in [0,1]$ and $z\in \R^n$.
So we can apply (1.11), which says that
$$
H^2(E \setminus F_1) \leq H^2(F_1 \setminus E)
+ 4 h(2).
\leqno (9.40)
$$
Since $E$ and $F_1$ coincide out $\overline B$, we can add
$H^2(E \cap F_1 \cap \overline B)$ to both sides and get that
$$
H^2(E \cap \overline B) \leq H^2(F_1 \cap \overline B) + 4 h(2).
\leqno (9.41)
$$
We shall introduce other competitors later, obtained by deformation
of $F_1$ inside $B$, so we shall not really use (9.41) directly, 
but even so it will be useful to estimate $H^2(F_1 \cap \overline B)$.

\ms\proclaim Lemma 9.42.
We have that 
$$
\limsup_{\xi \to 0_+} H^2(F_1 \cap \overline B) 
\leq H^2(\Sigma \cap B) + 2550
\int_{E \cap \partial B} \dist(z,\Gamma) \, dH^1(z).
\leqno (9.43)
$$
where $\Gamma = \cup_{j,k} \Gamma_{j,k}$ as in (9.3).

\ms
First observe that 
$$
f_1(E\cap B(0,1-\xi)) = p_2(E\cap B(0,1-\xi))
\i \Sigma \cap \overline B
\leqno (9.44)
$$ 
by (9.38), because $\psi(t) = 1$ for $t \leq 1-\xi$,
and by Lemma 9.4. Then set 
$$
A_\xi = B \sm B(0,1-\xi).
\leqno (9.45)
$$
If we prove that 
$$
\limsup_{\xi \to 0_+} H^2(f_1(E\cap A_\xi)) \leq 2550
\int_{E \cap \partial B} \dist(z,\Gamma) \, dH^1(z),
\leqno (9.46)
$$
Lemma 9.42 will follow at once, because 
$H^1(\Sigma \cap \partial B) < +\infty$ (and hence
$H^2(\Sigma \cap \partial B) = 0$) trivially, and also
$H^2(f_1(E \cap \partial B)) = 0$, since
$H^1(f_1(E \cap \partial B)) < +\infty$ because $f_1$ is 
Lipschitz and $H^1(E \cap \partial B) < +\infty$ by (6.2).

We shall use the area formula to compute $H^2(f_1(E\cap A_\xi))$, 
the coarea formula to estimate the result in terms of an integral 
on $A_\xi$, and then (9.46) and the lemma will follow from 
our additional precautionary assumptions (4.3) and (4.4). 

When the ambient dimension is $n=3$, we could avoid 
using most of this with a trick. First, we would modify 
$p_2$ slightly, so that $p_2(E\cap A_\xi) \i \partial B$, 
and then we would also modify (9.38) (by pushing radially
on $\partial B$) so that $f_1(E\cap A_\xi)$ is also contained 
in $\partial B$, and in fact in an arbitrary small neighborhood 
in $\partial B$ of the union of the arcs of geodesics from the points
$x\in E\cap \partial B$ to their projections $p_2(x) \in \Gamma$.
The area of this small piece of $\partial B$ around $\Gamma$ would 
be easier to estimate, just by covering it with balls centered on
$E\cap \partial B \sm \Gamma$, and (9.46) would follow.

Let us return to the general case and use the area theorem. 
Recall from Theorem~2.11 of [DS] 
(also see Section 2 of [D3] for the adaptation to the 
present setting) that $E$ is rectifiable; the area theorem
(Corollary 3.2.20 in [Fe]) says that 
$$ 
H^2(f_1(E\cap A_\xi))
\leq \int_{E\cap A_\xi} J_{f_1}(z) \, dH^2(z),
\leqno (9.47)
$$
where $J_{f_1}(z)$ denotes the approximate Jacobian of 
$f_1$ on $E$ at $z$.

We need to estimate Jacobians. Denote by $P(z)$ the 
approximate tangent plane to $E$ at $z$, which exists for $H^2$-almost 
every $z\in E$ (because $E$ is rectifiable). Incidentally, $P(z)$ is 
even a true tangent plane because $E$ is Ahlfors-regular (see for 
instance Exercise 41.21 on page 277 of [D2]). 
Pick an orthonormal basis $(v,w)$ of the vector space parallel to
$P(z)$, with $v$ orthogonal to the radial direction $(0,z)$.
Denote by $\theta(z)\in [0,\pi/2]$ the (unoriented) angle of $(0,z)$ 
with $w$, so that $|\langle w,z\rangle| = |z| \cos\theta(z)$.

We shall compute $J_{f_1}$ at a point 
$z\in E\cap A_\xi$ where $f_1$ (or equivalently $p_2$,
since $\psi(|z|)$ is smooth there, and by (9.38))
is differentiable in the direction of $P(z)$. Notice that 
this is the case almost-everywhere.

In the direction of $v$, $\psi(|\cdot|)$ has a vanishing
differential, and we get that 
$$
|Df_1(z)(v)| \leq \psi(|z|) \, |Dp_2(z)| 
+ [1-\psi(|z|)]  \leq 50 
\leqno (9.48)
$$
by (9.38), and because $p_2$ is $50$-Lipschitz by Lemma 9.4.

In the direction of $w$, we perform a similar computation, but
also have a term that comes from the derivative of $\psi$.
The derivative of $\psi(|\cdot|)$ in the 
direction of $w$ is bounded by $\xi^{-1}\cos\theta(z)$,
so we need to add
$$ 
[\xi^{-1} \cos\theta(z)] |p_2(z)-z| 
\leq 51 \xi^{-1} \cos\theta(z) \dist(z,\Sigma),
\leqno (9.49)
$$
because $p_2$ is $50$-Lipschitz and by (9.5). Thus  
$$  
|D f_1(z)(w)| 
\leq 50 + 51 \xi^{-1} \cos\theta(z) \dist(z,\Sigma)
\leqno (9.50)
$$
Hence
$$
J_{f_1}(z)  \leq |Df_1(z)(v)| \, |D f_1(z)(w)| 
\leq 2500+2550\xi^{-1}\cos\theta(z)\dist(z,\Sigma),
\leqno (9.51)
$$
and (9.47) yields
$$\eqalign{
H^2(f_1(E\cap A_\xi)) 
&\leq \int_{E\cap A_\xi} 
[2500+2550\xi^{-1} \cos\theta(z)\dist(z,\Sigma_F)] \, dH^2(z)
\cr& = 2500 H^2(E\cap A_\xi) + 2550 \xi^{-1}\int_{E\cap A_\xi} 
\cos\theta(z)\dist(z,\Sigma_F) \, dH^2(z).
}\leqno (9.52)
$$

Next let us apply the coarea formula (Theorem 3.2.22 in [Fe]) 
to the restriction to $E \cap A_\xi$ of the $C^1$ mapping $h: z \to |z|$, 
and integrate against the continuous function 
$z \to \dist(z,\Sigma)$; we get that
$$\eqalign{
\int_{E\cap A_\xi} \dist(z,\Sigma) &J_h(z) \, dH^2(z)
\cr&= \int_{t\in (1-\xi,1)} \int_{E\cap \partial B(0,t)}
\dist(z,\Sigma) \, dH^1(z) \, dt.
}\leqno (9.53)
$$
Here $J_h(z)$ is a one-dimensional jacobian; it is the largest size 
of the differential of $h$ applied to a vector of the tangent plane to
$E$ at $z$; if we keep the same basis $(v,w)$ of the vector space parallel 
to $P(z)$ to compute $J_h(z)$, we get that 
$J_h(z) = |Dh(z)(w)| = \cos\theta(z)$.
We compare to (9.52) and get that
$$\eqalign{
H^2(f_1(E\cap A_\xi)) 
&\leq 2500 H^2(E\cap A_\xi) 
\cr& \hskip 1cm + 2550 \xi^{-1}
\int_{t\in (1-\xi,1)} \int_{E\cap \partial B(0,t)}
\dist(z,\Sigma) \, dH^1(z) \, dt.
}\leqno (9.54)
$$
Now we use our assumptions (4.3) and (4.4)
(and recall that here $r=1$); (4.4) implies that 
$H^2(E\cap A_\xi)$ tends to $0$ when $\xi$ tends to $0$,
and (4.3) (applied to $f(z)=\dist(z,\Sigma)$) says that
$$
\lim_{\xi \to 0_+} \xi^{-1}
\int_{t\in (1-\xi,1)} \int_{E\cap \partial B(0,t)}
\dist(z,\Sigma) \, dH^1(z) \, dt
= \int_{E\cap \partial B} \dist(z,\Sigma) \, dH^1(z).
\leqno (9.55)
$$
So (9.54) implies that
$$
\limsup_{\xi \to 0_+} \xi^{-1} H^2(f_1(E\cap A_\xi)) 
\leq 2550 \int_{E\cap \partial B} \dist(z,\Sigma) \, dH^1(z).
\leqno (9.56)
$$
Finally observe that $\Gamma \i \Sigma$, because 
(9.1) and (9.3) say that, out of $B(0,19/20)$, 
$\Sigma$ coincides with the cone $\Gamma^\ast$ over the union $\Gamma$ 
of the curves $\Gamma_{j,k}$ constructed in Section 7. 
So $\dist(z,\Sigma) \leq \dist(z,\Gamma)$ for $z\in E\cap \partial B$,
and (9.46) and Lemma 9.42 follow from (9.56).
\qed

\ms
We shall end this section with estimates on the right-hand side of (9.43).
We start with $\int_{E \cap \partial B} \dist(z,\Gamma) \, dH^1(z)$,
so we want to estimate $\dist(z,\Gamma)$ for $z\in E \cap \partial B$,
and we shall start with $z \in g = \bigcup_{j,k} g_{j,k}$.
Denote by $x_{j,k}$ and $y_{j,k}$ the endpoints of $\C_{j,k}$
(the order will not matter). The endpoints of $g_{j,k}$ and 
$\Gamma_{j,k}$ are $\varphi(x_{j,k})$ and $\varphi(y_{j,k})$,
by (6.12) and (7.30). We denote by $\rho_{j,k}$ the geodesic in 
$\partial B$ with the same endpoints (that is, $\varphi(x_{j,k})$ 
and $\varphi(y_{j,k})$). By (7.31) and the definition (7.4),
$$
H^1(\Gamma_{j,k} \setminus g_{j,k}) 
\leq  H^1(g_{j,k} \setminus \Gamma_{j,k}) 
\leq C \eta^{-2} [\length(g_{j,k}) - \length(\rho_{j,k})].
\leqno (9.57)
$$
In addition,
$$
\length(g_{j,k}) = H^1(g_{j,k}) 
\leq \length(\C_{j,k}) + C_4 \tau
\leq \length(\rho_{j,k}) + C \tau
\leqno (9.58)
$$
because $g_{j,k}$ is simple, by (6.43), and because $\C_{j,k}$ and 
$\rho_{j,k}$ are geodesics of $\partial B$, with almost the same 
endpoints (recall that $\rho_{j,k}$ goes from $\varphi(x_{j,k})$ to 
$\varphi(y_{j,k})$ and use (6.14)). 
Thus $H^1(g_{j,k} \setminus \Gamma_{j,k}) \leq C \eta^{-2} \tau$,
and hence
$$
\dist(z,\Gamma_{j,k}) \leq C \eta^{-2} \tau \ \hbox{ for } z\in g_{j,k},
\leqno (9.59)
$$
because $g_{j,k}$ is connected and has common endpoints with
$\Gamma_{j,k}$. Let us also check that
$$
\dist(z,\Gamma) \leq C \eta^{-2} \tau 
\ \hbox{ for } z\in E\cap \partial B.
\leqno (9.60)
$$

Let $z\in E\cap \partial B$ be given. By (4.2), we can find 
$z_1 \in X$ such that $|z_1-z| \leq 100 \varepsilon$.
Then $z_2 = z_1/|z_1|$ lies in $X\cap \partial B$
(because $X$ is a cone), and $|z_2-z| \leq 200 \varepsilon r$
because $||z_1|-1| \leq 100 \varepsilon$. By definitions, 
$z_2$ lies in some $\C_{j,k}$ and it will be enough to
prove that $\dist(z_2,\Gamma_{j,k}) \leq C \eta^{-2}\tau$, or even,
in view of (9.59), that $\dist(z_2,g_{j,k}) \leq C \tau$.

By (6.15), $g_{j,k}$ is contained in a little tube of width $C\tau$ 
around $\C_{j,k}$. By (6.14), its endpoints lie within $C\tau$
of the endpoints of $\C_{j,k}$. If $z_2$ lies within $C\tau$
of an endpoint of $\C_{j,k}$, we are happy because it lies
close to an endpoint of $g_{j,k}$. Otherwise, we know that
$g_{j,k}$ crosses the little hyperdisk of radius $C \tau$
perpendicular to $\C_{j,k}$ at $z_2$, and we are again happy.
This proves that $\dist(z_2,g_{j,k}) \leq C \tau$, and (9.60) 
follows.

Let us now drop the dependence on $\eta$
from our notation. This does not matter, because we can 
choose $\tau$ and $\varepsilon$ very small, depending on $\eta$.
We obtain that
$$\eqalign{
\int_{E \cap \partial B} \dist(z,\Gamma) &\, dH^1(z)
\leq C \tau H^1(E \cap \partial B\setminus \Gamma)
\cr&
\leq C\tau H^1(E \cap \partial B \setminus g)
+ C\tau H^1(g \setminus \Gamma)
\cr& 
\leq C\tau H^1(E \cap \partial B \setminus g)
+ C\tau \sum_{j,k} [\length(g_{j,k}) - \length(\rho_{j,k})]
\cr&
= C\tau H^1(E \cap \partial B \setminus g)
+  C\tau [H^1(g)-H^1(\rho)], 
}\leqno (9.61)
$$
because $\dist(z,\Gamma) = 0$ on $\Gamma$, by (9.60),
(9.57), because the $g_{j,k}$ are essentially disjoint, 
and where we now set 
$$
\rho = \bigcup_{(j,k)\in \widetilde J} \rho_{j,k} \i \partial B.
\leqno (9.62)
$$ 
Thus $\rho$ is the union of the geodesics that connect the 
endpoints of the $g_{j,k}$ and the $\Gamma_{j,k}$, i.e., 
the $\varphi(x)$, $x\in V$. With the notation of Section 2,
$$
\rho = \varphi_\ast(K), \ \hbox{ with }
K = X\cap \partial B.
\leqno (9.63)
$$
\ms
The estimate (9.61) of what we lose in the deformation of
$E$ will be enough for our purposes. 
Let us also look at what we win in (9.43) with the term
$H^2(\Sigma \cap B(0,1))$.

Recall from the discussion above (9.3) that 
$\Sigma = \cup_{j,k} \Sigma_{j,k}$,
where $\Sigma_{j,k}$ is the graph of some function $G_{j,k}$
over a sector $\overline D_{j,k}$: see above (9.1).
The function $G_{j,k}$ was constructed as in Section~8, starting
with the cone $\Gamma_{j,k}^\ast$ over $\Gamma_{j,k}$, which 
corresponds to a (homogeneous) function $F_{j,k}$. The various 
$\Gamma_{j,k}^\ast$ are essentially disjoint, so (9.2)
says that 
$$\eqalign{
H^2(\Sigma \cap B)
&\leq H^{2}(\Gamma^\ast\cap B) - 10^{-4}
\sum_{j,k} [\length(\Gamma_{j,k}) - \length(\rho_{j,k})]
\cr&
\leq H^{2}(\Gamma^\ast\cap B) - 10^{-4} [H^1(\Gamma)-H^1(\rho)].
}\leqno (9.64)
$$
Set
$$
\delta_1 = H^1(\Gamma)-H^1(\rho) \geq 0
\leqno (9.65)
$$
(because the $\rho_{j,k}$ are geodesics with the same endpoints as 
the $\Gamma_{j,k}$),
$$
\delta_2 = H^1(g)-H^1(\Gamma) \geq 0
\leqno (9.66)
$$
(by the first part of (9.57) and because the curves are simple 
and essentially disjoint), and
$$
\delta_3 = H^1(E\cap \partial B)-H^1(g)
= H^1(E\cap \partial B \setminus g)
\leqno (9.67)
$$
(because $g \i E\cap \partial B$): thus
$$
H^1(E\cap \partial B) - H^1(\rho) = \delta_3+\delta_2+\delta_1.
\leqno (9.68)
$$

We are ready to combine the various estimates.
Recall from (9.41) that 
$H^2(E \cap \overline B) \leq H^2(F_1 \cap \overline B) + 4 h(2)$.
for every small $\xi >0$. Then (9.41) and (9.43) say that
$$\eqalign{
H^2(E \cap \overline B)
&\leq 4 h(2) + H^2(\Sigma \cap B) + 2550
\int_{E \cap \partial B} \dist(z,\Gamma) \, dH^1(z)
\cr&
\leq 4 h(2) + H^2(\Sigma \cap B) 
+ C \tau H^1(E \cap \partial B \setminus g)
+  C\tau [H^1(g)-H^1(\rho)]
\cr& 
= 4 h(2) + H^2(\Sigma \cap B) + C \tau \delta_3 + C \tau 
(\delta_1+\delta_2)
\cr&
\leq 4 h(2) + H^{2}(\Gamma^\ast\cap B) - 10^{-4} [H^1(\Gamma)-H^1(\rho)]
+ C \tau (\delta_1 + \delta_2 + \delta_3)
\cr&
= 4 h(2) + {1 \over 2} H^1(\Gamma) - 10^{-4} \delta_1
+ C \tau (\delta_1 + \delta_2 + \delta_3)
\cr&
= 4 h(2) + {1 \over 2} H^1(E\cap \partial B)
-{1 \over 2} \, (\delta_2 + \delta_3) - 10^{-4} \delta_1
+ C \tau (\delta_1 + \delta_2 + \delta_3)
\cr&
\leq {1 \over 2} H^1(E\cap \partial B)
- 10^{-5} (\delta_1 + \delta_2 + \delta_3) + 4 h(2) 
\cr&
= {1 \over 2} H^1(E\cap \partial B)
- 10^{-5} [H^1(E\cap \partial B) - H^1(\rho)] + 4 h(2) 
}\leqno (9.69)
$$
by (9.61), (9.64), and because $\Gamma^\ast$ is the cone over
$\Gamma$, if $\tau$ is small enough, and by (9.68).

\ms
Observe that (9.69) is the same thing as (4.7), because we reduced to 
$x=0$ and $r=1$, and by (9.63).
We also have (4.6), because here $H^1(\varphi_\ast(K)) = H^1(\rho)
= H^1(E\cap \partial B) - \delta_1 -\delta_2 -\delta_3$
by (9.68).
So it will be enough to establish (4.8), and Theorem 4.5 will follow.

Notice that if by chance $H^1(\rho) \leq 2d(0)$, 
then (9.69) says that
$$
H^2(E \cap B)
\leq {1 \over 2} H^1(E\cap \partial B)
- 10^{-5} [H^1(E\cap \partial B) - 2d(0)] + 4 h(2)
\leqno (9.70)
$$
so (4.8) holds with $\alpha = 10^{-5}$.
We cannot expect this to happen automatically, but we shall
see in the next section that otherwise, and if $X$ is a full minimal 
length cone such that $H^2(X\cap B) \leq d(0)$, the cone over $\rho$ 
is far from being minimal, and we can replace
our first competitor $F_1$ with a better one, improve (9.69),
and get (4.8).

\bigskip
\noindent {\bf 10. Our second competitor and the proof of Theorem 4.5}
\medskip

Let us first summarize the situation so far. We started form a reduced
almost minimal set $E$ in $U$ and a ball $B(x,r)$ centered on $E$ that 
satisfies the assumptions (4.1)-(4.4), reduced by translation and 
dilatation to $x=0$ and $r=1$, did some constructions, and eventually
obtained a competitor $F_1$ for $E$ in $B$, and proved (4.6), (4.7), 
and (9.69). In the special case when $H^1(\rho) \leq 2d(0)$, (9.69)
implies (4.8) and we are finished. Otherwise, when
$$
H^1(\rho) > 2d(0),
\leqno (10.1)
$$
we still need to do some work, and this is where we shall use the 
additional assumptions (about full length minimal cones) in Theorem 4.5.
But in the mean time, let us try to continue the construction, improve
$F_1$ neat its tip, and get a slightly stronger variant of (9.69).
We do this and state an intermediate result under the current assumptions 
(4.1)-(4.4), because we shall also use the result in Section 11.
We keep the notations of Section~9, and in particular
we continue to denote by $\varphi_\ast(X)$ the cone over $\rho$, as 
in (9.63) and near (2.8). The next lemma says that if we can find a deformation 
$\widetilde X$ of $\varphi_\ast(X)$ in $B$, with a smaller Hausdorff 
measure, we can plug it at the tip of $F_1$, get a better competitor
$F_2$, and improve on (9.69).

\ms \proclaim Lemma 10.2. Suppose we can find a Lipschitz function
$f : \R^n \to \R^n$ such that 
$$
f(x)=x \hbox{ for $x\in \R^n\sm B$ and $f(B) \i B$,}
\leqno (10.3)
$$
and such that if we set $\widetilde X = f(\varphi_\ast(X))$,
$$
H^2(\widetilde X \cap B) \leq H^2(\varphi_\ast(X)\cap B) - A
\leqno (10.4)
$$
for some $A>0$, then 
$$ 
H^2(E \cap \overline B) \leq {1 \over 2} H^1(E\cap \partial B)
- 10^{-5} [H^1(E\cap \partial B) - H^1(\rho)] - \kappa^2 A + 4 h(2).
\leqno (10.5)
$$

We shall not use (10.1)  or the full length assumption in the proof.
So we want to use $f$ and $\widetilde X$ to improve
our competitor $F_1$ from Section 9. Let us first check that
$$
F_1 \cap B(0,\kappa) \i \Sigma\cap B(0,\kappa) 
= \varphi_\ast(X) \cap B(0,\kappa).
\leqno (10.6)
$$
For the inclusion, recall that $F_1 = f_1(E)$. Since (9.44) says that 
$f_1(E \cap B(0,1-\xi)) \i \Sigma$, we just need to check that 
$f_1(z) \notin B(0,1/2)$ when $z \in E \sm B(0,1-\xi)$.
This is clear when $z \in E \sm B$, because then $f_1(z)=z$ 
by (9.39). Otherwise, $f_1(z)$ lies somewhere between $p_2(z)$ and 
$z$, by (9.38) and (9.37), and it is enough to check that 
$|p_2(z)-z| \leq 1/3$. By Lemma 9.4, $p_2$ is $50$-Lipschitz and
$p_2(z) = z$ for $z\in \Sigma \cap \overline B \sm B(0,1/3)$,
so it is enough to check that 
$\dist(z,\Sigma \cap \overline B) \leq 10^{-3}$.

Recall from the beginning of Section 9 that 
$\Sigma$ is the union of the $\Sigma_{j,k}$ (by (9.3)), and that
$\Sigma_{j,k}$ is the graph over $D_{j,k}$ of some $C \eta$-Lipschitz
function $G_{j,k}$ (see above (9.1)). In addition, $G_{j,k}$ vanishes on 
$\partial D_{j,k}\cap B$ (or equivalently, $\Sigma_{j,k}$ contains 
$\partial D_{j,k} \cap B$), by (8.6), (8.4), and the line above (8.2). 
Then every point of $D_{j,k}\cap B$ lies within $C\eta$ from $\Sigma \cap B$.

Now $z\in U$ by (9.36), and hence 
$\dist(z,D_{j,k}) =\dist(z,D_{j,k} \cap B) \leq C_7 \tau$ for some 
$(j,k) \in \widetilde J$, by (9.12) and the line below. Then 
$\dist(z,\Sigma_{j,k}\cap \overline B) \leq C_7 \tau + C \eta$, as 
needed for the inclusion in (10.6).

For the second part of (10.6), recall from (8.7) that 
$G_{j,k}(z) = 0$ when $|z| \leq 2 \kappa$; then
$$
\Sigma_{j,k} \cap B(0,\kappa) = D_{j,k} \cap B(0,\kappa)
\leqno (10.7)
$$ 
(recall that $G_{j,k}$ is $C \eta$-Lipschitz by (8.9)).
Denote by $a$ and $b$ the endpoints of $\C_{j,k}$; then
the endpoints of $\Gamma_{j,k}$ are $\varphi(a)$ and $\varphi(b)$
(see below (9.1), for instance), so $D_{j,k}$ coincides in $B(0,\kappa)$
with the cone over the geodesic from $\varphi(a)$ to $\varphi(b)$,
which is also denoted by $\varphi_\ast(\C_{j,k})$ in (2.8). 
Since $\varphi_\ast(X)$ is the cone over $\dsp\varphi_\ast(K) =
\bigcup_{j,k}\varphi_\ast(\C_{j,k})$, (10.7) implies that 
$\Sigma \cap B(0,\kappa) = \varphi_\ast(X) \cap B(0,\kappa)$. 
This completes our proof of (10.6).

Now we can glue the deformation $\widetilde X = f(\varphi_\ast(X))$
at the tip of $F_1$. Set 
$$
f^\kappa (z) = \kappa f(\kappa^{-1} z)
\ \hbox{ for } z \in \R^n,
\leqno (10.8)
$$
and notice that 
$$
f^\kappa(z) = z \hbox{ when } |z| \geq \kappa
\hbox{ and $f^\kappa(z) \in B(0,\kappa)$ otherwise;}
\leqno (10.9)
$$
then set $f_2 = f^\kappa \circ f_1$. Observe that $f_2$ is Lipschitz, 
and $f_2$ is the endpoint of a one-parameter family of functions that 
satisfy (1.5)-(1.9) with $\widehat W \i \overline B$, just by (10.9) 
and because $f_1$ has this property (see below (9.39)). 
So we may apply (1.11) with $F_2 = f_2(E)$, just as we did for (9.40) 
and (9.41), and get that 
$$
H^2(E \setminus F_2) \leq H^2(F_2 \setminus E) + 4 h(2)
\leqno (10.10)
$$
and, since $E$ and $F_2$ coincide out $\overline B$,
$$
H^2(E \cap \overline B) \leq H^2(F_2 \cap \overline B) + 4 h(2).
\leqno (10.11)
$$
We want to bound the right-hand side of (10.11) slightly better than what 
we did in Section~9. We do not need to modify our estimate of 
$H^2(F_1 \cap \overline B \sm B(0,\kappa))$, because (10.9) says that 
$F_2 =F_1$ out of $B(0,\kappa)$. On $B(0,\kappa)$, we simply used (9.44) 
to say that 
$$
H^2(F_1 \cap B(0,\kappa)) \leq H^2(\Sigma \cap B(0,\kappa))
= H^2(\varphi_\ast(X) \cap B(0,\kappa)) = \kappa^2 
H^2(\varphi_\ast(X) \cap B),
\leqno (10.12)
$$
by (10.6) and because $\varphi_\ast(X)$ is a cone. Here we say instead that
$$
\eqalign{
F_2 \cap B(0,\kappa)
&= f^\kappa (F_1) \cap B(0,\kappa)
= f^\kappa(F_1 \cap B(0,\kappa))
\i f^\kappa(\varphi_\ast(X) \cap B(0,\kappa))
\cr&= \kappa f(\kappa^{-1} [\varphi_\ast(X) \cap B(0,\kappa)])
=\kappa f(\varphi_\ast(X) \cap B)
\cr&=\kappa [f(\varphi_\ast(X)) \cap B]
= \kappa [\widetilde X \cap B]
}\leqno (10.13)
$$
by (10.6) and (10.9), because $\varphi_\ast(X)$ is a cone, by (10.3), and by 
definition of $\widetilde X$. This yields
$$
H^2(F_2 \cap B(0,\kappa)) \leq \kappa^2 H^2(\widetilde X \cap B)
\leq \kappa^2 H^2(\varphi_\ast(X)\cap B) - A \kappa^2
\leqno (10.14)
$$
by (10.4). In other words, compared with (10.12), we won
$A \kappa^2$. Recall that we obtained in (9.69) that
$$
H^2(E \cap \overline B) \leq 
{1 \over 2} H^1(E\cap \partial B)
- 10^{-5} [H^1(E\cap \partial B) - H^1(\rho)] + 4 h(2);
\leqno (10.15)
$$
then the extra $A \kappa^2$ yields (10.5), and Lemma 10.2
follows.
\qed

\ms
Return to the proof of Theorem 4.5. We need to prove (4.8)
in the remaining case when (10.1) holds, and now we shall
use the extra assumptions in Theorem  4.5, namely, that 
$H^2(X\cap B) \leq d(0)$
and that $X$ is a full length minimal 
cone with constants $C_1$ and $\eta_1$. Then
$$
H^1(K) = H^1(X\cap \partial B) = 2 H^2(X\cap B) 
\leq 2 d(0) < H^1(\rho) = H^1(\varphi_\ast(K))
\leqno (10.16)
$$
because $K = X\cap \partial B$, $X$ is a cone, and by (10.1) and
(9.63).
Also notice that $\varphi \in \Phi(\eta_1)$ (where $\Phi(\eta_1)$
is as in (2.7)), by (6.14) and if $\tau$ is small enough.
So Definition 2.10 says that there is a deformation $\widetilde X$ 
of $\varphi_\ast(X)$ in $B$ such that
$$
H^2(\widetilde X \cap B) \leq H^2(\varphi_\ast(X)\cap B)
- C_1^{-1} [H^1(\varphi_\ast(K)) - H^1(K)].
\leqno (10.17)
$$
That is, there is a Lipschitz function $f : \R^n \to \R^n$
such that (10.3) holds and for which $\widetilde X = f(\varphi_\ast(X))$ 
satisfies (10.17).

Thus  we can apply Lemma 10.2 and get (10.5) with 
$$
A = C_1^{-1} [H^1(\varphi_\ast(K)) - H^1(K)]
= C_1^{-1} [H^1(\rho) - H^1(K)]
\leqno (10.18)
$$
by (9.63). Set $\alpha = C_1^{-1} \kappa^2$; 
we can safely assume that $\alpha < 10^{-5}$, so (10.5) also holds 
with $10^{-5}$ replaced with $\alpha$, because 
$H^1(E\cap \partial B) - H^1(\rho) 
= \delta_1 + \delta_2 + \delta_3 \geq 0$ by (9.68).
Thus we get that
$$\eqalign{
H^2(E \cap \overline B) &\leq {1 \over 2} H^1(E\cap \partial B)
- \alpha [H^1(E\cap \partial B) - H^1(\rho)] - \kappa^2 A + 4 h(2)
\cr&=
{1 \over 2} H^1(E\cap \partial B)
- \alpha [H^1(E\cap \partial B) - H^1(K)] + 4 h(2)
\cr&\leq
{1 \over 2} H^1(E\cap \partial B)
- \alpha [H^1(E\cap \partial B) - 2d(0)] + 4 h(2)
}\leqno (10.19)
$$
by (10.5), (10.18), and (10.16).
Since here $r=1$ and $x=0$, we just established (4.8). This completes 
our proof of Theorem~4.5.
\qed

\ms
To  end this section, we now explain how we can construct
competitors for cones (to be used in Lemma 10.2 or just
to prove that a given cone has the full length property).

We keep the same notation as in Section 2. That is,
$X$ is a minimal cone, we set $K = X\cap \partial B$, and we decompose 
$K$ as a union of arcs $\C_{j,k}$, $(j,k) \in \widetilde J$, that 
connect vertices $x\in V$.
Then we pick $\varphi \in \Phi(\eta_1)$, with $\eta_1 < \eta_0/10$
(as in (2.7)), and define $\dsp\varphi_\ast(K) 
= \bigcup_{(j,k) \in \widetilde J}\varphi_\ast(\C_{j,k}) \i \partial B$ 
and the cone $\varphi_\ast(X)$ over $\varphi_\ast(K)$
as we did near (2.8). The set $\rho$ of (9.62) in Section 9 was a typical
example of this.
We are interested in finding better competitor
for $\varphi_\ast(X)$ in $B$, when this is easy to do.

For each vertex $x\in V$, we define a deviation $\alpha_\varphi(x)$ 
from the ideal angles at $\varphi(x)$ as follows.

First suppose that $x\in V_0$; thus there are three curves
$\C_{j,k}$ which have $x$ as one of their endpoints. Denote by
$w_1$, $w_2$, and $w_3$ the unit tangent vectors to the
three $\C_{j,k}$ (pointing away from $\varphi(x)$, and set
$$
\alpha_\varphi(x) = |w_1+w_2+w_3|
\leqno (10.20)
$$
Of course $\alpha_\varphi(x) =0$ when the three $\C_{j,k}$ make
$120^\circ$ angles at $\varphi(x)$; $\alpha_\varphi(x)$ is not 
exactly equivalent to the maximum deviation 
from $120^\circ$ of the mutual angles 
of the $w_j$, because when $w_1$ and $w_2$ make a $120^\circ$ angle
and $w_1+w_2+w_3$ is almost orthogonal to the plane that
contains $w_1$ and $w_2$, the angle deviation is of the order
of $|w_1+w_2+w_3|^{2}$. A definition with angles would not
be precise enough in this case.

When $x\in V_1$, there are only two arcs $\C_{j,k}$ leaving
from $x$, we denote by $w_1$ and $w_2$ their tangent direction
at $\varphi(x)$, and by $\theta \in (0,\pi]$ the angle of $w_1$ and $w_2$,
with the convention that $\theta$ is close to $\pi$, and we set
$$
\alpha_\varphi(x) = \pi- \theta.
\leqno (10.21)
$$
In this case, there is no difficulty, $\alpha_\varphi(x)$ is 
equivalent to $|w_1+w_2|$, and we decided to take the ``simpler" 
definition. Finally we set
$$
\alpha_+(\varphi) = \sup_{x\in V} \alpha_\varphi(x).
\leqno (10.22)
$$
Since we assumed that $\varphi \in \Phi(\eta_1)$ with
$\eta_1$ small, we still get that $\alpha_+(\varphi)$ is
small (because our description of minimal cones readily gives 
$\alpha_+(\varphi) = 0$ when $\varphi$ is the identical map).
The next lemma says that when $\alpha_+(\varphi) \neq 0$,
we can win $C^{-1} \alpha_+(\varphi)^2$ by deforming 
$\varphi_\ast(X)$ in $B$.

\ms \proclaim Lemma 10.23. Let $X$, $\varphi \in \Phi(\eta_1)$,
and $\alpha_+(\varphi)$ be as above. We can find a Lipschitz function
$f : \R^n \to \R^n$ such that (10.3) holds, and
$$
H^2(\widetilde X \cap B) \leq H^2(\varphi_\ast(X)\cap B) 
- C^{-1} \alpha_+(\varphi)^2
\leqno (10.24)
$$
for $\widetilde X = f(\varphi_\ast(X))$. Here $C$ depends only
on $n$ and the constant $\eta_0$ in Section 2 
(associated to the description of the minimal cones).

\ms
To prove the lemma, we first select $x\in V$ such that
$\alpha_\varphi(x) = \alpha_+(\varphi)$. We we shall first assume
that $x\in V_0$ (but the other case is simpler).

Choose coordinates so that the first axis is in the direction of
$\varphi(x)$, and hence $\varphi(x) = (1,0) \in \R \times \R^{n-1}$.
Our mapping $f$ will be of the form
$$
f(z) = z + \psi(z_1,\rho) \, v,
\leqno (10.25)
$$
where $\rho = (z_2^2 + \cdots + z_n^2)^{1/2}$, and
$v = (0,v_2, \cdots,v_n)$ is a small vector perpendicular 
to the direction of $\varphi(x)$. We choose the nonnegative
smooth bump function $\psi$ so that
$$
\psi \hbox{ is supported in }
V = [1/4, 1/2] \times [0,10^{-1} \eta_0],
\leqno (10.26)  
$$
where $\eta_0$ is as in (2.2) and (2.3),
$$
0 \leq \psi \leq 1 \ \hbox{  and } \ 
\partial_2 \psi =: {\partial \psi \over \partial r} \leq 0 
\, \hbox{ everywhere,}
\leqno (10.27)
$$
$$
\psi(z_1,0) = 1 \hbox{ for } |z_1-1/3| \leq 1/10, 
\leqno (10.28)
$$
and 
$$
|\nabla \psi| \leq 100 \eta_0^{-1} 
\hbox{ everywhere.}
\leqno (10.29)
$$
We shall choose $v$ such that $|v| \leq 200^{-1}\eta_0$.
Then (10.3) holds automatically, and also $f$ is a bijection, 
for instance because $|Df-Id| < 1/2$ everywhere.

We intend to compute $H^2(f(\varphi_\ast(X)\cap B))$
with the area formula. We may restrict the computation to
$W = \big\{ z\in \R^n \, ; \, (z_1,\rho) \in V \big\}$,
because $f(z)=z$ out of $W$. Let us check that in $W$, 
$\varphi_\ast(X)$ is composed of three faces $P_1$, $P_2$, and $P_3$, 
which are half-planes that make angles that are close to $120^\circ$,
and have a common boundary, the half line through $\varphi(x)$. 
Since $\varphi_\ast(X)$ is the cone over $\varphi_\ast(K)$, it is 
enough to  check that in $\partial B \cap B(\varphi(x),4 \eta_0/10)$,
$\varphi_\ast(K)$ is composed of three arcs of great circles that
leave from $\varphi(x)$ with angles that are close to $120^\circ$.
This last follows from the description of $X$ in Section 2
(and in particular (2.2) and (2.3)), and the fact 
that $\varphi \in \Phi(\eta_1)$ for some $\eta_1 < \eta_0/10$.
The vector plane parallel to $P_j$ is spanned by $\varphi(x)$ and the 
tangent vector $w_j$ at $\varphi(x)$ to the arc of $\varphi_\ast(K)$ 
that corresponds to $P_j$. Thus the three $w_j$ are orthogonal to 
$\varphi(x)$, and they make with each other angles that are close to 
$120^\circ$.

We shall do the area computation one face at a time, so let us
take care of $P_1 \cap W$. Let us chose our coordinates so that
$w_1 = (0,1,0 \cdots,0)$. Also set $e_1 = \varphi(x)$ for convenience.
We need to compute the differential of $f$ in directions of $e_1$
and $w_1$. Notice that $\rho = z_2$ on $F_1$, so 
$f(z) = z + \psi(z_1,z_2) \, v$, and the derivatives are 
$Df(e_1) = e_1 + \partial_1\psi(z_1,z_2) \, v$
and $Df(w_1) = w_1 + \partial_2\psi(z_1,z_2) \, v$.

Next write $v = \beta w_1 + v_3$, where the two first coordinates
of $v_3$ vanish (recall that $v$ is orthogonal to $e_1$). Then
$$\eqalign{
Df(e_1) \wedge &Df(w_1)
= [e_1 + \partial_1\psi(z_1,z_2) \, v] \wedge 
[w_1 + \partial_2\psi(z_1,z_2) \, v]
\cr& = e_1 \wedge w_1 + [\partial_2\psi(z_1,z_2)\, e_1 
- \partial_1\psi(z_1,z_2) \, w_1] \wedge v
\cr& = e_1 \wedge w_1 + [\partial_2\psi(z_1,z_2)\, e_1 
- \partial_1\psi(z_1,z_2) \, w_1] \wedge [\beta w_1 + v_3]
\cr& = [1+\beta \partial_2\psi(z_1,z_2)] \, e_1 \wedge w_1 
+ \partial_2\psi(z_1,z_2)\, e_1 \wedge v_3
- \partial_1\psi(z_1,z_2) \, w_1 \wedge v_3
}\leqno (10.30)
$$
and the jacobian determinant of the restriction of $f$ to
$F_1$ is
$$\eqalign{
J_1(z) &= |Df(e_1) \wedge Df(w_1)| 
\cr&= \big\{ [1+\beta \partial_2\psi(z_1,z_2)]^2 + 
\partial_2\psi(z_1,z_2)^2 |v_3|^2 
+ \partial_1\psi(z_1,z_2)^2 |v_3|^2 \big\}^{1/2}
\cr&
\leq 1 + \beta \partial_2\psi(z_1,z_2) + C |v|^2
}\leqno (10.31)
$$
by (10.29) and because $|v|^2 = \beta^2 + |v_3|^2$.
Notice also that $\beta = \langle v,w_1 \rangle$; hence, when
we apply the area formula to compare $H^2(f(P_1\cap W))$ to
$H^2(P_1\cap W)$, we get that
$$\eqalign{
H^2(f(P_1\cap W)) &- H^2(P_1\cap W)
= \int_{P_1\cap W} [J_1(z)-1] \, dz
\cr&\leq  \int_{P_1\cap W} \big[\beta \partial_2\psi(z_1,\rho) 
+ C |v|^2 \big] \, dz
\cr& \leq \langle v,w_1 \rangle \int_{P_1\cap W} \partial_2\psi(z_1,\rho) 
+ C |v|^2 H^2(P_1\cap W)
\cr&
\leq \langle v,w_1 \rangle \int_{P_1\cap W} \partial_2\psi(z_1,\rho) 
+ C |v|^2.
}\leqno (10.32)
$$ 
We used some coordinates to prove this, but the final result is 
stated in a more invariant way, and is also valid for the faces
$P_2$ and $P_3$. The integral
$a = - \int_{P_j\cap W} \partial_2\psi(z_1,\rho)$ has the same value 
for $j=1,2,3$ because of the radial symmetry of $\psi$, and then
$$\eqalign{
a &= - \int_{P_j\cap W} \partial_2\psi(z_1,\rho) \, dz
\geq - \int_{|z_1-1/3| \leq 1/10} 
\int_{0}^{10^{-1} \eta_0} 
\partial_2 \psi(z_1,\rho) d\rho dz_1
\cr&
= \int_{|z_1-1/3|\leq 1/10}  
[\psi(z_1,0)-\psi(z_1,10^{-1} \eta_0)] dz_1
= \int_{|z_1-1/3| \leq 1/10}dz_1 = 1/5
}\leqno (10.33)
$$
by (10.27), (10.26), and (10.28). The analogue of (10.32)
for the face $P_j$ is
$$
H^2(f(P_j\cap W)) - H^2(P_j\cap W)
\leq - a \langle v,w_j \rangle + C |v|^2,
\leqno (10.34)
$$
and when we sum over $j$ we get that
$$\eqalign{
H^2[f(&\varphi_\ast(X))\cap B] - H^2[\varphi_\ast(X)\cap B]
= H^2[f(\varphi_\ast(X))\cap W] - H^2[\varphi_\ast(X)\cap W]
\cr&
= \sum_j [H^2(f(P_j\cap W)) - H^2(P_j\cap W)]
\leq - a \,\langle v,w_1+w_2+w_3 \rangle + C |v|^2
}\leqno (10.35)
$$
because $f$ is bijective and $f(z)=z$ out of $W$ and by (10.34).
Set $s=w_1+w_2+w_3$, and choose $v=cs$, where $c$ is a small positive
constant. If $c$ is small enough, our constraint that
$|v| \leq 200^{-1}\eta_0 $ is satisfied. Also,
$a \langle v,w_1+w_2+w_3 \rangle = ac |s|^2$, while
$C |v|^2 = C c^2 |s|^2 < ac |s|^2/2$ (if $cC \leq a/2$), so (10.35)
says that
$$\eqalign{
H^2(f(\varphi_\ast(X)) \cap B)
- H^2(\varphi_\ast(X) \cap B)
&\leq -{ac \over 2} \, |s|^2 
\leq - {c \over 10} \, |s|^2 
\cr&= - {c \over 10} \alpha_\varphi(x)^2 = 
- {c \over 10} \alpha_+(\varphi)^2
}\leqno (10.36)
$$
by (10.20) and our choice of $x$. Thus (10.24)
holds in this case.

\ms
Let us now assume that $x\in V_1$. Denote by
$w_1$ and $w_2$ the unit tangent vectors to the two arcs of 
$\varphi_\ast(K)$ that leave from $\varphi(x)$. Thus 
$\alpha_+(\varphi) = a_\varphi(x) = \pi - \theta$,
where $\theta \in (0,\pi]$ denotes the angle between $w_1$ and $w_2$
(see (10.21)). As before, $\theta$ is close to $\pi$ because
$\varphi \in \Phi(\eta_1)$.

We may now repeat the same proof as above, with the same formula
for $f$, except that now we only have two faces $P_1$ and $P_2$ with 
the same boundary. We do the same computations as before, and we still
get that 
$$
H^2(f(\varphi_\ast(X)) \cap B) - H^2(\varphi_\ast(X) \cap B)
\leq - {c \over 10} \, |s|^2,
\leqno (10.37)
$$
but this time with $s=w_1+w_2$. It is very easy to see that 
$|s| \geq \alpha_+(\varphi)/2 = (\pi-\theta)/2$, and so (10.24) 
holds in this case too. 

Note that in this simpler case we could also have used the computations of 
Section 6, applied to a Lipschitz graph composed of two arcs of geodesics, to 
construct $f$.

So we proved (10.24) in both cases, and Lemma 10.23 follows.
\qed

\bigskip
\noindent {\bf 11. Approximation of $E$ by cones}
\medskip

So far we focused our energy on proving decay estimates for the density 
$\theta(r)$. Now we want to show that this density controls the 
geometric behavior of $E$. In this section, we show that it controls
the Hausdorff distance to minimal cones. In the next one, we shall 
explain how to deduce $C^{1+\alpha}$ estimates from this.

In this section too we are given a reduced almost-minimal set $E$ in 
$U \i \R^3$, with gauge function $h$, and a point $x\in E$. As before,
we set 
$$
\theta(r) = r^{-2} H^2(E\cap B(x,r)), \ 
d(x)=\lim_{r\to 0} \theta(r) \hbox{, and } \ f(r)=\theta(r)-d(x).
\leqno (11.1)
$$

Our main task will be to use the smallness of $f$ 
to show that $E$ is well approximated by certain types of cones.

Let us be more specific.
We start from a minimal cone $X$, and cut $K = X\cap \partial B(0,1)$
into arcs of circles $\C_{j,k}$, $(j,k) \in \widetilde J$,
as we did in Section 2. Thus the $\C_{j,k}$ are arcs of great circles.
The endpoints of these arcs lie in a finite set $V=V_0 \cup V_1$, where 
$V_0$ is the initial set of vertices where three arcs $\C_{j,k}$ meet with 
$120^\circ$ angles, and $V_1$ is the set of added vertices where only 
two $\C_{j,k}$ meet, with $180^\circ$ angles.

We shall approximate $E$ first by perturbations of $X$ by mappings 
$\varphi \in \Phi(\eta_1)$, with $\eta_1 \leq \eta_0/10$ (see near (2.7)), 
which means that we shall use the sets $\varphi_\ast(K)$ and $\varphi_\ast(X)$.
Recall that $\varphi_\ast(K)$ is obtained from $K$ by replacing
each arc $\C_{k,k}$ with endpoints $a$ and $b\in V$ with the geodesic
from $\varphi(a)$ to $\varphi(b)$, and that $\varphi_\ast(X)$
denotes the cone over $\varphi_\ast(K)$.

We shall denote by ${\cal Z}_0(X,\eta_1)$ the collection of
cones $\varphi_\ast(X)$ obtained this way, and by ${\cal Z}(X,\eta_1)$ 
the set of images of cones of ${\cal Z}_0(X,\eta_1)$ by translations.

For $Z \in {\cal Z}(X,\eta_1)$, we set $\alpha_+(Z) = 
\alpha_+(\varphi)$, where $\varphi\in \Phi(\eta_1)$ is such that
$Z$ is a translation of $\varphi_\ast(X)$, and $\alpha_+(\varphi)$
is defined in (10.20)-(10.22). Thus $\alpha_+(Z)$ measures the largest
difference between the position of the tangent vectors to 
$\varphi_\ast(K)$ at a vertex $\varphi(x)$ and the standard
position.

We denote by $\cal Z$ the union of the ${\cal Z}(X,\eta_1)$,
where $X$ is a minimal cone, and $\eta_1 < \eta_0/10$ (again with
$\eta_0$ as in (2.2)-(2.6)).
Finally we set
$$
\beta_{Z}(x,r) = \inf \big\{ d_{x,r}(E,X)+\alpha_+(Z)  
\, ; \, X \in {\cal Z} \hbox{ is centered at } x \big\}
\leqno (11.2)
$$
for $x\in E$ and $r>0$, where 
$$\eqalign{
d_{x,r} (E,Z) = r^{-1} \sup \big\{\dist(x,Z) \, &; \, 
x \in E \cap B(x,r) \big\}
\cr&
+ r^{-1}\sup \big\{\dist(x,E)  \, ; \, x \in Z \cap B(x,r) \big\}.
}\leqno (11.3)
$$

\medskip \proclaim Theorem 11.4. 
There is a positive constant $C$ such that if $E$ is a reduced 
almost minimal set in $U$ with gauge function $h$, $x \in E$, 
and $r_0>0$ are such that $B(x,110r_0) \i U$, then 
$$
\beta_{Z}(x,r_0) \leq C f(110r_0)^{1/3} + C h_1(110r_0)^{1/3}. 
\leqno (11.5) 
$$

\ms
See (1.16) for the definition of $h_1(t)$.
The power $1/3$ is probably not optimal; $1/2$ looks more plausible 
but seems to require more work too. Anyway, if $h_1$ is small enough, 
we proved that $f(r)$ decays like a power of $r$ when $r$ tends to $0$,
so Theorem 11.4 says that the normalized Hausdorff distance 
$\beta_{Z}(x,r)$ tends to $0$ like a power too. Later on, we shall
use this to control the normalized Hausdorff distance to minimal
cones, but for the moment we only care about cones in $\cal Z$.

\ms
The proof of Theorem 11.4 will keep us busy for the rest of this 
section. Let $x$ and $r_0$ be as in the statement. 
We can assume that $x=0$ (by translation invariance), and that
$$
f(110r_0) + h_1(110r_0)\leq \varepsilon_1,
\leqno (11.6)
$$ 
with $\varepsilon_1$ as small as we want, because otherwise
(11.5) holds with $C=2\varepsilon_1^{-1/3}$. 
Let us apply Proposition~7.24 in [D3]; 
observe that $\theta (110r_0) = d(0) + \varepsilon_1
\leq \inf_{0 < t < 11r_0/10} \theta(t) + C \varepsilon_1$
by (3.8) and (11.6), so we can apply the proposition, and find
a minimal cone $X$ centered at $0$ such that
$$
d_{0,100r_0}(E,X) \leq \varepsilon,
\leqno (11.7)
$$
where $\varepsilon > 0$ is any small constant given in advance,
and if $\varepsilon_1$ is small enough. We choose $\varepsilon$
as in Theorem 4.5 above, which will allow us to apply the 
construction of Sections 6-10.

Observe that 
$$
f(r) + h_1(r)\leq C\varepsilon_1
\ \hbox{ for } 0 < r \leq 60r_0,
\leqno (11.8)
$$
because $h_1$ is nondecreasing, and by (3.6). 
Set $v(r) = H^2(E \cap B(0,r))$, and 
denote by ${\cal R}$ the set of radii $r \in (0,2r_0)$
such that $\theta$ and $v$ are differentiable at $r$, 
$$
\theta'(r) = r^{-2} v'(r) - 2 r^{-3} v(r)
\ \hbox{ and } \ v'(r) \geq H^1(E\cap \partial B(0,r)),
\leqno (11.9) 
$$
and in addition (4.3) and (4.4) hold. Then   
$$
H^1((0,2r_0) \setminus {\cal R}) = 0
\leqno (11.10) 
$$
by Lemma 5.1, Lemma 5.5, (5.8), and Lemma 4.12.
We shall get better estimates when we consider radii
$r\in \cal R$ such that $j(r)$ is small, where we set
$$
j(r) = rf'(r) + f(r) + h_1(r) = r\theta'(r) + f(r) + h_1(r) + h(2r).
\leqno (11.11) 
$$
Our first task will be to fix $r\in \cal R$ and control the geometry of 
$E\cap \partial B(0,r)$ in terms of $j(r)$. Observe that
$$
H^1(E \cap \partial B(0,r)) \leq v'(r)
= r^2 \theta'(r) +2r^{-1} v(r)
= r^2 \theta'(r) +2r \theta(r) 
\leq 2r j(r) + 2 r d(0)
\leqno (11.12) 
$$ 
for $r \in {\cal R}$, by (11.9) and because $\theta(r)=f(r) + d(0)$.
A first consequence of this is that if we set
$$
{\cal R}_1 = \big\{ r \in {\cal R} \, ; \, j(r) \leq \tau d(0)
\big\};
\leqno (11.13) 
$$
then $H^1(E \cap B(0,r)) \leq 2 (1+\tau) r d(0)$ for $r\in {\cal R}_1$.
This is the same thing as (6.2), except that here we did not set $r=1$.
Thus we checked that for $r \in {\cal R}_1 \cap (r_0/10,2r_0)$,
the analogue when $r \neq 1$ of the assumptions (4.1)-(4.4), and (6.2) are 
satisfied, that were needed to construct the competitors $F_1$ and 
$F_2$ in Sections 6-10. This allows us to apply Lemma 10.2 to the set
$r^{-1} E$, with $A = C^{-1} \alpha_+(\varphi)^2$ coming from Lemma 10.23, 
and where $\varphi$ is still as in Lemma 6.11 and the definition of 
the cone $\rho = \varphi_\ast(K)$. 
Recall that Lemma~10.2 itself does not use (10.1) or the full length
assumption (see below its statement). Now (10.5) yields
$$\eqalign{
H^2(r^{-1} &E \cap B) 
\leq {1 \over 2} H^1(r^{-1} E\cap \partial B)
\cr&- 10^{-5} [H^1(r^{-1}E\cap \partial B) - H^1(\varphi_\ast(K))] 
- \kappa^2 A  
+ 4 h(2r)
}\leqno (11.14)
$$
or equivalently
$$ \leqalignno{
\hskip 1cm H^2&(E\cap B(0,r)) 
\leq {r \over 2} H^1(E\cap \partial B(0,r))
\cr&- 10^{-5} r [H^1(E\cap \partial B(0,r)) - r H^1(\varphi_\ast(K))] 
- C^{-1} \alpha_+(\varphi)^2 r^2 + 4 r^2 h(2r),
& (11.15)
}$$
where we drop the dependence on $\kappa$ from our notation. 
Define the $\delta_i$, $1 \leq i \leq 3$, as in (9.65)-(9.67), but for 
the set $r^{-1} E$. Then $H^1(E\cap \partial B(0,r)) - r H^1(\varphi_\ast(K))
= r(\delta_1 + \delta_2 + \delta_3)$ by (9.68) and (9.63), and 
(11.15) says that
$$ \leqalignno{
H^2(E\cap B(0,r)) 
&\leq {r \over 2} H^1(E\cap \partial B(0,r)) 
- 10^{-5} r^2 [\delta_1 + \delta_2 + \delta_3] 
- C^{-1} r^2 \alpha_+(\varphi)^2 + 4 r^2 h(2r)
\cr& \leq r^2 j(r) + r^2 d(0)
- 10^{-5} r^2 [\delta_1 + \delta_2 + \delta_3] 
- C^{-1} r^2 \alpha_+(\varphi)^2 + 4 r^2 h(2r)
& (11.16)
}$$
by (11.12). Thus
$$\leqalignno{
10^{-5} (\delta_1 + \delta_2 + \delta_3) 
+C^{-1} \alpha_+(\varphi)^2
&\leq j(r) + d(0) - r^{-2} H^2(E \cap B(0,r)) 
\cr&
= j(r) + d(0) - \theta(r) + 4 h(2r)
& (11.17)
}$$
and, since $\theta(r) \geq d(0) - C h_1(r)$ by (3.8), we get that
$$
\delta_1 + \delta_2 + \delta_3 + \alpha_+(\varphi)^2 \leq C j(r).
\leqno (11.18) 
$$

Let us first take care of the union $g = \cup_{j,k} g_{j,k}$,
where the $g_{j,k}$ are the simple arcs constructed in in Section 6
for $r^{-1}E$. For the moment we do this for a fixed $r \in {\cal R}_1$;
later on we shall connect the information coming from different radii, and 
get control on $E$ itself.

\medskip \proclaim Lemma 11.19.
There is a cone $X(r) \in {\cal Z}_0(X,C\tau)$ centered at $0$, 
such that $\alpha_+(Z) \leq C j(r)^{1/2}$ and
$$\eqalign{
\dist(x,g) \leq C &j(r)^{1/2}
\ \hbox{ for } x\in X(r) \cap \partial B
\ \hbox{ and }
\cr&
\dist(x,X(r) \cap \partial B) \leq C j(r)^{1/2}
\ \hbox{ for } x\in g.
}\leqno (11.20)
$$

\medskip 
Here $X$ is as in (11.7), and $\tau$ is chosen as in Section 6,
and in particular Lemma~6.12. For the proof, still denote by $\rho_{j,k}$ 
the geodesic of $\partial B$ with the same endpoints as $g_{j,k}$.
Also call $\rho$ the union of the $\rho_{j,k}$, as in (9.62). 
Obviously $\length(g_{j,k}) \geq \length(\rho_{j,k})$, but 
$$\eqalign{
\sum_{j} [\length(g_{j,k}) - \length(\rho_{j,k})]
&\leq H^1(g) - H^1(\rho) = \delta_1 + \delta_2
\cr& \leq \delta_1 + \delta_2 + \delta_3
\leq C j(r)
}\leqno (11.21)
$$
because the $g_{j,k}$ are simple and disjoint,
by (9.65) and (9.66), because $\delta_3 \geq 0$ by (9.67), 
and by (11.18).

It follows from (11.21) and a simple geometric argument 
(or brutally (7.9), the definitions (7.4) and (7.5), the fact that
$v(0)=0$, and Cauchy-Schwarz) that
$$\eqalign{
\dist(x,g) \leq C j(r)^{1/2}
&\ \hbox{ for } x \in \rho 
\ \hbox{ and }
\cr&
\dist(x,\rho) \leq C j(r)^{1/2}
\ \hbox{ for } x\in g.
}\leqno (11.22)
$$

Let $X(r)$ be the cone over $\rho$. 
Then (11.20) follows from (11.22). Also, $X(r)$ lies in 
${\cal Z}_0(X,C\tau)$ by Lemma 6.11 and the definition of
${\cal Z}_0(X,\tau)$ at the beginning of this section.
In addition, $\alpha_+(Z) = \alpha_+(\varphi)$ (see the definition 
above (11.2)), so $\alpha_+(Z) \leq C j(r)^{1/2}$ by (11.18).
Lemma 11.19 follows.
\qed

\medskip
Lemma 11.19 gives a good control on $g$ for each $r\in {\cal R}_1$.
We also control $E\cap \partial B(0,r) \setminus rg$, because
$$
H^1(r^{-1}E\cap \partial B \setminus g) 
= \delta_3 \leq C j(r)
\leqno (11.23)
$$
by (9.67) and (11.18) (recall that the three $\delta_j$ are 
nonnegative). 

\ms
Now we want to show that the minimal set $X(r)$ 
does not depend too much on $r$. For this we shall need to build transverse 
curves that are almost radial and meet $r g = rg(r)$ for different values 
of $r$. We shall first restrict to $r\in [r_0/3,r_0]$, 
where $r_0$ is as in the statement of Theorem 11.4. 
The construction of transverse curves will take some time.

We keep the minimal set $X$ of (11.7) and the decomposition of
$X\cap \partial B(0,1)$ into curves $\C_{j,k}$, and
denote by $\C_{j,k}'$ the arc of geodesic contained in $\C_{j,k}$ 
with the same center as $\C_{j,k}$, but only half the length. 
Then let $z \in \C_{j,k}'$ be given, and denote by $P_z$ the vector 
hyperplane through $z$ that is  perpendicular to $\C_{j,k}$ at $z \,$;
we want to construct a curve in $E \cap P_z$, that crosses 
the annulus $\overline B(0,r_0)\setminus B(0,r_0/3)$, but  let
us start with a connected set. The reader may start looking at 
Figure 11.1 for a picture.

\proclaim
Lemma 11.24. There is a compact connected $G_z \i E \cap P_z$,
which is contained in a $C \varepsilon r_0$-neighborhood
of the segment $I_z = [r_0 z/5, 3r_0 z/2]$,
and that contains points in $\partial B(0,r_0/3)$ and in
$\partial B(0,r_0)$.

\ms
We proceed as in Section 6, for the construction of the curves
$g_{j,k}$. By (11.7) and (11.6) or (11.8), and if $\varepsilon_1$ 
and $\varepsilon$ are chosen small enough, we can apply Lemma 16.19 
or Lemma~6.25 or Lemma~16.56 in [D3],     
and we get that $B(0,r_0)$ is a biH\"older ball for $E$.
[We used $60r_0$ instead of $50 r_0$ in the statement above so that 
we could bound $h(100r_0)$ in terms of $h_1(50r_0)$.]
This gives a biH\"older mapping $f : B(0,2r_0) \to \R^n$,
with the properties (6.16)-(6.19), with a small constant
$\tau > 0$ that will be chosen soon, but where we should adapt
the estimates to the scale $r_0$. We can keep the same cone
$X$ as in (11.7).

Denote by $H_0$ the cone over $\C_{j,k}$, and set
$H = H_0 \cap \overline B(0,5r_0/4) \setminus B(0,r_0/4)$. This is a
vaguely rectangular plane domain. Denote by $L_1$
and $L_2$ the two straight parts of its boundary, obtained
as the product by $[r_0/4,5r_0/4]$ of one of the two endpoints of
$\C_{j,k}$. Note that $f(L_1)$ and $f(L_2)$ lie on different
sides of the hyperplane $P_z$, by (6.16). 

Now set $F = H \cap f^{-1}(E \cap P_z) \i B(0,2r_0)$; then $F$ separates 
$L_1$ from $L_2$ in $H$, because if $\gamma$ is a path from 
$L_1$ to $L_2$ in $H$, $f(\gamma)$ must cross $P_z$. 
By 52.III.1 in [Ku] 
or the simpler Theorem~14.3 in [Ne] 
(see Remark 6.45), there is a connected subset of $F_1$ of 
$F$ that separates $L_1$ from $L_2$ in $H$. In particular,
it meets $\partial B(0,r_0/3)$, because otherwise we could use 
$H \cap \partial B(0,r_0/3)$ to connect $L_1$ to $L_2$
in $H\setminus F_1$. Similarly, $F_1$ meets $\partial B(0,4r_0/3)$.

Now set $G_z = f(F_1)$. By definition of $F$, it is contained
in $E \cap P_z$. Obviously, it is connected too. By what we just said 
and (6.16), it contains points in $\overline B(0,r_0/4 + \tau r_0)$ and 
out of $B(0,5r_0/4 - \tau r_0)$; by connectedness, it meets 
$\partial B(0,r_0/3)$ and $\partial B(0,r_0)$ too. 

We still need to check that $\dist(y,I_z) \leq C \varepsilon r_0$
for $y\in G_z$. Let $x\in F_1$ be such that $f(x)=y$.
First observe that $\dist(x,P_z) \leq \tau r_0$, by (6.16) and because 
$y \in P_z$. Since $x\in H$ and $P_z$ is the hyperplane perpendicular to 
$\C_{j,k}$ at $z$, this forces $x$ to be within $\tau r_0$ of the line segment
$[r_0z/4,5r_0z/4]$. Then $\dist(y,I_z) \leq 2 \tau r_0$, by (6.16).

Now recall that we took $z\in \C_{j,k}'$, i.e., far from the
extremities of $\C_{j,k}$; thus, if $\tau$ is small enough compared 
to the constant $\eta_0$ in (2.2) and (2.3), the description of $X$
in Section~2 says that, in the region where 
$\dist(y,I_z) \leq 3 \tau r_0$, $X$ coincides with the plane $H'$ 
that contains $\C_{j,k}$. Since 
(11.7) says that $\dist(y,X) \leq 100 \varepsilon r_0$,
we get that $\dist(y,H') \leq 100 \varepsilon r_0$,
and then $\dist(y,I_z) \leq 100 \varepsilon r_0$ because
$y\in P_z$. Lemma 11.24 follows.
\qed

\vskip 0.5cm  \hskip 2.5cm  
\epsfxsize = 7.5cm \epsffile{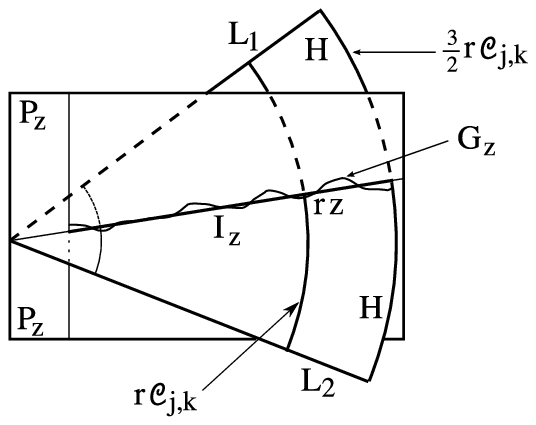}
\medskip
\noindent
\centerline{\bf Figure 11.1.} 
\medskip 

\ms
For most choices of $z\in \C_{j,k}'$, we want to find a nice simple
curve in $G_z$, show that it stays very close to a radius, and
at the same time that for many radii $r$ it meets the curve $rg(r)$
studied above (the same one as in Lemma 11.19). This will then make it easy 
to show that the cones $X(r)$ are close to each other. 
The acceptable choices of $z$ and $r$ will come from Chebyshev, 
and before we do this we need to control the averages of certain quantities.

\medskip\proclaim Lemma 11.25. Set $A = B(0,5r_0/3) \setminus B(0,r_0/6)$ 
and, for almost every $x\in E \cap  A$, denote by $\alpha(x) \in [0,\pi/2]$ 
the angle of the radius $[0,x]$ with the tangent plane to $E$ at $x$. Then
$$
\int_{x \in E \cap A} [1 - \cos\alpha(x)] \, dH^2(x) 
\leq C r_0^2 [f(2r_0)+ h_1(2r_0)].
\leqno (11.26)
$$

\medskip  
We want to prove this with the co-area theorem.
Recall that $E$ is rectifiable, and that it even has a true
tangent plane $P(x)$ at $x$ for $H^2$-almost-every $x\in E$.
(See the discussion after the statement of Lemma 4.12). 
So $\alpha(x)$ is defined almost-everywhere, and 
(4.13) and (4.14), applied with $f = {\bf 1}_A$, say that
$$
\int_{E\cap A} \cos\alpha(x) \, dH^2(x) 
= \int_{r_0/6 \leq r < 5r_0/3} H^1(E \cap \partial B(0,r)) \, dr.
\leqno (11.27)
$$
Let us check that
$$
H^1\big((r_0/6,5r_0/3)\setminus {\cal R}_1\big)
\leq C r_0 [f(2r_0)+ h_1(2r_0)],
\leqno (11.28)
$$
where $C$ may depend on $\tau$ (through (11.13)), 
but this does not matter any more because we shall no longer need to 
modify our choice of $\tau$. Since $\cal R$ has full measure, the 
definition (11.13) of ${\cal R}_1$ says that
$H^1\big((r_0/6,5r_0/3)\setminus {\cal R}_1\big)
\leq \tau^{-1} d(0)^{-1} \int_{[r_0/6,5r_0/3]} j(r) dr$, 
so it is enough to show that
$$
\int_{[r_0/6,5r_0/3]} j(r) dr \leq C r_0 [f(2r_0)+ h_1(2r_0)].
\leqno (11.29)
$$

Recall from (11.11) that $j(r) = rf'(r)+f(r)+h_1(r)+h(2r)$. 
There is  no problem with the contribution of $h_1(r)+h(2r)$,
because $h_1(r)+h(2r)\leq h_1(2r_0)$. Next consider $f(r)$.
Notice that for $r \leq 5r_0/6$,
$$
f(r) \leq f(2r_0) + C h_1(2r_0) 
\leqno (11.30)
$$
by (3.6) and because $B(0,r_0) \i U$, so the contribution of $f(r)$ 
in (11.29) is under control. 
We are left with with the contribution of $r f'(r)$. 
Observe that
$$\eqalign{
\int _{r_0/6 \leq r < 5r_0/3} &f'(r) \, dr
=\int _{r_0/6 \leq r < 5r_0/3} \theta'(r) \, dr
\leq \theta(2r_0) - \theta(r_0/6)
\cr&
= f(2r_0) - f(r_0/6) \leq f(2r_0) + C h_1(r_0/6)
\leq f(2r_0) + C h_1(2r_0)
}\leqno (11.31)
$$
by (1.3), (5.7), (3.8), and because $h_1$ is nondecreasing;
(11.29) and (11.28) follow at once.

\ms
Let us now focus on ${\cal R}_1$. For $r\in {\cal R}_1$, we can apply
to $r^{-1}E$ the construction (and notation) of the previous sections. Then
$$
r^2 \theta(r) = H^2(E \cap B(0,r)) 
\leq {r \over 2} \, H^1(E \cap \partial B(0,r))
+ 4 h(2r) \, r^2 
\leqno (11.32)
$$
by (9.69) and because $H^1(r^{-1}E\cap \partial B) - H^1(\rho)
= \delta_1 + \delta_2 +\delta_3 \geq 0$
(see (9.68) and the definitions above it). Equivalently,
$$
r H^1(E \cap \partial B(0,r))
\geq 2r^2 \theta(r) - 8 h(2r) \, r^2.
\leqno (11.33)
$$
But $\theta(r) \geq d(0) - C h_1(r) \geq
d(0) - C h_1(2r_0)$ by (3.8), and $h(2r) \leq h(10r_0/3) \leq C h_1(2r_0) = 
\int_0^{2r_0} h(2t) {dt \over t}$, because $h$ is nondecreasing and
by the definition (1.16) or (3.1) of $h_1$. Then 
$$
r H^1(E \cap \partial B(0,r)) \geq 2 d(0) \, r^2 - C h_1(2r_0) \, r^2.
\leqno (11.34)
$$
Thus (11.27) says that
$$\leqalignno{
\int_{E\cap A} \cos\alpha(x) &\, dH^2(x) 
\geq \int_{r\in {\cal R}_1 \cap (r_0/6,5r_0/3)} 
H^1(E \cap \partial B(0,r)) \, dr
\cr&
\geq \int_{r\in {\cal R}_1 \cap (r_0/6,5r_0/3)} 
[2 d(0) - C h_1(2r_0)] \, r dr
\cr& \geq d(0) \, \big({25 r_0^2 \over 9} -{r_0^2 \over 36}\big)  
- 2 d(0) \, r_0 H^{1}((r_0/6,5r_0/3) \setminus {\cal R}_1)
- C h_1(2r_0) \, r_0^2  &(11.35)
\cr& 
\geq d(0) \, \big({25 r_0^2 \over 9} -{r_0^2 \over 36}\big) 
- C  [f(2r_0)+ h_1(2r_0)] \, r_0^2
}
$$
by (11.28). On the other hand, 
$$\leqalignno{
H^2(E\cap A) 
&= {25 r_0^2 \over 9} \, \theta(5r_0/3) 
- {r_0^2 \over 36} \, \theta(r_0/3) 
= {25 r_0^2 \over 9} \, [d(0) + f(5r_0/3)] 
- {r_0^2 \over 36} \, \theta(r_0/3) 
\cr&
\leq {25 r_0^2 \over 9} \, [d(0) + f(5r_0/3)] 
- {r_0^2 \over 36} \, [d(0)-Ch_1(r_0)] 
\cr&
=  d(0) \, \big({25 r_0^2 \over 9} -{r_0^2 \over 36}\big)  
+ 3 r_0^2 f(2r_0) + C \, r_0^2 h_1(2r_0)   & (11.36)
\cr&
\leq \int_{E\cap A} \cos\alpha(x) \, dH^2(x) 
+ C [f(2r_0)+ h_1(2r_0)] \, r_0^2
}
$$
by (1.1), (3.8), (3.6), and (11.35). This  proves (11.26) and 
Lemma 11.25. \qed

\ms
Return to the compact connected set $G_z$ of Lemma 11.24.
Still fix $j$ and $k$, and for $z\in \C_{j,k}'$, choose $G_z$ as in 
Lemma 11.24. Also let $r_1$ and $r_2$
be given, with $r_0/3 \leq r_1 < r_2 \leq r_0$, and
pick points $y_1(z)\in G_z \cap \partial B(0,r_1)$ and 
$y_2(z)\in G_z \cap \partial B(0,r_2)$. Such points
exist, because $G_z$ is connected and contains points
of $\partial B(0,r_0/3)$ and $\partial B(0,r_0)$.
The following lemma says that $y_1$ and $y_2$ are 
often nearly aligned.

\medskip\proclaim Lemma 11.37.
There is a measurable function $q$ defined on $\C_{j,k}'$, such that
$$
\int_{\C_{j,k}'} q(z) dH^1(z) \leq C  [f(2r_0)+h_1(2r_0)]^{1/2}
\leqno (11.38)
$$
and
$$
\left| {y_1 \over |y_1|} - {y_2 \over |y_2|} \right|
\leq q(z)
\leqno (11.39)
$$
for all choices of $z\in \C_{j,k}'$, $r_1$, $r_2$, $y_1$, 
and $y_2$ as above.

\ms
Set $S_z = \big\{ y\in P_z \, ; \, \dist(y,I_z) \leq C \varepsilon r_0 
\big\}$, where $P_z$, $I_z$ and $C$ are as in Lemma 11.24; thus
$G_z \i S_z$ for all $z$. Notice also that $S_z \i A$, where 
$A$ is as in Lemma 11.25. 
Let us first check that
$$
H^1(G_z) \leq H^1(E \cap S_z) < +\infty \hbox{ for almost every choice of }
z\in \C_{j,k}'.
\leqno (11.40)
$$

Denote by $S$ the union of the $S_z$, $z\in \C_{j,k}'$, 
and define a Lipschitz mapping 
$\varphi : E \cap S \to \C_{j,k}'$ by $\varphi(x)=z$ when 
$x\in S_z$. 
We apply the co-area theorem (see 3.2.22 in [Fe]) to $\varphi$ 
(or its composition with the inverse of a parameterization of 
$\C_{j,k}'$ by arc-length, if you want a mapping to an interval) 
and get that for any Borel set $T \i E\cap S$,
$$
\int_{T} J_\varphi(x) \, dH^2(x) 
= \int_{z \in \C_{j,k}'} H^1(T \cap \varphi^{-1}(z)) \, dH^1(z)
= \int_{z \in \C_{j,k}'} H^1(T \cap S_z) \, dH^1(z),
\leqno (11.41)
$$
where $J_\varphi(x)$ is the appropriate directional Jacobian.
Since $\varphi$ is $6r_0^{-1}$-Lipschitz on $S$,
$J_\varphi(x) \leq 6r_0^{-1}$, and (11.41) with $T=E\cap S$ 
immediately yields
$$
 \int_{z \in \C_{j,k}'} H^1(E \cap S_z) \, dH^1(z) 
\leq 6 r_0^{-1} H^2(E \cap S) \leq C r_0
\leqno (11.42)
$$
(by local Ahlfors-regularity or more brutally because $f(60r_0) \leq 
\varepsilon$ by (11.8)).
This proves (11.40) (because we already knew its first part).
Thus we can safely take $q(z)=+\infty$ when
$H^1(E \cap S_z)=+\infty$ and restrict to the other case.

Since $G_z$ is connected and $H^1(G_z) < + \infty$,
there is a simple curve $\zeta \i G_z$ that goes from $y_{1}$ to $y_{2}$.
See for instance Proposition 30.14 in [D2].  
Denote by $\xi : I \to K$ a parameterization of $\zeta$ by arc-length;
such a parameterization exists, because 
$\length(\zeta) \leq H^1(G_z)< +\infty$.

Set $v(x)=x/|x|$ for $x\neq 0$. We are interested in the 
variations of  $v \circ \xi$. 
The differential of $v$ in the unit direction $a$ is
$Dv(x)\cdot a = |x|^{-1}\,\widetilde a$, where $\widetilde a$ is the 
projection of $a$ on the plane orthogonal to $x$. 
For almost every $s\in I$, $\xi'(s)$ exists, and is a unit vector. 
Then $|(v \circ \xi)'(s)| = |\xi(s)|^{-1} \sin\gamma(s)
\leq 6 r^{-1}_0 \sin\gamma(s)$, where $\gamma(s) \in [0,\pi/2]$ 
is the angle of $\xi'(s)$ with the radius $[0,\xi(s)]$. Hence
$$
\left| {y_1 \over |y_1|} - {y_2 \over |y_2|} \right|
=|v(y_{1})-v(y_{2})| \leq 6 r^{-1}_0 \int_{I} \sin\gamma(s) ds.
\leqno (11.43)  
$$

Next we want to bound the right-hand side of (11.43)  by integrals on
$E\cap S_z$, and then average over $z$. It will be more convenient to
replace $\sin\gamma(s)$ with quantities that depend only on 
$x=\xi(s)$.

We need to distinguish different types of points in $I$. 
We start with the good set $I_{1}$ of points $s\in I$ with the
following properties. First, $\xi'(s)$ exists and is a unit vector.
Next, $E$ has a (true) tangent $2$-plane $P(x)$ at $x=\xi(s)$. 
Also, $r = |x|$ lies in ${\cal R}_1$, which allows us to do 
the construction of Sections 6-9. 
Finally, we require the existence of a unit vector $v$
in the direction of $P(x)$, which is orthogonal to $e = x/|x|$,
and whose orthogonal projection on $P_z$ has a norm less than
$10^{-1}$. We shall see soon that this happens often.
Let us check that 
$$
\sin\gamma(s) \leq  10 \, \alpha(x) \hbox{ when $s \in I_{1}$,}
\leqno (11.44)
$$ 
where $\alpha(x)$ is the angle of $[0,x]$ with $P(x)$, as in 
(11.26) and (11.27). 

We may assume that $\alpha(x) \leq 10^{-1}$, because otherwise
(11.44) is trivial.
Denote by $e_2$ any of the two unit vectors
perpendicular to $P_z$. Set $V=P(x)-x$ (the vector plane 
parallel to $P(x)$). By definition of $\alpha(x)$,
we can find a unit vector $v_{1}$ in $V$, with 
$$
|v_{1}-e| \leq \alpha(x) \leq 10^{-1}.
\leqno (11.45)
$$ 
Notice that since $v_1$ is close to $e$ by (11.45), it is nearly 
orthogonal to $v$ (the unit vector whose existence was required 
in the definition of $E_1$). Thus $(v_1,v)$ is a basis of $V$, and the
projection in $V$ on the line through $v_1$ and parallel to $v$ 
has a norm smaller than $2$.

Recall that $\xi'(s)$ exists and is a unit vector, by definition
of $I_1$. Then $\xi'(s)\in V$, because $\xi(t) \in E$ for all $t$.
So we can write $\xi'(s) = a_1 v_1 + a_2 v$. Recall that 
the norm of the projection on the line through $v_1$ parallel to $v$ 
is less than $2$, and hence $|a_1| \leq 2$. 
We know that $\xi'(s) \cdot e_2=0$, because $\xi'(s) \in P_z$
(since the whole curve lies in $P_z$) and $e_2$ is orthogonal to $P_z$.
Then $a_1 v_1 \cdot e_2 + a_2 v \cdot e_2 = 0$. In addition,
$|v \cdot e_2| \geq 2/3$ because the orthogonal projection
of $v$ on $P_z$ has a norm smaller than $10^{-1}$, and $e_2$ spans the orthogonal
direction to $P_z$. Thus
$$\eqalign{
|a_2| &= |a_2 v \cdot e_2|/|v \cdot e_2|
\leq {3 \over 2} \, |a_2 v \cdot e_2|
= {3 \over 2} \, |a_1 v_1 \cdot e_2|
\cr&
\leq 3 |v_1 \cdot e_2| = 3 |(v_1-e)\cdot e_2| \leq 3 \alpha(x),
}\leqno (11.46)
$$
because $e \cdot e_2 = 0$ (since $x\in P_z$), and by (11.45). 
Now recall that $\gamma(s)$ is the angle of $\xi'(s)$ with 
the direction of $\xi(s) = x$  and $e = x/|x|$. Thus 
$\sin\gamma(s) = |\pi(\xi'(s))|$, where $\pi$ denotes the orthogonal
projection on $e^\perp$. Now
$$\eqalign{
\sin\gamma(s) &= |\pi(\xi'(s))| = |\pi(a_1 v_1 + a_2 v)|
\leq |a_1||\pi(v_1)| + |a_2|
\cr&  
\leq 2 |\pi(v_1)| + 3 \alpha(x)
= 2 |\pi(v_1- e)| + 3 \alpha(x)
\leq 5 \alpha(x)
}\leqno (11.47)
$$
because $\xi'(s) = a_1 v_1 + a_2 v$ and $|a_1|\leq 2$,
by (11.46), and by (11.45). Thus (11.44) holds in all cases.

\ms
There is another part of $I$ which will cause no problem;
this is the set $I_2$ of points $s\in I\setminus I_1$ such that 
$\alpha(\xi(s))$ is defined, and $\alpha(\xi(s)) \geq 10^{-1}$.
Notice that (11.44) (trivially) holds for $I_2$. Hence (11.43) yields
$$\eqalign{
\left| {y_1 \over |y_1|} - {y_2 \over |y_2|} \right|&
\leq 6 r_0^{-1} \int_{I} \sin\gamma(s) ds
\cr&
\leq 6 r_0^{-1} H^1(I\setminus [I_1\cup I_2]) 
+ 60 r_0^{-1} \int_{I_1 \cup I_2} \alpha(\xi(s)) ds.
}\leqno (11.48) 
$$

\ms
Now let us chase various exceptional sets. First, the set $I_3$ where
$\xi'(s)$ does not exist, or is not a unit vector, has vanishing 
measure, so we may drop it.

Next we need to control some exceptional sets in $E$,
before we return to $I$. Set $W_4 = \{ x \in E \cap S \, ; \, 
E \hbox{ has no tangent plane at } x \}$.
Notice that $H^2(W_4) = 0$, because $E$ is rectifiable 
and locally Ahlfors-regular. See Lemma~2.15 in [D3] 
and Exercise 41.21 in [D2].  

Set $W_5 = \{ x \in E \cap S \setminus W_4 \, ; \,
r \notin {\cal R}_1 \hbox{ and } \alpha(x) \leq 10^{-1}\}$. 
Notice that $\cos\alpha(x) \geq \cos(10^{-1}) \geq 1/2$ on $W_5$, so
$$\leqalignno{
H^2(W_5) &\leq 2 \int_{W_5} \cos\alpha(x) dH^2(x)
= 2 \int_{W_5} J_g(x) dH^2(x)
\cr&
\leq \int_{[r_0/6,5r_0/3] \setminus {\cal R}_1} 
H^1(E \cap \partial B(0,r)) \, dr
\leq \int_{[r_0/6,5r_0/3] \setminus {\cal R}_1} 2 r [j(r) + d(0)] \, dr
\cr&
\leq 2r_0 \, d(0) \, H^1([r_0/6,5r_0/3] \setminus {\cal R}_1)
+ 2r_0 \int_{[r_0/6,5r_0/3]} j(r)  \, dr &(11.49)
\cr&
\leq C r_0^2 [f(2r_0)+h_1(2r_0)],
}
$$
by (4.14) and the analogue of (4.13) for Borel subsets of $E$, 
because $r_0/6 \leq |x| \leq 5r_0/3$ when $x\in S_z$, by
the fact that $\cal R$ has full measure, and by (11.12), (11.28), and 
(11.29).

Finally let $W_6$ denote the set of points 
$x\in E \cap S\setminus W_4$ such that $\alpha(x) \leq 10^{-1}$
and $r=|x|\in {\cal R}_1$, but we cannot find a unit vector
$v$ in the direction $V$ of $P(x)$, as in the definition of
$I_1$. Then 
$$\eqalign{
H^2(W_6)  &
\leq 2\int_{W_6} \cos\alpha(x) dH^2(x) 
= 2 \int_{W_6} J_g(x) dH^2(x)
\cr&
\leq 2 \int_{[r_0/6,5r_0/3] \cap {\cal R}_1} H^1(W_6 \cap \partial B(0,r))
\, dr
}\leqno (11.50)
$$
because $\cos\alpha(x) \geq 1/2$ on $W_6$, and by (4.14) and (4.13).
Let us fix $r$, and prove that 
$$
H^1(W_6 \cap \partial B(0,r) \cap \Gamma_{j,k}') = 0
\leqno (11.51)
$$
where we set $\Gamma_{j,k}' = r \Gamma_{j,k}$ and
$\Gamma_{j,k}$ is the Lipschitz curve in $\partial B$
that was constructed in Section 7 (so that 
$\Gamma_{j,k}'\i \partial B(0,r)$). Indeed, for almost every point 
$x \in W_6 \cap \Gamma_{j,k}'$, $\Gamma_{j,k}'$ has a tangent line $L$ 
at $x$, and $x$ is a Lebesgue density point of $E \cap \Gamma_{j,k}'$
(recall that $W_6 \i E$).
Then $L$ is contained in $P(x)$ (which exists because 
$x\in W_6 \i E \cap S \setminus W_4$). Let $v$ be a unit
vector in the direction of $L$; we claim that it satisfies
the requirements in the definition of $I_1$ (and then $x\notin W_6$, 
a contradiction which proves (11.51)). We already know that 
$v\in V$, because $L \i P(x)$. Then $v$ is orthogonal to the
direction of $x$, because $\Gamma_{j,k}' \i \partial B(x,r)$.

We still need to check that $v$ is almost orthogonal to
$P_z$ (where $z\in \C_{j,k}'$ is the one for which $x\in S_z$). 
Without loss of generality, we can assume that $z=(1,0,\ldots 0)$
and $P_z$ is the vertical hyperplane plane $\{ x_2 = 0 \}$. 
This also means that our arc of geodesic $\C_{j,k}$ is contained 
in the horizontal $2$-plane $\{ x_3 = \ldots = x_{n} = 0 \}$.

By construction, $\Gamma_{j,k}'$ is a small Lipschitz graph over 
the horizontal $2$-plane which contains $\C_{j,k}$. This means that 
all the coordinates of $v$, except perhaps the first two, are 
smaller than $C\eta$ (see for instance (7.45)). The first
coordinate is less than $C\varepsilon$, because $v$ is orthogonal 
to $e = x/|x|$ and $e$ is very close to $z=(1,0,\ldots 0)$,
because $x \in S_z$ lies very close to $I_z = [r_0 z/6,5r_0z/3]$
(see Lemma 11.24 and the definition of $S_z$ above (11.40)).
So $v$ has a small projection on $P_z$, $x\notin W_6$,
and (11.51) holds. Hence (11.50) yields
$$
H^2(W_6) \leq 2\int_{[r_0/6,5r_0/3] \cap {\cal R}_1} 
H^1(W_6 \cap \partial B(0,r) \setminus \Gamma_{j,k}') \, dr.
\leqno (11.52)
$$
Let $x\in W_6 \cap \partial B(0,r) \setminus \Gamma_{j,k}'$
be given. Then $x \in S = \cup_{z\in \C_{j,k}'} S_z$, which 
is contained in the cone over a small neighborhood of $\C_{j,k}'$.
That is, $x/r$ lies in a small neighborhood of $\C_{j,k}'$.
Also recall that $\C_{j,k}'$ itself lies near the middle
of $\C_{j,k}$, hence reasonably far from the other $\C_{j',k'}$
(by (2.2) and (2.3)). Then $x/r$  does not lie in any other
$\Gamma_{j',k'}$, because $\Gamma_{j',k'}$ stays within
$C \eta$ of $\C_{j',k'}$ (because it is a small Lipschitz graph
with endpoints on $g_{j',k'}$, and by (6.15)), and if $\eta$
is small enough compared to $\eta_0$ in (2.2) and (2.3).
So $x \notin \Gamma' = \cup_{j',k'} \Gamma_{j',k'}'$. 
Recall that $W_6 \i E$, so
$$\eqalign{
H^1(W_6 \cap \partial &B(0,r) \setminus \Gamma_{j,k}') 
\leq H^1(E\cap \partial B(0,r) \setminus \Gamma')
= r H^1(r^{-1}E\cap \partial B \setminus \Gamma)
\cr&
\leq  r H^1(r^{-1} E\cap \partial B \setminus g)
+ r H^1(g \setminus \Gamma)
= r\delta_3 + r\sum_{j,k} H^1(g_{j,k} \setminus \Gamma)
\cr&
\leq r\delta_3 + r\sum_{j,k} H^1(g_{j,k} \setminus \Gamma_{j,k})
\cr&
\leq r\delta_3 + C r \sum_{j,k} [\length(g_{j,k})-\length(\rho_{j,k})]
\leq C r j(r),
}\leqno (11.53)
$$
where $g = \cup_{j,k} g_{j,k}$ and the $g_{j,k}$ are the curves
of Section 6, by (9.67), because the $g_{j,k}$
are simple and disjoint, by (9.57), and by (11.21).
Altogether
$$
H^2(W_6) 
\leq C \int_{[r_0/6,5r_0/3] \cap {\cal R}_1} r j(r) dr
\leq C r_0^2 \, [f(2r_0)+h_1(2r_0)]
\leqno (11.54)
$$
by (11.29).

Set $I_4 = I\setminus (I_1 \cup I_2 \cup I_3)$. If $s\in I_4$ and 
$x=\xi(s)$, $\xi'(s)$ exists and is a unit vector because $s\notin I_3$.
If $\xi(s) \notin W_4$, $P(x)$ is defined, $\alpha(x) \leq 10^{-1}$
because $s\notin I_2$, and then $\xi(s) \in W_5 \cup W_6$ because 
otherwise $s\in I_1$. Thus $\xi(I_4) \i W$, where $W=W_4 \cup W_5 \cup W_6$.

We may now return to (11.48); first observe that
$$
\int_{I_1 \cup I_2} \alpha(\xi(s)) ds
\leq \int_{E\cap S_z} \alpha(x) dH^1(x)
\leqno (11.55)
$$
because $\xi$ is a parameterization by arc-length of the simple 
curve $\zeta$, so the direct image of $ds$ on $I$ by $\xi$ is less 
than or equal to the restriction of $H^1$ to $E\cap S_z$.
Also, 
$$
H^1(I\setminus [I_1\cup I_2]) \leq H^1(I_3 \cup I_4) = H^1(I_4)
\leq H^1(S_z \cap W)
\leqno (11.56)
$$
by definition of $I_4$, because $H^1(I_3) = 0$,
because $\zeta(I_4) \i S_z \cap W$, and again because 
the image of $ds$ by $\xi$ is no greater than $H^1$.
Thus (11.48) yields
$$\eqalign{
\left| {y_1 \over |y_1|} - {y_2 \over |y_2|} \right|
\leq 6 r_0^{-1} H^1(S_z \cap W) 
+ 60 r_0^{-1} \int_{E\cap S_z} \alpha(x) dH^1(x).
}\leqno (11.57) 
$$

We take $q(z)=6 r_0^{-1} H^1(W \cap S_z) + 60 r_0^{-1} 
\int_{E\cap S_z} \alpha(x) dH^1(x)$; then (11.39) obviously holds,
and Lemma 11.37 will follow as soon as we prove (11.38).
First 
$$\eqalign{
\int_{z \in \C_{j,k}'} H^1(W \cap S_z) \, dH^1(z)
&= \int_{W} J_\varphi(x) \, dH^2(x) 
\cr&\leq 6 r_0^{-1} H^2(W) \leq C r_0 \, [f(r_0)+h_1(2r_0)]
}\leqno (11.58) 
$$
by (11.41) (with $T = W$), because $J_\varphi(x) \leq 6 r_0^{-1}$ 
almost-everywhere on $E$, because $H^2(W_4)=0$, and 
by (11.49) and (11.54) (recall that $W = W_4 \cup W_5 \cup W_6$). 
Similarly,
$$\leqalignno{
\int_{z \in \C_{j,k}'} &\int_{E\cap S_z} \alpha(x) \, dH^1(x) \, dH^1(z)
= \int_{E\cap S} J_\varphi(x) \, \alpha(x) \, dH^2(x) 
\cr&\leq 6 r_0^{-1} \int_{E\cap S} \alpha(x) \, dH^2(x)
\leq 6 r_0^{-1} H^2(E\cap S)^{1/2} \,
\Big\{ \int_{E\cap S} \alpha^2(x) \, dH^2(x) \Big\}^{1/2}
& (11.59) 
\cr&
\leq C \,\Big\{ \int_{E\cap S} [1- \cos\alpha(x)] \, dH^2(x) \Big\}^{1/2}
\leq C r_0 [f(2r_0)+ h_1(2r_0)]^{1/2}
}
$$
by the integral version of (11.41), Cauchy-Schwarz, (11.26), 
and because $S \i A$. Now 
$\int_{z \in \C_{j,k}'} q(z) \, dH^1(z) \leq C [f(2r_0)+ h_1(2r_0)]^{1/2}$,
by (11.58) and (11.59). This proves (11.38), and Lemma 11.37 follows.
\qed

\ms
We are now ready to show that the various cones
$X(r)$ are close to each other. We shall restrict our attention
to $r\in {\cal R}_2$, where 
$$ 
{\cal R}_2 = \big\{ r \in (r_0/3,r_0) \cap{\cal R} \, ; \, 
j(r) \leq [f(2r_0)+ h_1(2r_0)]^{2/3} \big\}.
\leqno (11.60) 
$$

\medskip \proclaim Lemma 11.61. We have that
$$
d_{0,1}(X(r_1),X(r_2)) \leq C [f(2r_0)+ h_1(2r_0)]^{1/3}
\ \hbox{ for } r_1, r_2 \in {\cal R}_2,
\leqno (11.62)
$$
where the local Hausdorff distance function $d_{0,1}$
is as in (11.3).

\ms
Let us first get rid of some exceptional sets.
Let $r \in {\cal R}_2$  be given; observe that
$r \in {\cal R}_1$ (if $\varepsilon$ is small enough;
see (11.8) and (11.13)), which allows the construction
of competitors the previous sections.
Set $E^\sharp = [r^{-1} E\cap \partial B \setminus g]
\cup [g \setminus \Gamma]$, 
where $g = \cup_{j,k} g_{j,k}$ is the union of the arcs
of Section 6, and $\Gamma = \cup_{j,k} \Gamma_{j,k}$ is  
the union of the Lipschitz arcs with the same endpoints
constructed in Section 7. Then
$$\eqalign{
H^1(E^\sharp) &
\leq \delta_3  + \sum_{j,k} H^1(g_{j,k} \setminus \Gamma)
\leq \delta_3  + \sum_{j,k} H^1(g_{j,k} \setminus \Gamma_{j,k})
\cr&
\leq \delta_3  + C\sum_{j,k} [\length(g_{j,k})-\length(\rho_{j,k})]
\leq \delta_3  + C [ H^1(g)-H^1(\rho)]
\cr&= \delta_3  + C [\delta_1 + \delta_2]
\leq C j(r)
\leq C [f(2r_0)+ h_1(2r_0)]^{2/3}
}\leqno (11.63)
$$
by (9.67), (9.57), (9.62) and the fact that the $\rho_{j,k}$
are disjoint, (9.65) and (9.66), and finally (11.18) and (11.60).

Denote by $Z_{j,k}(r)$ the set of points 
$z\in \C_{j,k}'$ such that the set $S_z \i P_z$ defined
above (11.40) meets $r E^\sharp$. 
Then $Z_{j,k}(r)$ is the image of $r E^\sharp\cap S_z$ 
by the $6r_0^{-1}$-Lipschitz projection $\varphi$ from $S$ to 
$\C_{j,k}'(r)$ (see below (11.40)), so
$$\eqalign{
H^1(Z_{j,k}(r)) \leq 6 r_0^{-1} H^1(r E^\sharp)
\leq C [f(2r_0)+ h_1(2r_0)]^{2/3}.
}\leqno (11.64)
$$

Let us prove that
for $z\in \C_{j,k}' \setminus Z_{j,k}(r)$,
$$
E \cap \partial B(0,r) \hbox{ meets $S_z$ exactly once, at a point
of $G_z$.} 
\leqno (11.65)
$$
Recall that $G_z$ is the connected subset of $E \cap P_z$
that shows up in Lemma 11.24. It is contained in $S_z$ (see below 
(11.39)) and meets $\partial B(0,r)$ because 
it is connected and meets $\partial B(0,r_0/3)$ and $\partial B(0,r_0)$.
So $E \cap \partial B(0,r)$ meets $S_z$ at least once (at a point of 
$G_z$), and we just need to check that the intersection is reduced 
to one point.  Let us first check that 
$$
E \cap \partial B(0,r) \cap S_z \i \Gamma'_{j,k},
\leqno (11.66)
$$
where again $\Gamma'_{j,k} = r \Gamma_{j,k}$.
Let $\xi' \in E \cap \partial B(0,r) \cap S_z$ be given
and set $\xi = \xi'/r$. 
Since $z\in \C_{j,k}' \setminus Z_{j,k}(r)$, $\xi$ lies out of
$E^\sharp$, and hence $\xi \in \Gamma=\cup_{(j',k') \in \widetilde J} 
\ \Gamma_{j',k'}$. 
We still need to check that $\xi$ lies in no other $\Gamma_{j',k'}$. 
Recall that the endpoints of $\Gamma_{j',k'}$ lie in
$g_{j',k'}$, which lies very close to $\C_{j',k'}$ by (6.15);
then the whole $\Gamma_{j',k'}$ lies close to 
$\C_{j',k'}$, because $\Gamma_{j',k'}$ is a Lipschitz graph
with small constant (see the first lines of Section 9, or
go to (7.45) and where $\pi$ is defined above (7.35)), 
and because we have a good control on the $2$-plane
that contains the origin and the endpoints of $\Gamma_{j',k'}$
(since the length of $\C_{j',k'}$ is at most $9\pi/10$ by (2.5)).
Now $\xi' \in S_z$, hence $\xi$ lies close to $\C_{j,k}'$, so $\xi$ lies 
far from the other $\C_{j',k'}$, because $\C_{j,k}'$ lies in the middle of 
$\C_{j,k}$ and by (2.2). We checked that $\xi \notin \Gamma_{j',k'}$ 
for $(j',k') \neq (j,k)$, and so (11.66) holds.

Now $\Gamma'_{j,k}$ is a Lipschitz graph with small constant  
over the $2$-plane that contains its endpoints (by (7.45)), so it cannot
cross the vertical hyperplane $P_z$ more than once; (11.65) follows.

\ms
Return to Lemma 11.61 and let $r_1, r_2 \in {\cal R}_2$ be given.
For each choice of $(j,k) \in \widetilde J$, we want choose three points $z_i$
in $\C_{j,k}' \sm [Z_{j,k}(r_1) \cup Z_{j,k}(r_2)]$. We cut
$\C_{j,k}'$ into three equal intervals $I_1$, $I_2$, and $I_3$,
with $I_1$ in the middle; thus $H^1(I_i) \geq \eta_0$ and
$\dist(I_1,I_3) \geq \eta_0$, where $\eta_0$ is still as in (2.2). 
By (11.64) (and (11.8)), we still have that
$H^{1}(I_i \setminus [Z_{j,k}(r_1) \cup Z_{j,k}(r_2)]) \geq 
\eta_0/2$. We choose 
$z_i \in I_i \setminus [Z_{j,k}(r_1) \cup Z_{j,k}(r_2)]$
such that
$$
q(z_i) \leq {2 \over \eta_0} \, \int_{\C_{j,k}'} q(z) \, dH^1(z) 
\leq C  [f(2r_0)+h_1(2r_0)]^{1/2},
\leqno (11.67)
$$
where $q$ is the measurable function of Lemma 11.37 and 
by (11.38).

Then (11.65), applied to each $z_i$, gives two points
$y_{1}(z_i) \in G_{z_i}\cap \partial B(0,r_1)$
and $y_{2}(z_i) \in G_{z_i}\cap \partial B(0,r_2)$.
In addition, Lemma 11.37 says that 
$$
\left| {y_{1}(z_i)\over |y_{1}(z_i)|} - 
{y_{2}(z_i) \over |y_{2}(z_i)|} \right|
\leq q(z_i) \leq C  [f(2r_0)+h_1(2r_0)]^{1/2}.
\leqno (11.68)
$$
By definition of $Z_{j,k}(r_1)$ and $E^\sharp$, the three points 
$y_{1}(z_i)/r_1$ lie in the set $g$ associated to $r_1$ by the construction 
of Section 6. Then (11.20) says that 
$$
\dist(y_{1}(z_i),X(r_1)) \leq C r_0 j(r_1)^{1/2}
\leq C r_0 [f(2r_0)+ h_1(2r_0)]^{1/3}
\leqno (11.69)
$$
because  $r_1 \in {\cal R}_2$ (see (11.60). Similarly, 
$\dist(y_{2}(z_i),X(r_2)) \leq C r_0 [f(2r_0)+ h_1(2r_0)]^{1/3}$.

For $i=1,2$, denote by $R^i_{j,k}$ the great circle of $\partial B$
that contains the arc of $X(r_i) \cap \partial B$ which
passes near $\C_{j,k}'$. With the notations of the definition of
${\cal Z}_0(X,\eta_1)$ at beginning of this section, 
$R^i_{j,k}$ is the great circle of $\partial B$
that contains $\varphi_\ast(\C_{j,k})$, where 
$\varphi: V  \to \partial B$ is the mapping that 
determines $X(r_i)$ as an element of ${\cal Z}_0(X,C\tau)$
(see Lemma 11.19).

By (11.69), the position of $y_{1}(z_1)$ and $y_{1}(z_3)$
determine $R^1_{j,k}$ within an error of at most 
$C [f(2r_0)+ h_1(2r_0)]^{1/3}$. [We do not need 
the middle point $y_{1}(z_2)$, but we use the fact that
$\dist(I_1,I_3) \geq \eta_0$.] 
Similarly, the position of $y_{2}(z_1)$ and $y_{2}(z_3)$
essentially says where $R^2_{j,k}$ is. Now (11.68) says that
the two $R^i_{j,k}$ are $C[f(2r_0)+ h_1(2r_0)]^{1/3}$-close
to each other. Since this holds for every index $(j,k)$, we get that
$X(r_1)$ is $C[f(2r_0)+ h_1(2r_0)]^{1/3}$-close
to $X(r_2)$ in $B(0,1)$, as needed for Lemma 11.61.
\qed

\medskip \proclaim Lemma 11.70.
There is a cone $Y \in {\cal Z}_0(X,C\tau)$ centered at $0$,
such that $\alpha_+(Y) \leq C [f(2r_0)+ h_1(2r_0)]^{1/3}$
and
$$
\dist(y,E) \leq C r_0 [f(2r_0)+ h_1(2r_0)]^{1/3} 
\ \hbox{ for } y\in Y\cap B(0,r_0).
\leqno (11.71)
$$

\ms
We want to take $Y = X(r)$ for some $r\in {\cal R}_2$,
so the first thing to do is check that ${\cal R}_2$ is not empty.
In fact,
$$\eqalign{
H^1((r_0/3,r_0)\setminus {\cal R}_2) 
&\leq C [f(2r_0)+ h_1(2r_0)]^{-2/3}
\int_{[r_0/3,r_0]} j(r) \, dr 
\cr&
\leq C r_0 [f(2r_0)+ h_1(2r_0)]^{1/3},
}\leqno (11.72)
$$
because $\cal R$ has full measure (by (11.10)), by the definition 
(11.60), and by (11.29). 

So we can pick $r_1\in {\cal R}_2$ and set $Y = X(r_1)$.
Then $\alpha_+(Y) \leq C j(r_1)^{1/2} \leq C [f(2r_0)+ h_1(2r_0)]^{1/3}$
by Lemma 11.19 and (11.60), and we just need to prove (11.71).

First consider $y\in Y\cap B(0,r_0)\setminus B(0,r_0/3)$.
By (11.72), ${\cal R}_2$ is $C r_0 [f(2r_0)+ h_1(2r_0)]^{1/3}$-dense,
so we can find $r\in {\cal R}_2$ such that $|r-|y|| \leq C r_0 
[f(2r_0)+ h_1(2r_0)]^{1/3}$. 

By Lemma 11.61, we can find $x\in X(r)\cap \partial B$
such that $\left| x- {y \over |y|} \right| 
\leq C [f(2r_0)+ h_1(2r_0)]^{1/3}$. Then $x' = rx$ lies in
$X(r)\cap \partial B(0,r)$, and
$|x'-y| \leq r \, \left| x- {y \over |y|} \right| 
+ \left| {r y \over |y|} - y \right| \leq 
C r_0 [f(2r_0)+ h_1(2r_0)]^{1/3}$.

Finally (11.20) gives a point $\xi \in g \i r^{-1} E\cap \partial B$ 
such that $|x-\xi| \leq j(r)^{1/2} \leq C [f(r_0)+ h_1(2r_0)]^{1/3}$;
thus $\dist(y,E) \leq |y-r\xi| \leq |y-x'|+r|x-\xi| \leq
C r_0 [f(2r_0)+ h_1(2r_0)]^{1/3}$, as needed.

\ms
To complete the proof of Lemma 11.70, we need to control
$Y \cap B(0,r_0/3)$ too. We'll do this by constructing new
cones $Y_l$, $l \geq 0$, corresponding to scales $2^l$ smaller.

By (11.8), $f(r) + h_1(r)\leq C\varepsilon_1$ for $0 < r \leq 60r_0$,
so the analogue of (11.7) holds with $r_0$ replaced with $2^{-l}r_0$;
thus we can apply the same argument as above, and in particular define 
sets ${\cal R}_2(l) \i (2^{-l} r_0/3,2^{-l} r_0)$, 
$l \geq 0$, and, for $r\in {\cal R}_2(l)$, a cone $X_l(r)$ as in
Lemmas~11.19 and 11.61. We also get a cone $Y_l$ centered at the 
origin, and such that 
$$\eqalign{
\dist(y,E) \leq C 2^{-l} r_0 [f(2^{1-l}r_0)+ &h_1(2^{1-l}r_0)]^{1/3}
\leq C' 2^{-l} r_0 [f(2r_0)+ h_1(2r_0)]^{1/3}
\cr&\ \hbox{ for } y\in Y_l\cap B(0,2^{-l}r_0)\sm B(0,2^{-l}r_0/3),
}\leqno (11.73)
$$
where the second inequality comes from (3.6). We need to show that
$$
d_{0,1}(Y_l,Y_{l+1}) \leq C [f(2^{1-l}r_0)+ h_1(2^{1-l}r_0)]^{1/3}
\leq C [f(2r_0)+ h_1(2r_0)]^{1/3}
\leqno (11.74)
$$
for $l \geq 0$. Notice that ${\cal R}_2(l)$ meets
${\cal R}_2(l+1)$, by the analogues of (11.72) for $l$ and $l+1$. 
Thanks to Lemma 11.61, we may replace $Y_{l+1}$ with any
$X_{l+1}(r)$, ${\cal R}_2(l+1)$ without really altering (11.74), so 
we can assume that $Y_l$ and $Y_{l+1}$ were both defined as sets 
$X_l(r)$ and $X_{l+1}(r)$ for a same 
$r \in {\cal R}_2(l) \cap {\cal R}_2(l+1)$.
We still need to say why we can even take $X_{l}(r)=X_{l+1}(r)$,
because the reader may be worried that, since the set $X_l$
that satisfies the analogue of (11.7) and is used in the construction of 
$X_{l}(r)$ is different from $X_{l+1}$, it could be impossible to use
$X_{l}(r)$ instead of $X_{l+1}(r)$, in the construction of
$X_{l+1}(r)$. This is not the case. The point is that
by (11.7) (and maybe at the cost of replacing $\varepsilon$ with
$\varepsilon/2$, we can safely use $X_{l}(r)$ in the construction
of $\rho$ and $X_{l+1}(r)$ as the cone over $\rho$
(see below (11.22)). We still have that 
$X_l(r) \in {\cal Z}_0(X_{l+1},C\tau)$ by the same argument
as for $X_{l}$ (and maybe with a slightly larger $C$), so the rest 
of the argument goes through essentially unchanged, and we can use 
$X_{l+1}$ or $X_{l}$ indifferently.

So we have (11.74), and it is now easy to deduce (11.71) from (11.73).
Indeed, let $y\in Y \cap B(0,r_0)$ be given, and let $l \geq 0$ be such that
$y\in Y\cap B(0,2^{-l}r_0)\sm B(0,2^{-l}r_0/3)$. If $l=0$, we
can use  (11.73) for $l=0$ and with $Y_0=Y$. Otherwise, (11.74) says 
that we can
find $y_1 \in Y_1\cap B(0,2^{-l}r_0)\sm B(0,2^{-l}r_0/3)$, with
$|y_1-y| \leq C 2^{-l} r_0 [f(2r_0)+ h_1(2r_0)]^{1/3}$. [Recall that
$Y$ and $Y_1$ are both cones, so $d_{0,2^{-l+1}}(Y_l,Y_{l+1})
= d_{0,1}(Y_l,Y_{l+1})$.]
If $l > 1$, we can even iterate this to find
$y_l \in Y_l\cap B(0,2^{-l}r_0)\sm B(0,2^{-l}r_0/3)$, with
$|y_1-y| \leq C l 2^{-l} r_0 [f(2r_0)+ h_1(2r_0)]^{1/3}$. 

Finally, (11.73) gives a point of $E$ at distance at most
$C 2^{-l} r_0 [f(2r_0)+ h_1(2r_0)]^{1/3}$ from $y_l$,
and altogether
$\dist(y,E) \leq C(l+1) 2^{-l} r_0 [f(2r_0)+ h_1(2r_0)]^{1/3}$.
This completes our proof of (11.71), and Lemma 11.70 follows.
\qed

\ms
\proclaim Lemma 11.75. We also have that
$$
\dist(x,Y) \leq C r_0 [f(2r_0)+ h_1(2r_0)]^{1/3}
\ \hbox{ for } x\in E \cap B(0,99r_0/100).
\leqno (11.76)
$$

\ms
Observe that (11.76) and (11.71) imply that 
$d_{0,99r_0/100}(E,Y) \leq C [f(2r_0)+ h_1(2r_0)]^{1/3}$.
Since $\alpha_+(Y) \leq C [f(2r_0)+ h_1(2r_0)]^{1/3}$ 
by Lemma 11.70 and $f(2r_0)+ h_1(2r_0) 
\leq f(110r_0)+ Ch_1(110r_0)]$ by (3.6), 
this is almost the same thing as (11.5);
the only difference is the extra $99/100$, but the reader will 
agree that we could easily have proved Lemma 11.70 with 
$100 r_0/99$, and then the proof of (11.76) below would have 
yielded (11.5). 

So  we just need to prove Lemma 11.75, and 
Theorem 11.4 will follow. Let $x\in E \cap B(0,99r_0/100)$ 
be given. We may assume that $x \neq 0$, because the origin 
lies in $Y$ anyway. Let $l \geq 0$ be such that 
$x\in B(0,2^{-l} \cdot 99r_0/100)\sm B(0,2^{-l-1}\cdot 99r_0/100)$,
We want to show that there is a geometric constant $C_1$ such that, 
for $0 < t < C_1^{-1} 2^{-l} r_0$, there is a minimal cone $Z(t)$ 
of type $\Bbb P$ or $\Bbb Y$, not necessarily centered at $x$, such that
$$
d_{x,t}(E,Z(t)) \leq \alpha,
\leqno (11.77)
$$
where the small constant $\alpha > 0$ will be chosen
later (depending on the constant $\eta_0$ of Section 2.
Recall that because of (11.8) and by the proof of (11.7),
there is a minimal cone $X_l$ (centered at the origin) such that 
$$
d_{0,2^{-l} r_0}(E,X_l) \leq \varepsilon,
\leqno (11.78)
$$
as in (11.6). We want to use the fact that since
$x$ is not too close to the center of $X_l$, there is a small
ball centered at $x$ where $X_l$ coincides with a cone of 
type $\Bbb Y$ or $\Bbb P$, and then we can use the results
of [D3] 
to control $E$ in smaller balls centered at $x$. Later on,
we shall compare our results with (11.71), and get that $x$
is close to $Y$.

We start with the case when the spine of $X_l$ meets
$B(x, 2^{-l-k} r_0)$, where the value of $k > 0$ will be
decided soon. By (11.78), we can find 
$z \in E \cap B(x,2^{-l-k+1} r_0)$, which lies within
$\varepsilon 2^{-l} r_0 \leq 2^{-l-k+1} r_0$ of the
spine of $X_l$. If $k$ is large enough, depending on
$\eta_0$ (in the description of minimal cones in 
Section 2), $X_l$ coincides with a cone of type $\Bbb Y$ 
in $B(z,2^{-l-k+11} r_0)$.
That is, $B(z,2^{-l-k+11} r_0)$ is so small compared to its
distance to the origin that it only meets the three faces
of $X_l$ that we already know about. We translate this cone 
of type $\Bbb Y$ so that its spine contains $z$, and we 
get a cone $Y'$ centered at $z$ such that 
$$
d_{z,2^{-l-k+10} r_0}(E,Y') \leq C\varepsilon.
\leqno (11.79)
$$
If $\varepsilon$ is small enough, we can apply Lemma 16.51 in [D3], 
and we get that $B(z,2^{-l-k+2} r_0)$ is a biH\"older ball
of type $\Bbb Y$ for $E$. To be precise, we are not interested in
the existence of a biH\"older parameterization this time,
but only in the fact that for every ball $B(y,t)$ centered on
$E$ and contained in $B(z,2^{-l-k+2} r_0)$, there is a 
minimal cone $Z(y,t)$, of type $\Bbb Y$ or $\Bbb P$,
such that $d_{y,t}(E,Z(y,t)) \leq \varepsilon'$. In [D3], 
this property is proved with an $\varepsilon'$ which is chosen 
so small that we can then apply Corollary 15.11 there,
and get the biH\"older parameterization. But here, we only need 
to take $\varepsilon' = \alpha$, $y=x$, and
$0 < t < 2^{-l-k+1} r_0$ (recall that 
$z \in E \cap B(x,2^{-l-k+1} r_0)$).
That is, we checked (11.77) in this first case.

Now suppose that the spine of $X_l$ does not meet $B(x, 2^{-l-k} r_0)$. 
If $k$ is large enough (our final condition on $k$), (2.3) says that the 
smaller ball $B(x, 2^{-l-k-1} r_0)$ only meets one face of $X_l$.
That is, $X_l$ coincides with a plane $P$ in 
$B(x, 2^{-l-k-1} r_0)$. We even know that 
$\dist(x,P) \leq \varepsilon 2^{-l} r_0$, because $x\in E$ 
and by (11.78). Translate $P$ so that it passes through $x$; 
we get a plane $P'$ such that
$$
d_{z,2^{-l-k-2} r_0}(E,P') \leq C\varepsilon.
\leqno (11.80)
$$
If $\varepsilon$ is small enough, we can apply Lemma 16.48 in [D3], 
and we get that $B(x,2^{-l-k-6} r_0)$ is a biH\"older ball
of type $\Bbb P$ for $E$. Again, we only need a simple consequence
of the proof, namely, the fact that for $t < 2^{-l-k-6} r_0$, 
there is a plane $P(t)$ such that $d_{x,t}(E,Z(t)) \leq \alpha$.
So (11.77) holds in this second case too, and we can take 
$C_1 = 2^{k+6}$.

\ms
Next we want to compare the sets $Z(t)$ to $Y_l$, where $Y_l$ 
is the analogue of $Y$ for the smaller ball $B(0,2^{-l}r_0)$,
that shows up in (11.73). 

Set $\delta = r_0 [f(2r_0)+ h_1(2r_0)]^{1/3}$ for convenience, and
let us show that 
$$
Y_l \hbox{ meets $B(x,t)$ for } C_2 \delta \leq t \leq 2^{-l} r_0,
\leqno (11.81)
$$
where $C_2$ will be chosen soon. When 
$(10C_1)^{-1} 2^{-l} r_0 \leq t \leq 2^{-l} r_0$,
we first use (11.78) to find $x_1 \in X_l \cap B(x,t/2)$;
then we use the fact that by construction, $Y_l$ is 
$C \tau 2^{-l} r_0$-close to $X_l$ in $B(0,2^{-l}r_0)$
(recall that  $Y_l$ was chosen via Lemma 11.70). 
This gives $y\in Y_l$ such that 
$|y-x_1| \leq C \tau 2^{-l} r_0 < t/2$
(if $\tau$ is small enough), and then $y\in Y_l \cap B(x,t)$,
as needed.

Now suppose that $Y_l$ meets $B(x,t)$ for some 
$t\leq (10C_1)^{-1} 2^{-l} r_0$. We want to show that
$Y_l$ also meets $B(x,t/2)$ if in addition $t\geq C_2 \delta$.
By (11.77), we have a cone $Z(10t)$ of type $\Bbb P$ or $\Bbb Y$, 
which is $10 \alpha rt$-close to $E$ in $B(x,10t)$. Let us check that
$$
\dist(y,Z(10t)) \leq 20 \alpha t
\hbox{ for } y \in Y_l \cap B(x,9t).
\leqno (11.82)
$$
Let $y \in Y_l \cap B(x,9t)$ be given. Observe that 
$x\in B(0,2^{-l} \cdot 99r_0/100)\sm B(0,2^{-l-1}\cdot 99r_0/100)$
by definition of $l$, so $y\in B(x,9t) \i B(0,2^{-l} r_0\sm B(0,2^{-l}r_0/3)$
(because $t \leq (10C_1)^{-1} 2^{-l} r_0$ and
$C_1 > 100$), and (11.73) says that we can find $z \in E$  
such that $|z-y| \leq C 2^{-l}\delta$. So $|z-y| \leq 10\alpha t$
if $t\geq C_2 \delta$ (and $C_2$ is large enough).
Then $z\in E \cap B(x,10t)$, and $\dist(z,Z(10t)) \leq 10\alpha t$
by definition of $Z(10t)$; (11.82) follows.

Now $Y_l$ is not exactly a minimal cone, but it lies in 
${\cal Z}_0(X_l,C\tau)$ (because it was chosen using Lemma 11.70),
and $X_l$ is a minimal cone. 
The fact that in $B(x,9t)$, it is contained in a thin tube around
$Z(10t)$ (by (11.82)), which itself is a minimal cone of type $\Bbb P$ or 
$\Bbb Y$ and contains points in $B(x,t/100)$ (by (11.77)), implies that 
$Y_l$ meets $B(x,t/2)$. Indeed, $Z(10t) \cap  \partial B(0,|x|) \cap B(x,10t)$
is an arc of geodesic in $\partial B(0,|x|)$ or a fork composed of 
three arcs of geodesics, and $Y_l \cap \partial B(0,|x|) \cap B(x,9t)$ 
is also composed of arcs of geodesics, meets $B(x,2t)$, and stays 
close to $Z(10t) \cap  \partial B(0,|x|) \cap B(x,10t)$, so it meets 
$B(x,t/2)$. So we are able to go from (11.81)  for $t$ to 
(11.81) for $t/2$ when $t\geq C_2 \delta$, and this completes the 
proof of (11.81).

If the range of $t$ in (11.81) is empty, then 
$\dist(x,Y) \leq |x| \leq 2^{-l} r_0 \leq C_2 \delta$
by definition of $l$; in this case the inequality in (11.76)
holds trivially. Otherwise, apply (11.82) to $t=C_2\delta$.
This gives a point $y_l\in Y_l$ such that $|y-x| \leq C_2 \delta$,
and we may even assume that $|y| \leq 2|x| \leq 2^{-l+1} r_0$.
Then we apply (11.74) many times, get successive points in the
$Y_m$, $0 \leq m \leq l$, and at the end a point $y_0 \in Y_0 = Y$
such that $|y_0-x| \leq C_2 \delta$. The argument is the same as for the
proof of (11.71) once we have (11.74). 

This completes our proof of (11.76). Lemma 11.75 follows, and so does 
Theorem 11.4 (see the comments below the statement of Lemma 11.75). 
\qed

\ms
\noindent{\bf Remark 11.83.}
In fact we could have proved the slightly more precise fact: 
if  $h$, $x$, and $r$ are as in the statement of Theorem 11.4, and if 
$f(110r_0) + h_1(110r_0)$ is small enough, as in (11.6),
then
$$
\beta_{Z}(x,r_0) \leq C f(2r_0)^{1/3} + C h_1(2r_0)^{1/3}. 
\leqno (11.84) 
$$
That is, we could replace $110$ in (11.5) with $2$.
Of course this looks like a better estimate, but it is unlikely that we shall 
ever  need the improvement, and so we shall stick with the slightly 
uglier $110$. But we shall use the following more specific version 
of Theorem~11.4.

\ms
\proclaim Corollary 11.85.
For each small enough $\tau > 0$, we can find $\varepsilon > 0$
and $C > 0$ such that if $E$ is a reduced almost minimal set in $U$ 
with gauge function $h$, $x \in E$, and $r_0>0$ are such that 
$B(x,110r_0) \i U$ and we can find
a minimal cone $X$ centered at $0$ such that 
$d_{x,100r_0}(E,X) \leq \varepsilon$, then we can find a cone
$Y \in {\cal Z}_0(X,\tau)$ such that
$$
d_{x,r}(E,Y) + \alpha_+(Y) \leq  C [f(110r_0) + h_1(110r_0)]^{1/3}.
\leqno (11.86)
$$

\ms
This is indeed what we proved once we found the minimal cone $X$ in (11.7), 
except that we just need to replace $\tau$ with a slightly smaller constant,
to account for the $C\tau$ in Lemma 11.70, for instance.
The interest compared Theorem  11.4 is that we say that we can keep 
the same minimal cone $X$ as a basis for our computations; this may be
useful in dimensions $n \geq 4$, for instance if we know that $X$
has the full length property. As before, (11.86) is only interesting
when $[f(110 r_0) + h_1(110 r_0)]^{1/3}$ is smaller than $\varepsilon$.
\qed

\bigskip
\noindent {\bf 12. Joint decay and a proof of Theorem 1.15}
\medskip

In this section we indicate how to deduce Theorem 1.15, or variants,
by combining Theorem 4.5, Corollary 11.85, and the simple estimates
from Section 5. 

As before, $E$ is a reduced almost minimal set in $U$ 
with gauge function $h$, we assume that the origin lies
in $E$, and we pick $r_0>0$ such that $B(0,110r_0) \i U$. 
In order to simplify the computations, we assume that 
$$
h(r) \leq C_0 r^{b}
\ \hbox{ for } 0 < r < 220 r_0
\leqno (12.1)
$$
for some choice of constants $C_0 \geq 0$ and $b\in (0,1]$. 
Note that then 
$$
h_1(r) \leq C C_0 r^{b}
\ \hbox{ for } 0 < r < 110 r_0 ;
\leqno (12.2)
$$
here and below, we denote by $C$ various constants that depend on $n$ 
and $b$, but not on $C_0$ or $\varepsilon$. We also assume that 
$$
d_{0,100r_0}(E,X) \leq \varepsilon/500
\leqno (12.3)
$$
for some cone $X$ centered at $0$. 
Recall that the normalized local Hausdorff distance $d_{x,r}$ is defined
in (11.3); here $\varepsilon$ is essentially chosen as in the previous 
sections (and in particular Section 6 and Corollary~11.85),
but we shall give a more precise account of quantifiers in the statement of 
Theorem 12.8. We also assume that
$$
H^2(X\cap B(0,1)) \leq d(0)
\leqno (12.4)
$$
and that
$$
\hbox{ $X$ is a full length minimal cone, with constants 
$\eta_1 \leq \eta_0/10$ and $C_1 \geq 1$.}
\leqno (12.5)
$$
See Definition 2.10 for full length minimal cones, and recall that
$\eta_0$ is the small geometric constant from (2.2) and (2.3).
We further assume that
$$
f(110r_0) + C_0 r_0^b \leq \varepsilon^3, 
\leqno (12.6)
$$
where $f(r) = \theta(r)-d(0)$ is as in (3.5). 

Our statement will also use the constant 
$a = {4 \alpha \over 1-2\alpha} > 0$ introduced in Lemma 5.11, where
$\alpha$ is the small positive in constant in (4.8). Thus $a$ depends
only on $C_1$ and the dimension $n$. In dimension $3$, all the 
minimal cones have the full length property with a fixed $C_1$
(see Section 14),  so we do not need to mention (12.5) 
and $a$ is an absolute constant. We choose an exponent $b_1$ such that
$$
0 < b_1 \leq b \ \hbox{ and } \  b_1<a.
\leqno (12.7)
$$

\ms
\proclaim Theorem 12.8.
For each choice of $b \in (0,1]$, $C_1 \geq 1$ , $\eta_1 \leq \eta_0/10$,
$b_1$ as in (12.7), and $\tau > 0$, there are constants $C \geq 0$ and
$\underline\varepsilon > 0$ such that the following statement holds for
$0 < \varepsilon \leq \underline\varepsilon$.
Let $E$ be a reduced almost minimal set in $U \i \R^n$,
with gauge function $h$, suppose that $0\in E$, and let $r_0 >0$ 
be such that $B(0,110 r_0) \i U$. Suppose that (12.1) holds 
for some $C_0 \geq 0$, and that we have (12.6). Finally
suppose that (12.3) holds for some minimal cone $X$ such that (12.4)
and (12.5) hold. Then there is a minimal cone $Y$ such that
$$
d_{0,r}(E,Y) \leq C \big( r/r_0 \big)^{b_1/3} 
\big\{ f(110r_0) + C_0 r_0^b) \big\}^{b_1/3} 
\leq C \big( r/r_0 \big)^{b_1/3}  \varepsilon
\leqno (12.9)
$$ 
for $0 < r \leq r_0/2$,
$$
Y \in {\cal Z}_0(X,\tau), 
\leqno (12.10)
$$
and
$$
\ d_{0,1}(Y,X) \leq C\varepsilon.
\leqno (12.11)
$$
In addition, $Y$ is the unique blow-up limit of $E$ at $0$, 
and $H^2(Y\cap B) = d(0)$.

The reader will find in [D3] 
a few examples of sufficient conditions that allow us to apply 
Theorem 12.8 (and Corollary 12.25 below).

Before we start with the proof of Theorem 12.8, observe that the 
quantifiers in  its statement are compatible with what we did so far. 
We shall keep the same choices as in the previous sections, except that 
in two occasions (near (12.18) and (12.20))
we shall again require $\varepsilon$ to be small enough.

We want to get $Y$ as the limit of a sequence of cones
$Y_l \in {\cal} Z_0(X,\tau)$, which we shall construct
by induction, so that
$$
d_{0,2^{-l}r_0}(E,Y_l) + \alpha_+(Y_l) 
\leq  C  2^{-b_1l/3} \mu ,
\leqno (12.12)
$$
where we set
$$
\mu = \big\{ f(110r_0) + C_0 r_0^b \big\}^{1/3} \leq \varepsilon
\leqno (12.13)
$$
(by (12.6)). Observe that the first cone $Y_0$ exists, directly by
Corollary 11.85, (12.2), (12.3), and (12.6).

Now let $l \geq 1$ be given, assume that we constructed 
$Y_0, \cdots, Y_{l-1}$, and let us find $Y_{l}$. 
First observe that for $1 \leq j < l$ (if $l > 1$; otherwise
we do not need the information), 
$$\eqalign{
d_{0,1}(Y_j,Y_{j-1}) &= d_{0,2^{-j-1}r_0}(Y_j,Y_{j-1})
\leq 2 d_{0,2^{-j}r_0}(Y_j,E) + 2 d_{0,2^{-j}r_0}(E,Y_{j-1})
\cr& 
\leq 2 d_{0,2^{-j}r_0}(Y_j,E) + 4 d_{0,2^{-j+1}r_0}(E,Y_{j-1})
\leq 6 C 2^{-b_1j/3} \mu
}\leqno (12.14)
$$
because $Y_j$ and $Y_{j-1}$ are cones, by the triangle inequality,
and by (12.12). Similarly,
$$\eqalign{
d_{0,1}(Y_0,X) & = d_{0,r_0/2}(Y_0,X)
\leq 2 d_{0,r_0}(Y_0,E) + 200 d_{0,100r_0}(E,X)
\leq 2C \mu  + 2 \varepsilon/5
}\leqno (12.15)
$$
by (12.12) and (12.3). Let us check that
$$
\hbox{almost-every $r\in (2^{-l}r_0,r_0)$ satisfies
the assumptions of Theorem 4.5.}
\leqno (12.16)
$$
The first assumption (4.1) is satisfied by (12.6), because
$h_1$ is nondecreasing, and by (3.6). Next (4.3) and (4.4)
hold almost-everywhere, by Lemma 4.12. So we just need to
prove (4.2), i.e., that 
$$
d_{0,100r}(E,X) \leq \varepsilon',
\leqno (12.17)
$$
where we denote by $\varepsilon'$ the small constant in (4.2)
(to avoid confusion), because the additional assumptions on $X$
are satisfied, by (12.4) and (12.5).
[This is why we want to keep the same $X$; incidentally the issue would
not arise in dimension 3.]

If $r \geq r_0/500$, (12.17) follows from (12.3). Otherwise,
let $j$ denote the largest integer such that 
$200r \leq 2^{-j} r_0$; then $0 \leq j < l$ 
because $2^{-l} r_0 < r \leq r_0/500$, and
$$
d_{0,100r}(E,X) \leq 2 d_{0,200r}(E,Y_j) + 2 d_{0,1}(Y_j,X)
\leq C' \mu  + 4 \varepsilon/5 \leq C\varepsilon
\leqno (12.18)
$$
again because $Y_j$ and $X$ are cones, by (12.15), iterations of 
(12.14), and (12.13). So (12.17) holds if $\varepsilon$ is small
enough. So (12.16) holds, we can apply Theorem 4.5 to almost every 
$r\in (2^{-l}r_0,r_0)$, and we get (4.8). This allows us to apply
Lemma 5.11, with $y=r_0$ and any $x\in (2^{-l}r_0,r_0)$. We get that
$$\eqalign{
f(x) &\leq (x/r_0)^a f(r_0) + C x^a \int_x^{r_0} r^{-a-1} h(2r) dr
\cr& = (x/r_0)^a f(r_0) + C C_0  x^a \int_x^{r_0} r^{b-a-1} dr
\cr& \leq (x/r_0)^a f(r_0) + C C_0 x^a 
\int_x^{r_0} r_0^{b-b_1} r^{b_1-a-1} dr
\cr& = (x/r_0)^a f(r_0) + C C_0 r_0^{b-b_1} x^{b_1} 
\cr& \leq (x/r_0)^a [f(110r_0) + C h_1(110 r_0)] 
+ C C_0 r_0^b \, (x/r_0)^{b_1}  
\cr& \leq (x/r_0)^{b_1} [f(110r_0) + C C_0 r_0^b] 
+ C C_0 r_0^b \, (x/r_0)^{b_1}  
\leq C (x/r_0)^{b_1} \mu^3 
}\leqno (12.19)
$$
by (5.13), (3.6), because $b_1 < a$, and by (12.2) and (12.13).
If $110 \cdot 2^{-l} < 1$,  we apply this to $x= 110 \cdot 2^{-l}r_0$, 
and we get that
$$
f(110 \cdot 2^{-l}r_0) \leq C 2^{-b_1l} \mu^3 \leq C 2^{-b_1l} \varepsilon^3 
\leqno (12.20)
$$
(by (12.13)). Otherwise, we simply observe that 
$f(110 \cdot 2^{-l}r_0) \leq f(110 r_0) + C h_1(110 r_0) 
\leq f(110 r_0) + C C_0 r_0^b )\leq  C \mu^3$ by (3.6), (12.2), 
and (12.13). This also yields  (12.20) (with a slightly larger constant), 
because  $2^l \leq 110$.

By (12.20) and (12.18) or (12.3), and (for the last time) if 
$\varepsilon$ is small enough, we may apply Corollary 11.85 to the 
radius $r= 2^{-l}r_0$; we get  a cone $Y_l \in {\cal Z}_0(X,\tau)$ 
such that
$$\eqalign{
d_{0,r}(E,Y_l) + \alpha_+(Y_l) &\leq  C [f(110 r) 
+ h_1(110 r)]^{1/3} 
\leq C 2^{-b_1l/3} \mu + C [C_0 r^b]^{1/3}
\cr&\leq C 2^{-b_1l/3} \mu + C (r/r_0)^{b/3} \mu
\leq C 2^{-b_1l/3} \mu,
}\leqno (12.21)
$$
by (12.20), (12.13), and because $b_1 \leq b$. 
This is the same as (12.12), and so we 
can construct the cones $Y_l$ by induction on $l$.

\ms
By (12.14), the $Y_l$ have a limit $Y$ when $l$ tends to $+\infty$,
and  
$$
d_{0,1}(Y,Y_l) \leq \sum_{j \geq l} d_{0,1}(Y_j,Y_{j+1})
\leq C  2^{-b_1l/3} \mu
\leqno (12.22)
$$
(where we use again the fact that the $Y_l$ are cones to simplify the
Hausdorff distance computations). Then
$$
d_{0,2^{-l-1}r_0}(E,Y) 
\leq d_{0,2^{-l}r_0}(E,Y_l) + 2 d_{0,1}(Y,Y_l)
\leq C  2^{-b_1l/3} \mu
\leqno (12.23)
$$
by (12.22) and (12.12). For (12.9) we also need the continuous
version of (12.23), so let $r\in (0 ,r_0/2]$ be given, and apply
(12.23) to the largest $l \geq 0$ such that $r \leq 2^{-l-1}r_0$;
we get that
$$\eqalign{
d_{0,r}(E,Y) &\leq 2d_{0,2^{-l-1}r_0}(E,Y) 
\leq C 2^{-b_1l/3} \mu  \leq  C \big( r/r_0 \big)^{b_1/3} \mu
\cr&= C \big( r/r_0 \big)^{b_1/3} \big\{ f(110r_0) + C_0 r_0^b \big\}^{1/3}
\leq C \big( r/r_0 \big)^{b_1/3} \varepsilon
}\leqno (12.24)
$$
by (12.13). This is just (12.9).

We already know that each $Y_l$ lies in ${\cal Z}_0(X,\tau)$
(see above (12.21)); (12.10) follows by taking limits. 
Next, (12.11) follows from (12.15), (12.14), and (12.13).
By (12.24), $r^{-1} E$ tends to $Y$ when $r$ tends to $0$, so
$Y$ is the only blow-up limit of $E$ at $0$. Then
Proposition~3.15 says that $Y$ is a minimal cone, 
with $H^2(Y\cap B) = d(0)$. This completes our proof of
Theorem 12.8.

\ms \proclaim Corollary 12.25.
For each choice of $n \geq 3$, $b \in (0,1]$, $C_1 \geq 1$, 
$\eta_1 > 0$ we can find $\beta > 0$ and $\varepsilon_1 > 0$ 
such that the following holds.
Let $U \i \R^n$ be open and let $E \i U$ be a reduced almost minimal 
set in $U \i \R^n$, with gauge function $h$. Suppose that
$0 \in E$, $r_0 > 0$ is such that $B(0,110r_0) \i U$, (12.1) holds 
for some choice of $b \in (0,1]$ and $C_0 \geq 0$,
$$
f(110 r_0) + C_0 r_0^b \leq \varepsilon_1,
\leqno (12.26)
$$
and 
$$
d_{0,100 r_0}(E,X) \leq \varepsilon_1
\leqno (12.27)
$$
for some full length minimal cone $X$ centered at $0$ 
such that (12.4) and (12.5) hold. Then for $0 < r \leq r_0$  
there is a $C^{1,\beta}$ diffeomorphism 
$\Phi : B(0,2r) \to \Phi(B(0,2r))$, such that 
$\Phi(0)=0$, $|\Phi(y)-y| \leq 10^{-2}r$ for $y\in B(0,2r)$,
and $E \cap B(0,r) = \Phi(X) \cap B(0,r)$.

\ms
As we shall see later, Theorem 1.15 follows fairly easily
from Corollary 12.25. In turn, Corollary 12.25 will be
easily deduced from the extension in [DDT] 
of Reifenberg's topological disk, once we have the following 
uniform control on the good approximation of $E$ by minimal 
cones in balls $B(x,r)$ near the origin, which will be proved first.

\ms
\proclaim Proposition 12.28.
Let $\varepsilon_0 > 0$ be given, and then let 
$E$ and $r_0$ be as in Corollary 12.25, with the same quantifiers,
except that now $\varepsilon_1$ may depend on $\varepsilon_0$ too. 
Then for $x\in E$ and $r > 0$ such that
$x\in E \cap B(0,10r_0)$ and $0 < r < 10 r_0$, we can find a minimal 
cone $Z(x,r)$, not necessarily centered at $x$ or at the origin, such that
$$
d_{x,r}(E,Z(x,r)) \leq (r/r_0)^{\beta} \varepsilon_0.
\leqno (12.29) 
$$

\ms
We choose the constant $b_1$ as in (12.7).
The minimal cones $Z(x,r)$ will be easier to find when
$x=0$, so we start with this case. 

\ms
\proclaim Lemma 12.30.
There is a minimal cone $Y$ such that 
$$
d_{0,r}(E,Y) \leq (r/r_0)^{b_1/3} \varepsilon_0
\hbox{ for $0 \leq r \leq 100 r_0$.}
\leqno (12.31) 
$$

\ms
We want to apply Theorem 12.8.
Let $\underline \varepsilon$ be as in the statement, and
choose $\varepsilon \leq \underline \varepsilon$ such that
$C \varepsilon$ in (12.11) is smaller that $10^{-3}\varepsilon_0$.
Next we check the assumptions; (12.1) and (12.4) hold by assumption,
(12.3) follows from (12.27) if $\varepsilon_1 < \varepsilon/500$
(and $X$ satisfies (12.4) and (12.5) by assumption), and finally
(12.6) follows from (12.26) if $\varepsilon_1 < \varepsilon^3$.

So Theorem 12.8 applies, and we get a minimal cone $Y$ as in 
its statement.
For $0 \leq r \leq r_0/2$, we apply (12.9) and get that
$$
d_{0,r}(E,Y) \leq C \big( r/r_0 \big)^{b_1/3} \varepsilon
\leqno (12.32) 
$$
as needed for (12.31). When $r_0/2 \leq r \leq 100 r_0$, 
$$\eqalign{
d_{0,r}(E,Y) 
&\leq d_{0,r}(E,X) + d_{0,1}(X,Y)
\leq (100 r_0/r) d_{0,100r_0}(E,X) + C\varepsilon
\cr& \leq (100 r_0/r) \varepsilon_1 + 10^{-3}\varepsilon_0 
\leq (r/r_0)^{b_1/3} \varepsilon_0
}\leqno (12.33) 
$$
because $X$ and $Y$ are cones, by (12.11), (12.27), 
our choice of $\varepsilon$, and if $\varepsilon_1$ is small enough. 
Lemma 12.30 follows.
Note that in addition we get (12.11) and (12.10), where $\tau >0$ 
is any small number given in advance.
\qed

\ms
We shall also apply Lemma 12.30 with smaller values of 
$\varepsilon_0$; this will be helpful to get additional properties,
and the only cost will be that we need to take $\varepsilon_1$
even smaller.

Next we want to take care of the small balls centered on
the set $E_Y$ of $Y$-points of $E \cap B(0,50 r_0)$.
Recall that $E_Y$ is the set of points $y\in E$ such that every 
blow-up limit of $E$ at $y$ is a cone of type $\Bbb Y$. 
By Proposition 3.14 and Lemma 14.12 in [D3], 
$$
E_Y = \big\{ y\in E \, ; \, d(y)=3 \pi/2 \big\}.
\leqno (12.34) 
$$

\ms
\proclaim Lemma 12.35.
For $y\in E \cap B(0,50 r_0) \sm \{ 0 \}$ there is a minimal cone 
$Y(y)$ of type $\Bbb Y$, centered at $y$,  such that 
$$
d_{y,r}(E,Y(y)) \leq (r/|y|)^{b_1/3} \varepsilon_0
\ \hbox{ for $0 \leq r \leq 10^{-2} \eta_0 |y|$,}
\leqno (12.36) 
$$
where $\eta_0$ is the constant from Section 2.

\ms
Let $y \in E_Y \cap B(0,50r_0)$ be given, set 
$\delta = 10^{-4} \eta_0 |y|$, and let us try to apply
Lemma~12.30 to the ball $B(y,\delta)$ (or to be exact,
to $E-y$ and $B(0,\delta)$). 

We shall prove the analogue of (12.27) with a minimal cone $Z$ of type
$\Bbb Y$; observe that $Z$ automatically satisfies (12.4), 
because $d(y) = 3 \pi/2$ by (12.34), and (12.5) with some constant 
that does not even depend on $n$, by Section 14. 

Denote by $\varepsilon'_1$ the value of  $\varepsilon_1$ in (12.26) 
and (12.27) that is needed for Lemma 12.30 to hold, with 
$\varepsilon_0$ replaced with the slightly smaller 
$[10^{-4} \eta_0]^{b_1/3} \varepsilon_0$, and let
us check the assumptions.

Obviously (12.1) holds, and $B(y,110\delta) \i B(0,100r_0) \i U$.
Then we need to check the analogue of (12.26) and (12.27).
For (12.27), we need to find a minimal cone $Z'$ of type $\Bbb Y$, 
centered at $y$, such that
$$
d_{y,100\delta}(E,Z') \leq \varepsilon'_1.
\leqno (12.37) 
$$
For (12.26), observe that since $d(y) = 3 \pi/2$ by (12.34),
we need to prove that
$$
(110\delta)^{-2} H^2(E \cap B(y,110\delta)) - {3 \pi \over 2} \, 
+ C_0 \delta^b \leq \varepsilon'_1.
\leqno (12.38) 
$$

In order to obtain both things, we apply Lemma 12.30 to $B(0,r_0)$, but 
with a smaller constant $\varepsilon'_0$ to be chosen soon.
The assumptions are still satisfied (if $\varepsilon_1$ is small
enough), so we get a minimal cone $Y$ such that (12.31) holds.
In particular, taking $r = 2|y|$, we get that
$$
d_{0,2|y|}(E,Y) \leq (2|y|/r_0)^{b_1/3} \varepsilon'_0 \leq \varepsilon'_0.
\leqno (12.39) 
$$
Let us check that there is a minimal cone $Z$, of type $\Bbb Y$ or 
$\Bbb P$, such that
$$
Y \cap B(y,200\delta) = Z \cap B(y,200\delta).
\leqno (12.40) 
$$
Recall from Section 2 that $K' = Y \cap \partial B(0,|y|)$
is composed of great circles or arcs of great circles $\C'_{j}$.

Set $D = B(y,\eta_0 |y|/10) \cap \partial B(0,|y|)$, and notice that 
$D$ meets $K'$ because $Y$ is a cone that contains point very near $y$ 
(by (12.39) and because $y\in E$). Let us first assume that $D$ meets
only one circle or arc of circle $\C'_{j}$. Denote by $Z$ the plane that 
contains $\C'_{j}$, and let us check (12.40). Observe that $D$ contains 
no endpoint of $\C'_{j}$, because each endpoint meets two other arcs of $K'$; 
thus $\C'_{j}$ crosses $D$, and $Z \cap D = \C'_{j} \cap D = K' \cap D$.
Now $Z$ is a cone, and $Y$ is the cone over $K'$, so $Z$ and $Y$ 
coincide in $B(y,\eta_0 |y|/20)$; (12.40) follows because 
$200 \delta = 2 \cdot 10^{-2} \eta_0 |y| < \eta_0 |y|/20$.

Now assume that $D$ meets at least two different circles or arcs $\C'_{i}$ 
and $\C'_{j}$. Pick $\xi \in \C'_{i}$; By (2.3), $\C'_{i}$ and $\C'_{j}$
have a common endpoint $x$ in $B(\xi,\eta_0 |y|/5) \i B(y,2\eta_0 |y|/5)$. 
Still by the description in Section 2, there is 
a third arc $\C'_{k}$ that ends at $x$, and the three arcs leave from $x$ 
with $120^\circ$ angles. 
Denote by $Z$ the minimal cone of type $\Bbb Y$ that contains
the beginning of these three arcs; we want to show that (12.40) 
holds for this $Z$. 

First observe that (2.2) says that the lengths of $\C'_{i}$, $\C'_{j}$,
and $\C'_{k}$ are at least $10 \eta_0 |y|$, so these arcs leave $D$
and $K' \cap D$ contains $Z \cap D$. We want to show that 
$K' \cap D \i Z$ as well.

Suppose $D$ meets some other $\C'_{l}$. As before, $\C'_{l}$
has a common extremity with $\C'_{i}$ in $B(y,2\eta_0 |y|/5)$,
which we call $y_i$. Note that $y_i \neq x$ because the three ends at
$x$ are already taken. Similarly, $\C'_{l}$ has a common endpoint
$y_j$ with $\C'_{j}$ and a common endpoint $y_k$ with 
$\C'_{k}$, both in $B(y,2\eta_0 |y|/5)$ and different from $x$. 
Since $\C'_{l}$ has only two ends, two of these points
are equal, for instance $y_j$ and $y_k$. But this is impossible,
as the only possible intersections of $\C'_{j}$ and $\C'_{k}$
are $x$ and $-x$ (they leave from $x$ with $120^\circ$ angles),
and we said that $y_j \in B(y,2\eta_0 |y|/5)$. So $D$ meets no other 
$\C'_{l}$. 

Another weird possibility would be that $\C'_{i}$, say,
is very long and returns to $D$ after leaving it. Then it would 
have another end in $D$ (because (2.3) says that it is not a full 
circle), and at this end it would meet two other arcs. We just proved
that these other arcs cannot be different from $\C'_{j}$ and 
$\C'_{k}$, but also that $\C'_{i}$ cannot share two ends in 
$B(y,2\eta_0 |y|/5)$ with $\C'_{j}$ or $\C'_{k}$. So $D \cap K'$
is reduced to the three arcs that leave from $x$, and hence 
$K' \cap D = Z \cap D$. As before, this implies that
$Y \cap B(y,\eta_0 |y|/20) = Z \cap B(y,\eta_0 |y|/20)$ because
$Y$ and $Z$ are cones, and then (12.40) holds in this second case as 
well.

Because of (12.40), (12.39) implies that
$$
d_{y,190 \delta}(E,Z) = d_{y,190 \delta}(E,Y) 
\leq {2|y| \over 190 \delta} \,  d_{0,2|y|}(E,Y)
\leq C \varepsilon'_0
\leqno (12.41) 
$$
because $\delta = 10^{-4} \eta_0 |y|$ (and so 
$B(0,190 \delta) \i B(0,2|y|)$). It does not matter that $C$ 
depends on $\eta_0$.

Let us apply Lemma~16.43 in [D3] to $B(0,190\delta)$, 
with $F = Z$ and $10^{-2}\varepsilon'_1$ playing the role of $\delta$ there.
Recall that $\varepsilon'_1$ was already chosen, a little before (12.37).
The assumptions are satisfied if $\varepsilon_1$ is small enough
(so that (12.26) controls $h(300 \delta)$) and 
$\varepsilon'_0$ is small enough (so that (12.41) controls 
the distance to $Z$). we get that
$$\eqalign{
H^2(E \cap B(y,&110\delta))
\leq H^2(Z \cap B(y,(1+ 10^{-2}\varepsilon'_1)110\delta)) 
+ 10^{-2} \varepsilon'_1 (190\delta)^2 
\cr&
\leq {3 \pi \over 2} [(1+ 10^{-2}\varepsilon'_1)110\delta)]^2
+ 10^{-2} \varepsilon'_1 (190\delta)^2 
\leq (110\delta)^2 \Big\{ {3 \pi \over 2} + {\varepsilon'_1 \over 2}
\Big\}
}\leqno (12.42)
$$
because $Z$ is of type $\Bbb Y$ or $\Bbb P$,  so its density in 
any ball is at most $3\pi/2$. We divide by $(110\delta)^2$,
add $C_0 \delta^b  \leq C_0 r_0^b \leq \varepsilon_1 
\leq \varepsilon'_1/2$ by (12.26) and
if $\varepsilon_1$ is small enough, and get (12.38).

Unfortunately, we cannot use $Z$ directly for (12.37),
because we do not know that it is centered at $y$. 
It would not matter much if $Z$ were a plane, but if we want to apply
Lemma 12.30 on $B(y,\delta)$, it is important to use a cone centered at 
$y$. So let us check that
$$
\hbox{$Z$ is of type $\Bbb Y$, and its spine meets
$B(y,10\varepsilon'_1 \delta)$.}
\leqno (12.43)
$$
Suppose not. Then $Z$ coincides with a plane in 
$B(y,2\varepsilon'_1 \delta)$. Let us apply 
Lemma~16.43 in [D3] again, 
this time in $B(y,3\varepsilon'_1 \delta)$, again with 
$F=Z$, and with $\delta = 10^{-2}$. We can still use (12.26)
to control $h(3\varepsilon'_1 \delta)$, and we observe that
$$
d_{0,4 \varepsilon'_1\delta}(E,Z) 
\leq {190 \over 4 \varepsilon'_1} \,  d_{0,190 \delta}(E,Z) 
\leq C \varepsilon'_0 / \varepsilon'_{1}
\leqno (12.44) 
$$
by (12.41). Again we take $\varepsilon'_0$  small enough,
depending on $\varepsilon'_{1}$, apply the lemma, and get that 
$$\eqalign{
H^2(E \cap B(y,\varepsilon'_1 \delta))
&\leq H^2(Z \cap B(y,(1+10^{-2})\varepsilon'_1 \delta)
+ 10^{-2}(3\varepsilon'_1 \delta)^2
\cr&\leq \pi \, [(1+10^{-2})\varepsilon'_1 \delta]^2 
+ 10^{-2}(3\varepsilon'_1 \delta)^2
\leq {4 \pi \over 3} \, (\varepsilon'_1 \delta)^2.
}\leqno (12.45) 
$$
On the other hand, $y \in E_Y$, so
$$
3\pi/2 = d(y) \leq 
(\varepsilon'_1 \delta)^{-2} H^2(E \cap B(y,\varepsilon'_1 \delta))
+ C h_1(\varepsilon'_1 \delta)
\leq {4 \pi \over 3} +  C h_1(\varepsilon'_1 \delta)
\leqno (12.46) 
$$
by (12.34), (3.8), and (12.45). This contradiction proves (12.43).

By (12.43), we can translate $Z$ by less than $10\varepsilon'_1 \delta$
to  get a set $Z'$ of type $\Bbb Y$ centered at $y$. Then 
$$
d_{y,100 \delta}(E,Z') \leq {190 \over 100} \, d_{0,190 \delta}(E,Z)
+ {10\varepsilon'_1 \delta \over 100 \delta}
\leq C \varepsilon'_0 + {\varepsilon'_1 \over 10}
< \varepsilon'_1
\leqno (12.47) 
$$
by (12.41), so we can use $Z'$ in (12.37).

We finally checked (12.37) and (12.38), which allows us to apply
Lemma 12.30 to $E-y$ in $B(0,\delta)$, and with respect to the
cone $Z'$. 
Recall that we even arranged to get a constant
$[10^{-4} \eta_0]^{b_1/3} \varepsilon_0$
instead of $\varepsilon_0$, so we get a minimal cone
$Y(y)$ centered at $y$, and such that
$$
d_{y,r}(E,Y(y)) \leq (r/\delta)^{b_1/3} \,  [10^{-4} \eta_0]^{b_1/3}
\varepsilon_0
\leqno (12.48) 
$$
for $0 \leq r \leq 100 \delta$. This is exactly the same thing as
(12.36), because $\delta = 10^{-4} \eta_0 |y|$.

To  complete the proof of Lemma 12.35, we still need to check
that $Y(y)$ is of type $\Bbb Y$. But we observed in the proof
of Lemma 12.30 (after (12.33)) that we can take $Y(y)$ in
${\cal Z}_0(Z',\tau)$ for some small $\tau$, as in (12.10), and then 
it is clear that $Y(y)$ is of type $\Bbb Y$, just like  $Z'$.
Lemma 12.35 follows.
\qed

\ms
We shall also need to control $E$ in small balls centered
on $E \sm [E_Y \cup \{ 0 \}]$, so we state a third lemma.
Set
$$
\rho(x) = \dist(x,E_Y \cup \{0 \}) \leq |x| \leq 10r_0
\leqno (12.49)
$$
for $x\in E \cap B(x,10r_0)$.

\ms
\proclaim Lemma 12.50.
For $x\in E \cap B(x,10r_0) \sm [E_Y \cup \{ 0 \}]$,
there is a plane $P(x)$ through $x$ such that 
$$
d_{x,r}(E,P(x)) \leq (r/\rho(x))^{b_1/3} \varepsilon_0
\ \hbox{ for $0 \leq r \leq 10^{-5} \eta_0 \rho(x)$.}
\leqno (12.51) 
$$

\ms
We shall proceed as in Lemma 12.35, and apply Lemma 12.35
or Lemma 12.30 to a small ball centered at $x$, where the choice
of the lemma will depend on whether $x$ is significantly closer to 
$E_Y$ than to the origin. The manipulation of quantifiers will be 
similar.

We start with the case when 
$$
\rho(x) \geq 10^{-4}\eta_0 |x|. 
\leqno (12.52) 
$$
Then we set $\delta = 10^{-7}\eta_0 |x|$, and we want to 
apply Lemma 12.30 to $E-y$ in $B(0,\delta)$, and with $\varepsilon_0$ 
replaced with $[10^{-7} \eta_0]^{b_1/3}\varepsilon_0$ . 
We denote by $\varepsilon'_1$ the constant in (12.26) and (12.27)
that will be enough to do this. This time we want to find a plane
$P$ through $x$ such that
$$
d_{x,100\delta}(E,P) \leq \varepsilon'_1
\leqno (12.53) 
$$
(as before, in (12.37)), and also prove that
$$
(110\delta)^{-2} H^2(E \cap B(y,110\delta)) - d(x) 
+ C_0 \delta^b \leq \varepsilon'_1.
\leqno (12.54) 
$$
(as for (12.38)). Note that (12.1) holds and $B(x,110\delta) \i U$
as before, and $P'$ automatically satisfies (12.4) (because $\pi$ is 
the smallest possible density), and (12.5) (see Section 14). 
So the other assumptions for Lemma 12.30 are satisfied.

Let us again apply Lemma 12.30 to $B(0,r_0)$, but with $\varepsilon_0$
replaced with a smaller $\varepsilon'_0$ to be chosen soon. We take 
$r = 2|x|$, and (12.31) says that
$$
d_{0,2|x|}(E,Y) \leq (2|x|/r_0)^{b_1/3} \varepsilon'_0 \leq \varepsilon'_0.
\leqno (12.55) 
$$
 
The same geometric argument as for (12.40) gives a 
minimal cone $Z$ of type $\Bbb Y$ or $\Bbb P$, such that 
$$
Y \cap B(x,200\delta) = Z \cap B(x,200\delta),
\leqno (12.56) 
$$
and so 
$$
d_{x,190\delta}(E,Z) = d_{y,190\delta}(E,Y) 
\leq {2|x| \over 190\delta} \,  d_{0,2|x|}(E,Y)
\leq C \varepsilon'_0
\leqno (12.57) 
$$
as for (12.41). Let us check that 
$$
\hbox{$Z$ is a plane, or
its spine does not meet $B(x,170\delta)$.}
\leqno (12.58) 
$$
Otherwise, $Z$ is of type $\Bbb Y$, and its spine 
meets $B(x,170\delta)$ at some point $\xi$. We want
to apply Proposition~16.24 in [D3] 
to the ball $B(\xi,10\delta)$ and get a point of $E_Y$
near $\xi$. The usual assumptions that $B(\xi,20\delta) \i U$ 
and $h_1(20\delta)$ is small
enough follow because $\delta = 10^{-7}\eta_0 |x| \leq r_0$ 
and $x\in B(0,10r_0)$, and by (12.26), and the main assumption
that $d_{x,10\delta}(E,Z)$ for some cone $Z$ of type $\Bbb Y$
centered at $\xi$ is small enough comes directly from (12.57)
if $\varepsilon'_2$ is small enough.
So  Proposition~16.24 in [D3] applies, 
and gives a point of type $\Bbb Y$ in $B(\xi,\delta)$.
But then $\dist(x,E_Y) \leq |x-\xi| + \delta \leq 171 \delta
= 171 \cdot 10^{-7}\eta_0 |x| < 10^{-4}\eta_0 |x|$. 
This contradiction with (12.52) or (12.49) proves (12.58).

By (12.57) and (12.58), there is a plane $P$ through $x$ such that
$$
d_{x,160\delta}(E,P) \leq C \varepsilon'_0
\leqno (12.59) 
$$
(take a plane $P'$ which coincides with $Z$ in $B(x,170\delta)$,
and then translate it by $\dist(x,P')  \leq 190 \delta d_{x,190\delta}(E,Z)
\leq  C \varepsilon'_0 \delta$ to make it go through $x$).
In particular, (12.53) holds if $\varepsilon'_0$ is small enough.

Next we apply Lemma~16.43 in [D3] to $B(0,140\delta)$, 
with $F = P$ and $\delta = 10^{-2}\varepsilon'_1$, where 
$\varepsilon'_1$ was chosen near (12.52).
As before, the assumptions are satisfied if $\varepsilon_1$ is 
small enough (so that (12.26) controls $h(200 \delta)$) and 
$\varepsilon'_0$ is small enough (so that (12.59) controls 
the distance to $P$). We get that
$$\eqalign{
H^2(E \cap B(y,&110\delta))
\leq H^2(P \cap B(y,(1+ 10^{-2}\varepsilon'_1)110\delta)) 
+ 10^{-2} \varepsilon'_1 (140\delta)^2 
\cr&
\leq \pi \, [(1+ 10^{-2}\varepsilon'_1)110\delta)]^2
+ 10^{-2} \varepsilon'_1 (140\delta)^2 
\leq (110\delta)^2 \Big\{ \pi + {\varepsilon'_1 \over 2}
\Big\}
}\leqno (12.60)
$$
as for (12.42), and because $P$ is a plane.
Thus
$$
(110\delta)^{-2} H^2(E \cap B(y,110\delta)) + C_0 \delta^b 
\leq \pi + {\varepsilon'_1 \over 2} + \varepsilon_1
\leq \pi + \varepsilon'_1 \leq d(x) + \varepsilon'_1
\leqno (12.61)
$$
by (12.26), if $\varepsilon_1$ is small enough, and 
because $d(x) \geq \pi$ on $E$. Then (12.54) holds.

So we completed the verification of (12.53) and (12.54),
Lemma 12.30 applies to $E-x$ in $B(0,\delta)$, and with $X=P$.
We get a minimal cone $P(x)$ centered at $x$ such that
$$\eqalign{
d_{x,r}(E,P(x)) &\leq (r/\delta)^{b_1/3} \,  [10^{-7} \eta_0]^{b_1/3}\varepsilon_0
\cr&
= (r/\rho(x))^{b_1/3} \,  (\rho(x)/\delta)^{b_1/3} \, [10^{-7} \eta_0]^{b_1/3}
\varepsilon_0
\cr&
\leq (r/\rho(x))^{b_1/3}\,  (|x|/\delta)^{b_1/3} \, [10^{-7} \eta_0]^{b_1/3}
\varepsilon_0
\leq (r/\rho(x))^{b_1/3} \varepsilon_0
}\leqno (12.62) 
$$
for $0 \leq r \leq 100 \delta$, because we applied Lemma 12.30 with 
$[10^{-7} \eta_0]^{b_1/3}\varepsilon_0$ instead of $\varepsilon_0$,
because $\rho(x) \leq |x|$ by (12.49), and since 
$\delta = 10^{-7} \eta_0 |x|$.
In addition, $P(x)$ is a plane, for instance because it lies in
${\cal Z}_0(P,\tau)$ for some small $\tau$. So  we proved
(12.51), and the range is large enough, because 
$100\delta =  10^{-5} \eta_0 |x| \geq  10^{-5} \eta_0 \rho(x)$.

\ms
We are left with the case when (12.52) fails, and  so
$$
\rho(x) < 10^{-4}\eta_0 |x|. 
\leqno (12.63) 
$$
Obviously $\rho(x) = \dist(x,E_Y)$, because the origin is much 
further.  Choose $y\in E_Y$ such that 
$$
|x-y| \leq 2 \rho(x) \leq 2 \cdot 10^{-4}\eta_0 |x|;
\leqno (12.64) 
$$
obviously and $0 < |y| \leq 2|x| \leq 20r_0$, so we can apply 
Lemma 12.35 to $y$. We do this with $\varepsilon_0$ replaced
with a smaller $\varepsilon'_0$, to be chosen soon; we get a minimal
cone $Y(y)$ of type $\Bbb Y$, centered at $y$, such that (12.36)
holds. We can take  $r = 3\rho(x)$ in (12.36), because 
$\rho(x) \leq 10^{-4}\eta_0 |x| < 2 \cdot 10^{-4}\eta_0 |y|$
by (12.63) and (12.64), and we get that
$$
d_{y,3\rho(x)}(E,Y(y)) \leq \varepsilon'_0 \, .
\leqno (12.65) 
$$
We claim that 
$$
\hbox{the spine  of $Y(y)$ does not meet $B(x,\rho(x)/2)$.}
\leqno (12.66) 
$$
Indeed, otherwise we can find $\xi \in B(x,\rho(x)/2)$ in the spine,
and we can apply Proposition~16.24 in [D3] 
to the ball $B(\xi, \rho)$, as we did near (12.58), to find a
point $z$ of type $\Bbb Y$ in $B(\xi, \rho(x)/10)$. 
Then $\rho(x) \leq |x-z| \leq |x-\xi| + \rho(x)/10) \leq 6\rho(x)/10$,
a contradiction which proves (12.66). 

By (12.65) and (12.66), there is a plane $P$ through $x$ such that 
$$
d_{y,\rho(x)/3}(E,P) \leq 18\varepsilon'_0 \, .
\leqno (12.67) 
$$
(take a plane $P'$ that coincides with $Y(y)$ in $B(x,\rho(x)/2)$,
and translate it by $\dist(P',x) \leq 3 \rho(x) d_{y,3\rho(x)}(E,Y(y)) 
\leq 3\varepsilon'_0 \rho(x)$  to make it go through $x$).

This is a good analogue of (12.59). Set $\delta = 10^{-3}\rho(x)$.
If $\varepsilon'_0$ is small enough, we can apply Lemma 12.30 to $E-x$, 
in the ball $B(0,\delta)$, with the cone $X=P$, and where
we replace $\varepsilon_2$ with $10^{-b_1}\varepsilon_2$.
The analogue of (12.27) comes directly from (12.67), and
the analogue of (12.28) is deduced from (12.67) and 
Lemma~16.43 in [D3], as we did for (12.54) in (12.60)-(12.61). 

We get a plane $P(x)$ through $x$ such that 
$$\eqalign{
d_{x,r}(E,P(x)) &\leq (r/\delta)^{b_1/3} \,  10^{-b_1}\varepsilon_0
= (r/\rho(x))^{b_1/3} \,  (\rho(x)/\delta)^{b_1/3} \, 10^{-b_1}
\varepsilon_0
\cr&
= (r/\rho(x))^{b_1/3} \varepsilon_0
}\leqno (12.68) 
$$
for $0 \leq r \leq 100 \delta = 10^{-1}\rho(x)$, as for (12.62). 
This  is better than (12.51); Lemma~12.50 follows.
\qed

\ms\noindent{\bf Proof of Proposition 12.28.}

The proof that follows is not very efficient, but uses
a minimal amount of information. See Remark 12.81, and
13.1 concerning possible improvements.   
Let us first say how we can compute an exponent $\beta > 0$
that works. Set $t =  b_1/3$, and define a function $F$ on
$V = \big\{ (u,v,w) \in [0,1]^3 \, ; \, u+v+w =1 \big\}$ by
$$
F(u,v,w) = {\rm Max} \big\{ tu, tv-u,tw-u-v \big\}. 
\leqno (12.69) 
$$
we take
$$
\beta = \inf \big\{ F(u,v,w) \, ; \, (u,v,w) \in V \big\}.
\leqno (12.70) 
$$
It is easy to check that $\beta > 0$ (for instance, by compactness),
but the result of estimates in terms of the already small $t>0$ is so 
disappointingly small that we shall not bother.

Let $x\in E \cap  B(0,10r_0)$ and $0 < r_0 \leq 10 r_0$ be given.
We may try to find the desired cone $Z(x,r)$ in for different ways,
and we shall decide about which one depending on the relative values 
of $r$, $r+\rho(x)$, and $r+\rho(x)+|x|$. Set 
$\lambda = r /(30r_0) < 1$, and define  $(u,v,w) \in V$ by
$$
{r \over r+\rho(x)} = \lambda^u \, , \, 
{r+\rho(x) \over r+\rho(x)+|x|} = \lambda^v  \hbox{ , and } \,
{r+\rho(x) +|x| \over 30 r_0} = \lambda^w. 
\leqno (12.71) 
$$

There is a trivial attempt, which will allow us to take care of
the case when $r$ is not too small: we can try $Z(x,r) = X$,
where $X$ is as in (12.27), This yields 
$$\eqalign{
d_{x,r}(E,Z(x,r)) &= d_{x,r}(E,X) 
\leq (100r_0/r) \, d_{x,100r_0}(E,X) 
\cr& \leq 100 \varepsilon_1 r_0/r
= {100 \varepsilon_1 \over \varepsilon_0 } \, 
(r_0/r)^{1+\beta} \, \varepsilon_0 (r/r_0)^{\beta},
}\leqno (12.72) 
$$
which is enough for (12.29) if 
$$
(r/r_0)^{1+\beta}  \geq 100 \varepsilon_1/\varepsilon_0.
\leqno (12.73) 
$$

\ms
We start with a first case when $F(u,v,w) = tu$ (and hence
$\beta \leq ut$, by (12.70)). We want to apply Lemma 12.50,
with $\varepsilon'_0 = \varepsilon_0/60$.
We can only do this if $r \leq 10^{-5} \eta_0 \rho(x)$
(and $\varepsilon_1$ is small enough), but 
if this is the case, we get a plane $P(x)$ such that
$$\eqalign{
d_{x,r}(E,P(x)) &\leq (r/\rho)^{b_1/3} \varepsilon'_0
\leq 2(r/(r + \rho(x)))^{b_1/3} \varepsilon'_0
= 2 \lambda^{tu} \varepsilon'_0 
\leq 2 \lambda^{\beta} \varepsilon'_0
\cr& = 2 (r/30r_0)^{\beta} \varepsilon'_0
\leq 60 (r/r_0)^{\beta} \varepsilon'_0
= (r/r_0)^{\beta} \varepsilon_0
}\leqno (12.74) 
$$
because $r \leq \rho(x)$ and $b_1/3 = t$, and by (12.71).
This is enough for (12.29).

\ms
If instead $F(u,v,w) = tu$ but 
$$
r > 10^{-5} \eta_0 \rho(x),
\leqno (12.75) 
$$
we want to use (12.72), so we shall check that $r$ is not 
too small because $r/\rho$
represents a good part of $\lambda$.

Since $F(u,v,w) = tu$, we get that $tu \geq tv-u$ and so
$v \leq t^{-1} (1+t) u$. Similarly, $tu \geq tw-u-v$, so
$w \leq t^{-1} (1+t) u + t^{-1} v$, and then
$1 = u+v+w \leq A(t) u$, where we don't need to compute
$A(t)$ explicitly. Hence 
$$\eqalign{
r/r_0 &= 30\lambda \geq 30 \, \lambda^{A(t) u} 
= 30\, \Big( {r \over r+\rho(x)} \Big)^{A(t)}
\cr&
= 30\, \Big( {1 \over 1+\rho(x)/r} \Big)^{A(t)}
\geq 30\, \Big( {1 \over 1+10^{5} \eta_0^{-1}} \Big)^{A(t)}
}\leqno (12.76) 
$$
by (12.75). Thus (12.73) holds (if $\varepsilon_1$ is small enough),
and (12.72) allows us to conclude.

\ms
Our next case is when $F(u,v,w) = tv-u$. This time, we 
want to apply Lemma 12.35 to $B(y,2\rho(x)+r)$, 
with the constant $\varepsilon'_0 = \varepsilon_0/5$, 
and where $y\in E_Y \cup \{ 0 \}$ is such that 
$|y-x| < 2 \rho(x)$. For this we shall need to assume that
$$
2\rho(x)+r \leq 10^{-2} \eta_0 |y|
\leqno (12.77) 
$$
Observe that if (12.77) holds, $\rho(x)$ is much smaller than 
$|y|$ and $|x|$, so $y \in E_Y \cap B(0,50 r_0)$ and we can
indeed apply Lemma 12.35 and get a minimal cone $Y(y)$ such that
(12.36) holds, and so
$$\leqalignno{
d_{x,r}(E,Y(y)) &
\leq {2\rho(x)+r\over r} \, d_{y,2\rho(x)+r}(E,Y(y))
\leq {2\rho(x)+r\over r} \, 
\Big( { 2\rho(x)+r \over |y|}\Big)^{b_1/3} \varepsilon'_0 
\cr&
\leq 5 \, {\rho(x)+r\over r} \, 
\Big( { \rho(x)+r \over |x|+\rho(x)+r}\Big)^{b_1/3} \varepsilon'_0 
= 5 \lambda^{-u} \lambda^{b_1 v/3} \varepsilon'_0 
= 5 \lambda^{tv-u} \varepsilon'_0 
&(12.78) 
\cr& 
= 5 \lambda^{F(u,v,w)} \varepsilon'_0
\leq 5 \lambda^{\beta} \varepsilon'_0
\leq 5 (r/r_0)^{\beta} \varepsilon'_0
= (r/r_0)^{\beta} \varepsilon_0
}
$$
because $\big| |y| - (|x|+\rho(x)+r) \big| \leq |y-x| + \rho(x)+r$ 
is much smaller than $|y|$ (by (12.77)), and by (12.71) and various 
definitions. So (12.29) holds in this case.

\ms
Next assume that $F(u,v,w) = tv-u$, but (12.77) fails.
Observe that $tu \leq tv-u$ and $tw-u-v \leq tv-u$ because 
$F(u,v,w) = tv-u$, so $u \leq t(1+t)^{-1}v$ and 
$tw \leq (1+t)v$, and finally $1 = u+v+w \leq B(t) v$
for some positive $B(t)$. Thus
$$
r/r_0 = 30\lambda \geq 30 \lambda^{B(t) v} 
= 30\, \Big( {r+\rho(x) \over r+\rho(x)+|x|}\Big)^{B(t)}.
\leqno (12.79) 
$$
In addition, $|x| \leq |y| + |x-y| \leq |y| + 2 \rho(x)
\leq C (x+\rho(x))$ because (12.77) fails, so 
$r/r_0 \geq C^{-1}$, (12.73) holds, and (12.72) yields (12.29).

\ms
Our last case is when $F(u,v,w) = tw-u-v$. Then we apply Lemma 12.30, 
with $\varepsilon'_0 = \varepsilon_0/30$
and to $B(0,r+|x|)$ and get that
$$\eqalign{
d_{x,r}(E,Y) &\leq {r+|x| \over r} \, d_{x,r+|x|}(E,Y) 
\leq {r+|x| \over r} \Big( {r+|x| \over r_0}\Big)^{b_1/3} \varepsilon'_0
\cr&
\leq 30 \, {r + \rho(x) + |x| \over r}\, 
\Big( {r + \rho(x) +|x| \over 30 r_0}\Big)^{b_1/3} \varepsilon'_0
= 30 \lambda^{-u-v}\lambda^{tw} \varepsilon'_0
\cr&
= 30 \lambda^{F(u,v,w)} \varepsilon'_0
\leq 30 \lambda^{\beta} \varepsilon'_0
\leq 30 (r/r_0)^{\beta} \varepsilon'_0
= (r/r_0)^{\beta} \varepsilon_0 \, .
}\leqno (12.80) 
$$
(12.71), (12.70), and the usual computation. This proves (12.29) in 
our last case, and establishes Proposition 12.28.
\qed

\ms\noindent{\bf Remark 12.81.} It is clear that the power $\beta$
that we get here is far from optimal. A typical place where we lose 
a lot of information is in the proof of Lemma 12.35, when we apply 
Lemma 12.30, obtain (12.39), and drop the possibly very small 
$(2|y|/r_0)^{b_1/3}$ from our estimate. The effect is that we have 
to start decay estimates for smaller balls from scratch. We also 
do such cultural revolutions in (12.55) and (12.65).

The author decided not to try to write a more efficient proof in detail;
see 13.1 in the next section for suggestions about how to do this.
Note that even with this additional work, our final estimates will
not look too good, because of $\alpha$ in (4.8).

\ms\noindent{\bf Proof of Corollary 12.25.}
Corollary 12.25 follows from Proposition 12.28 and 
the generalization in [DDT] 
of Reifenberg's topological disk. 
More precisely, Section 10 in [DDT] 
gives sufficient conditions weaker than the conclusions
of Proposition 12.28 for the existence of $\Phi$ as in the
conclusion of Corollary 12.25, but only of class $C^1$.
See in particular the discussion that starts after (10.22)
there, up  to the end of the paper.
The proof applies and gives $C^{1,\beta}$ estimates
(apply the Whitney extension theorem to H\"older-continuous functions).

The reader should not pay too much attention to the existence
of the mapping $\Phi$; we stated it this way to avoid long 
descriptions, but what seems important to the author
is the decomposition of $E \cap B(0,2r)$ into $C^{1+\beta}$
faces $F_{j}$, which make $120^\circ$ angles with each other.
Here the decomposition into faces comes from the
H\"older description in [D3], 
and both the fact that they are $C^{1+\beta}$
and that they make the right angles along the edges comes from 
Proposition 12.28. The point is that we have a control
on the variations of the tangent plane, because (12.29) also
gives the existence of a tangent cone $Z$ to $E$  at
$x$, with $d_{0,1}(Z,Z(x,r)) \leq C (r/r_0)^\beta \varepsilon_0$
(compare the $Z(x,2^{-k}r)$ and sum a geometric series).
\qed

\ms\noindent{\bf Proof of Theorem 1.15.}
Let $E$ be as in the theorem, and let $x\in E$ be given.
Without loss of generality, we may assume that $x=0$.
Let $X$ be a blow-up limit of $E$ at $x$, with the full length
property.  [See Definition 3.11 for blow-up limits.]
When $n=3$, we simply need to pick any blow-up limit
(and such a limit exists, see (3.13)), because every minimal cone 
has the full length property (see Section 14).  
When $n > 3$, we don't know that much, so we put the existence 
of $X$ in the assumptions.

Then there are arbitrarily small radii $r_0$ that satisfy the 
assumptions of Theorem~12.8 (notice that $f(r)$ tends to $0$ 
by its definition (3.5), and that (12.5) holds because 
$H^2(X \cap B(0,1)) = d(x)$ by Proposition 3.14). 
Theorem 12.8 says that there is a unique minimal cone $Y$
(which is then equal to $X$), and (12.9) even says how fast
$d_{0,r}(E,X) = d_{0,r}(E,Y)$  tends to $0$. Then we can apply
Corollary~12.25  to $r_0$ for every $r_0$ small enough, and
Theorem~1.15 follows.
\qed

\bigskip
\noindent {\bf 13. Small generalizations and improvements}
\medskip

In this section we are concerned about various ways to take 
the same proofs as above and get slightly better results.
The comments in 13.32 and 13.41 are independent from the previous ones,
but 13.22 depends on 13.1.

\ms\noindent {\bf 13.1 Better estimates for Proposition 12.28.}

We observed in Remark 12.81 that in the proof of Proposition 12.28,
we may lose a lot of information when we apply Lemma~12.30 or 12.35
to estimate the analogue of $f(r)$ for a small ball $B$ centered away
from the origin, and where we have a good approximation of $E$ by a 
minimal cone near $B$. For instance, in (12.39) we had a point
$y\in E_Y$ and we obtained that
$$
d_{0,2|y|}(E,Y) \leq (2|y|/r_0)^{b_1/3} \varepsilon'_0 
\leq \varepsilon'_0,
\leqno (13.2) 
$$
and we decided to use the last inequality alone.
If instead we keep the whole estimate, the proof of 
(12.41) gives a minimal cone $Z$ of type $\Bbb Y$ such that
$$
d_{y,190\delta}(E,Z) \leq {2 |y| \over 190\delta} \, d_{0,2|y|}(E,Y)
\leq C (2|y|/r_0)^{b_1/3} \varepsilon'_0 \, , 
\leqno (13.3) 
$$
because $|y|$ and $\delta = 10^{-4} \eta_0 |y|$ are comparable.
Then we decided to use Lemma 16.43 in [D3] 
to get a mild control on
$$
f(110 \delta)  = (110 \delta)^{-2} H^2(E\cap B(y,110 \delta))
-{3\pi \over 2}
\leqno (13.4) 
$$
(see (12.42)). When $|y|/r_0$ is very small, we should do something 
more clever instead. 

First observe that the spine of $Z$ should pass within 
$C\delta d_{y,190\delta}(E,Z)$ from $y$, by the same argument as 
for (12.43); this allows us to replace $Z$ with a cone $Y$ of type
$\Bbb Y$, parallel to $Z$, whose spine contains $y$, and such that
$$
d_{y,180\delta}(E,Y) \leq C (|y|/r_0)^{b_1/3} \varepsilon'_0 \, .
\leqno (13.5) 
$$
We claim that in such circumstances,
$$
f(110 \delta) \leq C (|y|/r_0)^{b_1/3} \varepsilon'_0 \, .
\leqno (13.6)
$$
Let us even formulate a slightly more general statement to this effect.

Let $E$ be a reduced almost minimal set in $U$, with gauge function 
$h$, suppose that $0 \in E$ and $B(0,2r_0) \i U$, and that 
$$
f(3r_0/2) \hbox{ and $h_1(2r_0)$ are small enough,}
\leqno (13.7)
$$
to set the stage. Let $Y$ is a minimal cone centered at the 
origin, and such that 
$$
H^2(Y\cap B(0,1)) = d(0).  
\leqno (13.8)
$$
Then
$$
f(r_0) \leq C d_{0,2r_0}(E,Y) + 9 h(3r_0).
\leqno (13.9)
$$

When $n=3$, this is rather easy. We want to use $Y$ to
construct a competitor, but we shall need to add a small rim
around $Y \cap \partial B(y,r_0)$. Set $B = B(y,r_0)$
and
$$
Y_\eta = \big\{ z\in B \, ; \, \dist(z,Y) \leq \eta \big\},
\leqno (13.10)
$$
where we take $\eta = 2r_0 d_{y,2r_0}(E,Y)$.
Thus $E \cap B \i Y_\eta$ by definition of 
$d_{0,2r_0}(E,Y)$. 

Set $Y' = [Y \cap B] \cup [Y_\eta \cap \partial B]$. 
By ``elementary geometry", there is a Lipschitz mapping $h$ 
such that $h(z)=z$ out of $B$, $h(z)\in \overline B$ for $z\in B$, 
and which maps $Y_\eta \cap B$ to a subset of $Y'$. Then set $F = h(E)$; 
it is easy to see that $F$ is a competitor for $E$ in $B$, an so
$$\eqalign{
H^2(E\cap B) &\leq H^2(F \cap \overline B) + 4r_0^2 h(2r_0)
\leq H^2(Y') + 4r_0^2 h(2r_0)
\cr&\leq H^2(Y\cap B(0,r_0)) + H^2(Y_\eta \cap \partial B) 
+ 4r_0^2 h(2r_0)
\cr&
\leq r_0^2 H^2(Y\cap B(0,1)) + C r_0^2 d_{0,2r_0}(E,Y) 
+ 4r_0^2 h(2r_0),
}\leqno (13.11)
$$
as needed for (13.9).

When $n>3$, we need to construct another competitor, because 
$H^2(Y_\eta \cap \partial B) = +\infty$. 
What seems the easiest at this point is to use the construction 
of Sections 6-9. Let us check that we can choose $r \in (r_0,3r_0/2)$ 
such that (4.3) and (4.4) hold, and also
$$
H^1(E \cap \partial B(0,r))  \leq 2(1+\tau) \, d(0) \, r
\leqno (13.12)
$$
(the analogue of (6.2) here). We don't nee to worry about 
(4.3) and (4.4), because Lemma~4.12 says that they hold almost 
everywhere, so let us assume that (13.12) fails 
for almost every $r \in (r_0,3r_0/2)$. Set 
$\theta(r) = r^{-2} H^2(E\cap B(0,r))$ as usual.
By (5.7), (5.6), and (5.8),
$$\leqalignno{
f(3r_0/2) - f(r_0) &\geq \int_{r_0}^{3r_0/2} \theta'(r) dr
\geq \int_{r_0}^{3r_0/2}  
[r^{-2} H^1(E\cap \partial B(0,r)) -2r^{-1}\theta(r)] dr
\cr&
\geq  \int_{r_0}^{3r_0/2} 
r^{-1} [2(1+\tau) \, d(0) -2\theta(r)] dr
& (13.13)
}$$
because (13.12) fails a.e. Recall from (3.5) that
$\theta(r) = d(0) + f(r)$; thus
$$
\eqalign{
f(3r_0/2) - f(r_0) &
\geq  \int_{r_0}^{3r_0/2} 2r^{-1} [\tau \, d(0) -2f(r)] dr
}\leqno (13.14)
$$
This is impossible if $f(3r_0/2)$ and $h_1(2r_0)$ are small enough, because 
(3.6) says that then $f(r_0)$ and $f(r)$ are very small too. 

So we can find $r$ such that (13.12) holds too, and this allows
us to follow the construct of Sections 6-9, with $X=Y$. 
Most of our estimates are probably useless here, but anyway the 
first line of (9.69), suitably normalized because here $r \neq 1$, says that
$$
H^2(E \cap B)
\leq (3r_0)^2 h(3r_0) + H^2(\Sigma \cap B(0,r)) + 2550
\int_{E \cap \partial B(0,r)} \dist(z,\Gamma) \, dH^1(z)
\leqno (13.15)
$$
for some nice union $\Sigma$ of Lipschitz graphs, and some
union $\Gamma$ of small Lipschitz graphs in $\partial B(0,r)$.
Notice that $\dist(z,\Gamma) \leq 4r_0 d_{0,2r_0}(E,Y)$ on
$E\cap \partial B(0,r)$, by definition of $d_{0,2r_0}(E,Y)$ and
construction of $\Gamma$; The construction of $\Sigma$, where
we tend to minimize areas of surfaces over triangular regions,
already gives $H^2(\Sigma \cap B(0,r)) \leq H^2(Y\cap B(0,r))
+ C r_0^2 d_{0,2r_0}(E,Y)$, but if this were not the case, 
we could improve on it, or even work directly on the cone over $\Gamma$, 
to get a deformation $\Sigma'$ (something that looks like
$Y'$ above, but with tilted walls), that satisfies this estimate.
The point is that it is much easier to deform a small Lipschitz
graph, rather than the set $E$, when we do not know the topology of 
$E$. So (13.15) yields 
$$
H^2(E \cap B) \leq H^2(Y\cap B(0,r)) + C r_0^2 d_{0,2r_0}(E,Y)
+ (3r_0)^2 h(3r_0),
\leqno (13.16)
$$
as needed for (13.9).

\ms
Return to Proposition 12.28. Let us also modify our estimate
near (12.55) and (12.65), as suggested above, and otherwise 
follow the proof of Proposition 12.28. The worse estimate is
when $B(x,r)$ is fairly far from $E_Y$, and even further
from the origin. In this case we choose $y\in E_Y$ such that
$|y-x| \leq 2 \rho(x)$, and first apply Theorem 12.8 
(as in Lemma~12.30) to get $Y$ with the property (12.9). 
Thus 
$$
d_{0,2|y|}(E,Y) \leq C (|y|/r_0)^{b_1/3}\varepsilon'_0.
\leqno (13.17)
$$
We apply (13.9) and get that 
$$
f(C^{-1} |y|) \leq C (|y|/r_0)^{b_1/3}\varepsilon'_0,
\leqno (13.18)
$$
where $f$ is computed with respect to $y$. By Theorem 12.8,
we get a cone $Y(y)$ of type $\Bbb Y$, with 
$$
d_{y,\rho}(E,Y(y)) \leq C (\rho/|y|)^{b_1/3} 
\big[ (|y|/r_0)^{b_1/3} \varepsilon'_0\big]^{b_1/3}.
\leqno (13.19)
$$
We take $\rho$ a little larger than $|x-y|$, and use $Y(y)$
to get a plane $P$ through $x$ such that
$$
d_{x,c\rho}(E,P) \leq C d_{y,\rho}(E,Y(y))
\leq C (\rho/|y|)^{b_1/3} 
\big[ (|y|/r_0)^{b_1/3} \varepsilon'_0\big]^{b_1/3}
\leqno (13.20)
$$
for some small $c>0$. By (13.9), we get a similar estimate
for $f(c\rho/2)$, where $f$ is now computed with respect to 
the center $x$. We apply Theorem 12.8 one last time, to get 
a plane $P(x)$ through $x$ such that
$$\eqalign{
d_{x,r}(E,P(x)) &\leq C (r/\rho)^{b_1/3} \, (\rho/|y|)^{{b_1^2}/9}
\,(|y|/r_0)^{{b_1^3}/27} (\varepsilon'_0)^{{b_1^3}/27}
\leq \varepsilon''_0 (r/r_0)^{{b_1^3}/27},
}\leqno (13.21)
$$
where $C''_0 = C (\varepsilon'_0)^{{b_1^3}/27}$ is still as small
as we want.

Thus, with a little more work (not entirely written down here) 
we could take $\beta = b_1^3/27$ in Proposition 12.28, and then
in Corollary 12.25.

\ms\noindent {\bf 13.22 Larger gauge functions $h$.}

When the gauge function $h$ is larger than a power of $r$,
but not too small, some part of Theorem 1.15 and Corollary 12.25 
stays true. Obviously, when $h$ becomes larger, we expect less 
regularity from the corresponding almost minimal sets $E$.
Trivial counterexamples consist in taking your favorite minimal 
set $E$ (for instance a line or a plane), distorting it so that
some spiraling occurs, and computing a gauge function
for the distorted set.

Here we only want to say that if $h$ is not too large, we can keep 
the same proof as above, and get $C^1$ estimates (instead of 
$C^{1+\beta}$), with some control on the modulus of continuity.
Let us only take an example, and assume that 
$$
h(r) \leq C [\log(A/r)]^{-b}
\leqno (13.23)
$$
for some constants $A, b > 0$ and $r$ small enough.
Our proof appears to give that $E$ is locally $C^1$-equivalent
to a minimal cone (under the same other assumptions as in
Theorem~1.15 and Corollary 12.25 ) when $b > 30$; however,
we shall not check every detail here, and the estimate is unlikely 
to be optimal.

Notice that the results from Sections 4-11 (and in particular 
Theorem 4.5) only require the Dini condition (1.16), so we just need 
to worry about Sections 12 and 13.

We start with the decay rate for $f$ that follows from Theorem 4.5,
which was computed in Example 5.21. 
If $f$ satisfies the differential inequality that comes from 
Theorem 4.5 and $h$ satisfies (13.23), we get that $f$ decays at 
least like $C [\log(A/r)]^{-b}$ again. So far, we require $b>1$,
and only because we always assume the Dini condition (1.16).
Thus, instead of (12.19), we get that 
$$
f(x) \leq (x/r_0)^{a} f(r_0) + C (x/A)^{a/2} 
+ C \big[\log\big({A \over 2x}\big)\big]^{-b}  
\leq C_{A,r_0} \, \big[\log\big({A \over 2x}\big)\big]^{-b} 
\leqno (13.24)
$$
for $x\in (2^{-l}r_0,r_0)$. Then we apply Corollary 11.85
to get that
$$\leqalignno{
d_{0,2^{-l}r_0}(E,Y_l) + \alpha_+(Y_l)
&\leq C 2^{-al/3} f(r_0)^{1/3} + C (2^{-l}r_0/A)^{a/6} 
+ C \big[\log\big({A \over 2^{-l+1}r_0}\big)\big]^{-b/3}  
\cr&
\leq C_{A,r_0} \, \big[\log\big({A \over 2^{-l+1}r_0}\big)\big]^{-b/3}
& (13.25)
}
$$
as in (12.21). We get a similar estimate for the relative distances 
between the successive cones $Y_l$, and for the existence of a 
limit $Y$ and (12.22), we simply demand that $b>3$ to get a converging
series, and then sum over $l$ such that $2^{-l}r_0 \leq 2r$ to get that
$$\eqalign{
d_{0,r}(E,Y) + \alpha_+(Y)
&\leq C (r/r_0)^{a/3} f(r_0)^{1/3} + C (r/A)^{a/6} 
+ C \big[\log\big({A \over 2r}\big)\big]^{1-{b\over3}}  
\cr&
\leq C_{A,r_0} \, \big[\log\big({A \over 2r}\big)\big]^{1-{b\over3}}
}\leqno (13.26)
$$
for $0 < r < r_0/2$, which is supposed to replace (12.9).

For the analogue of Theorem 1.15 and Corollary 12.25, we continue the
argument as suggested in 13.1. The analogue of the estimates 
(13.17)-(13.21) (the worst case scenario) seems (modulo computation 
mistakes) to be
$$
d_{0,2|y|}(E,Y) \leq C 
(|y|/r_0)^{a/3} (\varepsilon'_0)^{1/3} + C (|y|/A)^{a/6} 
+ C \big[\log\big({A \over 4|y|}\big)\big]^{1-{b\over3}},
\leqno (13.27)
$$
$$
f(C^{-1} |y|) \leq C (|y|/r_0)^{a/3} (\varepsilon'_0)^{1/3} 
+ C (|y|/A)^{a/6} 
+ C \big[\log\big({A \over 4|y|}\big)\big]^{1-{b\over3}},
\leqno (13.28)
$$
$$\leqalignno{
d_{y,\rho}(E,Y(y)) &\leq C (\rho/|y|)^{a/3} f(C^{-1} |y|)^{1/3}
+ C (\rho/A)^{a/6} 
+ C \big[\log\big({A \over 2\rho}\big)\big]^{1-{b\over3}}
\cr&
\leq C (\rho/|y|)^{a/3}(|y|/r_0)^{a/9}(\varepsilon'_0)^{1/9}
+  C (\rho/|y|)^{a/3} (|y|/A)^{a/18} 
\cr& \hskip0.7cm
+ C (\rho/|y|)^{a/3} 
\big[\log\big({A \over 4|y|}\big)\big]^{{1\over 3}-{b\over 9}}
+ C (|y|/A)^{a/6} 
+ C \big[\log\big({A \over 2\rho}\big)\big]^{1-{b\over3}}
&(13.29)
\cr&
\leq C (\rho/r_0)^{a/9} (\varepsilon'_0)^{1/9}
+ C (\rho/A)^{a/18} + C (|y|/A)^{a/6}
\cr& \hskip0.7cm
+ C (\rho/|y|)^{a/3} 
\big[\log\big({A \over 4|y|}\big)\big]^{{1\over 3}-{b\over 9}}
+C \big[\log\big({A \over 2\rho}\big)\big]^{1-{b\over3}},
}
$$
$$\eqalign{
d_{x,c\rho}(E,P) &\leq C d_{y,\rho}(E,Y(y))
\leq C (\rho/r_0)^{a/9} (\varepsilon'_0)^{1/9}
+ C (\rho/A)^{a/18} + C (|y|/A)^{a/6}
\cr& \hskip1cm
+ C (\rho/|y|)^{a/3} 
\big[\log\big({A \over 4|y|}\big)\big]^{{1\over 3}-{b\over 9}}
+C \big[\log\big({A \over 2\rho}\big)\big]^{1-{b\over3}},
}
\leqno (13.30)
$$
a similar estimate for $f(c\rho/2)$, and finally
$$\eqalign{
d_{x,r}(E,P(x)) &\leq 
C (r/\rho)^{a/3} f(c\rho/2)^{1/3}
+ C (r/A)^{a/6} 
+ C \big[\log\big({A \over 2r}\big)\big]^{1-{b\over3}}
\cr&
\leq C (r/r_0)^{a/27} (\varepsilon'_0)^{1/27}
+  C (r/A)^{a/54} + C (r/\rho)^{a/3} (|y|/A)^{a/18}
\cr&\hskip0.5cm
+ C (r/A)^{a/6} 
+ C (r/\rho)^{a/3} (\rho/|y|)^{a/9}
\big[\log\big({A \over 4|y|}\big)\big]^{{1\over 9}-{b\over 27}}
\cr&\hskip1cm
+ C (r/\rho)^{a/3} 
\big[\log\big({A \over 2\rho}\big)\big]^{1-{b\over3}}
+  C \big[\log\big({A \over 2r}\big)\big]^{{1\over 3}-{b\over 9}}
\cr&\leq C_{A,r_0} 
\big[\log\big({A \over 2r}\big)\big]^{{1\over 9}-{b\over 27}}.
}\leqno (13.31)
$$

Recall that we want to apply Section 10 in [DDT] 
to get $C^1$ estimates on $E$ and, at least if we do this
brutally, we want upper bounds $\varepsilon_k$ for 
$d_{x,r}(E,Z(x,r))$ when $r \sim 2^{-k}r_0$ such that 
$\sum_k \varepsilon_k < +\infty$. So we want the power in
(13.31) to be less than $-1$, i.e., that 
${b\over 27}-{1\over 9} > 1$, or equivalently $b>30$.
Again, all this is subject to verification and improvement.

\ms\noindent {\bf 13.32 Weaker full length conditions.}

The full length condition in Definition 2.10 can be seen
as a strange non-degeneracy condition concerning the 
the length function $H^1(\varphi_\ast(K))$ on deformations
of the cone $X$ through mappings $\varphi \in \Phi(\eta_1)$.
If we have less precise estimates, some part of Theorem 4.5
will still hold, but unfortunately not enough to prove the
local $C^1$-equivalence of $E$ to a minimal cone.

Let us say that the (reduced) minimal cone $X$ satisfies
a weak full length condition of order $N > 1$ when
there is a standard decomposition of $K = X\cap \partial B(0,1)$ 
as in Section 2, an $\eta_1 < \eta_0/10$,
and a constant $C_1 \geq 1$, such that if 
$\varphi \in \Phi(\eta_1)$ is such that
$$
H^1(\varphi_\ast(K)) > H^1(K),
\leqno (13.33)
$$
then there is a deformation $\widetilde X$ of $\varphi_\ast(X)$ 
in $B = B(0,1)$ such that
$$
H^2(\widetilde X \cap B) \leq H^2(\varphi_\ast(X)\cap B)
- C_1^{-1} [H^1(\varphi_\ast(K)) - H^1(K)]^N.
\leqno (13.34)
$$
[The notations are the same as in Definition 2.10, and we only
added the exponent $N$.]

In Theorem 4.5, when we replace the full length condition 
above (4.8) with a weak full length condition of order $N$,
we get the following weaker form of (4.8):
$$\eqalign{
H^2(E \cap B(x,r))
\leq {r \over 2} &\, H^1(E\cap \partial B(x,r))
\cr&- \alpha r^2 \, [r^{-1} H^1(E\cap \partial B(x,r)) - 2 d(x)]^N_+
+ 4 r^2 h(2r),
}\leqno (13.35)
$$
where the $A_+$ denotes the nonnegative part of $A$.

When $H^1(E\cap \partial B(x,r)) \leq 2 d(x) r$, we claim
no gain, and (13.35) follows from (4.7) because (4.6) says that
the part inside the brackets is nonnegative. In the other case,
notice that (13.35) is weaker than (4.8) (which corresponds 
to $N=1$). 

We proceed as before; there is no need to change anything before
the first time we used the full length condition, which happens
near (10.16). Here we need to replace (10.17) with
$$
H^2(\widetilde X \cap B) \leq H^2(\varphi_\ast(X)\cap B)
- C_1^{-1} [H^1(\varphi_\ast(K)) - H^1(K)]_+^N.
\leqno (13.36)
$$
Then we apply Lemma 10.2, as before, and get (10.5) with
$$
A = C_1^{-1} [H^1(\varphi_\ast(K)) - H^1(K)]_+^N
= C_1^{-1} [H^1(\rho) - H^1(K)]_+^N
\leqno (13.37)
$$
That is,
$$\eqalign{
H^2(E \cap \overline B) &\leq {1 \over 2} H^1(E\cap \partial B)
- 10^{-5} [H^1(E\cap \partial B) - H^1(\rho)] 
\cr& \hskip 2.2cm
- C_1^{-1}\kappa^2 [H^1(\rho) - H^1(K)]_+^N + 4 h(2)
\cr& \leq {1 \over 2} H^1(E\cap \partial B)
- C_1^{-1}\kappa^2  [H^1(E\cap \partial B) - H^1(\rho)]_+^N 
\cr& \hskip 2.2cm
- C_1^{-1}\kappa^2 [H^1(\rho) - H^1(K)]_+^N + 4 h(2)
\cr& \leq {1 \over 2} H^1(E\cap \partial B)
- \alpha [H^1(E\cap \partial B) - H^1(K)]_+^N + 4 h(2)
\cr& \leq {1 \over 2} H^1(E\cap \partial B)
- \alpha [H^1(E\cap \partial B) - 2d(0)]_+^N + 4 h(2)
}\leqno (13.38)
$$
because $H^1(E\cap \partial B) - H^1(\rho)$
is small and nonnegative (recall from (9.68) that
$H^1(E\cap \partial B) - H^1(\rho) 
= \delta_1 + \delta_2 + \delta_3 \geq 0$), 
with $\alpha = C_1^{-1}\kappa^2$, and then by 
(10.16). This gives (13.35) after re-scaling (recall that for 
(13.38) we assumed that $x=0$ and $r=1$).

\ms
The new estimate (13.35) yields some decay for the function
$f$, but not as much as before. Notice that (13.35) is the same
as (5.27), except that we used $\alpha' = 2^N \alpha$ there.  
Lemma 5.28 and Remark 5.32 then say that if we use a gauge function
$h$ such that
$$
h(r) \leq C [{\rm Log}({1 \over r})]^{-{N \over N-1}}
\ \hbox{ for $r$ small}
\leqno (13.39) 
$$
(we are allowed larger functions $h$, but the final estimate is then worse),
and as long as we can use (12.35), we get that
$$
f(x) \leq C_1 \Big[{\rm Log}\Big({C_2 \over x}\Big)\Big]^{-{1 \over N-1}}
\ \hbox{ for $x$ small,}  
\leqno (13.40) 
$$
as in (5.33) and (5.34). However, two related things happen.
First, (13.40) is not enough to get $C^1$ estimates, or the 
uniqueness of the tangent cone to $E$ at $0$, as in Section 12.
But also, we cannot sum up our estimates for smaller balls as we did 
for Theorem 12.8, and as a consequence we cannot use the same cone
$X$ to approximate $E$ in $B(0,r)$ for $r$ small. Thus, if we want to 
obtain (13.40), we have to find other cones with the weaker full 
length property, so that we can apply the analogue of Theorem 4.5 and
get (13.35) at smaller scales. The simplest way to do this is to assume 
that every blow up limit of $E$ at $0$ is a minimal cone with the weaker full 
length property, with uniform constants $N$, $C_1$, and $\eta_1$.
Then we get (13.40) for $r$ small, and some estimates on the 
$\beta_Z(0,r)$ that we can derive from Theorem 11.4.

\ms\noindent {\bf 13.41. Minimal-looking cones.}

Recall from Remark 2.14 that a minimal-looking cone is a 
(reduced) cone $E$ such that $K=E\cap \partial B(0,1)$ is as in 
the description of minimal cones given at the beginning of Section~2
(up to (2.3)). This comes with a small constant $\eta_0 > 0$
attached (the constant in (2.2) and (2.3)).

Our main results do not use the fact that we deal with true minimal
cones. That is, the minimality is only used through the description
of Section 2. Thus, for instance, Theorem 4.5 is still true if
$X$ in the assumptions is a (full length for (4.8)) minimal-looking
cone. Then we can use the differential inequality as before
(if we can find an $X$ for almost every radius).
Of course, all the standard ways which give us a cone $X$  
as in the statement will give a minimal cone.

Similarly, the results of Section 12 go through. For instance,
if in Theorem 12.8, we only know that the cone $X$ in (12.3)
is a full length minimal-looking cone, we still get the same 
conclusion. At the beginning of the argument, we only know that
the cone $Y$ in the conclusion is minimal-looking, but then it
turns out that (12.9) holds, so $Y$ is a blow-up limit
of $E$ at $x$, which forces $Y$ to be minimal. [Again, this also
makes it unlikely that we will find $X$ in (12.3) without knowing 
that it is minimal.] The same remarks apply to Corollary 12.25,
Proposition 12.28, and Lemmas 12.30, 12.35, and 12.50.

\bigskip
\noindent {\bf 14. Examples of full length minimal cones}
\medskip

The main purpose of this section is to verify that the
standard minimal cones of dimension $2$ have the full length
property of Definition 2.10. We shall systematically use
Lemma 10.23 to do this, so let us review the notation.

We consider a minimal cone $X$, take a standard decomposition 
of $K = X \cap \partial B$ (where $B$ is the unit ball),
and consider the deformation of $\varphi_\ast(K)$ and $\varphi_\ast(X)$ 
of $K$ and $X$ with a mapping $\varphi \in \Phi(\eta_1)$
(see Section 2). We are allowed to choose $\eta_1$,
and we shall do it so that the main term in some expansion
is larger than the errors. 

Lemma 10.23 says that there is a competitor $\widetilde X$
for $\varphi_\ast(X)$ in $B$ such that 
$$
H^2(\widetilde X \cap B) \leq H^2(\varphi_\ast(X)\cap B) 
- C^{-1} \alpha_+(\varphi)^2,
\leqno (14.1) 
$$
where $\alpha_+(\varphi)$ is defined in (10.20)-(10.22),
and is the maximum deviation (from the standard position 
that would occur in a minimal cone) of the position of two 
or three tangent vectors to $\varphi_\ast(K)$ at a vertex. 
So (2.12) and the full length property will follow as soon 
as we check that
$$
H^1(\varphi_\ast(K)) - H^1(K) \leq C \alpha_+(\varphi)^2
\leqno (14.2)
$$
when
$$
H^1(\varphi_\ast(K)) > H^1(K).
\leqno (14.3)
$$

\ms\proclaim Lemma 14.4.
The planes in $\R^n$ have the full length property.

\ms
Let $P$ be a plane, choose three vertices $z_j$, $0 \leq j \leq 2$,
in $K = P  \cap \partial B$, for instance at equal distances 
from each other, to define the standard decomposition of $K$.
Then let $\varphi \in \Phi(\eta_1)$ be given and set 
$w_j = \varphi(z_j)$ for $0 \leq j \leq 2$.

We may assume that $w_1$ and $w_2$ lie in the horizontal plane
$P_0 = \big\{ x\in \R^n \, ; \, x_3 = \cdots = x_n = 0 \big\}$,
and even, if we identify $P_0$ with $\Bbb C$ for convenience, 
that for $j=1,2$, $w_i = e^{i \theta_j}$ for some $\theta_j$ 
which is close to $2\pi j/3$. Denote by $\xi$ the point
of $P_0 \cap \partial B$ that lies closest to $w_0$, and set
$\delta = |\xi-w_0| = \dist(w_0,P_0 \cap \partial B)$. 
Also  denote by 
$d$ the geodesic distance on $\partial B$. Then
$$\leqalignno{
H^1(\varphi_\ast(K)) &- H^1(K) = H^1(\varphi_\ast(K)) - 2\pi
= d(w_0,w_1) + d(w_0,w_2) + d(w_1,w_2) - 2\pi
\cr&= d(w_0,w_1) + d(w_0,w_2) + d(w_1,w_2) - 
[d(\xi,w_1) + d(\xi,w_2) + d(w_1,w_2)]
\cr& = [d(w_0,w_1)-d(\xi,w_1)] + [d(w_0,w_2)-d(\xi,w_2)]
\leq C \delta^2 \leq C \alpha_+(\varphi)^2.
& (14.5)
}
$$
Thus (14.2) holds, and Lemma 14.4 follows.
\qed

\ms\proclaim Lemma 14.6.
The cones of type $\Bbb Y$ in $\R^n$ have the full length property.

\ms
The proof will be a little unnerving (at least for the author),
because we cannot really trust the pictures in $\Bbb R^3$.
Let $Y$ be a cone of type $\Bbb Y$ and $\varphi \in \Phi(\eta_1)$
be given, denote by $z^+$ and $z^-$ the two vertices of $K$,
and chose additional vertices $z_1$, $z_2$, and $z_3$ near the middle
of each arc of $K$, to define a standard decomposition.
Then set $w^\pm = \varphi(z^\pm)$ and $w_j = \varphi(z_j)$;
these are the vertices of $\varphi_\ast(K)$.

Choose coordinates of $\R^n$ so that $w^+$ and $w^-$ lie in
a vertical plane in $\R^3$, and even 
$$
w^\pm = (\sin\rho,0,\pm\cos\rho,0)
\leqno (14.7)
$$
for some $\rho \in [0,10^{-2}]$. (We just rotate the coordinates
in $\R^3$ so that the midpoint of the $w^\pm$ lies on the 
nonnegative first axis; then $\rho$ is small because the $w^\pm$
are almost antipodal.)

Denote by $\Gamma_j^\pm$ the arc of geodesic between $w^\pm$
and $w_j$, and set $\Gamma_j = \Gamma_j^+ \cup \Gamma_j^-$.
Thus $\varphi_\ast(K)$ is the union of the $\Gamma_j$.
Denote by $v_j$ the point of $\Gamma_j$ whose third coordinate
vanishes. Let us check that
$$
H^1(\Gamma_j) \leq d(v_j,w^+) + d(v_j,w^-) + C\alpha_+(\varphi)^2,
\leqno (14.8)
$$
where we still denote by $d$ the geodesic distance on 
$\partial B$.  Suppose, for the sake of definiteness, that the 
third coordinate of $w_j$ is nonnegative, so that $v_j$ lies on 
$\Gamma_j^-$ and the the geodesic $g^-$ from $v_j$ to $w^-$ is 
contained in $\Gamma_j^-$. Set $G=\Gamma_j \sm g^-$; thus $G$ is 
composed of a small arc $g$ of $\Gamma_j^-$ that goes from $v_j$
to $w_j$, followed by $\Gamma_j^+$ that goes from $w_j$ to
$w^+$. The angle between the two is different from $\pi$
by at most $\alpha_+(\varphi)$, by (10.21), so
$H^1(G) \leq d(v_j,w^+) + C \alpha_+(\varphi)^2$.
We add $H^1(g^-) = d(v_j,w^-)$ to both sides and get (14.8).

Next we evaluate $d(v_j,w^\pm)$. 
Write $v_j = (a_j,b_j,0,\xi_j)$ with $\xi_j \in \R^{n-3}$.
Set $|\xi_j| = \sin\alpha_j$, so that $a_j^2 + b_j^2 = \cos^2\alpha_j$,
and choose $\theta_j$ so that $a_j = \cos\theta_j\cos\alpha_j$ 
and $b_j = \sin\theta_j\cos\alpha_j$. Thus
$$
v_j = (\cos\theta_j\cos\alpha_j,\sin\theta_j\cos\alpha_j, 0, \xi_j).
\leqno (14.9)
$$
Since $d(v_j,w^\pm)$ is the length of an arc of great circle
that goes from $v_j$ to $w^\pm$, $|w^\pm-v_j| = 2 \sin(d(v_j,w^\pm)/2)$.
Then 
$$
\cos(d(v_j,w^\pm)) = 1-2\sin^2(d(v_j,w^\pm)/2) 
= 1 - {1 \over 2} \, |w^\pm-v_j|^2 .
\leqno (14.10)
$$
Next (14.7) and (14.9) yield
$$\leqalignno{
|w^\pm-v_j|^2 &
= [\cos\theta_j \cos\alpha_j-\sin\rho]^2 
+ \sin^2\theta_j \cos^2\alpha_j + \cos^2\rho + \sin^2\alpha_j
\cr&
= 1 + \cos^2\theta_j \cos^2\alpha_j - 2 \cos\theta_j \cos\alpha_j\sin\rho
+ \sin^2\theta_j \cos^2\alpha_j + \sin^2\alpha_j
&(14.11)
\cr&
= 1 + \cos^2\alpha_j - 2 \cos\theta_j \cos\alpha_j\sin\rho
+ \sin^2\alpha_j
= 2 - 2 \cos\theta_j \cos\alpha_j\sin\rho
}
$$
and (14.10) says that 
$\cos(d(v_j,w^\pm)) = \cos\theta_j \cos\alpha_j\sin\rho$.
Set $t_j^\pm = \pi/2-d(v_j,w^\pm)$; these will be easier to
manipulate because they are small. Record that
$$
\sin t_j^\pm = \cos(d(v_j,w^\pm)) = \cos\theta_j \cos\alpha_j\sin\rho,
\leqno (14.12)
$$
so $|\sin t_j^\pm| \leq \sin\rho \leq 10^{-2}$, hence
$|t_j| \leq {\rm Arcsin}(10^{-2})$ and then
$$
|t_j^\pm - \sin t_j^\pm| \leq {1 \over 6} \, |t_j^\pm|^3 
\leq {1 \over 3} \, |\sin^3t_j^\pm| \leq {1 \over 3} \, \sin^3\rho.
\leqno (14.13)
$$
Also, 
$$
3 \pi - \sum_{j \, ; \, \pm} d(v_j,w^\pm) 
= \sum_{j \, ; \, \pm} \big[{\pi \over 2}-d(v_j,w^\pm) \big] 
= \sum_{j \, ; \, \pm} t_j^\pm
\leqno (14.14)
$$
so 
$$
\Big| 3 \pi - \sum_{j \, ; \, } d(v_j,w^\pm) 
- 2 \sin\rho \sum_{j} \cos\theta_j \cos\alpha_j \Big|
\leq 2\sin^3\rho,
\leqno (14.15)
$$
by (14.12) and (14.13).

Next we want to evaluate angles. Let us first prove that
$$
\rho \leq 2\alpha_+(\varphi).
\leqno (14.16)
$$
Set $s_j^\pm = w^\pm - \langle w_j,w^\pm \rangle w_j$; 
observe that $s_j^\pm$ lies in the plane through $w^\pm$ and $w_j$, 
and is orthogonal to $w_j$, so it is parallel to the tangent of
$\Gamma_{j}^\pm$ at $w_j$. Also, $|s_j^\pm| \geq 9/10$ because
$\langle w_j,w^\pm \rangle$ is small (recall that $\varphi$
does not move the points very far, and that $z_j$ is orthogonal 
to $z^\pm$). Let $\pi - \beta_j$ denote the angle of $\Gamma_j^+$ 
and $\Gamma_j^-$ at $w_j$, or equivalently of $s_j^+$ and $s_j^-$. 
Then 
$$
|\beta_j| \leq \alpha_+(\varphi)
\leqno (14.17)
$$
by definition of $\alpha_+(\varphi)$. Denote by $p$
the orthogonal projection of $s_j^+$ on the direction
of $s_j^-$, and set $p^\perp = s_j^+ -p$. 
Then $|p^\perp| = |s_j^+||\sin\beta_j|$.

We want to use a determinant to evaluate this, so
we consider the vector space $V$ spanned by 
$w_j$, $w^-$, and $w^+$, with the Euclidean structure
inherited from $\R^n$ and any choice of orientation.
[If the three vectors are not independent, pick any
$V$ of dimension $3$ that contains them.]
Then compute $D_j = \det(w_j,w^-,w+)$ in $V$. Since
we may remove a linear combination of two vectors to the third
one, $D_j =\det(w_j,s_j^-,s_j^+) = \det(w_j,s_j^-,p^\perp)$,
and so $|D_j|=|s_j^-||p^\perp| = |s_j^-||s_j^+||\sin\beta_j|$
because the three vectors are orthogonal and by the computation
above. So 
$$
|D_j| \leq |\sin\beta_j| \leq \sin(\alpha_+(\varphi)), 
\leqno (14.18)
$$
by (14.17). We can also compute $D_j$ brutally. Denote
by $e_1$ and $e_3$ the first and third elements of the
canonical basis, complete the basis of $V$, and denote
by $x$, $z$, and $y$ the coordinates of $w_j$.
Thus $|y| = \dist(w_j,P_v)$, where $P_v$ is the vertical 
plane spanned by $e_1$ and $e_3$. A simple computation
using (14.7) says that 
$$
|D_j| = 2|y|\sin\rho \cos\rho = \sin(2\rho)\dist(w_j,P_v).
\leqno (14.19)
$$
Now we want to select $j$ so that $\dist(w_j,P_v)$ is not too small.
Recall that $w_j = \varphi(z_j)$, and that the three $z_j$ 
lie in a $2$-plane orthogonal to the line through $z^+$,
where they make $120^\circ$ degree angles with each other.
At least one $z_j$ lies at distance $\geq 1/2$ from the 
vector space $W$ through $e_1$ and $z^+$. 
Recall that $z^+$ is quite close to $w^+$, which is close to
$e_3$ because $\rho \leq 10^{-2}$, and since these vectors are 
almost orthogonal to $e_1$, $W$ makes a small angle
with $P_v$. Then $\dist(z_j,P_v) \geq 4/10$, and 
$\dist(w_j,P_v) \geq 1/3$ for the corresponding $w_j$. So
$$
\sin(2\rho) \leq 3 \sin(2\rho)\dist(w_j,P_v)
\leq 3|D_j| \leq 3\sin(\alpha_+(\varphi))
\leqno (14.20) 
$$
by (14.19) and (14.18), which proves (14.16).

\ms
Next we want to control some other term of (14.15),
namely $S =  \sum_{j} \cos\theta_j \cos\alpha_j$.
Recall from (14.9) and the definition above (14.8) that
$\cos\theta_j \cos\alpha_j = \langle v_j,e_1 \rangle$,
where $v_j$  is the point of $\Gamma_j$ such that
$\langle v_j,e_3 \rangle = 0$.

Denote by $\widehat \Gamma_j^+$ the geodesic that
leaves from $w^+$ in the same direction as $\Gamma_j^+$,
but which we continue a little further, so that its length
is $2\pi/3$, say. Let $h_j: [0,2\pi/3] \to \partial B$ denote 
the parameterization of $\widehat \Gamma_j^+$ by arc-length. 
We want to  check that
$$
|v_j - h_j(\pi/2)| \leq 3 \alpha_+(\varphi).
\leqno (14.21) 
$$
First notice that $\langle h_j(t),w^+ \rangle = \cos t$, so 
$$
\cos t - \rho \leq \langle h_j(t),e_3 \rangle 
\leq \cos t + \rho
\leqno (14.22) 
$$
because $|w^+-e_3| \leq \rho$ by (14.7).
Let $t_j\in[0,2\pi/3]$ be such that $\langle h_j(t_j),e_3 \rangle=0$.
Then $|\cos t_j| \leq \rho$, or equivalently, 
$|t_j-\pi/2|\leq {\rm Arcsin}(\rho)$. 

If $v_j$ lies in $\Gamma_j^+$, then $v_j = h(t_j)$ and
$|v_j - h_j(\pi/2)| \leq |t_j-\pi/2| \leq {\rm Arcsin}(\rho)
\leq 3\alpha_+(\varphi)$, by (14.16) and as needed for (14.21).

Otherwise, if $v_j \in \Gamma_j^-$, we go from $w^+$ to $v_j$
as follows. First we follow $\Gamma_j^+$ up to $w_j$, and
then we turn by at most $\alpha_+(\varphi)$ (by definition
of $\alpha_+(\varphi)$) to follow $\Gamma_j^+$. Let $T$
be such that $h_j(T) = w_j \,$; the continuation of the trip
can be parameterized by a new function $h'_j$, with
$h'_j(T)=h_j(T)=w_j$, and then 
$|h'_j(t)-h_j(t)| \leq 2(t-T)\,\alpha_+(\varphi)$
for $T \leq t \leq T+10^{-1}$. We cross the hyperplane plane $x_3=0$
much before that, because for $t=T+10^{-1}$, $h_j(t)$ lies
far under the plane, by (14.22).

Now $v_j = h'_j(t)$, where $t\in [T,T+10^{-1}]$ is such that
$\langle h'_j(t),e_3 \rangle = 0$. Then $|\langle h_j(t),e_3 \rangle| 
\leq |h'_j(t)-h_j(t)|\leq 2 (t-T) \,\alpha_+(\varphi) \leq \alpha_+(\varphi)/5$, so 
$|\cos t| \leq \rho + \alpha_+(\varphi)/5$ by (14.22), and
$|t-\pi/2| \leq 11\alpha_+(\varphi)/5$ by (14.16). Finally
$$\eqalign{
|v_j-h_j(\pi/2)| & = |h'_j(t)-h_j(\pi/2)|
\leq |h_j(t)-h_j(\pi/2)| + \alpha_+(\varphi)/5
\cr& \leq |t-\pi/2| + \alpha_+(\varphi)/5
\leq 3 \alpha_+(\varphi),
}\leqno (14.23) 
$$
which proves (14.21) in our second case.

We prefer the $h_j(\pi/2)$ because they are easy to localize.
In fact, $h_j(\pi/2)$ is precisely the unit direction of $\Gamma_j$
at $w^+$, so $\big|\sum_j h_j(\pi/2)\big| \leq \alpha_+(\varphi)$
by the definition (10.20) and (10.22).
Then $\big|\sum_j v_j \big| \leq 10 \alpha_+(\varphi)$ by (14.23),
and 
$$
\Big|\sum_{j} \cos\theta_j \cos\alpha_j \Big|
= \Big|\sum_{j} \langle v_j,e_1 \rangle \Big|
= \Big|\langle \sum_{j} v_j,e_1 \rangle \Big|
\leq \Big|\sum_j v_j \Big| \leq 10 \alpha_+(\varphi)
\leqno (14.24) 
$$
by (14.9). Thus
$$\eqalign{
\big| 3 \pi - \sum_{j \, ; \, \pm} &d(v_j,w^\pm) \big| 
\leq 2 \sin\rho \Big|\sum_{j} \cos\theta_j \cos\alpha_j \Big|
+ 2\sin^3\rho
\cr& \leq 20 \sin\rho \  \alpha_+(\varphi) + 2\sin^3\rho
\leq 40 \alpha_+(\varphi)^2 + 8 \alpha_+(\varphi)^3
\leq 41 \alpha_+(\varphi)^2
}\leqno (14.25) 
$$
by (14.15), (14.24), and (14.16). We now compare to (14.8)
and get that 
$$\eqalign{
H^1(\varphi_\ast(K)) = \sum_j H^1(\Gamma_j)
&\leq \sum_j [d(v_j,w^+) + d(v_j,w^-) + C\alpha_+(\varphi)^2]
\cr&\leq 3 \pi + C\alpha_+(\varphi)^2 = H^1(K) + C\alpha_+(\varphi)^2
}\leqno (14.26) 
$$
because $K$ is a cone of type $\Bbb Y$. Thus
(14.2) holds, and Lemma 14.6 follows.
\qed
\ms
We turn to the cone $T$ that was described in Section 1,
near Figure 1.2. That is,  $T\i \R^3$ is the cone over the 
union of the edges of a regular tetrahedron centered at the origin.

\ms\proclaim Lemma 14.27.
The isometric images in $\R^n$ of the cone $T$ 
have the full length property.

\ms
We shall denote by  $A'_j$, $1 \leq j \leq 4$, the vertices of the
regular tetrahedron used to define $T$; in the present situation,
the regular decomposition of $K = T \cap \partial B$ constructed in 
Section 2 is just its decomposition into six arcs of great circles, i.e., 
we don't need to add new vertices.

Let $\varphi \in \Phi(\eta_1)$ be given, with $\eta_1$ small;
then $\varphi_\ast(K)$ is a tetrahedron in $\R^n$, with the vertices 
$A_j = \varphi(A'_j)\in \partial B$. We want to show first that there is a
regular tetrahedron centered at the origin, with vertices $a_j \in 
\partial B$, $1 \leq j \leq 4$, and such that
$$
|a_j-A_j| \leq C \alpha_+(\varphi)
\ \hbox{ for } 1 \leq j \leq 4.
\leqno (14.28)
$$

We shall first reduce to the simpler situation when 
$\varphi_\ast(K)$ lies in $\R^3$. We may assume that
$A_1$, $A_2$, and $A_3$ lie in $\R^3$. Denote by $w_{i,j}$
the unit direction of the geodesic from $A_i$ to $A_j$,
at the point $A_i$ (and pointing away from $A_i$).
We know from the definitions (10.20) and (10.22) that
$|w_{1,2}+w_{1,3}+w_{1,4}| \leq \alpha_+(\varphi)$, so
$\dist(w_{1,4},\R^3) \leq \alpha_+(\varphi)$ because the
other vectors lie in $\R^3$. We follow the geodesic for less
than two units of length to go from $A_1$ to $A_4$, so
$\dist(A_4,\R^3) \leq 3\alpha_+(\varphi)$ (compare with the 
end of a nearby geodesic with the same length and contained in 
$\Bbb R^3$). Let $A' \in \Bbb R^3$ be such that 
$|A'-A_4| \leq 3 \alpha_+(\varphi)$, and then set
$A = A'/|A'|$. Thus $A \in \R^3 \cap \partial B$ and 
$|A-A_4| \leq 6 \alpha_+(\varphi)$.

When we replace $A_4$ with $A$, we modify the $w_{i,j}$
by less than $C \alpha_+(\varphi)$, so if we modify 
$\varphi$ by setting $\widetilde\varphi(a_4) = A$, we still have 
that $\alpha_+(\widetilde\varphi) \leq C\alpha_+(\varphi)$.
If we can prove (14.28) for tetrahedra in $\R^3$, we apply it to
$\widetilde\varphi_\ast(K)$ and get a regular tetrahedron 
which also works for $\varphi_\ast(K)$ (because 
$|A-A_4| \leq 6 \alpha_+(\varphi)$).

So we may assume that the $A_j$ lie in $\Bbb R^3$. This will be 
convenient, because this will allow us to use some identities 
on spherical triangles. 

Consider the spherical triangle $(A_1,A_2,A_3)$ for a moment,
and denote by $\alpha_j$ its angle at $A_j$, and $l_j$ the length
of the edge opposite to $A_j$. We shall use the following two formulae, 
which we take from [Be]:

$$
{\sin l_{1} \over \sin \alpha_{1}}
= {\sin l_{2} \over \sin \alpha_{2}}
= {\sin l_{3} \over \sin \alpha_{3}} 
\leqno (14.29)
$$
and
$$
\cos\alpha_1 = {\cos l_1 - \cos l_2 \, \cos l_3 \over \sin l_2 \sin l_3},
\leqno (14.30)
$$
which are respectively 18.6.13.4 and 18.6.13.7 in [Be]. 
Note that (14.30) looks a little strange at first sight, because
the denominator on the right-hand side could be small, but if this 
happens, the cosines of the numerator are all close to $1$, and the
numerator is small too. Let us check that
$$
\big|\alpha_{j} - {2 \pi \over 3}\big|
\leq 5 \alpha_+(\varphi).
\leqno (14.31)
$$
Denote by $w^1$, $w^2$, and $w^3$ the unit directions of the three
geodesics of $\varphi_\ast(K)$ that leave from $A_j$. Observe that
$$
|w^1+w^2+w^3| \leq \alpha_+(\varphi),
\leqno (14.32)
$$
by the definitions (10.20) and (10.22). 
Denote by $\alpha^1$,  $\alpha^2$, $\alpha^3$ the angles of the
$w^k$,  with $\alpha^k$ opposite to $w^k$, and choose the names so that
$\alpha_j = \alpha^1$.

Since we now work in $\R^3$, the three vectors $w^k$ lie in a same 
plane (the plane orthogonal to the direction of $A_j$), so 
$\alpha^1+ \alpha^2+\alpha^3 = 2\pi$. 
Let  us assume that $\alpha^2$ lies on the other side of 
${2 \pi \over 3}$ than $\alpha_j = \alpha^1$ (if not, $\alpha^3$
has this property); then  $|\alpha^1 - \alpha^2|  = |\alpha_j - \alpha^2| 
\geq \big|\alpha_{j} - {2 \pi \over 3}\big|$.
Let $h$ denote the coordinate of $w^1+w^2+w^3$ on the direction orthogonal
to $w^3$; then
$$\eqalign{
|h| = &\big|\sin\alpha^1-\sin\alpha^2 \big| =  \Big|2\sin\big((\alpha^1-\alpha^2)/2\big)
\cos\big((\alpha^1+\alpha^2)/2\big)\Big|
\cr&
\geq {1 \over 2} \Big| \sin\big((\alpha^1-\alpha^2)/2\big) \Big|
\geq {1 \over 5}\, |\alpha^1-\alpha^2| 
\geq {1 \over 5}\,\big|\alpha_{j} - {2 \pi \over 3}\big|
}\leqno (14.33)
$$
because only $w_1$ and $w_2$ contribute, and
all the $\alpha^j$ are fairly close to $2\pi/3$.
Now $|h| \leq |w^1+w^2+w^3|$, and (14.31) follows from 
(14.33) and  (14.32).

Return to (14.29) and (14.30). The $\sin\alpha_j$ in the formula
are all $5\alpha_+(\varphi)$-close to $\sqrt 3/2$, so (14.29)
says that the three $\sin l_j$ are all close to a same value. 
That is, set $y = \sin l_1$; then, for $j\neq 1$,
$$
|\sin l_j - y| = \Big|{y \sin\alpha_j \over \sin\alpha_1} - y \Big|
= {  y |\sin\alpha_j - \sin\alpha_1| \over \sin\alpha_1}
\leq 2 |\alpha_j - \alpha_1| \leq 20 \alpha_+(\varphi)
\leqno (14.34)
$$
by (14.29) and (14.31).

Denote by $L$ the the length of each arc of $K$. We can compute
$L$, because it is the side-length of an equilateral spherical triangle 
with $2\pi/3$ angles, and (14.30) says that $X = \cos L$ is a solution of
${ X- X^2 \over 1-X^2} = \cos(2\pi/3) = -1/2$, 
so ${X \over 1+X} = -1/2$ and $X =-1/3$. Now 
$$
|l_j-L| \leq 10^{-1}  \ \hbox{ for } 1 \leq j \leq 3
\leqno (14.35)
$$
because $\varphi \in \Phi(\eta_1)$ for some small $\eta_1$ that we 
can choose. In the range where (14.35)  holds, the right-hand
side of (14.30) is a Lipschitz function of the three $\sin l_j$.
When the three arguments are equal to $y=\sin l_1$, we get the
value ${x-x^2 \over 1-x^2} = {x \over 1+x}$, where we set 
$x = \cos l_1$. Thus
$$\eqalign{
\big| {x \over 1+x} + {1 \over 2} \big|
&\leq \big| {x \over 1+x}  - \cos \alpha_1 \big| 
+ \big| \cos \alpha_1 + {1 \over 2} \big|
\cr&
\leq C \sum_{j=2}^3 |\sin l_j - y| 
+\big| \alpha_1 - {3\pi \over 2} \big|
 \leq C \alpha_+(\varphi)
}\leqno (14.36)
$$
by (14.30), (14.34), and (14.31). We solve for $x$ and get 
that $|x+1/2| \leq C \alpha_+(\varphi)$, which in turn shows that
$|l_1- L| \leq C \alpha_+(\varphi)$ (because $\cos l_1 = x$
and $\cos L = -1/2$). This implies that
$$
|l_j-L| \leq C \alpha_+(\varphi)  \ \hbox{ for } 1 \leq j \leq 3
\leqno (14.37)
$$
because (14.34) gives some control on $|l_j-l_1|$.

\ms
Of course what we did for the spherical triangle with vertices 
$A_1$, $A_2$ and $A_3$ can also be done for the other triangles,
and we get that each angle in $\varphi_\ast(K)$ is
$5\alpha_+(\varphi)$-close to $2\pi/3$, and the length of each arc
of $\varphi_\ast(K)$ is $C\alpha_+(\varphi)$-close to $L$.

The regular tetrahedron promised in (14.28) can now be constructed
by hand: we take $a_1 = A_1$, then $a_2 \in \partial B$ in the plane 
that contains $A_1$ and $A_2$, close to $A_2$ and so that 
$d(a_1,a_2) = L$. It is automatically $C\alpha_+(\varphi)$-close 
to $A_2$. The two last points $a_3$ and $a_4$ are the only
two points of $\partial B$ that lie at geodesic distance $L$
from $a_1$ and $a_2$, and it is easy to see that they also
lie within $C\alpha_+(\varphi)$ of $A_3$ and $A_4$.

So we get (14.28), and now we need to prove (14.2).
Define a function $F$ on $(\partial B)^4$ by 
$F(a_1,a_2,a_3,a_4) = 
\sum_{1 \leq j < k \leq 4} d_{\partial B} (a_j,a_k)$
(the sum of the length of the geodesic arcs that connect the 
vertices). For the tetrahedron of (14.28), we just get
$F(a_1,a_2,a_3,a_4) = 6 L$. 

Now $F$ is smooth near $(a_1,a_2,a_3,a_4)$ , with a bounded second derivative, 
and  $(a_1,a_2,a_3,a_4)$ is a critical point of $F$; this is easy to check directly,
because the geodesic arcs that connect the vertices make angles of 
$120^\circ$ when they meet. Then 
$$\eqalign{
|H^1(\varphi_\ast(K)) - H^1(K)| 
&= |F(A_1,A_2,A_3,A_4)-F(a_1,a_2,a_3,a_4)|
\cr&\leq C \sum_{1 \leq i \leq 4} |A_i-a_i|^2 \leq C \alpha_+(\varphi)^2
}\leqno (14.38)
$$
by Taylor's formula and (14.28). This proves (14.2), and
Lemma 14.27 follows.
\qed

\ms\noindent{\bf Remark 14.39.}
The full length constants in Lemmas 14.4, 14.6, and 14.27
do not depend on $n$, for instance because the $\varphi_\ast(K)$
have at most five points, so that we can always assume that they lie in 
$\R^5$.

\ms\noindent{\bf Remark 14.40.}
The sufficient condition (14.2) that implies the 
full length condition can be checked independently on the
connected components of $K = E\cap \partial B$. That is,
if $E$ is a minimal (or minimal-looking) cone, then
$K$ satisfies (14.2) if and only if every connected component
of $K$ satisfies it. 

This is easy to check; the only point is that by the description
of Section 2, $K$ is the finite union of its connected components 
$K_j$, which lie at distances at least $\eta_0$ from each other. Then 
$\varphi_\ast(K)$ is the disjoint union of the $\varphi_\ast(K_j)$,
and (14.2) can be verified independently for each $K_j$.

Thus, when $E$ is a finite union of planes, or sets of type $\Bbb Y$
or $\Bbb T$ that only meet at the origin, $E$ automatically has the 
full length property. But even in this case the conditions under 
which $E$ is a minimal cone are not known. In fact, the only known 
example seems to be the union of minimal cones that lie in 
orthogonal spaces. 

It is a little less obvious that the full length condition itself
can be checked independently on the connected components.
Let us merely suggest a strategy that we could use to prove that
if every connected component of $K$ satisfies the full length 
condition, then $K$ itself satisfies it.

Let $\varphi \in \Phi(\eta_1)$ be such that 
$H^1(\varphi_\ast(K)) > H^1(K)$;
we can find a component $K_1$ of $K$ and a deformation 
$\widetilde X_1 = f(X_1)$ of the cone $X_1$ over $\varphi_\ast(K_1)$ 
that allows to save some amount $A$ of area, as in (2.12). 

Set $K_2 = K \sm K_1$, and denote by $X_2$ the cone over
$\varphi_\ast(K_2)$. We would like to take $f(x) = x$ on $X_2$,
and then say that
$H^2(f(X_1 \cup X_2) \cap B(0,1))
\leq H^2(f(X_1) \cap B(0,1)) + H^2(f(X_2) \cap B(0,1))
= H^2( \widetilde X_1 \cap B(0,1)) + H^2(X_2 \cap B(0,1))
\leq H^2(\varphi_\ast(K))- A$ and conclude.

However, $X_1$ and $X_2$ meet at the origin, so we cannot take
two different definitions of $f$ there, and some surgery near
the origin is needed.
Notice that the question only arises when $n \geq 4$, 
because we know the minimal cones when $n=3$. 
Then $\varphi_\ast(K_1)$ and $\varphi_\ast(K_2)$
are not only disjoint, but (by general position)
they are not linked in $\partial B$. That is, we can deform 
$\varphi_\ast(K_1)$ to a point $z_1$ inside 
$\partial B \sm \varphi_\ast(K_2)$, and then deform $\varphi_\ast(K_2)$ 
to a point in $\partial B \sm \{ z_1 \}$.
This should allow us to modify $\varphi_\ast(E) = X_1 \cup X_2$ 
in a small ball near the origin, so that the images of $X_1$ and $X_2$ are 
mostly far apart, except for a small line segment that connects them. 
[Follow the deformation above, where the distance to the origin
plays the role of the time parameter, on a small annulus.]
Then the initial plan of taking different 
definitions of the deformation on the two pieces can be realized:
what happens on the line segment does not matter, because
its image will have finite $H^1$-measure.

\bigskip
REFERENCES

\smallskip 
\item {[Al]} F. J. Almgren, Existence and regularity almost everywhere 
of solutions to elliptic variational problems with constraints, 
Memoirs of the Amer. Math. Soc. 165, volume 4 (1976), i-199.
\smallskip
\item {[Be]} M. Berger, \underbar{G\'{e}om\'{e}trie, Vol. 5}. La sph\`{e}re pour 
elle-m\^eme, g\'{e}om\'{e}trie hyperbolique, l'espace des sph\`{e}res. 
CEDIC, Paris; Fernand Nathan, Paris, 1977. 
\smallskip 
\item {[D1]} G. David, Limits of Almgren-quasiminimal sets, 
Proceedings of the conference on Harmonic Analysis, 
Mount Holyoke, A.M.S. Contemporary Mathematics series, Vol. 320 
(2003), 119-145.
\smallskip
\item {[D2]} G. David, Singular sets of minimizers for 
the Mumford-Shah functional,
Progress in Mathematics 233 (581p.), Birkh\"auser 2005.
\smallskip
\item {[D3]} G. David, Low regularity for almost-minimal sets in 
$\R^3$, preprint, Universit\'{e} de Paris-Sud, 2005, and submitted
in 2007.
\item {[DDT]} G. David, T. De Pauw, and T. Toro,
A generalization of Reifenberg's theorem in $\Bbb R^3$, 
to appear, Geometric And Functional Analysis.
\smallskip
\item {[DS]} G. David and S. Semmes, Uniform rectifiability and 
quasiminimizing sets of arbitrary codimension, 
Memoirs of the A.M.S. Number 687, volume 144,  2000.
\item {[Fe]} H. Federer, \underbar{Geometric measure theory}, 
Grundlehren der Mathematishen Wissenschaf-ten 
153, Springer Verlag 1969.
\smallskip
\item {[He]} A Heppes, Isogonal sph\"arischen netze, Ann. Univ. Sci. Budapest
E\"otv\"os Sect. Math. 7 (1964), 41-48.
\smallskip
\item {[Ku]} C. Kuratowski, \underbar{Topologie}, vol. II,  troisi\`{e}me \'{e}dition,
Monografie Matematyczne, Tom XX. Polskie Towarzystwo 
Matematyczne, Warsawa, 1952, xii+450 pp, or reprinted by
ƒditions Jacques Gabay, Sceaux, 1992, iv+266 pp.
\smallskip
\item {[La]} E. Lamarle, Sur la stabilit\'{e} des syst\`{e}mes liquides en 
lames minces, M\'{e}m. Acad. R. Belg. 35 (1864), 3-104.
\item {[Le]} A. Lemenant, Sur la r\'{e}gularit\'{e} des minimiseurs 
de Mumford-Shah en dimension 3 et sup\'{e}rieure,
Thesis, Universit\'{e} de Paris-sud, June 2008.
\smallskip  
\item {[Ma]}  P. Mattila, \underbar{Geometry of sets and 
measures in Euclidean space}, Cambridge Studies in
Advanced Mathematics 44, Cambridge University Press l995.
\smallskip
\item {[Mo1]} F. Morgan, Size-minimizing rectifiable currents, 
Invent. Math. 96 (1989), no. 2, 333-348.
\smallskip
\item {[Ne]} M. H. A. Newman,  \underbar{Elements of the topology of plane 
sets of points}, 
Second edition, reprinted, Cambridge University Press, New York 1961.
\smallskip
\item {[R1]} E. R. Reifenberg, Solution of the Plateau Problem for 
$m$-dimensional surfaces of varying topological type,
Acta Math. 104, 1960, 1--92.
\smallskip
\item {[R2]} E. R. Reifenberg, An epiperimetric inequality related to 
the analyticity of minimal surfaces, Ann. of Math. (2) 80, 1964, 1--14.
\smallskip
\item {[St]}	E. M. Stein, \underbar{Singular integrals and 
differentiability properties of functions},
Princeton university press 1970.
\smallskip
\item {[Ta]} J. Taylor, The structure of singularities in 
soap-bubble-like 
and soap-film-like minimal surfaces, 
Ann. of Math. (2) 103 (1976), no. 3, 489--539.
\smallskip

 \bigskip
\vfill \vfill \vfill\vfill
\noindent Guy David,  
\smallskip\noindent 
Math\'{e}matiques, B\^atiment 425,
\smallskip\noindent 
Universit\'{e} de Paris-Sud, 
\smallskip\noindent 
91405 Orsay Cedex, France
\smallskip\noindent 
guy.david@math.u-psud.fr

\bye